\documentclass{ws-book9x6}
\usepackage{ws-book-thm}           
\usepackage{ws-book-har}           

\title{The hyperboloidal foliation method}

\usepackage{mathabx}

\usepackage{ulem}
\usepackage{pdflscape}

\usepackage{amsfonts}
\usepackage[mathcal]{euscript}
\usepackage{mathrsfs}

\usepackage{ulem}
\usepackage{graphicx}

\numberwithin{equation}{section}

\numberwithin{theorem}{section}

\numberwithin{proposition}{section}

\numberwithin{lemma}{section}

\numberwithin{remark}{section}

\newcommand \Kcal {\mathcal K}
\newcommand \Hcal {\mathcal H}

\newcommand \bei {\begin{itemize}}
\newcommand \eei {\end{itemize}}

\newcommand \Zscr {\mathscr Z}

\newcommand \mbar {\underline m}

\newcommand \eps {\epsilon}

\newcommand \be {\begin{equation}}
\newcommand \ee {\end{equation}}
\newcommand \bel {\be \label}

\newcommand \RR {\mathbb R}

\newcommand \ih {\widehat {\imath}}

\newcommand \ic {\widecheck {\imath}}

\newcommand \jh {\widehat{\jmath}}

\newcommand \jc {\widecheck{\jmath}}

\newcommand \kh {\widehat{k}}

\newcommand \kc {\widecheck{k}}

\newcommand \lh {\widehat{l}}

\newcommand \lc {\widecheck{l}}

\newcommand \alphar {{{\alpha}}}
\newcommand \betar {{{\beta}}}
\newcommand \gammar{{\gamma}}
\newcommand \ar {{{a}}}
\newcommand \br {{{b}}}

\newcommand \jr {{{j}}}
\newcommand \kr {{{k}}}

\newcommand \Is {I^{\sharp}}
\newcommand \If {I^{\flat}}
\newcommand \Id {I^{\dag}}
\newcommand \Js {J^{\sharp}}
\newcommand \Jd {J^{\dag}}

\newcommand \delu {\underline{\del}}
\newcommand \Au {\underline{A}}
\newcommand \Bu {\underline{B}}
\newcommand \Gu {\underline{G}}

\newcommand \Tu {\underline{T}}

\newcommand \xb{\bar x}
\newcommand \delb{\bar \del}
\newcommand \Gbar{\bar G}

\newcommand \thetau{\underline{\theta}}
\newcommand \Gammau{\underline{\Gamma}}
\newcommand \Thetau{\underline{\Theta}}

\newcommand \del \partial
\newcommand \la \langle
\newcommand \ra \rangle

\newcommand \Hf {{\textbf H}}

\newcommand \Gb {{\bf G}}
\newcommand \Ib {{\bf I}}
\let\oldmarginpar\marginpar
\renewcommand\marginpar[1]{\-\oldmarginpar[\raggedleft\footnotesize #1]%
{\raggedright\footnotesize #1}}

\newcommand \Kcoef {\kappa_1}

 \usepackage[perpage]{footmisc}

\usepackage{rotating}

\DeclareMathSymbol{*}{\mathbin}{symbols}{"03} 
\DeclareMathSymbol{\ast}{\mathbin}{symbols}{"03}

\begin{document}

\centerline{\Large \bf  The hyperboloidal foliation method}

\vskip.3cm  


\

\

\

\centerline{\Large Philippe G. LeFloch\footnote{
Laboratoire Jacques-Louis Lions,
 Centre National de la Recherche Scientifique (CNRS), 
Universit\'e Pierre et Marie Curie (Paris 6), 4 Place Jussieu, 75252 Paris, France.
\newline
Email: {\it contact@philippelefloch.org, ma@ann.jussieu.fr}
\,  Blog: {\it http://philippelefloch.org} 
}
and Yue Ma$^\ast$
}

\

\

\

\

\

\titlepages                        



\begin{preface}

\vskip-1.6cm 

The Hyperboloidal Foliation Method presented in this monograph 
is based on a  $(3+1)$--foliation of Minkowski spacetime by hyperboloidal hypersurfaces.
It allows us to establish global-in-time existence results for systems of nonlinear wave equations posed on a curved spacetime and 
to derive uniform energy bounds and optimal rates of decay in time. We are also able to encompass the wave equation and the Klein-Gordon equation in a unified framework and to establish a well-posedness theory for nonlinear wave-Klein-Gordon systems and a large class of nonlinear interactions. 

The hyperboidal foliation of Minkowski spactime we rely upon in this book has the advantage of being {\sl geometric in nature} and, especially, {\sl invariant} under Lorentz transformations. As stated, our theory applies to many systems arising in mathematical physics and involving a massive scalar field, such as the Dirac-Klein-Gordon system. 
As it provides uniform energy bounds and optimal rates of decay in time, our method appears to be very robust and should extend to even more general systems. 

We have built upon many earlier studies of nonlinear wave equations or Klein-Gordon equations, especially by 
Sergiu Klainerman, Demetri Christodoulou, Jalal Shatah, Alain Bachelot, and many others. 
The coupling of  nonlinear wave-Klein-Gordon systems was first understood by  
Soichiro Katayama who succeeded to establish an existence theory of such systems.

Importantly, in developing the Hyperboloidal Foliation Method, we were inspired by earlier work on the Einstein equations of general relativity by Helmut Friedrich, Vincent Moncrief, and Anil Zenginoglu. 

We are very grateful Soichiro Katayama for observations he made to the authors on a preliminary version of this monograph. 

Last but not least, the authors are very grateful to their respective families for their strong support. 

\vskip.4cm 

\begin{flushright}
{\it Philippe G. LeFloch and Yue Ma}

Paris, September 2014
\end{flushright}

\end{preface}



\tableofcontents



\setcounter{page}{1}


\chapter[Introduction]{Introduction\label{cha:1}}

\section{Background and main objective}
\label{sec:11}

We are interested in nonlinear wave equations posed on the $(3+1)$-dimensional
Minkowski spacetime and, especially, in models arising in mathematical physics and involving a nonlinear coupling with the Klein-Gordon equation. 
A typical example is provided by the Dirac-Klein-Gordon equation. 
(Cf.~Section~\ref{sec:14}, below.) Our study is also motivated by the Einstein equations of general relativity when the matter model is a massive scalar field. The Klein-Gordon equation also describes nonlinear waves that propagate in fluids or in elastic materials.

The so-called `vector field method'  was introduced by 
\citet{Klainerman80,Klainerman85,Klainerman86,Klainerman87}.   
It is based on weighted norms defined from the conformal Killing fields of Minkowski spacetime
and Sobolev-type arguments,
and yields a global-in-time, well-posedness theory for the initial value problem for  nonlinear wave equations,
when the initial data have small amplitude and 
are ``localized'', that is, compactly supported.

This method relies on a bootstrap argument and on an analysis of the time decay of solutions. It applies to quadratic nonlinearities satisfying the so-called `null condition' introduced by Klainerman and Christodoulou.
A vast literature is now available on nonlinear wave equations. (Cf.~Section~\ref{sec:13} for further references.) As far as coupled systems of wave and Klein-Gordon equations are concerned, the current state-of-the-art 
is given by a recent work by \citet{Katayama12a} who succeded to extend the vector field method. 

In this monograph, building upon these earlier works, we introduce a novel approach, which we refer to as the `hyperboloidal foliation method' and we establish a global-in-time existence theory for a broad class of nonlinear wave-Klein-Gordon systems. 
In short, by working with suitably weigthed spacetime norms, we are able to cover the large class of quadratic nonlinearities, also recently treated by \citet{Katayama12a}), while our method demands limited regularity on the initial data and provides sharp bounds on the asymptotic profile of solutions. The hyperboloidal foliation method introduces a novel methodology, which takes its roots in an observation\footnote{recalled in Section~\ref{sec:21}, below} by \citet{Hormander97} for the Klein-Gordon equation.

Let us denote by $\Box$ the wave operator\footnote{Our convention here is $\Box := \del_t\del_t - \sum_{a=1}^3\del_a\del_a$.}
 in Minkowski space. For our purpose in this monograph, a simple (yet challenging) model of interest is provided by the following two equations which  
couple a wave equation with a Klein-Gordon equation
\bel{eq:WKG}
\aligned
\Box u = P(\del u, \del v),
\\
\Box v + v = Q(\del u, \del v),
\endaligned
\ee
where the unknowns are the two scalar fields $u, v$.
This model describes the nonlinear interactions between a massless scalar field and a massive one. Here, the nonlinear terms
$P=P(\del u, \del v)$ and $Q=Q(\del u, \del v)$ are quadratic forms in the first-order spacetime derivatives $\del u, \del v$, and account for self-interactions as well as interactions between the two fields.

Recall that global-in-time existence results for nonlinear wave equations (without Klein-Gordon components) is established when the nonlinearities satisfy the null condition. (Cf.~\eqref{main structure c}, below.) 
On the other hand, the vector field method applies also to the nonlinear Klein-Gordon equation, as shown by \citet{Klainerman85}. Recall also that the global existence problem for the nonlinear Klein-Gordon equation
was also solved independently by \citet{Shatah85} with a different method. 

However, when one attempts to tackle {\sl coupled} systems of wave and Klein-Gordon equations like \eqref{eq:WKG},
one faces a major challenge due to the fact that one of the conformal Killing fields associated with the wave equation (the scaling vector field denoted below by $t\del_t + r\del_r$) is {\sl not} a conformal Killing field for the Klein-Gordon equation and, therefore, can no longer be used in the vector field analysis.
\citet{Katayama12a} succeeded to circumvent this difficulty and established a global existence theory for a class of coupled systems which includes \eqref{eq:WKG}. 
His method relies on a novel $L^\infty$-$L^\infty$ estimate. (See Section~\ref{sec:13} for further historical background.)


\section{Statement of the main result}

The new method we provide in the present monograph relies on a fully geometric foliation. It is `robust' as it is expected to be applicable to large classes of curved spacetimes and nonlinear hyperbolic equations. 

We state here the main result that we will establish in this monograph. 
We are interested in the Cauchy problem for the following large class of {\bf nonlinear systems of wave--Klein-Gordon equations}:
\bel{main eq main}
\aligned
&\Box w_i + G_i^{\jr\alphar\betar}(w,\del w)\del_{\alphar}\del_{\betar}w_{\jr} + c_i^2 w_i = F_i(w,\del w),
\\
&w_i(B+1,x) =  {w_i}_0,
\\
 & \del_t w_i(B+1,x) =  {w_i}_1,
\endaligned
\ee
in which the unknowns are the functions $w_i$ ($1\leq i \leq n_0$) defined on Minkowski space $\RR^{3+1}$
and ${w_i}_0,  {w_i}_1$ are prescribed initial data. 
Throughout, Latin indices $a,b,c$ will
take values within $1,2,3$, while Greek indices $\alpha,\beta,\gamma$ 
take values within $0,1,2,3$.
Einstein's convention on repeated indices is in order throughout.

We assume the {\bf symmetry conditions}
\bel{pre condition symmetry}
G_i^{j\alpha\beta} = G_j^{i\alpha\beta},\quad G_i^{j\alpha\beta} = G_i^{j\beta\alpha}
\ee
and, for definiteness,
the {\bf wave-Klein-Gordon structure}
\bel{eq:lesc}
c_i \,
\begin{cases}
= 0, \qquad & \hskip.75cm 1\leq i\leq j_0,
\\
 \geq \sigma,   \qquad & j_0+1 \leq i \leq n_0,
\end{cases}
\ee
where $\sigma>0$ is a (constant, positive) lower bound for the mass coefficients of Klein-Gordon equations. 
\begin{subequations}\label{main structure}
We decompose the {\bf ``curved metric''} coefficients 
$G_i^{j\alpha\beta}(w,\del w)$ and the {\bf interaction terms} $F_i(w,\del w)$ in the form
\bel{main structure a}
G_i^{j\alpha\beta}(w,\del w) = A_i^{j\alpha\beta\gammar \kr}\del_{\gammar} w_{\kr} + B_i^{j\alpha\beta \kr}w_{\kr},
\ee
\bel{main structure b}
F_i(w,\del w) = P_i^{\alphar\betar\jr\kr}\del_{\alphar}w_{\jr} \del_{\betar}w_{\kr} + Q_i^{\alphar\jr\kr} w_{\kr}\del_{\alphar}w_{\jr} + R_i^{\jr\kr}w_{\jr}w_{\kr},  
\ee
in which we can restrict 
the summation in \eqref{main structure a} and \eqref{main structure b}
 to the range $j\leq k$.
For simplicity in the presentation of the method and without genuine loss of generality,
we focus on quadratic nonlinearities and assume that the coefficients $A_i^{j\alpha\beta\gamma k}$, $B_i^{j\alpha\beta k}$, $P_i^{\alpha\beta jk}$, $Q_i^{\alpha jk}$, and $R_i^{jk}$ are constants.

In order to simplify the notation, we adopt the following index convention: 
$$
\aligned
& \text{all indices $i,j,k,l$ take the values $1,\ldots,n_0$,}
\\
&\text{all indices $\ih,\jh,\kh,\lh$ take the values $1,\ldots, j_0$,}
\\
&\text{all indices $\ic,\jc,\kc,\lc$ take the values $j_0+1,\ldots,n_0$.}
\endaligned
$$
It is also convenient to write
$$
u_{\ih} := w_{\ih}
$$
for the components satisfying wave equations, or {\bf wave components} for short,
 and, analogously, 
$$
v_{\ic} := w_{\ic}
$$
for the components satisfying Klein-Gordon equations, or {\bf Klein-Gordon components}.

Our main assumptions are the {\bf null condition for wave components}
\index{null}
\bel{main structure c}
\aligned
& A_{\ih}^{\jh\alphar\betar\gammar \kh}\xi_{\alphar}\xi_{\betar}\xi_{\gammar}
= B_{\ih}^{\jh\alphar\betar \kh}\xi_{\alphar}\xi_{\betar}
= P_{\ih}^{\alphar\betar\jh\kh}\xi_{\alphar}\xi_{\betar} = 0
 \\
&\text{whenever } (\xi_0)^2 - \sum_a (\xi_a)^2 = 0,
\endaligned
\ee
and the {\bf non-blow-up condition} 
\bel{main structure d}
B_i^{\jc\alpha\beta\kh}=R_i^{j\kh} = R_i^{\jh k} = 0.
\ee
Moreover, we impose that  
\bel{main structure d2}
Q_i^{\alpha j\kh} = 0,
\ee
\end{subequations}
which is our only genuine restriction\footnote{When this condition is violated, solutions may not have the time decay and asymptotics of solutions to linear wave or Klein-Gordon equations in Minkowski space.}  
 required for the present implementation of the hyperboloidal foliation method.  
We emphasize that 
$$
\aligned
& \text{the null condition \index{null}
 in \eqref{main structure c} is imposed for the quadratic forms}
\\
& \text{associated with the {\sl wave components} only,}
\endaligned
$$
 and that no such restriction is required for the Klein-Gordon components.

We observe\footnote{as pointed out to the authors by S. Katayama}
 that \eqref{main structure d} combined with the symmetric condition \eqref{pre condition symmetry} leads to the following restriction on $B_{i}^{j\alpha\beta k}$:
\bel{main structure f}
B_{\jc}^{\ih\alpha\beta\kh} = 0. 
\ee
This class of systems was also studied in \citet{Katayama12a} by a completely different approach. 

We are now in a position to state the main result that we establish in the present monograph with the Hyperboloidal Foliation Method.  

\begin{theorem}[Global well-posedness theory]
\label{main thm main}
Consider the initial value problem \eqref{main eq main} with smooth initial data
posed on the spacelike hypersurface $\{t= B+1\}$ of constant time
and compactly supported in the ball $\{t= B+1; \, |x|\leq B\}$.
Under the conditions \eqref{pre condition symmetry}--\eqref{main structure},
there exists a real $\eps_0>0$ such that, for all initial data
${w_i}_0, {w_i}_1: \RR^3 \to \RR$ satisfying the smallness condition
\be
\sum_i \| {w_i}_0 \|_{\Hf^6(\RR^3)} + \| {w_i}_1 \|_{\Hf^5(\RR^3)} < \eps_0,
\ee
the Cauchy problem  \eqref{main eq main} admits a unique, smooth global-in-time solution. In addition, the energy of the wave components
--that is, $\sum_{|I|\leq 3} E_{m,c_i}(s,Z^I u_{\ih})$ defined in Section~\ref{sec:21} below--
remains globally bounded in time. 
\end{theorem}

In the special case $n_0=j_0$, the system \eqref{main eq main} contains only wave equations
and the statement in Theorem~\ref{main thm main} reduces to the classical existence result for quasilinear wave equations satisfying the null condition.
Our method is somewhat simpler than the classical proof in this case, as we will show in Chapter~\ref{cha:10}. 
In the opposite direction, if we take 
$j_0 = 0$, the system under consideration contains Klein-Gordon components only, and our result reduces to the classical existence result for quasilinear Klein-Gordon equations.

An outline of this monograph is as follows.
In Chapter~\ref{cha:2}, we introduce some basic notations on the hyperboloidal foliation and the associated energy, and we formulate our bootstrap assumptions.
Chapter~\ref{cha:3} is devoted to the derivation of fundamental properties of vector fields and commutators and their decompositions, which we will use throughout this book.

In Chapter~\ref{cha:4}, we discuss the null condition in the proposed semi-hyperboloidal frame
and we derive preliminary estimates on first- and second-order derivatives of the solutions.
In Chapter~\ref{cha:5}, we present some technical tools and, especially, a Sobolev inequality on hyperboloids
and a Hardy inequality along the hyperboloidal foliation.

At this juncture, in Chapter~\ref{cha:10} we apply our method
to scalar semilinear wave equations and we provide a new proof of the standard existence result for such equations 

We pursue our analysis, in Chapter~\ref{cha:6}, with fundamental estimates in the $L^\infty$ and $L^2$ norms, which follow from our bootstrap assumptions.
Chapter~\ref{cha:7} is devoted to controlling certain `second-order derivatives' of 
the wave components (in a sense explained therein), while
Chapter~\ref{cha:8} deals with quadratic terms satisfying the null condition and also discusses additional estimates that mainly rely on the time decay of solutions. 

Next, in Chapter~\ref{cha:9}, we derive $L^2$ estimates for nonlinear interaction terms and, therefore, complete our bootstrap argument. 

Finally, for the sake of completness, 
in Chapter~\ref{cha:11} we sketch the local-in-time existence theory.
Furthermore, the bibliography at the end of this book provides the reader with further material of interest.


\section{General strategy of proof}
\label{sec:12}

In this section, we present several key features of our method, while
 referring to Chapter~\ref{cha:2} for the relevant notions (hyperboloidal foliation, bootstrap estimates, etc.).

\bei

\item {\bf Hyperboloidal foliation.}

Most importantly, in this book
we propose to work with the {\sl family of hyperboloids} which generate a foliation of the interior of the light cone in Minkowski spacetime. (Cf.~\eqref{eq:cone}, below.) 
In contrast with other foliations of Minkowski space which are adopted in the literature, the hyperboidal foliation
has the advantage of being {\sl geometric in nature} and {\sl invariant} under Lorentz transformations.
It is therefore quite natural to search for an estimate of the energy defined on these hypersurfaces, rather than the energy defined on flat hypersurfaces of constant time, as is classically done.

\item {\bf The semi-hyperboloidal frame.}

Furthermore, to the hyperboloidal foliation we attach a {\sl semi-hyperboloidal frame} (as we call it), which consists of three vectors tangent to the hyperboids plus a timelike vector.   
This  frame has several advantages in the analysis in comparison with, for instance, the `null frame' which is often used in the literature and, instead, is defined from vectors tangent to the light cone. Importantly, the semi-hyperboloidal frame is {\sl regular} (in the interior of the light cone, which is the region of interest), while the null frame is singular at the center ($\{r=0\}$, say).

\item {\bf Decomposition of the wave operator.}

We also introduce a decomposition of the wave operator $\Box$ with respect to the semi-hyperboloidal frame, which yields us an expression of the second-order time derivative $\del_t \del_t$ of the wave components 
in terms of better-behaved derivatives. (Cf.~Proposition~\ref{decompW}, below.)

\item {\bf The hyperboloidal energy.}

Our method takes advantage of the {\sl full} expression\footnote{In constrast, \citet{Hormander97} only sought for a control of the zero-order term of the energy.}
 of the energy flux induced
on the hyperboloids 
in order to estimate certain {\sl weigthed derivatives on the hyperboloids}. This appears to be essential in order to encompass wave equations and Klein-Gordon equations in a single framework.

\item {\bf  Sobolev inequality on hyperboloids.}

In order to establish decay estimates (or $L^\infty$ estimates) on the solutions, we must control
 various commutators of fields and of operators and, next, apply suitable embedding theorems. To this purpose, we rely on a Sobolev inequality on hyperboloids, first derived by \citet{Hormander97}. 

\item {\bf  Hardy inequality on hyperboloids.}

Furthermore, we also need a new embedding estimate, that is, a Hardy inequality on hyperboloids, which we
establish in this book and is essential in eventually deriving an $L^2$ estimate on the `metric coefficients'.  Section~\ref{sec:52} for a precise statement.

\item {\bf Bootstrap strategy and hierarchy of energy bounds}

Our bootstrap formulation below consists of a {\sl hierarchy of energy bounds}, involving several levels of regularity of the wave components and the Klein-Gordon components.  
This rather involved bootstrap argument is necessary (and natural) in order to handle the 
coupling of wave equations and Klein-Gordon equations: the derivatives of different order of Klein-Gordon components enjoy different decay behaviors and different energy bounds.
 
\eei

We refer to Chapter~\ref{cha:2} for further details and continue with several observations concerning the scope of  Theorem~\ref{main thm main}.

1. As stated in the theorem, the energy of the wave components, 
that is, the quantity
$E_{m,c_i}(s,Z^I u_{\ih})$ (defined in Section~\ref{sec:21}) 
remains globally bounded for all $|I|\leq 3$, that is, up to fourth-order derivatives. Hence, the wave components have not only `small' amplitude but also `small' energy.
At the end of Chapter~\ref{cha:2}, we also establish that the standard flat energy
(that is, the quantity $\|\del_{\alpha}Z^Iu_{\ih}(t,\cdot)\|_{L^2(\RR^3)}$ defined at the end of Chapter 2, below) is also uniformly bounded for all times.
Standard methods of proof lead to possibly {\sl unbounded} high-order energies.

2. We also emphasize that, in  Theorem~\ref{main thm main}, the initial data belong
 to $H^6$,
while \citet{Katayama12a} assumes a very high regularity on the initial data  (that is, a bound on the first $19$ derivatives). In this latter paper however, the initial data 
need not be compactly supported, 
but have sufficient decay in all spatial directions. 

3. In principle, $H^4$ would be the optimal regularity in order to work with a uniformly bounded metric and to 
apply the vector field technique: namely, 
to guarantee the coercivity of, both, the flat and the hyperboloidal energy functionals, 
we need a sup-norm bound of the second-order derivatives 
the `curved metric' terms  $G_i^{j\alpha\beta}(w,\del w)$. 
In spatial dimension three,
Sobolev's embedding theorem $H^m \subset L^\infty$ holds, provided $m> 3/2$. Allowing only integer exponents, we see that $H^4$ would be optimal.

4. Certain nonlinear interaction terms may lead to a {\sl finite time blow-up} of the solutions, especially 
$$
u_{\ih} u_{\jh}, \quad u_{\ih} \del_\alpha u_{\jh}, 
\quad u_{\ih}  v_{\jc}, \quad u_{\ih}  \del_\alpha v_{\jc}
$$
and are, therefore, naturally excluded in Theorem~\ref{main thm main}.

5. In short, by denoting by $Q$ an arbitrary quadratic nonlinearity and by $N$ an arbitrary quadratic null form
and by using the notation $u,v$ for arbitrary wave/Klein-Gordon components,
the terms allowed in Theorem~\ref{main thm main} for wave equations are
$$
\aligned
& Q(v,v),\quad Q(v,\del v),\quad Q(v,\del u), \quad Q(v,\del\del v)\quad Q(v, \del\del u), \quad
\\
& Q(\del v,\del v),\quad Q(\del v, \del u), \quad Q(\del v,\del\del v), \quad Q(\del v,\del\del u),
\\
& Q(\del u,\del\del v),
\\
&N(u,\del\del u),\quad N(\del u,\del u),\quad N(\del u, \del\del u), 
\endaligned
$$ 
while, in Klein-Gordon equations, we can include the terms
$$
\aligned
& Q(v,v),\quad Q(v,\del v),\quad Q(v,\del u),\quad 
Q(v,\del\del v), \quad Q(v,\del\del u), 
\\
& Q(\del v,\del v),\quad Q(\del v, \del u), \quad Q(\del v,\del\del v),\quad Q(\del v,\del\del u),
\\
& Q(\del u, \del\del v), \quad Q(\del u,\del u),\quad Q(\del u, \del\del u).
\endaligned
$$


\section{Further references on earlier works}
\label{sec:13}

For a background on the vector field method and the global well-posedness for
nonlinear wave equations, in additional to the references cited earlier, especially 
the pioneering paper \citet{Klainerman80},  
we refer to the textbooks \citet{Hormander97} and \citet{Sogge08}.
Additional background on nonlinear wave equations is found in \citet{Strauss89}.
We do not attempt to review the large literature and only mention some selected results,
while referring the reader to the bibliography at the end of this monograph.

As already mentioned, the first results of global existence for nonlinear wave equations in three spatial dimensions
were established by \citet{Klainerman86} and \citet{Christodoulou86} under the assumption
that the nonlinearities satisfy the null condition
\index{null}
and when the equation is posed in the flat Minkowski space.
A (very) large literature is available for equations posed on curved spaces and, once more, we do not try to be exhaustive.
We refer to  \citet{Lindblad90}, \citet{KlainermanSideris96}, \citet{KlainermanSelberg97},
\citet{KlainermanMachedon97}, \citet{BahouriChemin99},
 \citet{Tataru00,Tataru01,Tataru02},
\citet{Alinhac04,Alinhac06}, \citet{LindbladRodnianski}, and \citet{LNS12}.

The Klein-Gordon equation on curved spaces was also studied in \citet{Bachelot94,Bachelot11}.

Since the decay of solutions to
 the (linear) Klein-Gordon equation is $t^{-d/2}$ in dimension $d \geq 1$, the decay function $t^{-d/2}$ is not integrable in dimension two and specific arguments are required in two dimensions:
\citet{DFX04} have treated quadratic quasilinear Klein-Gordon systems in two space dimensions and, more precisely,
 coupled systems of two equations with masses satisfying $m_1 \neq 2 m_2$ and $m_2 \neq 2m_1$ with general nonlinearities.
Furthermore, \citet{DFX04} could treat the case of equality when the null condition is assumed. This work simplified and generalized (by including resonant cases)  the earlier works by \citet{OTT95,OTT96}, \citet{Tsutsumi03a,Tsutsumi03b} and \citet{Sungawa03,Sungawa04}. See also \citet{KOS} for the algebraic characterization of the null condition 
and the asymptotic behavior of solutions, as well 
as \citet{KawaharaSunagawa} for a condition weaker than the null condition. 

More recently, \citet{Germain10} revisited the global existence theory in dimension three
for coupled Klein-Gordon equations with different speeds and systematically analyzed resonance effects. Systems of wave equations for different speeds  were
 also studied by \citet{Yokoyama} and \citet{SiderisTu01} under the null condition.
See also \citet{HoshigaKubo00}, \citet{Katayama12c}, \citet{KatayamaYokoyama}, 
and \citet{KubotaYokoyama}. 

On the other hand, Klein-Gordon equations in one space dimension are treated by quite different methods. (Cf.~\citet{Delort01,Sungawa03,Candy13}.) 

The problem of global existence for coupled systems of wave and Klein-Gordon equations have attracted much
less attention so far in the literature. In addition to the references already quoted, let us mention
 \citet{Bachelot88} who first treated the Dirac-Klein-Gordon system.
Furthermore, results on the blow-up of solutions were established by \citet{John79,John81} and, more recently,  \citet{Alinhac00}.
 
More recenty, a novel method to study nonlinear wave equations (applicable also to other dispersive systems) 
was introduced in \citet{Shatah10}, \citet{GMS12}, 
\citet{PS13}, 
\citet{P13}, and \citet{BG14}, which is based on an analysis of space-time resonances. 
See also the review by \citet{Lannes13}. 

Hyperboloidal foliations were used first by \citet{Friedrich81,Friedrich83,Friedrich02} in order to establish a global existence result for the Einstein equations.  His proof was based on a conformal transformation of the Einstein equations and an  analysis of the regularity of its solutions at infinity.
This was motivated by earlier work by \citet{Penrose} on the compactification of spacetimes.
This idea was later developed by \citet{Frauendiener98,Frauendiener02,Frauendiener04} and \citet{RinneMoncrief}.
The importance of hyperboloidal foliations for general hyperbolic systems was emphasized in \citet{Zenginoglu08,Zenginoglu11} in order to numerically compute solutions within an unbounded domain. The general standpoint in these works is that, by compactification of the spacetime, one can conveniently formulate an `artificial' outer boundary and numerically compute asymptotic properties of interest.


\section{Examples and applications}

\label{sec:14}

\subsection*{The Maxwell-Klein-Gordon system}

The theory presented in this book applies to many systems arising in mathematical physics, and we present here a few of them. For instance, the Maxwell-Klein-Gordon system in Coulomb gauge takes the form of a system of nonlinear wave equations for real-valued unknown $A_j$ and a complex-valued field $\phi$
\be
\aligned
-\Box A_j &= - \Im(\phi \overline {\del_j \phi}) + |\phi|^2 \, A_j - \del_j \del_t A_0,
\\
-\Box \phi &= 2 \sqrt{-1} \big( - A^j \del_j \phi + A_0 \del_t \phi \big) + \sqrt{-1} \del_t A_0 \phi +
\big(A_\alpha A^\alpha + m^2 \big) \phi,
\endaligned
\ee
with auxillary unknown $A_0$ given by
\be
\Delta A_0 = - \Im(\phi \overline {\del_t \phi}) + |\phi|^2 \, A_0,
\ee
supplemented with an elliptic constraint equation imposed on the initial data
\be
\aligned
\del^j A_j= 0.
\endaligned
\ee


\subsection*{The Dirac-Klein-Gordon system}

Consider next the following coupling between the Dirac equation and the wave or Klein-Gordon equation
(with $m, \sigma \geq 0$):
\bel{eq:DiracKG}
\aligned
& - \sqrt{-1} \sum_{\alpha=0}^3\Gamma_\alpha \, \del_\alpha \psi+ m \psi = \beta \, v \Gamma_0 \Gamma_1 \Gamma_2 \Gamma_3 \psi,
\\
& \Box v + \sigma^2 v = \psi^\dag K \psi,
\endaligned
\ee
in which the unknown are the ($\mathbb{C}^4$-valued) spinor field $\psi$ and the (real-valued) scalar field $v$.
We have denoted by $\psi^\dag$ the complex conjugate transpose of $\psi$.
Here, $\beta$ is a coupling constant and $K$ a constant $4 \times 4$ matrix,
while the $4 \times 4$ matrices $\Gamma_\alpha$ are the so-called Dirac matrices which are essentially characterized by the commutation conditions
\be
\Gamma_\alpha \Gamma_\beta + \Gamma_\beta \Gamma_\alpha = - 2 I \, m_{\alpha \beta},
\ee
where $m_{\alpha\beta}$ is the Minkowski metric $\text{diag} \big(1, -1, -1,- 1 \big)$
and $I$ denotes the $4 \times 4$ identity matrix. 
From the Dirac equation, on can deduce second-order equations for real-valued unknowns (so that our theory applies): this is done 
by composing the Dirac operator with itself
(since, roughly speaking, the Dirac operator is the ``square-root'' of the wave operator) and then considering the real and imaginary parts of $\psi$.


\subsection*{The Einstein equations}

Although our theory in its present form does not directly apply to the 
Einstein equations of general relativity
\bel{eq:Eins}
G_{\alpha\beta} = T_{\alpha\beta},
\ee
it is nonetheless motivated by this system and we expect a suitable extension of 
our method to apply to \eqref{eq:Eins}. 
The left-hand side $G_{\alpha\beta}$
of \eqref{eq:Eins} is the Einstein tensor of a spacetime $(M, g)$, that is, 
a Lorentzian $(3+1)$-dimensional manifold,
while the right-hand side  $T_{\alpha\beta}$ denotes the energy momentum tensor of a matter field, which in our context can be assumed to be a set of massless and massive scalar fields.


\chapter[The hyperboloidal foliation and the bootstrap strategy]{The hyperboloidal foliation and the bootstrap strategy \label{cha:2}}

\section{The hyperboloidal foliation and the Lorentz boosts}
\label{sec:21}

We will work with the foliation of the interior of the light cone in Minkowski spacetime $\RR^{3+1}$, defined as below.

\begin{figure}
  \centering
  \includegraphics[width=1.0\textwidth]{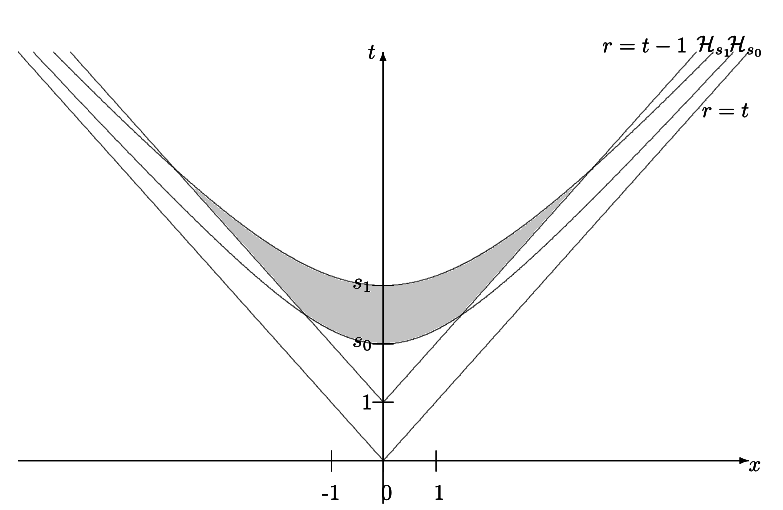}\\
  \caption{The set $\Kcal_{[s_0,s_1]}$.}\label{fig K-2-3}
\end{figure}

We introduce the {\bf hyperboloidal hypersurfaces}
\be
\Hcal_s := \big\{ (t,x) \, \big/ \, t>0; \,  t^2 - |x|^2 =s^2 \big\}
\ee
with hyperbolic radius $s>0$, where $(t,x)=(t, x^a)= (t, x^1, x^2, x^3)$ denote  Cartesian coordinates, and we write $r^2 := |x|^2 = \sum_a (x^a)^2$.
We then consider the interior of the (future) light-cone
\be
\Kcal := \big\{ (t,x) \, / \, |x|< t-1 \big\}
\ee
and, with $s_1 > s_0> 1$, 
as well as the truncated conical region
\be
\aligned
\Kcal_{[s_0,s_1]}
:=& \big\{(t,x) \, / \,  |x| < t-1,\quad (s_0)^2 \leq t^2- |x|^2\leq (s_1)^2,{t>0} \big\}
\\
=&\bigcup_{s_0 \leq s \leq s_1} (\Hcal_s \cap \Kcal).
\endaligned
\ee
This set is thus limited by two hyperboloids and is naturally foliated by hyperboloids. See Fig.~\ref{fig K-2-3} for a display of the set $\Kcal_{[s_0,s_1]}$.

\begin{figure}
  \centering
  \includegraphics[width=1.0\textwidth]{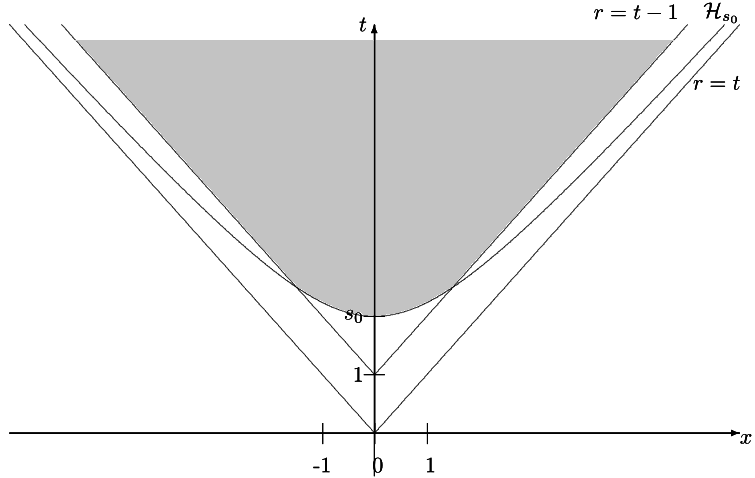}\\
  \caption{The set $\Kcal_{[s_0,+\infty)}$.}\label{fig K-2-infty}
\end{figure}

Taking now $s_1=+\infty$, we will use the notation
\bel{eq:cone}
\aligned
\Kcal_{[s_0,+\infty)}
&:= \big\{(t,x) \, \big/ \, |x|< t-1,\,(s_0)^2\leq t^2 - |x|^2 \big\}
\\
&= \bigcup_{s \geq s_0} (\Hcal_s\cap \Kcal).
\endaligned
\ee
We refer to Fig.~\ref{fig K-2-infty} for a display of this set. 

In the following, we will be interested in functions supported in the conical region
 $\Kcal_{[s_0,+\infty)}$.
Observe that the set $\Kcal_{[s_0,+\infty)}$ is neither closed nor open
and
 a function supported in this set, by definition, vanishes near the future 
light cone $\{r = t-1\}$ but 
can be non-vanishing on the surface $\Hcal_{s_0}$. 

The region $\Kcal \cap \{|x|\leq t/2\}$ will be also of interest in our estimates below, when we will investigate the behavior of solutions away from the light cone. Note in passing that the uniform estimate $t\leq \frac{2 \sqrt{3}}{3}s$ holds in $\Kcal \cap \{|x|\leq t/2\}$. 


We consider first the Klein-Gordon equation
\bel{intro eq linear}
\aligned
&\Box u + \sigma^2 \, u = f,
\\
&u(t,x)|_{\Hcal_{s_0}} = u_0(t,x),\quad u_t(t,x)|_{\Hcal_{s_0}} = u_1(t,x),
\endaligned
\ee
with given $\sigma, s_0=B+1>1$. Recall that the symbol $\Box$ denotes the wave operator in Minkowski spacetime whose metric has the signature $(1,-1,-1,-1)$.
In \eqref{intro eq linear}, the initial data $u_0, u_1$ are prescribed and compactly supported in the ball
$\big\{ |x| \leq B \big\}$ of radius $B= s_0-1$.
The source-term function $f$ is supported in $\Kcal_{[s_0,+\infty)}$,
so that, by the principle of propagation at finite speed, the solution $u=u(t,x)$ to \eqref{intro eq linear} is also
supported in $\Kcal_{[s_0,+\infty)}$. 
In the same manner, the solution of the main system \eqref{main eq main} studied in this monograph
is also supported in $\Kcal_{[s_0,+\infty)}$ and vanishes 
in a neighborhood of the light cone
$\big\{ r=t-1 \big\}$.

Next, let us introduce the {\bf hyperbolic rotations} or {\bf Lorentz boosts} (by rising and lowering the indices with the Minkowski metric with signature $(+,-,-,-)$)
\be
\aligned
&L_a := - x_a\del_0  + x_0\del_a = x^a\del_t + t\del_a,
\endaligned
\ee
which are tangent vectors to the hyperboloids. 
Denote by $\Zscr$ the family of {\bf admissible vector fields} consisting
of all vectors
\bel{vectorsZ}
Z_\alpha  := \del_\alpha,
\qquad
Z_{3+a} := L_a.
\ee
Observe that, for any $Z,\,Z'\in \Zscr$, the Lie bracket
$[Z,Z']$ also belongs to $\Zscr$, so that this set is a Lie algebra. For any multi-index $I = (\alpha_1,\alpha_2,\ldots, \alpha_m)$ of length $|I|:=m$, we denote by $Z^I$ the $m$-th order differential operator
$Z^I := Z_{\alpha_1} \ldots Z_{\alpha_m}$. We also denote by $\del^I$ the $m$-th order derivative operator $\del^I := \del_{\alpha_1}\del_{\alpha_2}\ldots\del_{\alpha_m}$ (here $0\leq \alpha_i\leq 3$) and $L^I$ the $m$-th order derivative operator $L^I := L_{\alpha_1}L_{\alpha_2}\dots L_{\alpha_m}$ (here $4\leq \alpha_j\leq 6$).

Since we will be working within $\Kcal$, we have $|x^a/t| \leq 1$ in $\Kcal$ and, therefore, the
{\bf spatial rotations}
\bel{eq:srot}
\Omega_{ab} := x^a\del_b - x^b\del_a
\ee
need not be included explicitly in our analysis,
since these fields can be recovered from $\Zscr$ via the identities 
\be
\Omega_{ab} = {x^a \over t} L_b -  {x^b \over t} L_a.
\ee
Within the cone $\Kcal$, the coefficients $x^a/t$ are smooth, bounded, and homogeneous of degree zero (see also Lemma \ref{pre lem homo}).

Now we study the energy associated with the hyperboloidal foliation. Using $\del_t u$ as multiplier for the equation \eqref{intro eq linear},
it is easy to derive the following  {\bf energy inequality} for all $s_1 \geq s_0$:
\be
\big(E_{m,\sigma}(s_1,u) \big)^{1/2} \leq \big( E_{m,\sigma}(s_0,u) \big)^{1/2} + \int_{s_0}^{s_1} \bigg(\int_{\Hcal_s}f^2 \, dx\bigg)^{1/2} ds,
\ee
where
the {\bf energy on the hyperboloids} is 
 defined as\footnote{The subscript refers to the Minkowski metric.}
\bel{exp14}
E_{m,\sigma}(s_1,u) := \int_{\Hcal_{s_1}}  \Big( \sum_{a=1}^{3} \big((x^a/t)\del_t u + \del_a u\big)^2 + ((s_1/t)\del_t u)^2 + \sigma^2 u^2 \Big) \, dx
\ee
with $dx=dx^1 dx^2 dx^3$. 
 \citet{Hormander97}) established a Sobolev-type estimate adapted to this inequality (cf.~Lemma 7.6.1 therein, or refer to \eqref{pre ineq sobolev} in Chapter \ref{cha:5}, below) and arrived at the $L^\infty$ estimate
\bel{eq:ZI}
\aligned
&\sup_{\Hcal_{s_1}} t^{3/2}|u|
\\
& \leq C\sum_{|I|\leq 2} E_{m,\sigma}(s_1,Z^Iu)^{1/2}
\\
& \leq C\sum_{|I|\leq 2}E_{m,\sigma}(s_0, Z^I u)^{1/2} + C\sum_{|I|\leq 2}\int_{s_0}^{s_1} \bigg(\int_{\Hcal_s} (Z^I f)^2 \, dx\bigg)^{1/2} ds, 
\endaligned
\ee 
where the summations are over all admissible vector fields in $\Zscr$. 
Importantly, this argument yields the (optimal) rate of decay $t^{-3/2}$
enjoyed by solutions to the Klein-Gordon equation.

\section{Semi-hyperboloidal frame}
\label{sec:22}

Our analysis in the present work is based on the {\bf semi-hyperboloidal frame} ---as we call it--- defined by\footnote{In contrast, a 
standard method relies on the so-called null frame containing two null vectors tangent to the light cone.}
\be
\delu_0 := \del_t, \quad \qquad
\delu_a := t^{-1}L_a = \big(x^a/t\big)\del_t + \del_a. 
\ee
The transition matrices between the semi-hyperboloidal frame and the natural frame $\del_\beta$, that is, 
$\delu_\alpha  = \Phi_\alpha ^{\betar}\del_{\betar}$
and $\del_\alpha  = \Psi_\alpha ^{\betar}\delu_{\betar}$
are found to be 
\bel{pre Phi}
\Phi:=
\begin{pmatrix}
&1        &0         &0         &0        \\
&x^1/t    &1         &0         &0        \\
&x^2/t    &0         &1         &0        \\
&x^3/t    &0         &0         &1        \\
\end{pmatrix},
\qquad
\Psi:= \Phi^{-1} =
\begin{pmatrix}
&1        &0         &0         &0        \\
&-x^1/t   &1         &0         &0        \\
&-x^2/t   &0         &1         &0        \\
&-x^3/t   &0         &0         &1        \\
\end{pmatrix}.
\ee
With our choice of frame, the matrices $\Phi, \Psi$ are smooth within the cone $\Kcal$.

We adopt the following notation and convention.  In order to express the components of a tensor in a frame, we always use Roman font with upper and lower indices for its components in the natural frame,
while we use {\sl underlined} Roman font for its components in the semi-hyperboloidal frame. Hence, a tensor is expressed as
$T = T^{\alpha\beta}\del_\alpha \otimes \del_\beta$
in natural frame, and as 
$T = \Tu^{\alpha\beta}\delu_\alpha \otimes \delu_\beta$
in the semi-hyperboloidal frame.  For example, the Minkowski metric is expressed in the semi-hyperboloidal frame as 
$$
m = \underline{m}^{\alpha\beta}\delu_\alpha \otimes \delu_\beta
$$
with 
\be
\big(\underline{m}^{\alpha\beta}\big)
 =
\left(
\begin{array}{cccc}
s^2/t^2 &x^1/t &x^2/t &x^3/t
\\
x^1/t &-1 &0 &0
\\
x^2/t &0 &-1 &0
\\
x^3/t &0 &0 &-1
\end{array}
\right)
\ee
and 
\be
\big(\underline{m}_{\alpha\beta}\big)
 =
\left(
\begin{array}{cccc}
1 &x^1/t &x^2/t &x^3/t
\\
x^1/t &(x^1/t)^2-1 &x^1x^2/t^2 &x^1x^3/t^2
\\
x^2/t &x^2x^1/t^2 &(x^2/t)^2-1 &x^2x^3/t^2
\\
x^3/t &x^3x^1/t^2 &x^3x^2/t^2 &(x^3/t)^2-1
\end{array}
\right).
\ee

Similarly, a second-order differential operator $T^{\alphar\betar}\del_{\alphar}\del_{\betar}$ can be written in the semi-hyperboloidal frame, so that
by writing 
$T^{\alphar\betar} = \Tu^{\alphar'\betar'} \Phi_{\betar'}^{\betar} \Phi_{\alphar'}^{\alphar}$ (with obvious notation), we obtain the following
 decomposition formula for any function $u$  
\bel{pre frame change of frame}
T^{\alphar\betar}\del_{\alphar}\del_{\betar}u
= \Tu^{\alphar\betar}\delu_{\alphar}\delu_{\betar} u
 +T^{\alpha\beta} (\del_\alpha \Psi_\beta^{\beta'})\delu_{\beta'}u.
\ee
In particular, for the wave operator, we obtain
\be
\Box u = \mbar^{\alphar\betar}\delu_{\alphar}\delu_{\betar} u
+ m^{\alphar\betar}\big(\del_\alpha \Psi_{\betar'}^{\betar}\big)\delu_{\betar'}u,
\ee
where  
we know that $\mbar^{00} = (t^2 - r^2)/t^2 = s^2/t^2$.
This leads us immediately to the following key identity.

\begin{proposition}[Semi-hyperboloidal decomposition of $\Box$]
\label{decompW}
In the semi-hyperboloidal frame, the wave operator admits the decomposition
\bel{pre expression of wave under oneframe}
\aligned
(s/t)^2\delu_0\delu_0 u
= \,
\Box u 
& - \mbar^{0\ar}\delu_0 \delu_{\ar} u - \mbar^{\ar 0}\delu_{\ar} \delu_0 u -\mbar^{\ar\br}\delu_{\ar} \delu_{\br}u
\\
& - m^{\alphar\betar}\big(\del_{\alphar}\Psi_{\betar'}^{\betar}\big)\delu_{\betar'}u.
\endaligned
\ee 
\end{proposition}

We will also need the following property.  

\begin{proposition}\label{pre lem frame}
For any  two-tensor $T$ defined in the cone $\Kcal$ and for all indices $\alpha, \beta, I$ and admissible field $Z$, one has
$$
\big|Z^I \Tu^{\alpha\beta}\big| \lesssim \sum_{\alpha',\beta'\atop |I'|\leq |I|} \big|Z^{I'}T^{\alpha'\beta'}\big|
\qquad \text{ in } \Kcal.
$$
\end{proposition}

The proof of this result (given below) will rely on the following two lemmas.

\begin{lemma}[Homogeneity lemma]
\label{pre lem homo}
Let $f=f(t,x)$ be a smooth function defined in the
closed region $\{t\geq 1,|x|\leq t\}$
and assumed to be homogeneous of degree $\eta$ in the sense that:
\be
f(pt,px) = p^{\eta}f(t,x) \quad \text{ for all }  
 p \geq 1/t.
\ee
For all multi-indices $I_1, I_2$, the following estimate holds for some positive constant $C(n,|I_1|,|I_2|,f)$:
\bel{pre homogeneity estimates}
\aligned
&\big|\del^{I_1}Z^{I_2}f(t,x)\big|\leq  C(n,|I_1|,|I_2|,f) \, t^{-|I_1|+\eta}
\qquad \text{ in } \Kcal = \big\{|x| < t-1 \big\}.
\endaligned
\ee
\end{lemma}

\begin{proof}
We observe that $\del_a f$ is homogeneous of degree $\eta -1$ and $L_a f$ is homogeneous of degree $\eta$. We also observe that if some functions $f_i$ are homogeneous of degree $\eta_i$, then the product $\prod_{i} f_i$ is homogeneous of degree $\sum_{i}\eta_i$ and then any $Z^I f$ is homogeneous.

We claim that the degree of homogeneity of $Z^If$, denoted by $\eta'$, 
is
not higher than $\eta$. This can be checked by induction, as follows. Namely, this is clear when $|I|=1$. Moreover, assume that for all $|I|\leq m$, $Z^If$ is homogeneous of degree $\eta'\leq \eta$, then we now check the same property for all $|I|=m+1$. Namely, assume that $Z^I = Z_1 \,  Z^{I'}$ with $|I'|=m$.  When $Z_1=\del_{\alpha}$, then $Z^If = \del_{\alpha}\big(Z^{I'}f\big)$ and we observe that $Z^{I'}f$ is homogeneous of degree $\eta'\leq \eta$. 
We see that $\del_{\alpha}\big(Z^{I'}f\big)$ is again homogeneous of degree $\eta'-1<\eta'\leq\eta$. When $Z_1=L_a$, then $L_a\big(Z^{I'}f\big)$ is homogeneous of degree $\eta'\leq \eta$. This completes the induction argument that $Z^If$ is homogeneous of degree $\eta$ at most.

Furthermore, observe that if $f$ is homogeneous of degree $\eta$, then $\del^{I_1}f$ is homogeneous of degree $\eta-|I_1|$. This is so since  $\del_{\alpha}f$ is homogeneous of degree $\eta-1$. 
Hence, 
$\del^{I_1}Z^{I_2}f$ is homogeneous of degree $\eta-|I_1|$ and is smooth within the cone $\Kcal$.

Next, in order to establish \eqref{pre homogeneity estimates}, we take $(t,x)\in \Kcal$
and compute 
\bel{e:300}
\del^{I_1}Z^{I_2}f(t,x) = t^{-|I_1|+\eta'}\del^{I_1}Z^{I_2}f(1,x/t)
\qquad \text{ with } \eta' \leq \eta.
\ee
By a continuity argument in the compact set $\{t=1,|x|\leq 1\}$, there exists a positive constant $C(n,|I_1|,|I_2|)$ such that
$$
\sup_{|x|\leq 1}\big|\del^{I_1}Z^{I_2}f(1,x)\big|\leq C(n,|I_1|,|I_2|), 
$$
so that, in view of \eqref{e:300}, the desired result is proven.
\end{proof}

We now estimate the coefficients of the matrices $\Phi$ and $\Psi$. First, we observe that $x^a/t$ is homogeneous of degree $0$ and smooth in the closed region $\{t\geq 1, |x|\leq t\}$.

\begin{lemma}[Changes of frame]
With the notation above, the following two estimates hold for all multi-indices $I_1, I_2$:
\bel{pre lem commutator pr5}
\big|\del^{I_1}Z^{I_2}\Phi_{\alpha}^{\beta}\big|\leq C(n,|I_1|,|I_2|) \, t^{-|I_1|}
\quad \text{ in } \Kcal.
\ee
\bel{pre lem commutator pr6}
\big|\del^{I_1}Z^{I_2}\Psi_{\alpha}^{\beta}\big|\leq C(n,|I_1|,|I_2|) \, t^{-|I_1|}
\quad \text{ in } \Kcal.
\ee
\end{lemma}

\begin{proof}[Proof of Proposition \ref{pre lem frame}]
From 
$\Tu^{\alpha\beta} = T^{\alpha'\beta'}\Psi_{\alpha'}^{\alpha}\Psi_{\beta'}^{\alpha}$
we find 
$$
Z^I\Tu^{\alpha\beta}
= Z^I\big(T^{\alpha'\beta'}\Psi_{\alpha'}^{\alpha}\Psi_{\beta'}^{\alpha}\big)
= \sum_{I_1+I_2=I}Z^{I_1}\big(\Psi_{\alpha'}^{\alpha}\Psi_{\beta'}^{\alpha}\big)Z^{I_2}T^{\alpha'\beta'}
$$
and, therefore, 
$$
\aligned
\big|Z^I\Tu^{\alpha\beta}\big|
\leq& \sum_{I_1+I_2=I}\big|Z^{I_1}\big(\Psi_{\alpha'}^{\alpha}\Psi_{\beta'}^{\alpha}\big)\big|\,\big|Z^{I_2}T^{\alpha'\beta'}\big|
\\
\leq& C(|I|)\sum_{|I_2|\leq |I|}\big|Z^{I_2}T^{\alpha'\beta'}\big|. 
\endaligned
$$
Here, we have used that 
$\big|Z^{I_1}\big(\Psi_{\alpha'}^{\alpha}\Psi_{\beta'}^{\alpha}\big)\big|\leq C(|I_1|)$
which follows from \eqref{pre lem commutator pr6}.
\end{proof}



\section{Energy estimate for the hyperboloidal foliation}
\label{sec:23}

As presented in Chapter \ref{cha:1}, we are interested in the following class of wave-Klein-Gordon type systems
\bel{pre eq energy}
\aligned
&\Box w_i + G_i^{\jr\alphar\betar}\del_{\alphar}\del_{\betar} w_{\jr} + c_i^2w_i = F_i,
\\
&w_i|_{\Hcal_{s_0}} = {w_i}_0,\quad \del_t w_i|_{\Hcal_{s_0}} = {w_i}_1,
\endaligned
\ee
with unknowns $w_i$ ($1 \leq i \leq n_0$), where
the metric $G_i^{j\alpha\beta}$ (with some abuse of notation) 
and the source-terms $F_i$ are  
supported in the cone $\Kcal_{[s_0,+\infty)}$,
and ${w_i}_0$, ${w_i}_1$ are supported on the initial hypersurface $\Hcal_{s_0}\cap \Kcal$.
To guarantee the hyperbolicity, we assume the symmetry conditions \eqref{pre condition symmetry}.
For definiteness, we may also assume \eqref{eq:lesc} on the constants $c_i$.

We introduce the following energy associated with the Minkowski metric on each hyperboloid $\Hcal_s$:
\bel{pre pre expression of energy}
\aligned
&E_{m,c_i}(s,w_i)
\\
&:= \int_{\Hcal_s}\bigg(
(\del_t w_i)^2 + \sum_a (\del_a w_i)^2 + (2x^{\ar}/t)\del_tw_i \del_{\ar} w_i + c_i^2 w_i^2 \bigg) \, dx
\\
&=\int_{\Hcal_s}\bigg( 
\sum_a (\delu_a w_i)^2 + \big((s/t)\del_t w_i\big)^2 + c_i^2 w_i^2 \bigg)\, dx
\\
&= \int_{\Hcal_s}\bigg( 
\sum_a \big((s/t) \del_a w_i\big)^2 + t^{-2}(Sw_i)^2
+ t^{-2}\sum_{a < b}\big(\Omega_{ab}w_i\big)^2 
+ c_i^2 w_i^2 
\bigg)\, dx,
\endaligned
\ee
where we use the notation \eqref{eq:srot} and
\be
\aligned
&t^{-1}S := \del_t + \sum_a(x^a/t)\del_a.
\endaligned
\ee
When $c_i=0$, we may also write $E_m(s,w_i): = E_{m,0}(s,w_i)$ for short.

On the other hand, the curved energy which is naturally
 associated with the principal part of \eqref{pre eq energy} is defined as
\bel{pre expression of curved energy}
\aligned
E_{G,c_i}(s,w_i):= E_{m,c_i}(s,w_i)
& + 2\int_{\Hcal_s} \big(\del_t w_i \del_{\betar}w_{\jr} G_i^{\jr\alpha\betar}\big)_{0\leq \alpha\leq 3} \cdot (1,-x/t) dx
\\
& - \int_{\Hcal_s} \big(\del_{\alphar}w_i\del_{\betar}w_{\jr} G_i^{\jr\alphar\betar}\big) \, dx,
\endaligned
\ee
where the second term in the right-hand side involves the Euclidian inner product of the vectors
$\big(\del_t w_i \del_{\betar}w_{\jr} G_i^{\jr\alpha\betar}\big)_{0\leq \alpha\leq 3}$ and $(1,-x/t)$.

\begin{proposition}[Hyperboloidal energy estimate]
\label{pre lem energy}
Given a constant $\Kcoef>1$ and some locally integrable functions $L, M \geq 0$, the following property holds. 
Let $(w_i)_{1 \leq i \leq n_0}$ be the (local-in-time) solution to \eqref{pre eq energy} defined on some hyperbolic time interval $[s_0, s_1]$ and
suppose that the metric and source satisfy the following coercivity conditions 
\bel{pre lem energy curved energy is big}
\Kcoef^{-2}\sum_i E_{m,c_i}(s,w_i)\leq \sum_i E_{G,c_i}(s,w_i) \leq \Kcoef^2 \sum_i E_{m,c_i}(s,w_i),
\ee
\bel{pre lem energy curveterm is small}
\aligned
&
\bigg|\int_{\Hcal_s}\frac{s}{t}\bigg( \del_{\alphar}G_i^{\jr\alphar\betar}\del_tw_i \del_{\betar}w_{\jr}
- \frac{1}{2}\del_tG_i^{\jr\alphar\betar}\del_{\alphar}w_i\del_{\betar}w_{\jr}\bigg)  \,dx\bigg|
\\
& \leq M(s) E_{m,c_i}(s,w_i)^{1/2},
\endaligned
\ee
and the bound
\bel{pre lem energy source}
\sum_i\|F_i\|_{L^2(\Hcal_s)} \leq L(s), 
\qquad s\in [s_0,s_1].
\ee
Recall also that the symmetry conditions \eqref{pre condition symmetry} are assumed, that is,
$G_i^{j\alpha\beta} = G_j^{i\alpha\beta} = G_i^{j\beta\alpha}$, 
and that $F_i$, $G_i^{j\alpha\beta}$ are supported in $\Kcal_{[s_0,s_1]}$.
The following energy estimate holds (for all $s \in [s_0, s_1]$):
\bel{eq:500}
\aligned
&\bigg(\sum_i E_{m,c_i}(s,w_i) \bigg)^{1/2}
\\
&\leq \Kcoef^2 \bigg(\sum_i E_{m,c_i}(s_0,w_i)\bigg)^{1/2}
+ C\Kcoef^2\ \int_{s_0}^s \big(L(\tau) + \, M(\tau)\big) \, d\tau,
\endaligned
\ee 
where the constant $C>0$ depends on the the structure of the system \eqref{pre eq energy} only. 
\end{proposition}

The following remarks are in order:
\begin{itemize}

\item In view of \eqref{pre pre expression of energy} and \eqref{eq:500}, the $L^2$ norm of $\delu_a w_i$ and $(s/t)\del_\alpha w_i$ are uniformly controlled on each hypersurface $\Hcal_s$.
It is expected that these weighted expressions
 enjoy better decay than $\del_\alpha w_i$ itself.

\item In view of 
$\sum_a (x^a/r) \big((r/t)\del_a + (x^a/r) \del_t\big)w_i = t^{-1}Sw_i$, 
we see that the weighted scaling derivative 
$t^{-1}Sw_i$ is also controlled.  

\item In Chapter~\ref{cha:10}, we will see that the energy $E_m(s_0,w_i)$
on the initial hyperboloid  is
controlled by the $\Hf^1$ norm of the initial data on the initial slice $t=s_0$.
\end{itemize}

\begin{proof} 
In view of the symmetry \eqref{pre condition symmetry} and by using $\del_t w_i$ as a multiplier,
we easily derive the energy identity
$$
\aligned
& \sum_i  \frac{1}{2} \del_t \bigg( \sum_\alpha(\del_\alpha  w_i)^2 + c_i^2 w_i^2 \bigg) - \sum_i  \sum_a\del_a\big(\del_a w_i \del_t w_i\big)
\\
& 
 + \sum_i  \del_{\alphar}\big(G_i^{\jr\alphar\betar}\del_t w_i\del_{\betar}w_{\jr}\big)
- \sum_i  \frac{1}{2}\del_t\big(G_i^{\jr\alphar\betar}\del_{\alphar}w_i\del_{\betar}w_{\jr}\big)
\\
&= \sum_i \del_t w_i F_i
+ \sum_i\bigg(\del_{\alphar}G_i^{\jr\alphar\betar} \del_tw_i \del_{\betar}w_{\jr}
- \frac{1}{2}\del_tG_i^{\jr\alphar\betar}\del_{\alphar}w_i\del_{\betar}w_{\jr}\bigg).
\endaligned
$$
We integrate this identity over the region $\Kcal_{[s_0,s]}$ and use Stokes' formula. Note that
by the property of propagation at finite speed, the solution $(w_i)$ is defined in $\Kcal_{[s_0,+\infty)}$
and vanishing in a neighborhood of the light cone. 
By our assumption on the support of
 $G_i^{j\alpha\beta}$, we obtain 
$$
\aligned
&\frac{1}{2}\sum_i \big(E_{G,c_i}(s,w_i) - E_{G,c_i}(s_0,w_i)\big)
\\
&=\sum_i
\int_{\Kcal_{[s_0,s]}} \Big(
\del_t w_i F_i
+ \del_{\alphar}G_i^{\jr\alphar\betar} \del_tw_i \del_{\betar}w_{\jr} -
\frac{1}{2}\del_tG_i^{\jr\alphar\betar}\del_{\alphar}w_i\del_{\betar}w_{\jr} \Big) \, dtdx
\\
&=\sum_i
\int_{s_0}^s \bigg(\int_{\Hcal_\tau} (\tau/t) \del_t w_i F_i \, dx\bigg) d\tau
\\
& \quad 
+ \sum_i\int_{s_0}^s \bigg(\int_{\Hcal_\tau} (\tau/t)
\Big( \del_{\alphar}G_i^{\jr\alphar\betar} \del_tw_i \del_{\betar}w_{\jr} -
\frac{1}{2}\del_tG_i^{\jr\alphar\betar}\del_{\alphar}w_i\del_{\betar}w_{\jr}
\Big) \, dx\bigg) d\tau,
\endaligned
$$
which leads us to
$$
\aligned
\frac{d}{ds}\sum_i E_{G,c_i}(s,w_i)
 = 
&2\sum_i\int_{\Hcal_s}
\Big(
(s/t)\del_{\alphar}G_i^{\jr\alphar\betar} \del_tw_i \del_{\betar}w_{\jr}
\\
&\qquad\qquad
-(s/2t)\del_tG_i^{\jr\alphar\betar}\del_{\alphar}w_i\del_{\betar}w_{\jr}
+ (s/t) \del_t w_i F_i \Big) \, dx.
\endaligned
$$
So, with the assumptions \eqref{pre lem energy curveterm is small} and \eqref{pre lem energy source}, we get 
$$
\aligned
&\bigg(\sum_i E_{G,c_i}(s,w_i)\bigg)^{1/2}\frac{d}{ds}\bigg(\sum_i E_{G,c_i}(s,w_i)\bigg)^{1/2}
\\
&\leq \sum_i
\bigg(
M(s) + \| F_i \|_{L^2(\Hcal_s)}  \bigg) E_{m,c_i}(s,w_i)^{1/2}
\\
&\leq C M(s) \Big(  \sum_i  E_{m,c_i}(s,w_i) \Big)^{1/2} 
      + C  \, \Big( \sum_i \| F_i \|_{L^2(\Hcal_s)}^2 \Big)^{1/2}  \Big( \sum_i E_{m,c_i}(s,w_i) \Big)^{1/2}
\\
&= C \, \Big( 
M(s)  + L(s) \Big) \Big( \sum_i  E_{m,c_i}(s,w_i) \Big)^{1/2}  
\\
&\leq C\kappa_1 \big( M(s)+ L(s) \Big) \, \Big( \sum_i E_{G,c_i}(s,w_i)\Big)^{1/2},
\endaligned
$$
which yields
$$
\quad\frac{d}{ds}\bigg(\sum_i E_{G,c_i}(s,w_i)\bigg)^{1/2}  \leq C\Kcoef \big(L(s) + M(s)\big).
$$
We integrate this inequality over the interval $[s_0,s]$ and obtain
$$
\aligned
&
\bigg(\sum_{i=1}^{n_0} E_{G,c_i}(s,w_i)\bigg)^{1/2}
\\
&\leq
 \bigg(\sum_{i=1}^{n_0} E_{G,c_i}(s_0,w_i)\bigg)^{1/2} + C\Kcoef\int_{s_0}^s L(\tau)d\tau + C\Kcoef \int_{s_0}^sM(\tau)d\tau.
\endaligned
$$
By using the condition \eqref{pre lem energy curved energy is big},
this completes the proof.
\end{proof}


\section{The bootstrap strategy}
\label{sec:24}

We are now in a position to outline our method of proof of Theorem~\ref{main thm main}. From now on, we assume that the assumptions therein are satisfied. 
It will be necessary to distinguish between three levels of regularity and, in order to describe this {\sl scale of regularity,}
 we will use the following convention on the indices in use:
\be\label{index-I}
\aligned
&\text{$\Is$:  multi-index of order $\leq 5$,} 
\\
&\text{$I^{\dag}$: multi-index of order $\leq 4$,} 
\\
&\text{$I$: multi-index of order $\leq 3$,} 
\endaligned
\ee
which we call {\bf admissible indices.}
Throughout, $C,C_0,C_1,C^*, \ldots$ are constants depending only on the structure of the system \eqref{main eq main}, such as $j_0,\,k_0,\,B, c_i$.

We will use certain norms on the hyperboloids and, so,
 if $u$ is a function supported in $\Kcal_{[s_0,+\infty)}$, we set 
(for $s\geq s_0$): 
\bel{eq:noteNorm}
\|u\|_{L^p(\Hcal_s)} := \bigg(\int_{\Hcal_s}\big|u(t,x)\big|^pdx\bigg)^{1/p} = \bigg(\int_{\mathbb{R}^3}\big|u\big(\sqrt{s^2+|x|^2},x\big)\big|^pdx\bigg)^{1/p}.
\ee
The proof of Theorem~\ref{main thm main} relies on 
Propositions \ref{main prop1}--\ref{main prop2}, stated now.
The first proposition below concerns the construction of initial data on the initial hyperboloid $\Hcal_{B+1}$, and will be established in Chapter~\ref{cha:11}.

\begin{proposition}[Initialization of the argument]
\label{main prop1}
For any sufficiently large constant $C_0>0$, there exists a positive constant $\eps'_0 \in (0,1)$ depending only on $B$ and $C_0$ such that, for every initial data satisfying
\bel{eq:datainit}
\sum_i  \big(\|{w_i}_0\|_{H^6(\RR^3)}+ \| {w_i}_1 \|_{H^5(\RR^3)}\big) \leq \eps \leq \eps_0',
\ee
the local-in-time solution to \eqref{main eq main} associated with this
 initial data extends to the region limited by the constant time hypersurface $t=B+1$
and the hyperboloid $\Hcal_{B+1}$. 
Furthermore, it satisfies the uniform bound
\bel{eq:inegEm}
\sum_j
E_m(s_0,Z^{\Is}w_j)^{1/2}\leq C_0\eps.
\ee
for all admissible vector fields $Z$ and all admissible indices $\Is$.
\end{proposition}

Given some constants $C_1, \eps>0$ and $\delta\in (0, 1/6)$ and a hyperbolic time interval $[s_0, s_1]$, we  call {\bf hierarchy of energy bounds} with parameters $(C_1,\eps, \delta)$ the following five inequalities (for all $s\in [s_0,s_1]$ and all admissible fields and admissible indices):
\begin{subequations}\label{proof energy assumption}
\bel{proof energy assumption b}
E_m(s,Z^{\Is} u_{\ih})^{1/2} \leq C_1\eps s^{\delta} \quad \text{for } 1\leq \ih\leq j_0,
\ee
\bel{proof energy assumption a}
E_{m,\sigma}(s,Z^{\Is} v_{\jc})^{1/2} \leq C_1\eps s^{\delta} \quad \text{ for } j_0+1\leq \jc\leq n_0,
\ee
\bel{proof energy assumption d}
E_m(s,Z^{\Id} u_{\ih})^{1/2} \leq C_1\eps s^{\delta/2} \quad \text{for } 1\leq \ih\leq j_0,
\ee
\bel{proof energy assumption c}
E_{m,\sigma}(s,Z^{\Id} v_{\jc})^{1/2} \leq C_1\eps s^{\delta/2} \quad \text{ for } j_0+1\leq \jc\leq n_0,
\ee
\bel{proof energy assumption e}
E_m(s,Z^I u_{\ih})^{1/2} \leq C_1\eps \quad \text{for } 1\leq \ih\leq j_0.
\ee 
\end{subequations}
Observe that \eqref{proof energy assumption e} concerns the wave components only and that the upper bound is independent of time. 

At this juncture, since
$c_{\ic}\geq \sigma>0$ for all $j_0+1\leq \ic\leq n_0$, we have 
\begin{equation}\label{equiv c_i-sigma}
E_{m,\sigma}(s,Z^{\Is} v_{\ic})\leq E_{m,c_{\ic}}(s,Z^{\Is} v_{\ic}) \leq (c_{\ic}/\sigma)^2E_{m,\sigma}(s,Z^{\Is} v_{\ic}),
\end{equation}
so that the two energy expressions are equivalent.

Given some constants  $\Kcoef>1$, $C^*, C_1, \eps>0$ and $\delta \in (0, 1/6)$ and a hyperbolic time interval $[s_0, s_1]$, we call {\bf hierarchy of metric-source bounds}  with parameters $(\Kcoef, C^*, C_1, \eps,\delta)$. 
 the following three sets of estimates:
\bei
\begin{subequations}
\item For all $|\Is|\leq 5$,
\bel{main bootstrap 1}
\Kcoef^{-2} \sum_i E_{G,c_i}(s,Z^{\Is} w_i) \leq \sum_i E_{m,c_i}(s,Z^{\Is}w_i)\leq \Kcoef^2 \sum_i E_{G,c_i}(s,Z^{\Is} w_i),
\ee
\bel{main bootstrap 2}
\aligned
&\bigg|\int_{\Hcal_s}\frac{s}{t}\bigg( \del_{\alphar}G_i^{\jr\alphar\betar}\del_tZ^{\Is}w_i \del_{\betar}Z^{\Is}w_{\jr}
- \frac{1}{2}\del_tG_i^{\jr\alphar\betar}\del_{\alphar}Z^{\Is}w_i\del_{\betar}Z^{\Is}w_{\jr}\bigg) dx\bigg|
\\
&\leq C^*(C_1\eps)^2s^{-1+\delta}E_{m,c_i}(s,Z^{\Is}w_i)^{1/2} 
\\
& =: M(\Is,s)E_{m,c_i}(s,Z^{\Is}w_i)^{1/2}, 
\endaligned
\ee
and\footnote{with the notation \eqref{eq:noteNorm}}
\bel{main bootstrap 3}
\aligned
&\sum_{i=1}^{n_0} \big{\|}[G_i^{\jr\alphar\betar}\del_{\alphar}\del_{\betar},Z^{\Is}]w_{\jr}\big{\|}_{L^2(\Hcal_s)}
 + \sum_{i=1}^{n_0}\big{\|}Z^{\Is} F_i\big{\|}_{L^2(\Hcal_s)}
\\
&\leq C^*(C_1\eps)^2s^{-1+\delta} =: L(\Is,s).
\endaligned
\ee
\end{subequations}

\begin{subequations}
\item For all $|\Id|\leq 4$,
\bel{main bootstrap 2'}
\aligned
&\bigg|\int_{\Hcal_s}\frac{s}{t}\bigg( \del_{\alphar}G_i^{\jr\alphar\betar}\del_tZ^{\Id}w_i \del_{\betar}Z^{\Id}w_{\jh}
- \frac{1}{2}\del_tG_i^{\jr\alphar\betar}\del_{\alphar}Z^{\Id}w_i\del_{\betar}Z^{\Id}w_{\jr}\bigg) dx\bigg|
\\
&\leq C^*(C_1\eps)^2s^{-1+\delta/2}E_{m,c_i}(s,Z^{\Id}w_i)^{1/2} 
\\
&=: M(\Id,s)E_{m,c_i}(s,Z^{\Id}w_i)^{1/2}
\endaligned
\ee
and
\bel{main bootstrap 3'}
\aligned
&
\sum_{i=1}^{n_0} \big{\|}[G_i^{\jr\alphar\betar}\del_{\alphar}\del_{\betar},Z^{\Id}]w_{\jr}\big{\|}_{L^2(\Hcal_s)}
 + \sum_{i=1}^{n_0}\big{\|}Z^{\Id} F_i\big{\|}_{L^2(\Hcal_s)}
\\
&\leq C^*(C_1\eps)^2s^{-1+\delta/2} =:L(\Id,s).
\endaligned
\ee
\end{subequations}

\begin{subequations}
\label{main bootstrap 4 5}
\item For all $|I|\leq 3$,
\bel{main bootstrap 4}
\aligned
&\bigg|\int_{\Hcal_s}\frac{s}{t}\bigg( \del_{\alphar}G_{\ih}^{\jr\alphar\betar}\del_tZ^Iu_{\ih} \del_{\betar}Z^Iu_{\jh}
- \frac{1}{2}\del_tG_{\ih}^{\jr\alphar\betar}\del_{\alphar}Z^Iu_{\ih}\del_{\betar}Z^Iu_{\jh}\bigg) dx\bigg|
\\
&\leq C^*(C_1\eps)^2s^{-3/2+2\delta}E_m(s,Z^{I}u_{\ih})^{1/2} 
\\
&=: M(I,s)E_m(s,Z^{I}u_{\ih})^{1/2}
\endaligned
\ee
and 
\bel{main bootstrap 5}
\aligned
&\sum_{i=1}^{n_0} \big{\|}[G_i^{\jr\alphar\betar}\del_{\alphar}\del_{\betar},Z^I]w_{\jr}\big{\|}_{L^2(\Hcal_s)}
\\
&+ \sum_{i=1}^{n_0}\big{\|}Z^I F_i\big{\|}_{L^2(\Hcal_s)}
+ \| Z^I\big(G_{\ih}^{\jc\alpha\beta}  \del_\alpha\del_\beta v_{\jc}\big) \|_{L^2(\Hcal_s)}
\\
&\leq C^*(C_1\eps)^2s^{-3/2+2\delta} =:L(I,s).
\endaligned
\ee
\end{subequations}
\eei
Observe that \eqref{main bootstrap 4} concerns the wave components, only.

\begin{proposition}
\label{theclaim} 
Fix some $\delta \in (0, 1/6)$. There exists a positive constant $\eps_0''$ such that for all $\eps,C_1>0$ with $C_1\eps \leq 1$ and $C_1\eps\leq \eps_0''$, there exists a constant $\Kcoef>1$ and a constant $C^*>0$ (both determined by the structure of the system \eqref{main eq main}) such that the hierarchy of energy bounds \eqref{proof energy assumption} with parameters $C_1,\eps, \delta$ implies the hierarchy of metric--source bounds \eqref{main bootstrap 1}---\eqref{main bootstrap 5} with parameters $\Kcoef, C^*,C_1,\eps,\delta$.
\end{proposition}

The proof of this proposition will occupy a major part of this monograph, especially Chapters~\ref{cha:6} to~\ref{cha:9}. Now we admit this proposition and give the proof of the main result, which is going to be essentially based on the following observation.

\begin{proposition}[Enhancing the hierarchy of energy bounds]
\label{main prop2}
Let $\delta \in (0,1/6)$ and let $C_0>0$ be a constant and $C_1>C_0$ be a sufficiently large constant. Then, there exists $\eps_1>0$ such that the following enhancing property holds. 
For any solution $(w_i)$ to \eqref{main eq main} defined
in $\Kcal_{[s_0,s_1]}$ and provided, for $\eps\in [0, \eps_1]$, 
\bei 

\item[(1)]  $\sum_{i=1}^{n_0}E_{m,c_i}(s_0,Z^{\Is}w_i)^{1/2}\leq C_0\eps$ for all admissible $Z$ and $\Is$, 

\item[(2)] the hierarchy 
\eqref{proof energy assumption} holds with parameters $(C_1,\eps,\delta)$
 in the time interval $[s_0,s_1]$, 
\eei
then 
the following {\bf improved energy estimates} also hold  
\begin{subequations}\label{eq:enhanced}
\begin{equation}
\label{eq:enhanced2}
E_m(s,Z^{\Is} u_{\ih})^{1/2} \leq \frac{1}{2} C_1\eps s^{\delta} \quad \text{for } 1\leq \ih\leq j_0,
\end{equation}
\begin{equation}
\label{eq:enhanced1}
E_{m,\sigma}(s,Z^{\Is} v_{\jc})^{1/2} \leq \frac{1}{2}C_1\eps s^{\delta} \quad \text{ for } j_0+1\leq \jc\leq n_0,
\end{equation}
\begin{equation}
\label{eq:enhanced4}
E_m(s,Z^{\Id} u_{\ih})^{1/2} \leq \frac{1}{2} C_1\eps s^{\delta/2} \quad \text{for } 1\leq \ih\leq j_0,
\end{equation}
\begin{equation}
\label{eq:enhanced3}
E_{m,\sigma}(s,Z^{\Id} v_{\jc})^{1/2} \leq \frac{1}{2}C_1\eps s^{\delta/2} \quad \text{ for } j_0+1\leq \jc\leq n_0,
\end{equation}

\begin{equation}
\label{eq:enhanced5}
E_m(s,Z^I u_{\ih})^{1/2} \leq \frac{1}{2} C_1\eps \quad \text{for } 1\leq \ih\leq j_0.
\end{equation}
\end{subequations}
\end{proposition}

\begin{proof} We will use here the conclusion of Proposition~\ref{theclaim}. 
We assume that $\eps_1$ is sufficiently small with 
$\eps_1\leq C_1^{-1}\min\{1, \eps''_0\}$ such 
 that, for all $\eps\leq \eps_1$, the conclusion in Proposition \ref{theclaim} holds.

\vskip.15cm 

\noindent{\bf Step I. High-order energy estimates.} To the equations \eqref{main eq main}, we apply the operator $Z^{\Is}$ (with $|\Is|\leq 5$ throughout) and obtain  
$$
\aligned
& \Box (Z^{\Is}w_i) + G_i^{\jr\alphar\betar}\del_{\alphar}\del_{\betar} (Z^{\Is}w_{\jr}) + c_i^2 Z^{\Is} w_i
= [G_i^{\jr\alphar\betar}\del_{\alphar}\del_{\betar},Z^{\Is}]w_{\jr} + Z^{\Is}F_i.
\endaligned
$$
In view of  Proposition \ref{pre lem energy} and thanks to the conditions
 \eqref{main bootstrap 1}--\eqref{main bootstrap 3} implied by 
 Proposition~\ref{theclaim}, 
we find
$$
\aligned
& \bigg(\sum_i E_{m,c_i}(s,Z^{\Is}w_i)\bigg)^{1/2}
\\
&\leq
 \Kcoef^2 \bigg(\sum_i E_{m,c_i}(s_0,Z^{\Is}w_i)\bigg)^{1/2} +  C\Kcoef^2 \int_{s_0}^s\big(L(\Is,\tau) + M(\Is,\tau)\big)d\tau.
\endaligned
$$
In view of Proposition \ref{theclaim} combined with 
the condition (1) in the proposition, we then have 
$$
\bigg(\sum_i E_{m,c_i}(s,Z^{\Is}w_i)\bigg)^{1/2} \leq \Kcoef^2C_0\eps + C\Kcoef^2 C^*(C_1\eps)^2\int_{s_0}^s \tau^{-1+\delta}d\tau.
$$ 
Choosing  
$\eps_1 \leq \frac{\delta(C_1 - 2\Kcoef^2C_0)}{2CC^*\Kcoef^2C_1^2}$, we obtain
$$
\bigg(\sum_i E_{m,c_i}(s,Z^{\Is}w_i)\bigg)^{1/2}\leq 
{1 \over 2} C_1 \eps s^{\delta},
$$
which leads to the enhanced energy bound
\be
E_{m,\sigma}(s,Z^{\Is}w_i)^{1/2}\leq E_{m,c_i}(s,Z^{\Is}w_i)^{1/2} \leq {1 \over 2}  C_1\eps s^{\delta}.
\ee
This establishes \eqref{eq:enhanced1}-\eqref{eq:enhanced2}.


\vskip.15cm 

\noindent{\bf Step II. Intermediate energy estimates.} We will next rely on the metric--source bounds \eqref{main bootstrap 2'}-\eqref{main bootstrap 3'} (with $|\Id|\leq 4$ throughout) given by the conclusion
of  Proposition~\ref{theclaim}.
To the equation \eqref{main eq main}, we apply the operator $Z^{\Id}$ and obtain
$$
\aligned
& \Box (Z^{\Id}w_i)
 + G_i^{\jr\alphar\betar}\del_{\alphar}\del_{\betar} (Z^{\Id}w_{\jr}) + c_i^2 Z^{\Id} w_i
= [G_i^{\jr\alphar\betar}\del_{\alphar}\del_{\betar},Z^{\Id}]w_{\jr} + Z^{\Id}F_i.
\endaligned
$$
From Proposition \ref{pre lem energy} and thanks to the conditions \eqref{main bootstrap 1}
and \eqref{main bootstrap 2'}-\eqref{main bootstrap 3'}, we find
$$
\aligned
& \bigg(\sum_i E_{m,c_i}(s,Z^{\Id}w_i)\bigg)^{1/2}
\\
&\leq \Kcoef^2\bigg(\sum_i E_{m,c_i}(s_0,Z^{\Id}w_i)\bigg)^{1/2} +  C\Kcoef^2\int_{s_0}^s \big(L(\Id,\tau) + M(\Id,\tau)\big)d\tau.
\endaligned
$$
We thus have 
$$
\bigg(\sum_i E_{m,c_i}(s,Z^{\Id}w_i)\bigg)^{1/2} \leq \Kcoef^2 C_0\eps + C\Kcoef^2 C^*(C_1\eps)^2\int_{s_0}^s \tau^{-1+\delta/2}d\tau.
$$ 
Provided $\eps$ is sufficiently small so that  
$\eps_1 \leq \frac{\delta(C_1-2\Kcoef^2C_0)}{4CC^*\Kcoef^2C_1^2}$, 
we find
$$
\bigg(\sum_i E_{m,c_i}(s,Z^{\Id}w_i)\bigg)^{1/2}\leq {1 \over 2} C_1\eps s^{\delta/2},
$$
which leads to the enhanced energy bound
\be
E_{m,\sigma}(s,Z^{\Id}w_i)^{1/2}\leq E_{m,c_i}(s,Z^{\Id}w_i)^{1/2} \leq {1 \over 2} C_1\eps s^{\delta/2}.
\ee
This establishes \eqref{eq:enhanced3}-\eqref{eq:enhanced4}.


\vskip.15cm 

\noindent{\bf Step III. Low-order energy estimates.} We will now rely on the metric--source bounds \eqref{main bootstrap 4}-\eqref{main bootstrap 5} (with $|I|\leq 3$ throughout). We apply $Z^I$ to the wave equations in \eqref{main eq main} and obtain 
$$
\aligned
& \Box (Z^I u_{\ih}) + G_{\ih}^{\jh\alphar\betar}\del_{\alphar\betar} (Z^Iu_{\jh}) 
= [G_{\ih}^{\jh\alphar\betar}\del_{\alphar}\del_{\betar},Z^I]u_{\jh} - Z^I\big(G_{\ih}^{\jc\alphar\betar}\del_{\alphar}\del_{\betar}v_{\jc}\big) + Z^IF_{\ih}.
\endaligned
$$
In view of Proposition~\ref{pre lem energy} 
and thanks to 
the conditions \eqref{main bootstrap 1} and  \eqref{main bootstrap 4}-\eqref{main bootstrap 5}, we obtain 
$$
\aligned
&\bigg(\sum_{\ih} E_m(s,Z^I u_{\ih} )\bigg)^{1/2}
\\
&\leq \Kcoef^2 \bigg(\sum_{\ih} E_m(s_0,Z^I u_{\ih} )\bigg)^{1/2} + C\Kcoef^2\int_{s_0}^s \big(L(I,\tau) + M(I,\tau)\big) \, d\tau
\endaligned
$$
and, therefore, 
$$
\bigg(\sum_{\ih} E_m(s,Z^I u_{\ih} )\bigg)^{1/2} 
\leq \Kcoef^2 C_0\eps + C\Kcoef^2C^*(C_1\eps)^2\int_{s_0}^s\tau^{-3/2+2\delta}d\tau.
$$
By recalling that $\delta<1/6$, it follows that
$$
\bigg(\sum_{\ih} E_m(s,Z^I u_{\ih} )\bigg)^{1/2}
\leq \Kcoef^2 C_0\eps 
+\frac{CC^*\Kcoef^2C_1^2\eps^2(B+1)^{2\delta-1/2}}{(1/2) - 2\delta}.
$$
Provided
$\eps_1 \leq \frac{(1-4\delta)(C_1 - 2\Kcoef^2C_0)}{4CC^*\Kcoef^2C_1^2(B+1)^{2\delta-1/2}}$,
we obtain the enhanced energy bound
for the wave components
$$
E_m(s,Z^I u_{\ih} )^{1/2} \leq {1 \over 2} C_1\eps
$$
and this establishes \eqref{eq:enhanced5}.

\vskip.15cm 

\noindent{\bf Step IV.} In conclusion, by assuming that  
$$
\eps_1 \leq \min\Big(\frac{\delta(C_1-2\Kcoef^2C_0)}{4CC^*\Kcoef^2C_1^2},\,\frac{(1-4\delta)(C_1 - 2\Kcoef^2C_0)}{4CC^*\Kcoef^2C_1^2(B+1)^{2\delta-1/2}},\,\eps_0'',\, C_1^{-1}\Big),
$$ 
all conditions \eqref{eq:enhanced} are satisfied and the proof is completed.
\end{proof}

Observe that  \eqref{main bootstrap 1}, \eqref{main bootstrap 2}, 
\eqref{main bootstrap 2'}, and \eqref{main bootstrap 4} concern only the ``metric'' $G_i^{j\alpha\beta}$ and depend mainly on $L^\infty$ bounds on the solution and its derivatives: these
 bounds will be established in Chapter ~\ref{cha:8}.
On the other hand, the inequalities \eqref{main bootstrap 3}, \eqref{main bootstrap 3'}, and \eqref{main bootstrap 5} are $L^2$ type estimates and will be derived later in Chapter~\ref{cha:9}.

We complete this chapter with a proof of our main result ---by assuming that Propositions \ref{main prop1} and \ref{theclaim}  
are established.

\begin{proof}[Proof of Theorem \ref{main thm main}]
Let $(w_i)$ be the unique local-in-time solution to \eqref{main eq main} associated with the initial data $({w_i}_0, {w_i}_1)$. 
Given our assumption on the support of $({w_i}_0, {w_i}_1)$
and according to the property of propagation at finite speed, 
this solution $(w_i)$ must be supported in the region $K_{[s_0,+\infty)}$.

As guaranteed by Proposition~\ref{main prop1}, there exist 
constants  $C_0,\eps'_0$ such that, provided
$\|{w_i}_0 \|_{H^6(\RR^3)} + \| {w_i}_1 \|_{H^5(\RR^3)} \leq \eps \leq \eps'_0$,
we have the energy bound $E_G(s_0,Z^{\Is}w_j)^{1/2}\leq C_0\eps$.
 
Let $[s_0,s^*]$ be the largest time interval (containing $s_0$) on which \eqref{proof energy assumption} holds with some parameters $C_1,\eps>0$ 
and 
$\delta \in (0,1/6)$ fixed and for a sufficiently large constant $C_1>C_0$. 
By continuity, we have $s^*>s_0$.

Proceeding by contradiction, let us assume that $s^*< +\infty$, so that at the time $s=s^*$, one (at least) of the inequalities
 \eqref{proof energy assumption} must be an equality. That is, at least one of the following conditions holds:
\bel{proof energy assumption last}
\aligned
&E_{m,\sigma}(s^*,Z^{\Is} v_{\jc})^{1/2} = C_1\eps {s^*}^{\delta} \quad \text{ for } j_0+1\leq \jc\leq n_0,
\\
&E_m(s^*,Z^{\Is} u_{\ih})^{1/2} = C_1\eps {s^*}^{\delta} \quad \text{for }   1\leq \ih\leq j_0,
\\
&E_{m,\sigma}(s^*,Z^{\Id} v_{\jc})^{1/2} = C_1\eps {s^*}^{\delta/2} \quad \text{ for } j_0+1\leq \jc\leq n_0,
\\
&E_m(s^*,Z^{\Id} u_{\ih})^{1/2} = C_1\eps {s^*}^{\delta/2} \quad \text{for } 1\leq \ih\leq j_0,
\\
&E_m(s^*,Z^I u_{\ih})^{1/2} = C_1\eps,\quad \text{for } 1\leq \ih\leq j_0.
\endaligned
\ee
Yet, according to Proposition \ref{main prop2} there exists $\eps_1>0$ such that, for  $\eps\leq \eps_1$, the following enhanced energy bounds
also hold:
$$
\aligned
&E_{m,\sigma}(s^*,Z^{\Is} v_{\jc})^{1/2} \leq \frac{1}{2}C_1\eps {s^*}^{\delta} \quad \text{ for } j_0+1\leq \jc\leq n_0,
\\
&E_m(s^*,Z^{\Is} u_{\ih})^{1/2} \leq \frac{1}{2} C_1\eps {s^*}^{\delta} \quad \text{for } 1\leq \ih\leq j_0,
\\
&E_{m,\sigma}(s^*,Z^{\Id} v_{\jc})^{1/2} \leq \frac{1}{2}C_1\eps {s^*}^{\delta/2} \quad \text{ for } j_0+1\leq \jc\leq n_0,
\\
&E_m(s^*,Z^{\Id} u_{\ih})^{1/2} \leq \frac{1}{2} C_1\eps {s^*}^{\delta/2} \quad \text{for } 1\leq \ih\leq j_0,
\\
&E_m(s^*,Z^I u_{\ih})^{1/2} \leq \frac{1}{2} C_1\eps \quad \text{for } 1\leq \ih\leq j_0.
\endaligned
$$
Clearly, this is impossible, unless $s^* = +\infty$. 

Hence, we have proven that the inequalities \eqref{proof energy assumption} hold for all $s\in [s_0,+\infty)$ and, by the local existence criteria in Theorem~\ref{local semi-hyper} (combined with Sobolev's inequality), this local solution extends for all times. By taking 
$\eps_0 := \min\big(\eps_0',\,\eps_1\big)$, 
this completes the proof of Theorem~\ref{main thm main}.
\end{proof}

\section{Energy on the hypersurfaces of constant time}

To end this chapter, we prove that the wave components of the global-in-time solutions which we have constructed for \eqref{main eq main} enjoy a uniform energy estimate also on the standard hypersurfaces of constant time $t$.

\begin{proposition}\label{prop1 energy-flat}
Consider the wave components $u_i$ of a global-in-time solution $w_i$ satisfying \eqref{proof energy assumption} and \eqref{main bootstrap 4 5} (with $\eps_0$ sufficiently small and $C_1\eps<1$). The following estimate also holds for all 
time $t \geq B+1$:
\begin{equation}
\aligned
& \|\del_t Z^Iu(t,\cdot)\|_{L^2(\RR^3)} + \sum_{a}\|\del_aZ^I u(t,\cdot)\|_{L^2(\RR^3)}
\\
& \leq C_3(B,\delta,\eps_0)C_1\eps,
\endaligned
\end{equation}
where $C_3$ depends upon $\delta, \eps_0$ and the structure of the system only.
\end{proposition}

The proof of this result relies on a modified version of the fundamental energy estimate, as now stated.

\begin{lemma}
Let $(u_i)$ be the solution to the Cauchy problem
\begin{equation}
\aligned
&\Box u_i + G_i^{j\alpha\beta}\del_{\alpha\beta}u_i = F_i,
\\
&u_i|_{\Hcal_{s_0}} = {u_i}_0,\quad \del_t u_i|_{\Hcal_{s_0}} = {u_i}_0,
\endaligned
\end{equation}
where $F_i$ and $G_i^{j\alpha\beta}$ are defined   
in $\Kcal_{[s_0,+\infty)}$, 
while 
$u_i|_{\Hcal_{s_0}},\del_t u_i|_{\Hcal_{s_0}}$ are supported in
$\Hcal_{s_0}\cap \Kcal$.
There exists $\varepsilon_0>0$ such that, provided
$$
\max_{i,j,\alpha,\beta}|G_i^{j\alpha\beta}|\leq \varepsilon_0,
$$
 there exists a positive constant $C=C(\eps_0)$ 
such that 
\begin{equation}\label{eq lem1 energy-flat}
\aligned
&\bigg(\|\del_t u(t_0,\cdot)\|_{L^2(\RR^3)} + \sum_a\|\del_a u(t_0,\cdot)\|_{L^2(\RR^3)}\bigg)^2
\\
& \leq \, C\sum_i E_G(t,u_i) +C\sum_i \int_{t_0}^{(t_0^2+1)/2}
\int_{\Hcal_s}(s/t)\big|\del_t u_i F_i\big| \, dxds
\\
&\quad + C\int_{t_0}^{(t_0^2+1)/2}
\int_{\Hcal_s}(s/t)\big|\del_{\alphar}G_i^{\jr\alphar\betar} \del_tu_i \del_{\betar}u_{\jr}
- \frac{1}{2}\del_tG_i^{\jr\alphar\betar}\del_{\alphar}u_i\del_{\betar}u_{\jr}\big|\,dx ds.
\endaligned
\end{equation}
\end{lemma}

\begin{figure}
  \centering
  \includegraphics[width=1.0\textwidth]{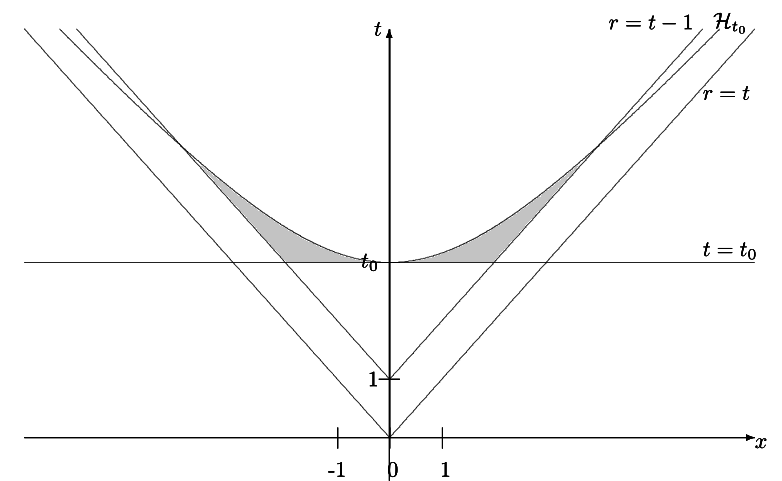}\\
  \caption{The region $\Kcal_{t_0}$}
  \label{fig kt}
\end{figure}

\begin{proof}
The proof relies on an energy estimate within the region $\Kcal_{t_0}:=\{(t,x)\big||x|\leq t-1, t_0 \leq t \leq \sqrt{|x|^2+t_0^2}\}$, which we display in 
Fig. \ref{fig kt}. 
We begin from the general identity
$$
\aligned
&\sum_i\bigg(\frac{1}{2} \del_t \sum_\alpha \big(
\del_\alpha  u_i \big)^2 + \sum_a\del_a\big(\del_a u_i \del_t u_i\big)
\\
& \qquad 
 + \del_{\alphar}\big(G_i^{\jr\alphar\betar}\del_t u_i\del_{\betar}u_{\jr}\big)
- \frac{1}{2}\del_t\big(G_i^{\jr\alphar\betar}\del_{\alphar}u_i\del_{\betar}u_{\jr}\big)\bigg)
\\
&= 
\sum_i \del_t u_i F_i
+ \sum_i\bigg(\del_{\alphar}G_i^{\jr\alphar\betar} \del_tu_i \del_{\betar}u_{\jr}
- \frac{1}{2}\del_tG_i^{\jr\alphar\betar}\del_{\alphar}u_i\del_{\betar}u_{\jr}\bigg)
\endaligned
$$
and we integrate it within $\Kcal_{t_0}$ with respect to the volume form $dtdx$. 

Applying Stokes' formula, we find 
$$
\aligned
&\frac{1}{2}\sum_i E_G(t,u_i)
\\
&-\frac{1}{2}\sum_i \int_{\RR^3}\bigg(\sum_{\alpha}|\del_\alpha u_i|^2 + 2G_i^{j0\beta}\del_tu_i\del_\beta u_j - G_i^{j\alpha\beta}\del_{\alpha}u_i\del_{\beta}u_j\bigg)(t,\cdot) dx
\\
&=
\int_{\Kcal_{t_0}} \sum_i\bigg(\del_t u_i F_i + \del_{\alphar}G_i^{\jr\alphar\betar} \del_tu_i \del_{\betar}u_{\jr}
- \frac{1}{2}\del_tG_i^{\jr\alphar\betar}\del_{\alphar}u_i\del_{\betar}u_{\jr}\bigg)\,dxdt
\\
&=\int_{\Kcal_{t_0}}(s/t) \sum_i\bigg(\del_t u_i F_i + \del_{\alphar}G_i^{\jr\alphar\betar} \del_tu_i \del_{\betar}u_{\jr}
- \frac{1}{2}\del_tG_i^{\jr\alphar\betar}\del_{\alphar}u_i\del_{\betar}u_{\jr}\bigg)\,dx ds,
\endaligned
$$
where we recall that $s = \sqrt{t^2 - |x|^2}$. This implies
$$
\aligned
&\sum_i \int_{\RR^3}\bigg(\sum_{\alpha}|\del_\alpha u_i|^2 + 2G_i^{j0\beta}\del_tu_i\del_\beta u_j - G_i^{j\alpha\beta}\del_{\alpha}u_i\del_{\beta}u_j\bigg)(t,\cdot) dx
\\
&\leq  \sum_i E_G(t,u_i)
\\
& \quad + \sum_i\int_{\Kcal_{t_0}}(s/t)\big|\del_t u_i F_i + \del_{\alphar}G_i^{\jr\alphar\betar} \del_tu_i \del_{\betar}u_{\jr}
- \frac{1}{2}\del_tG_i^{\jr\alphar\betar}\del_{\alphar}u_i\del_{\betar}u_{\jr}\big|\,dx ds,
\endaligned
$$
which is bounded by 
$$
\aligned
&\sum_i E_G(t,u_i)
 +\sum_i \int_{t_0}^{(t_0^2+1)/2}
\int_{\Hcal_s}(s/t)\big|\del_t u_i F_i\big| \, dxds
\\
&+ \int_{t_0}^{(t_0^2+1)/2}
\int_{\Hcal_s}(s/t)\big|\del_{\alphar}G_i^{\jr\alphar\betar} \del_tu_i \del_{\betar}u_{\jr}
- \frac{1}{2}\del_tG_i^{\jr\alphar\betar}\del_{\alphar}u_i\del_{\betar}u_{\jr}\big|\,dx ds.
\endaligned
$$

We then observe that
$$
\aligned
& \sum_i\int_{\RR^3}\big|2G_i^{j0\beta}\del_tu_i\del_\beta u_j - G_i^{j\alpha\beta}\del_{\alpha}u_i\del_{\beta}u_j\big| \, dx
\\
& \leq C\max_{i,j,\alpha,\beta}|G_i^{j\alpha\beta}|\sum_{i,\alpha}\|\del_{\alpha}u_i(t,\cdot)\|_{L^2(\RR^3)}
\endaligned
$$
and thus, provided $\max_{i,j,\alpha,\beta}|G_i^{j\alpha\beta}|$ is sufficiently small,
$$
\sum_i\int_{\RR^3}\big|2G_i^{j0\beta}\del_tu_i\del_\beta u_j - G_i^{j\alpha\beta}\del_{\alpha}u_i\del_{\beta}u_j\big| \, dx\leq \frac{1}{2}\sum_{i,\alpha}\|\del_{\alpha}u_i(t,\cdot)\|_{L^2(\RR^3)}.
$$
We thus have 
$$
\aligned
& \frac{1}{2}\sum_{i,\alpha}\|\del_{\alpha}u(t,\cdot)\|_{L^2(\RR^3)}
\\
& \leq\sum_i \int_{\RR^3}\bigg(\sum_{\alpha}|\del_\alpha u_i|^2 + 2G_i^{j0\beta}\del_tu_i\del_\beta u_j - G_i^{j\alpha\beta}\del_{\alpha}u_i\del_{\beta}u_j\bigg) dx
\endaligned
$$
and the desired result is proven.
\end{proof}

\begin{proof}[Proof of Proposition \ref{prop1 energy-flat}]
We recall the equations satisfied by the wave components:
$$
\Box (Z^I u_{\ih}) + G_{\ih}^{\jh\alphar\betar}\del_{\alphar\betar} (Z^Iu_{\jh}) 
= [G_{\ih}^{\jh\alphar\betar}\del_{\alphar}\del_{\betar},Z^I]u_{\jh} - Z^I\big(G_{\ih}^{\jc\alphar\betar}\del_{\alphar}\del_{\betar}v_{\jc}\big) + Z^IF_{\ih}
$$
and apply \eqref{eq lem1 energy-flat}, so that 
$$
\aligned
&\sum_{\ih,\alpha}\|\del_\alpha Z^Iu_{\ih}(t,\cdot)\|_{L^2(\RR^3)}
\\
&\leq C\sum_{\ih} E_G(t,Z^Iu_{\ih})
+C\sum_{\ih} \int_{t_0}^{(t_0^2+1)/2}
\int_{\Hcal_s}(s/t)\big|\del_t Z^Iu_{\ih} Z^IF_{\ih}\big| \, dxds
\\
&\quad+\sum_{\ih} \int_{t_0}^{(t_0^2+1)/2}
\int_{\Hcal_s}(s/t)\big|\del_t Z^Iu_{\ih} [G_{\ih}^{\jh\alphar\betar}\del_{\alphar}\del_{\betar},Z^I]u_{\jh}\big| \, dxds
\\
&\quad+\sum_{\ih} \int_{t_0}^{(t_0^2+1)/2}
\int_{\Hcal_s}(s/t)\big|\del_t Z^Iu_{\ih}
Z^I\big(G_{\ih}^{\jc\alphar\betar}\del_{\alphar}\del_{\betar}v_{\jc}\big)\big| \, dxds
\\
&+ C\sum_{\ih}\int_{t_0}^{(t_0^2+1)/2} \hskip-.15cm 
\int_{\Hcal_s}
\hskip-.15cm
(s/t)\big|\del_{\alphar}G_{\ih}^{\jr\alphar\betar} \del_t Z^Iu_{\ih} \del_{\betar}Z^Iu_{\jr}
- \frac{1}{2}\del_tG_{\ih}^{\jr\alphar\betar}\del_{\alphar}Z^Iu_{\ih}\del_{\betar}Z^Iu_{\jr}\big| dx ds.
\endaligned
$$

By \eqref{proof energy assumption e} (which was established within our bootstrap argument and holds on the time interval $[B+1,+\infty)$), we have
$$
\sum_{\ih} E_G(t,Z^Iu_{\ih}) \leq (C_1\eps)^2. 
$$
Also, the second term in the right-hand side is uniformly bounded, as follows:
$$
\aligned
&\int_{t_0}^{(t_0^2+1)/2}
\int_{\Hcal_s}(s/t)\big|\del_t Z^Iu_{\ih} Z^IF_{\ih}\big| \, dxds
\\
&\leq \int_{t_0}^{(t_0^2+1)/2}\|(s/t)\del_tZ^Iu_{\ih}\|_{L^2(\Hcal_s)} \|Z^IF_{\ih}\|_{L^2(\Hcal_s)}ds
\\
&\leq C_1\eps \int_{t_0}^{(t_0^2+1)/2}\|Z^IF_{\ih}\|_{L^2(\Hcal_s)}ds
\\
&\leq C_1\eps \int_{t_0}^{(t_0^2+1)/2}CC^*(C_1\eps)^2s^{-3/2+2\delta}ds
\\
&\leq CC^*(C_1\eps)^3(1/2 - 2\delta)^{-1}(B+1)^{-1/2+2\delta}.
\endaligned
$$
The third and fourth terms are bounded in the same manner, that is, 
$$
\aligned
&\int_{t_0}^{(t_0^2+1)/2}
\int_{\Hcal_s}(s/t)\big|\del_t Z^Iu_{\ih} [G_{\ih}^{\jh\alphar\betar}\del_{\alphar}\del_{\betar},Z^I]u_{\jh}\big| \, dxds
\\
& \leq CC^*(C_1\eps)^3(1/2 - 2\delta)^{-1}(B+1)^{-1/2+2\delta}
\endaligned
$$
and 
$$
\aligned
&
\int_{t_0}^{(t_0^2+1)/2}
\int_{\Hcal_s}(s/t)\big|\del_t Z^Iu_{\ih}
Z^I\big(G_{\ih}^{\jc\alphar\betar}\del_{\alphar}\del_{\betar}v_{\jc}\big)\big| \, dx ds
\\
&\leq CC^*(C_1\eps)^3(1/2 - 2\delta)^{-1}(B+1)^{-1/2+2\delta}.
\endaligned
$$
The last term is bounded by applying \eqref{main bootstrap 4}:
$$
\aligned
&\int_{t_0}^{(t_0^2+1)/2}
\int_{\Hcal_s}(s/t)\big|\del_{\alphar}G_{\ih}^{\jr\alphar\betar} \del_t Z^Iu_{\ih} \del_{\betar}Z^Iu_{\jr}
- \frac{1}{2}\del_tG_{\ih}^{\jr\alphar\betar}\del_{\alphar}Z^Iu_{\ih}\del_{\betar}Z^Iu_{\jr}\big|\,dx ds
\\
&\leq CC^*(C_1\eps)^3(1/2 - 2\delta)^{-1}(B+1)^{-1/2+2\delta}.
\endaligned
$$
\end{proof}

\chapter[Decompositions and estimates for the commutators]{Decompositions and estimates for the commutators \label{cha:3}}
 
\section{Decompositions of commutators. I}
\label{sec:31}

In this chapter, we present technical results concerning the commutators $ [X,Y]u := X(Yu)- Y(Xu)$ of certain operators $X,Y$ associated with the set of admissible vector fields $Z$ (cf.~Chapter~\ref{cha:2}) 
and 
applied to functions $u$ defined in the cone 
$\Kcal=\{|x|< t-1\}$.
In order to derive uniform bounds, we will rely on homogeneity arguments and 
on the observation that all the coefficients in the following decomposition 
are smooth within $\Kcal$.

First of all, all admissible vector fields $\del_\alpha, L_a$ under consideration 
are Killing fields for the flat wave operator $\Box$, so that
 the following commutation relations hold: 
\be
[\del_\alpha, \,\Box]=0, \qquad [L_a,\, \Box] =0. 
\ee
By introducing the notation
\begin{subequations}
\label{pre commutator base}
\bel{pre commutator base L-P}
[L_a,\del_{\beta}]u =: \Theta_{a\beta}^{\gamma}\del_{\gamma}u,
\ee
\bel{pre commutator base P-B}
[\del_\alpha,\delu_{\beta}]u =: t^{-1}\Gammau_{\alpha\beta}^{\gamma}\del_{\gamma}u,
\ee
\bel{pre commutator base L-Bs}
[L_a,\delu_{\beta}]u =: \Thetau_{a\beta}^{\gamma}\delu_{\gamma}u,
\ee
\end{subequations}
we find easily that 
\bel{pre commutator base'}
\aligned
&\Theta_{ab}^{\gamma} = -\delta_{ab}\delta_0^{\gamma}, \quad
&&\Theta_{a0}^{\gamma} = -\delta_a^{\gamma},
\\
&\Gammau_{ab}^{\gamma}= \delta_{ab}\delta_0^{\gamma},\quad
&&\Gammau_{0b}^{\gamma} = -\frac{x^b}{t}\delta_0^{\gamma} = \Psi^0_b\delta_0^{\gamma},\quad&&&
 \Gammau_{\alpha0}^{\gamma}  = 0,
\\
&\Thetau_{ab}^{\gamma} = -\frac{x^b}{t}\delta^{\gamma}_a = \Psi^0_b\delta^{\gamma}_a,\quad
&&\Thetau_{a0}^{\gamma} = -\delta^{\gamma}_a + \frac{x^a}{t}\delta^{\gamma}_0 =  -\delta^{\gamma}_a + \Phi_0^a\delta^{\gamma}_0,
\endaligned
\ee
where $\Phi$ and $\Psi$ were defined in \eqref{pre Phi}. 
All of these coefficients are smooth in the cone $\Kcal$ and for each index $I$,
the functions $Z^I\Pi$ are bounded for all $\Pi \in \big\{ \Theta, \Thetau, \Gammau \big\}$.
Furthermore, we can also check that 
\be\label{pre commutator base''}
\Thetau_{ab}^0 = 0 \quad \text{ so that } \quad  
[L_a,\delu_b] = \Thetau_{ab}^{\beta}\delu_{\beta} = \Thetau_{ab}^c\delu_c.
\ee

\begin{lemma}[Algebraic decompositions of commutators. I]
\label{lem-com1}
There exist {\rm constants} $\theta_{\alpha J}^{I\beta}$ such that, for all sufficiently regular function $u$ defined in the cone $\Kcal$, the following identities hold for all $\alpha$ and $I$:
\bel{pre lem commutator pr1}
[Z^I,\del_\alpha ]u = \sum_{\atop |J|< |I|}\theta_{\alpha J}^{I\beta}\del_{\beta}Z^Ju.
\ee
\end{lemma}

\begin{proof}
The proof is done by induction on $|I|$. When $|I|=1$, \eqref{pre lem commutator pr1} is implied by \eqref{pre commutator base L-P}. Assume next that \eqref{pre lem commutator pr1} is valid for $|I|\leq k$ and let us derive this property for all $|I|= k+1$. Let $Z^I$ be a product with index satisfying $|I|=k+1$,
 where $Z^I = Z_1 Z^{I'}$ with $|I'|=k$ and $Z_1$ is one of the fields $\del_{\gamma}, L_a$.
In view of 
$$
[Z^I,\del_\alpha ]u = [Z_1Z^{I'},\del_\alpha ]u = Z_1\big([Z^{I'},\del_\alpha ]u\big) + [Z_1,\del_\alpha ]Z^{I'}u
$$
and
$$
Z_1\big([Z^{I'},\del_\alpha ]u\big)
= Z_1\bigg(\sum_{|J|\leq k-1}\theta_{\alpha J}^{I'\gamma}\del_{\gamma}Z^{J}u\bigg)
=\sum_{|J|\leq k-1}\theta_{\alpha J}^{I'\gamma}Z_1\del_{\gamma}Z^{J}u, 
$$
we have
$$
[Z^I,\del_\alpha ]u = \sum_{|J|\leq k-1}\theta_{\alpha J}^{I'\gamma}Z_1\del_{\gamma}Z^{J}u + [Z_1,\del_\alpha ]Z^{I'}u.
$$

First of all, if $Z_1=\del_{\beta}$, we have the commutation property $[Z_1,\del_\alpha ]Z^{I'}u = 0$ and
$$
\aligned
\sum_{|J|\leq k-1}\theta_{\alpha J}^{I'\gamma}Z_1\del_{\gamma}Z^{J}u
=& \sum_{|J|\leq k-1}\theta_{\alpha J}^{I'\gamma}\del_{\gamma}Z_1\,  Z^{J}u,
\endaligned
$$
so that \eqref{pre lem commutator pr1} is established in this case.

Second, if $Z_1 = L_a$, we apply \eqref{pre commutator base L-P} and write
$$
\aligned
\sum_{|J|\leq k-1}\theta_{\alpha J}^{I_1\gamma}L_a\del_{\gamma}Z^{J}u
=&\sum_{|J|\leq k-1}\theta_{\alpha J}^{I'\gamma}\del_{\gamma}L_a Z^{J}u
+ \sum_{|J|\leq k-1}\theta_{\alpha J}^{I'\gamma}[L_a,\del_{\gamma}]Z^{J}u
\\
=&\sum_{|J|\leq k-1}\theta_{\alpha J}^{I'\gamma}\del_{\gamma}L_a Z^{J}u
+ \sum_{|J|\leq k-1}\theta_{\alpha J}^{I'\beta}\Theta_{a\beta}^{\gamma}\del_{\gamma}Z^{J}u.
\endaligned
$$
So, \eqref{pre lem commutator pr1} holds for $|I|=k+1$, which completes the induction argument.
\end{proof}


\section{Decompositions of commutators. II and III}
\label{sec:32} 

\begin{lemma}[Algebraic decompositions of commutators. II]
\label{lem-com2}
For all sufficiently regular functions $u$ defined in the cone $\Kcal$,
 the following identity holds:
\bel{pre lem commutator pr2}
[Z^I,\delu_{\beta}]u = \sum_{\atop |J|< |I|}\thetau_{\beta J}^{I \gamma}\del_{\gamma}Z^Ju, 
\ee
where the functions $\thetau_{\beta J}^{I\gamma}$ are smooth and satisfy in $\Kcal$ :
\bel{pre lem commutator pr4a}
\big|\del^{I_1}Z^{I_2}\thetau_{\beta J}^{I\gamma}\big| \leq C\big(n,|I|,|I_1|,|I_2|\big) \, t^{-|I_1|}.
\ee
\end{lemma}

\begin{proof} Consider the  identity
$$
\aligned
\,[Z^I,\delu_{\beta}]u
= [Z^I, \Phi_{\beta}^{\gamma}\del_{\gamma}]u
=& \sum_{I_1+I_2=I\atop |I_2|<|I|}
Z^{I_1}\Phi_{\beta}^{\gamma}Z^{I_2}\del_{\gamma}u
+\Phi_{\beta}^{\gamma}[Z^I,\del_{\gamma}]u.
\endaligned
$$
In the first sum, we commute $Z^{I_2}$ and $\del_{\gamma}$ and obtain
$$
\aligned
\,[Z^I,\delu_{\beta}]u
=& \sum_{I_1+I_2=I\atop |I_2|<|I|}
Z^{I_1}\Phi_{\beta}^{\gamma}\del_{\gamma}Z^{I_2}u
 + \sum_{I_1+I_2=I\atop |I_2|<|I|}
Z^{I_1}\Phi_{\beta}^{\gamma}[Z^{I_2},\del_{\gamma}]u
+ \Phi_{\beta}^{\gamma}[Z^I,\del_{\gamma}]u
\\
=& \sum_{I_1+I_2=I\atop |I_2|<|I|}
Z^{I_1}\Phi_{\beta}^{\gamma}\del_{\gamma}Z^{I_2}u
 + \sum_{I_1+I_2=I}
Z^{I_1}\Phi_{\beta}^{\gamma}[Z^{I_2},\del_{\gamma}]u
\\
=& \sum_{I_1+I_2=I\atop |I_2|<|I|}
Z^{I_1}\Phi_{\beta}^{\gamma}\del_{\gamma}Z^{I_2}u
 + \sum_{I_1+I_2=I\atop |J|<|I_2|}
\big(Z^{I_1}\Phi_{\beta}^{\gamma}\big)\theta_{\gamma J}^{I_2\alpha}\del_\alpha Z^Ju.
\endaligned
$$
Hence, $\thetau_{\gamma J}^{I\alpha}$ are linear combinations of $Z^{I_1}\Phi_{\beta}^{\gamma}$ and
$\theta_{\gamma J}^{I_2\alpha}Z^{I_1}\Phi_{\beta}^{\gamma}$ with $|I_1|\leq |I|$ and $|I_2|\leq |I|$,
which yields \eqref{pre lem commutator pr2}.
Note that $\theta_{\gamma J}^{I_2\alpha}$ are constants, so that
$$
\del^{I_3}Z^{I_4}\big(\theta_{\gamma J}^{I_2\alpha}Z^{I_1}\Phi_{\beta}^{\gamma}\big)
= \theta_{\gamma J}^{I_2\alpha}\del^{I_3}Z^{I_4}Z^{I_1}\Phi_{\beta}^{\gamma}
$$
and, by \eqref{pre lem commutator pr5}, we arrive at \eqref{pre lem commutator pr4a}.
\end{proof}

\begin{lemma}[Algebraic decompositions of commutators. III]
For all sufficiently regular functions $u$ defined in the cone $\Kcal$, the following identities hold:
\bel{pre lem commutator pr3}
[Z^I,\delu_c]u = \sum_{|J|<|I|}\sigma^{Ib}_{cJ}\delu_bZ^Ju + t^{-1}\sum_{|J'|<|I|}\rho_{cJ'}^{I\gamma}\del_{\gamma}Z^{J'} u,
\ee
where the functions $\sigma_{c J}^{I\gamma}$ and $\rho_{cJ}^{I\gamma}$ are 
smooth and 
satisfy in $\Kcal$:
\begin{subequations}\label{pre lem commutator pr4}
\bel{pre lem commutator pr4b}
\big|\del^{I_1}Z^{I_2}\sigma_{\beta J}^{I\gamma}\big| \leq C(n,|I|,|I_1|,|I_2|)t^{-|I_1|},
\ee
\bel{pre lem commutator pr4c}
\big|\del^{I_1}Z^{I_2}\rho_{\beta J}^{I\gamma}\big| \leq C(n,|I|,|I_1|,|I_2|)t^{-|I_1|}.
\ee
\end{subequations}
\end{lemma}

\begin{proof}
We proceed by induction. The case $|I|=1$ is easily checked from \eqref{pre commutator base P-B}, \eqref{pre commutator base L-Bs} and \eqref{pre commutator base''}. 
Assuming that \eqref{pre lem commutator pr3} with \eqref{pre lem commutator pr4b} and \eqref{pre lem commutator pr4c} hold for $|I|\leq k$, we consider indices $|I|=k+1$. To do this, let $Z^I$ be a product of operators with index $|I|=k+1$. Then, $Z^I = Z_1Z^{I'}$ with $|I'|=k$ and the following holds:
$$
[Z^I,\delu_c]u = [Z_1Z^{I'},\delu_c] = Z_1\big([Z^{I'},\delu_c]u\big) + [Z_1,\delu_c]Z^{I'}u.
$$

\vskip.3cm

\noindent {\bf Case $Z_1 = \del_{\alpha}$.} In this case, we write
$$
\aligned
&Z_1\big([Z^{I'},\delu_c]u\big)
\\
&=  \del_\alpha \bigg(\sum_{|J|<|I'|}
\sigma_{cJ}^{I'b}\delu_bZ^Ju
+ t^{-1}\sum_{|J'|<|I'|}\rho_{cJ'}^{I'\gamma}\del_{\gamma}Z^{J'}u\bigg)
\\
&= \sum_{|J|<|I'|}\del_\alpha \big(\sigma_{cJ}^{I'b}\big)\delu_bZ^Ju
 + \sum_{|J|<|I'|}\sigma_{cJ}^{I'b}\delu_b\del_\alpha Z^Ju
 + \sum_{|J|<|I'|}\sigma_{cJ}^{I'b}[\del_\alpha,\delu_b]Z^Ju
\\
&\quad
 + \del_\alpha t^{-1}\sum_{|J'|<|I'|}\big(\rho_{cJ'}^{I'\gamma}\big)\del_{\gamma}Z^{J'}u
+ t^{-1}\sum_{|J'|<|I'|}\del_\alpha \big(\rho_{cJ'}^{I'\gamma}\big)\del_{\gamma}Z^{J'}u
\\
 & \quad
  + t^{-1}\sum_{|J'|<|I'|}\rho_{cJ'}^{I'\gamma}\del_{\gamma}\del_\alpha Z^{J'}u.
\endaligned
$$ 
For the third term, we recall \eqref{pre commutator base P-B} and write
$$
\sigma_{cJ}^{I'b}[\del_{\alpha},\delu_b]Z^Ju = t^{-1}\sigma_{cJ}^{I'b}\Gammau_{\alpha b}^\gamma \del_{\gamma}Z^Ju.
$$
Hence, we find 
$$
\aligned
&Z_1\big([Z^{I'},\delu_c]u\big)
\\
&= \sum_{|J|<|I'|}\del_\alpha \big(\sigma_{cJ}^{I'b}\big)\delu_bZ^Ju
 + \sum_{|J|<|I'|}\sigma_{cJ}^{I'b}\delu_b\del_\alpha Z^Ju
 + t^{-1}\sum_{|J|<|I'|}\sigma_{cJ}^{I'b}\Gammau_{\alpha b}^\gamma\del_{\gamma}Z^Iu
\\
&\quad 
+ \del_\alpha t^{-1}\sum_{|J'|<|I'|}\big(\rho_{cJ'}^{I'\gamma}\big)\del_{\gamma}Z^{J'}u
  + t^{-1}\sum_{|J'|<|I'|}\del_\alpha \big(\rho_{cJ'}^{I'\gamma}\big)\del_{\gamma}Z^{J'}u
\\
  &\quad + t^{-1}\sum_{|J'|<|I'|}\rho_{cJ'}^{I'\gamma}\del_{\gamma}\del_\alpha Z^{J'}u.
\endaligned
$$

We conclude that, when $Z_1 = \del_\alpha $,  $\sigma_{\alpha J}^{I\gamma}$ are linear combinations of $\del_{\alpha}\sigma_{cJ}^{I'b}$ and $\sigma_{cJ}^{I'b}$ with $|J|<|I'|$. For all $\del^{I_1}Z^{I_2}$, we have
$$
\aligned
\del^{I_1}Z^{I_2}\big(\del_\alpha \sigma_{cJ}^{I'b}\big)
=&\del_{\alpha}\del^{I_1}Z^{I_2}\sigma_{cJ}^{I'b} + \del^{I_1}\big([Z^{I_2},\del_{\alpha}]\sigma_{cJ}^{I'b}\big)
\\
=&\del_{\alpha}\del^{I_1}Z^{I_2}\sigma_{cJ}^{I'b}
 +\sum_{|I_2'|<|I_2|}\theta_{\alpha I_2'}^{I_2\beta}\del_{\beta} \del^{I_1}Z^{I_2'}\sigma_{cJ}^{I_2\beta}. 
\endaligned
$$
By applying the induction assumption, we find 
$$
\aligned
\big|\del_{\alpha}\del^{I_1}Z^{I_2}\sigma_{cJ}^{I'b}\big|\leq C(n,I_1,I_2,I')t^{-|I_1|-1},
\endaligned
$$
so that 
$$
\big|\del^{I_1}Z^{I_2}\big(\del_{\alpha}\sigma_{cJ}^{I'b}\big)\big|\leq C(n,I_1,I_2,I')t^{-|I_1|}.
$$
Similarly, $\sigma_{cJ}^{I'b}$ satisfies the same estimate and we see that, in this case, $\sigma_{\alpha J}^{I\gamma}$ satisfies the desired estimate.

Note that $\rho_{\alpha J}^{I\gamma}$ are linear combinations of the following terms:
$$
\aligned
&\sigma_{cJ}^{I'b}\Gammau_{\alpha b}^\gamma,\quad t\del_{\alpha}t^{-1}\rho_{cJ}^{I'\gamma},\quad \del_{\alpha}\rho_{cJ}^{I\gamma},\quad \rho_{cJ}^{I'\gamma}.
\endaligned
$$
Observe that $\Gammau_{\alpha\beta}^\gamma$ are linear combinations of $\Phi_{\alpha}^\beta$ and recall the estimate \eqref{pre lem commutator pr5}, we see that each term can be controlled by $C(n,I_1,I_2,I)t^{-|I_1|}$.

\vskip.3cm 

\noindent{\bf Case $Z_1 = L_a$.} In this case,  we write
$$
\aligned
&Z_1\big([Z^{I'},\delu_c]u\big)
\\
&= L_a\big(\sum_{|J|<|I'|}\sigma_{cJ}^{I'b}\delu_bZ^Ju
  + t^{-1}\sum_{|J'|<|I'|}\rho_{cJ'}^{I'\gamma}\del_{\gamma}Z^{J'}u\big)
\\
&= \sum_{|J|<|I'|}\big(L_a\sigma_{cJ}^{I'b}\big)\delu_bZ^Ju
 +\sum_{|J|<|I'|}\sigma_{cJ}^{I'b}\delu_bL_aZ^{J}u
 +\sum_{|J|<|I'|}\sigma_{cJ}^{I'b}[L_a,\delu_b]Z^Ju
\\
&\quad + t^{-1}(tL_at^{-1})\sum_{\gamma\atop |J'|<|I'|}\rho_{cJ'}^{I'\gamma}\del_{\gamma}Z^{J'}u
\\
&\quad+ t^{-1}\sum_{|J'|<|I'|}\big(L_a\rho_{cJ'}^{I'\gamma}\big)\del_{\gamma}Z^{J'}u
 + t^{-1}\sum_{|J'|<|I'|}\rho_{cJ'}^{I'\gamma}\del_{\gamma}L_aZ^{J'}u
\\
& \quad
 + t^{-1}\sum_{|J'|<|I'|}\rho_{cJ'}^{I'\gamma}[L_a,\del_{\gamma}]Z^{J'}u.
\endaligned
$$
In the right-hand side, we apply \eqref{pre commutator base L-Bs} on the third term and \eqref{pre commutator base L-P} on the last term, and observe that the coefficient of the fourth term 
$$
tL_a(t^{-1}) = -\frac{x^a}{t} = \Psi^0_a.
$$
We have 
$$
\aligned
&Z_1\big([Z^{I'},\delu_c]u\big)
\\
&=\sum_{|J|<|I'|}\big(L_a\sigma_{cJ}^{I'b}\big)\delu_bZ^Ju
 +\sum_{|J|<|I'|}\sigma_{cJ}^{I'b}\delu_bL_aZ^{J}u
 +\sum_{|J|<|I'|}\sigma_{cJ}^{I'b}\Thetau_{ab}^c\delu_c Z^Ju
\\
&\quad+ t^{-1}\Psi_a^0\sum_{|J'|<|I'|}\rho_{cJ'}^{I'\gamma}\del_{\gamma}Z^{J'}u
\\
&\quad+ t^{-1}\sum_{|J'|<|I'|}\big(L_a\rho_{cJ'}^{I'\gamma}\big)\del_{\gamma}Z^{J'}u
 + t^{-1}\sum_{|J'|<|I'|}\rho_{cJ'}^{I'\gamma}\del_{\gamma}L_aZ^{J'}u
\\
&\quad
 + t^{-1}\sum_{|J'|<|I'|}\rho_{cJ'}^{I'\gamma}\Theta_{a\gamma}^{\beta}\del_{\beta} Z^{J'}u.
\endaligned
$$ 
We also recall that
$$
[Z_1,\delu_c]Z^{I'}u = [L_a,\del_c]Z^{I'}u = \Thetau_{ac}^b\delu_b Z^{I'}u.
$$ 
We observe that $\sigma_{cJ}^{Ib}$ are linear combinations of
$L_a\sigma_{cJ}^{I'b}$, $\sigma_{cJ}^{I'b}$, $\sigma_{cJ}^{I'b}\Thetau_{ab}^c$ and $\Thetau_{ac}^b$.
The estimate on the first two terms is a direct consequence of the induction assumption \eqref{pre lem commutator pr4b}. On the other hand, the third term is estimated as follows:
$$
\del^{I_1}Z^{I_2}\sigma_{cJ}^{I'b}\Thetau_{ab}^c
= \sum_{I_3+I_4=I_1\atop I_5+I_6=I_2}\del^{I_3}Z^{I_5}\sigma_{cJ}^{I'b}\,  \del^{I_4}Z^{I_6}\Thetau_{ab}^c.
$$
The first factor is bounded by the induction assumption \eqref{pre lem commutator pr4b}, and the second factor can be bounded by \eqref{pre lem commutator pr6}, by observing that  $\Thetau_{ab}^c$ is a linear combination of $\Psi_{\gamma}^{\delta}$. In the same manner, we get the desired estimate on $\Thetau_{ac}^b$. So, \eqref{pre lem commutator pr4b} is established in the case $Z_1=L_a,\, |I|=k+1$. 

In the same way, $\rho_{cJ}^{I\gamma}$ are linear combinations of
$\Psi_a^0\rho_{cJ'}^{I'\gamma}$,
$L_a\rho_{cJ'}^{I'\gamma}$,
$\rho_{cJ'}^{I'\gamma}$,
$\rho_{cJ'}^{I'\gamma}\Theta_{a\gamma}^{\beta}$.
The second and the third terms are bounded by the induction assumption \eqref{pre lem commutator pr4c}. The first term is estimated by applying \eqref{pre lem commutator pr6} and the induction assumption \eqref{pre lem commutator pr4c}. 
For the last term, note that $\Theta_{a\gamma}^{\beta}$ are linear combinations of $\Psi_{\delta}^{\delta'}$ and constants. Hence,
by applying \eqref{pre lem commutator pr6} and the induction assumption \eqref{pre lem commutator pr4c}, the desired bound is reached. Finally, by combining these estimates together, \eqref{pre lem commutator pr4c} is established in the case $|I|=k+1$, which completes the argument
and, consequently, the proof of Lemma \ref{lem-com2}.
\end{proof}


\section{Estimates of commutators}
\label{sec:33}

The following statement is now immediate in view of \eqref{pre lem commutator pr1}, \eqref{pre lem commutator pr2}, and \eqref{pre lem commutator pr3}.

\begin{lemma}\label{pre lem commutator1}
For all sufficiently regular functions $u$ defined in the cone $\Kcal$, the following estimates hold:
\bel{pre lem commutator trivial}
\big|[Z^I,\delu_{\alpha}] u\big| + \big|[Z^I,\del_\alpha ] u\big|\leq C(n,|I|)\sum_{\beta,|J|<I}\big|\del_{\beta}Z^J u\big|,
\ee
\bel{pre lem commutator bar}
\big|[Z^I,\delu_a]u\big|\leq  C(n,|I|)\sum_{b,|J_1|<|I|}\big|\delu_b Z^{J_1} u\big|
                             +C(n,|I|)\sum_{\gamma,J_2|<|I|}\big|t^{-1}\del_{\gamma}Z^{J_2}u\big|.
\ee
\end{lemma}

Furthermore, from \eqref{pre lem commutator pr1}, \eqref{pre lem commutator pr2}, and \eqref{pre lem commutator pr3}, we deduce the following estimates. Recall our convention $Z^I = 0$ in case we write $|I|<0$.

\begin{lemma}\label{pre lem commutator}
For all sufficiently regular functions $u$ defined in the cone $\Kcal$, the following estimates hold:
\bel{pre lem commutator second-order}
\big|[Z^I,\del_\alpha \del_{\beta}] u \big| 
\leq  C(n,|I|)\sum_{\gamma,\gamma'\atop |J|<|I|} \big|\del_{\gamma}\del_{\gamma'}Z^J u\big|,
\ee
\bel{pre lem commutator second-order bar}
\aligned
\big|[Z^I, \delu_a\delu_{\beta}] u\big| + \big|[Z^I, \delu_\alpha \delu_b] u\big|
\leq& C(n,|I|) \sum_{c,\gamma\atop|J_1|\leq |I|}\big|\delu_c \delu_{\gamma} Z^{J_1}u\big|
\\
&+ C(n,|I|)t^{-1}\sum_{\gamma\atop |J_2|\leq|I|}\big|\del_{\gamma}Z^{J_2}u\big|.
\endaligned
\ee
\end{lemma}

\begin{proof} 
1. To derive \eqref{pre lem commutator second-order}, we write
$$
[Z^I,\del_\alpha \del_{\beta}]u =\del_\alpha \big([Z^I,\del_{\beta}]u\big) + [Z^I,\del_\alpha ]\del_{\beta}u
$$
and, by applying \eqref{pre lem commutator pr1}, the first term in the right-hand 
side can be written as
$$
\del_\alpha \big([Z^I,\del_{\beta}]u\big)
= \del_\alpha \bigg(\sum_{|J|<|I|}\theta_{\beta J}^{I \gamma}\del_{\gamma}Z^Ju\bigg) =\sum_{|J|<|I|}\theta_{\beta J}^{I \gamma}\del_\alpha \del_{\gamma}Z^Ju, 
$$
which is bounded by
$
C(n,|I|)\sum_{\alpha,\beta\atop |J|<|I|}\big|\del_\alpha \del_{\beta}Z^Ju\big|.
$
The second term is estimated as follows:
$$
\aligned
\,[Z^I,\del_\alpha ]\del_{\beta}u
=&\sum_{|J|<|I|}\theta_{\alpha J}^{I\gamma}\del_{\gamma}Z^J\del_{\beta}u
\\
=&\sum_{|J|<|I|}\theta_{\alpha J}^{I\gamma}\del_{\gamma}\del_{\beta} Z^Ju
+ \sum_{|J|<|I|}\theta_{\alpha J}^{I\gamma}\del_{\gamma}[Z^J,\del_{\beta}]u
\\
=&\sum_{|J|<|I|}\theta_{\alpha J}^{I\gamma}\del_{\gamma}\del_{\beta} Z^Ju
+ \sum_{|J|<|I|}\theta_{\alpha J}^{I\gamma}\del_{\gamma}\bigg(\sum_{|J'|<|J|}\theta_{\beta J'}^{J\delta}\del_{\delta}Z^{J'}u\bigg)
\\
=&\sum_{|J|<|I|}\theta_{\alpha J}^{I\gamma}\del_{\gamma}\del_{\beta} Z^Ju
+ \sum_{|J|<|I|\atop |J'|<|J|}\theta_{\alpha J}^{I\gamma}\theta_{\beta J'}^{J\delta}\del_{\gamma}\del_{\delta}Z^{J'}u.
\endaligned
$$
This latter expression is bounded by
$
C(n,|I|)\sum_{\alpha,\beta\atop |J|<|I|}\big|\del_\alpha \del_{\beta}Z^Ju\big|, 
$
and \eqref{pre lem commutator second-order} is established.

 \ 

2. We now derive \eqref{pre lem commutator second-order bar}
and we begin by considering $[Z^I,\delu_a\delu_\beta]u$.
By \eqref{pre lem commutator pr2} and \eqref{pre lem commutator pr3}, we have 
$$
\aligned
&[Z^I,\delu_a\delu_{\beta}]u
\\
&=\delu_a[Z^I,\delu_{\beta}]u + [Z^I,\delu_a]\delu_{\beta}u
\\
&=\delu_a\bigg(\sum_{|J|<|I|}\thetau_{\beta J}^{I\gamma}\del_{\gamma}Z^J u\bigg)
 + \sum_{|J|<|I|} \sigma_{a J}^{I c}\delu_c Z^J\delu_{\beta}u
 + t^{-1}\sum_{|J|<|I|}\rho_{a J}^{I\gamma}\del_{\gamma}Z^J\delu_{\beta}u
\\
&=\sum_{|J|<|I|}\delu_a\thetau_{\beta J}^{I\gamma}\del_{\gamma}Z^J u
 +\sum_{|J|<|I|}\thetau_{\beta J}^{I\gamma}\delu_a\del_{\gamma}Z^J u
\\
&\quad
 + \sum_{|J|<|I|} \sigma_{a J}^{I c}\delu_c Z^J\delu_{\beta}u
 + t^{-1}\sum_{|J|<|I|}\rho_{a J}^{I\gamma}\del_{\gamma}Z^J\delu_{\beta}u,
\endaligned
$$
thus
$$
\aligned
&[Z^I,\delu_a\delu_{\beta}]u
\\
&=\sum_{|J|<|I|}\delu_a\thetau_{\beta J}^{I\gamma}\del_{\gamma}Z^J u
 +\sum_{|J|<|I|}\thetau_{\beta J}^{I\gamma}\delu_a\del_{\gamma}Z^J u
\\&
\quad
 + \sum_{|J|<|I|} \sigma_{a J}^{I c}\delu_c\delu_{\beta} Z^Ju
 + \sum_{|J|<|I|} \sigma_{a J}^{I c}\delu_c\big([Z^J,\delu_{\beta}]u\big)
 + t^{-1}\sum_{|J|<|I|}\rho_{a J}^{I\gamma}\del_{\gamma}Z^J\delu_{\beta}u
\endaligned
$$
and, therefore,
\bel{eq pr1 pre lem commutator}
\aligned
&[Z^I,\delu_a\delu_{\beta}]u
\\
&=\sum_{|J|<|I|}\delu_a\thetau_{\beta J}^{I\gamma}\del_{\gamma}Z^J u
 +\sum_{|J|<|I|}\thetau_{\beta J}^{I\gamma}\delu_a\big(\Psi_{\gamma}^{\gamma'}\delu_{\gamma'}Z^J u\big)
\\
&\quad
 + \sum_{|J|<|I|} \sigma_{a J}^{I c}\delu_c\delu_{\beta} Z^Ju
 + \sum_{|J|<|I|} \sigma_{a J}^{I c}\delu_c\big([Z^J,\delu_{\beta}]u\big)
 + t^{-1}\sum_{|J|<|I|}\rho_{a J}^{I\gamma}\del_{\gamma}Z^J\delu_{\beta}u
\\
&=\sum_{|J|<|I|}\thetau_{\beta J}^{I\gamma}\Psi_{\gamma}^{\gamma'}\delu_a\delu_{\gamma'}Z^J u
 + \sum_{|J|<|I|} \sigma_{a J}^{I c}\delu_c\delu_{\beta} Z^Ju
\\
&\quad
 +\sum_{|J|<|I|}\thetau_{\beta J}^{I\gamma}\delu_a\big(\Psi_{\gamma}^{\gamma'}\big)\delu_{\gamma'}Z^J u
 +\sum_{|J|<|I|}\delu_a\thetau_{\beta J}^{I\gamma}\del_{\gamma}Z^J u
\\
&\quad + \sum_{|J|<|I|} \sigma_{a J}^{I c}\delu_c\big([Z^J,\delu_{\beta}]u\big)
 + t^{-1}\sum_{|J|<|I|}\rho_{a J}^{I\gamma}\del_{\gamma}Z^J\delu_{\beta}u.
\endaligned
\ee
Here, the first and second terms can be bounded by
$\sum_{a,\beta,|J|<|I|}|\delu_a\delu_{\beta}Z^Ju|$. 
By recalling that
$\big|\delu_a\Psi_{\gamma}^{\gamma'}\big|\leq C(n)t^{-1}$ and $\big|\delu_a\thetau_{\beta J}^{I\gamma}\big|\leq C(n,|I|)t^{-1}$, the third and fourth can be bounded by $C(n,|I|)t^{-1}\sum_{\gamma, |J_2|\leq|I|}|\del_{\gamma}Z^{J_2}u|$.

So, we focus on the fifth term:
$$
\aligned
&\sum_{|J|<|I|} \sigma_{a J}^{I c}\delu_c\big([Z^J,\delu_{\beta}]u\big)
\\
&=\sum_{|J|<|I|} \sigma_{a J}^{I c}\delu_c\bigg(\sum_{|J'|<|J|}\thetau_{\beta J'}^{J\gamma}\del_{\gamma}u\bigg)
\\
&=\sum_{|J|<|I|\atop |J'|<|J|}\sigma_{a J}^{I c}\delu_c\big(\thetau_{\beta J'}^{J\gamma}\big)\del_{\gamma}u
+ \sum_{|J|<|I|\atop |J'|<|J|}\sigma_{a J}^{I c}\thetau_{\beta J'}^{J\gamma}\delu_c\del_{\gamma}u
\\
&=\sum_{|J|<|I|\atop |J'|<|J|}\sigma_{a J}^{I c}\delu_c\big(\thetau_{\beta J'}^{J\gamma}\big)\del_{\gamma}u
+ \sum_{|J|<|I|\atop |J'|<|J|}
\sigma_{a J}^{I c}\thetau_{\beta J'}^{J\gamma}\delu_c\big(\Psi_{\gamma}^{\gamma'}\delu_{\gamma'}u\big), 
\endaligned
$$
so
$$
\aligned
&\sum_{|J|<|I|} \sigma_{\beta J}^{I c}\delu_c\big([Z^J,\delu_{\beta}]u\big)
\\
&=\sum_{|J|<|I|\atop |J'|<|J|}\sigma_{a J}^{I c}\delu_c\big(\thetau_{\beta J'}^{J\gamma}\big)\del_{\gamma}u
+ \sum_{|J|<|I|\atop |J'|<|J|}
\sigma_{a J}^{I c}\thetau_{\beta J'}^{J\gamma}\delu_c\big(\Psi_{\gamma}^{\gamma'}\big)\delu_{\gamma'}u
\\
& \quad
+ \sum_{|J|<|I|\atop |J'|<|J|}
\sigma_{a J}^{I c}\thetau_{\beta J'}^{J\gamma}\Psi_{\gamma}^{\gamma'}\delu_c\delu_{\gamma'}u.
\endaligned
$$
Similarly, the first two terms can be bounded by 
$$
C(n,|I|)t^{-1}\sum_{\gamma, |J_2|\leq|I|}|\del_{\gamma}Z^{J_2}u|
$$ 
and the last term is bounded by $\sum_{\alpha\beta,|J|<|I|}|\delu_a\delu_\beta Z^J u|$.

For the last term in \eqref{eq pr1 pre lem commutator}, we perform the following calculation:
$$
\aligned
&t^{-1}\sum_{|J|<|I|}\rho_{aJ}^{I\gamma}\del_\gamma Z^J\delu_{\beta}u 
\\
&= t^{-1}\sum_{|J|<|I|}\rho_{aJ}^{I\gamma}\del_\gamma \delu_{\beta} Z^Ju 
+ t^{-1}\sum_{|J|<|I|}\rho_{aJ}^{I\gamma}\del_\gamma \big([Z^J,\delu_{\beta}]u\big)
\\
&= t^{-1}\sum_{|J|<|I|}\rho_{aJ}^{I\gamma}\del_\gamma\big(\Psi_\beta^{\beta'}\del_{\beta'}Z^J u\big)
+
  t^{-1}\sum_{|J|<|I|}\rho_{aJ}^{I\gamma}\del_\gamma \big(\thetau_{\beta J'}^{J\delta}\del_{\delta}Z^{J'}u \big)
\\
&= t^{-1}\sum_{|J|<|I|}\rho_{aJ}^{I\gamma}\Psi_\beta^{\beta'}\del_\gamma\del_{\beta'}Z^J u
  +t^{-1}\sum_{|J|<|I|}\rho_{aJ}^{I\gamma}\del_\gamma\Psi_\beta^{\beta'}\,\del_{\beta'}Z^J u
\\
& \quad +t^{-1}\sum_{|J|<|I|}\rho_{aJ}^{I\gamma}\thetau_{\beta J'}^{J\delta}\del_\gamma\del_{\delta}Z^{J'}u
  +t^{-1}\sum_{|J|<|I|}\rho_{aJ}^{I\gamma}\del_\gamma\thetau_{\beta J'}^{J\delta}\,\del_{\delta}Z^{J'}u.
\endaligned 
$$
The desired estimates on these terms are immediate thanks to
\eqref{pre lem commutator pr4c},  
\eqref{pre lem commutator pr4a}, and the fact that $|\del_{\alpha}\Psi_{\beta}^{\beta'}|\leq Ct^{-1}$.

Now, in order to control $[Z^I,\delu_\alpha\delu_b]u$, we consider the identity 
$$
[Z^I,\delu_\alpha\delu_b] u = [Z^I,\delu_b\delu_\alpha]u + [Z^I,[\delu_\alpha,\delu_b]]u.
$$
The estimate on the first term in the right-hand side is already done. We concentrate on the second term, via the following calculation:
$$
\aligned
& [Z^I,[\delu_\alpha,\delu_b]] u
\\
&= [Z^I,\Phi_{\alpha}^\beta[\del_\beta,\delu_b]]u - [Z^I,\delu_b\Phi_\alpha^\beta\del_\beta]u
\\
&=[Z^I,t^{-1}\Phi_{\alpha}^\beta\Gammau_{\beta b}^\gamma\del_\gamma]u - [Z^I,\delu_b\Phi_\alpha^\beta\del_\beta]u
\\
&=\sum_{I_1+I_2=I\atop |I_2|<|I|}Z^{I_1}\big(t^{-1}\Phi_{\alpha}^\beta\Gammau_{\beta b}^\gamma\big)Z^{I_2}\del_\gamma u
 + t^{-1}\Phi_{\alpha}^\beta\Gammau_{\beta b}^\gamma[Z^I,\del_\gamma]u
\\
 &\quad -\sum_{I_1+I_2=I\atop |I_2|<|I|}Z^{I_1}\big(\delu_b\Phi_{\alpha}^\beta\big) Z^{I_2}\del_\beta u
 -\delu_b\phi_{\alpha}^\beta[Z^I,\del_\beta]u
\endaligned
$$
and thus 
$$
\aligned
& [Z^I,[\delu_\alpha,\delu_b]] u
\\
&=\sum_{I_1+I_2=I\atop |I_2|<|I|}Z^{I_1}\big(t^{-1}\Phi_{\alpha}^\beta\Gammau_{\beta b}^\gamma\big)\del_\gamma Z^{I_2}u
\\
 & \quad +\sum_{I_1+I_2=I\atop |I_2|<|I|}Z^{I_1}\big(t^{-1}\Phi_{\alpha}^\beta\Gammau_{\beta b}^\gamma\big)[Z^{I_2},\del_\gamma]u
 + t^{-1}\Phi_{\alpha}^\beta\Gammau_{\beta b}^\gamma[Z^I,\del_\gamma]u
\\
 & \quad -\sum_{I_1+I_2=I\atop |I_2|<|I|}Z^{I_1}\big(\delu_b\Phi_{\alpha}^\beta\big)\del_\beta Z^{I_2} u
\\
 & \quad -\sum_{I_1+I_2=I\atop |I_2|<|I|}Z^{I_1}\big(\delu_b\Phi_{\alpha}^\beta\big)[Z^{I_2},\del_\beta]u
 -\delu_b\phi_{\alpha}^\beta[Z^I,\del_\beta]u.
\endaligned
$$
So, we have
$$
\aligned
&\,[Z^I,[\delu_\alpha,\delu_b]] u
\\
&=\sum_{I_1+I_2=I\atop |I_2|<|I|}Z^{I_1}\big(t^{-1}\Phi_{\alpha}^\beta\Gammau_{\beta b}^\gamma\big)\del_\gamma Z^{I_2}u
+\sum_{I_1+I_2=I}Z^{I_1}\big(t^{-1}\Phi_{\alpha}^\beta\Gammau_{\beta b}^\gamma\big)[Z^{I_2},\del_\gamma]u
\\
&\quad -\sum_{I_1+I_2=I\atop |I_2|<|I|}Z^{I_1}\big(\delu_b\Phi_{\alpha}^\beta\big)\del_\beta Z^{I_2} u
-\sum_{I_1+I_2=I}Z^{I_1}\big(\delu_b\Phi_{\alpha}^\beta\big)[Z^{I_2},\del_\beta]u.
\endaligned
$$
We observe that the function $t^{-1}\Phi_{\alpha}^\beta\Gammau_{\beta b}^\gamma$ and $\delu_b\Phi_{\alpha}^\beta$ are homogeneous of degree $-1$ and 
are smooth in the closed region $\{t\geq 1,|x|\leq t\}$. 

By Lemma \ref{pre lem homo}, we have
$$
\aligned
\big|Z^{I_1}\big(t^{-1}\Phi_{\alpha}^\beta\Gammau_{\beta b}^\gamma\big)\big|\leq& C(n,|I|)t^{-1},
\\
\big|Z^{I_1}\big(\delu_b\Phi_{\alpha}^\beta\big)\big|\leq& C(n,|I|)t^{-1}.
\endaligned
$$
We conclude with 
$$
\big|[Z^I,[\delu_\alpha,\delu_b]] u\big|\leq C(n,|I|)t^{-1}\sum_{\gamma\atop |J|<|I|}\big|\del_{\gamma}Z^Ju\big|,
$$
and this completes the proof of \eqref{pre lem commutator second-order bar}.
\end{proof}

\begin{lemma}
\label{pre lem commutator s/t}
For all sufficiently regular functions $u$ defined in the cone $\Kcal$, the following estimates hold:
\bel{pre lem commutator T/t'}
\big|Z^I\big((s/t)\del_\alpha u\big)\big| \leq \big|(s/t)\del_\alpha Z^Iu\big| + C(n,|I|)\sum_{\beta,|J|<|I|}\big|(s/t)\del_{\beta}Z^Ju\big|.
\ee
\end{lemma}

The proof of this lemma will rely on the following technical remark.

\begin{lemma}
\label{pre lem lem commutator s/t}
For all index $I$, the function 
\bel{pre lem s/t}
\Xi^I  := (t/s)Z^I(s/t)
\ee
defined in the closed cone $\overline{\Kcal} = \{|x|\leq t-1\}$, is smooth and all of its derivatives (of any order) 
are bounded in $\overline{\Kcal}$. Furthermore, it is homogeneous of degree $\eta$ with $\eta\leq 0$.
\end{lemma}

We admit here this result and give the proof of \eqref{pre lem commutator T/t'}, while 
the proof of \eqref{pre lem s/t} is given afterwards.

\begin{proof}[Proof of Lemma~\ref{pre lem commutator s/t}]
We observe that
$$
\aligned
\,[Z^I, (s/t)\del_\alpha ] u
=& \sum_{I_1+I_2=I\atop |I_1|<|I|}Z^{I_2}(s/t)\, Z^{I_1}\del_\alpha u
 + (s/t)[Z^I,\del_\alpha ]u
\\
=&\sum_{I_1+I_2=I\atop |I_1|<|I|}Z^{I_2}(s/t)\, [Z^{I_1},\del_\alpha ]u
+ \sum_{I_1+I_2=I\atop |I_1|<|I|}Z^{I_2}(s/t)\,\del_\alpha Z^{I_1}u
\\
&\quad + (s/t)[Z^I,\del_\alpha ]u
\\
=&\sum_{I_1+I_2=I}Z^{I_2}(s/t)\, [Z^{I_1},\del_\alpha ]u
+ \sum_{I_1+I_2=I\atop |I_1|<|I|}Z^{I_2}(s/t)\,\del_\alpha Z^{I_1}.
\endaligned
$$
By applying \eqref{pre lem s/t},
we find
$$
\aligned
\bigg|\sum_{I_1+I_2=I}Z^{I_2}(s/t)\, [Z^{I_1},\del_\alpha ]u\bigg|
\leq & \sum_{I_1+I_2=I}\big|Z^{I_2}(s/t)\big|\, \big|[Z^{I_1},\del_\alpha ]u\big|
\\
\leq& C(n,|I|)\sum_{|J|<|I|}\big|(s/t)\del_\alpha Z^{J}u\big|
\endaligned
$$
and
$$
\aligned
\bigg|\sum_{I_1+I_2=I\atop |I_1|<|I|}Z^{I_2}(s/t)\,\del_\alpha Z^{I_1}u\bigg|
& \leq C(n,|I|)\sum_{|J|<|I|}\big|(s/t)\del_\alpha Z^Ju\big|.
\endaligned
$$ 
\end{proof}

\begin{proof}[Proof of Lemma~\ref{pre lem lem commutator s/t}]
We consider the identities
\bel{pre lem lem commutator s/t 1}
\aligned
&(t/s)L_a(s/t) = -x^a/t,
\\
&(t/s)\del_a(s/t) = -x^a/s^2,\quad (t/s)\del_t(s/t) =  |x|^2s^{-2}t^{-1}.
\endaligned
\ee
In the cone $\Kcal = \{|x|<t-1\}$, the functions $x^a/s^2$, $|x|^2s^{-2}t^{-1}$,  and $x^a/t$ are smooth and bounded, while
$x^a/s^2$ and $|x|^2s^{-2}t^{-1}$ are homogeneous of degree $-1$, and 
 $x^a/t$ is homogeneous of degree $0$. 
All of their derivatives (of any order) are bounded in $\overline\Kcal$.  We have proved \eqref{pre lem s/t} in the case $|I| = 1$.

For the case $|I|>1$, we use an induction on $|I|$. Assume that for $|I|\leq k$, \eqref{pre lem s/t} holds. For an operator $Z^I$ with $|I| = k+1$, we suppose that $Z^I = Z_1\,Z^{I'}$ where $|I'|=k$ and we have 
$$
(t/s)Z^{I}(s/t) = (t/s)Z_1Z^{I'}(s/t) = Z_1\big((t/s)Z^{I'}(s/t)\big) - Z_1(s/t)Z^{I'}(s/t),
$$
where $Z_1$ can be $\del_\alpha$ or $L_a$. By the
 induction assumption 
and \eqref{pre lem lem commutator s/t 1}, the second term is a smooth function and is homogeneous of non-positive degree, while 
all of its derivatives are bounded in $\overline\Kcal$. 
We focus on the first term and, by the induction 
assumption, $(t/s)Z^{I'}(s/t)$ is smooth, homogeneous of non-positive degree, and is bounded in $\overline\Kcal$. We need to distinguish between two different cases, as follows.

\vskip.3cm

\noindent {\bf Case $Z_1 = \del_{\alpha}$.} In this case, $Z_1\big((t/s)Z^{I'}(s/t)\big)$ is homogeneous of degree less than or equal to $-1$, and by the induction assumption, all of its derivatives are bounded in $\overline \Kcal$.

\vskip.3cm

\noindent {\bf Case $Z_1 =  L_a$.} In this case, $Z_1\big((t/s)Z^{I'}(s/t)\big)$ is homogeneous of degree less than or equal to $0$.
$$
Z_1\big((t/s)Z^{I'}(s/t)\big) = x^a\del_t\big((t/s)Z^{I'}(s/t)\big) + t\del_a\big((t/s)Z^{I'}(s/t)\big).
$$
Denote by $f(t,x) = \del_{\alpha}\big((t/s)Z^{I'}(s/t)\big)$ and recall that $f(t,x)$ is homogeneous of degree $\eta$ with $\eta\leq -1$ so
$$
f(t,x) = (t/2)^{\eta}f(2,2x/t).
$$
Recall that $f(2,x)$ is bounded when $|x|\leq 2$. Recalling that $t\geq 1$, we get 
$$
\big|\del_{\alpha}\big((t/s)Z^{I'}(s/t)\big)\big| = |f(t,x)|\leq C(n,I')t^{\eta}\leq C(n,I')t^{-1} 
$$
thus
$$
\big|Z_1\big((t/s)Z^{I'}(s/t)\big)\big|\leq C(n,I')(1+|x^a|/t).
$$
Taking into consideration the fact that in $\overline\Kcal$, $|x^a|\leq t$, the desired result is proven.
\end{proof}

\chapter[The null structure in the semi-hyperboidal frame]{The null structure in the semi-hyperboidal frame \label{cha:4}}

\section{Estimating first-order derivatives}
\label{sec:41}

In this chapter, we derive various estimates on null quadratic or cubic forms.
Recall that a quadratic form $T^{\alpha\beta}\del_\alpha  u\del_{\beta} v$
acting on the gradient of functions $u, v$ is said to satisfy the {\bf null condition} if, for all null vectors $\xi\in \RR^4$, i.e.~all vectors satisfying
$(\xi_0)^2 - \sum_a (\xi_a)^2  =0$, one has
\be
T^{\alphar\betar} \xi_{\alphar}\xi_{\betar}= 0.
\ee
Similarly, a cubic form $A^{\alphar\betar\gammar}\xi_{\alphar}\xi_{\betar}\xi_{\gammar}$ is said to satisfy the null condition if, for all null vectors, one has
\be
A^{\alphar\betar\gammar} \xi_{\alphar}\xi_{\betar}\xi_{\gammar}= 0.
\ee
All the terms of interest will be linear combinations of factors $\del_{\alpha}u\,\del_{\beta}v$, $\del_{\gamma}u\del_{\alpha}\del_{\beta}v$, and $u\del_{\alpha}\del_{\beta}v$.
Throughout, the notation $s^2 = t^2 - r^2$ is in order.

\begin{proposition}
\label{pre lem null quadratic}
If $T^{\alphar\betar}\xi_{\alphar}\xi_{\betar}$ is a quadratic form satisfying the null condition, then for every index $I$
there exists a constant $C(I)>0$ such that
\be
\big|Z^I\Tu^{00}\big| \leq C(I) \, \frac{t^2-r^2}{t^2} = C(I) \, (s/t)^2 \quad \text{ in the cone } \Kcal.
\ee
Similarly, if a cubic form $A^{\alphar\betar\gammar}\xi_{\alphar}\xi_{\betar}\xi_{\gammar}$ satisfies
 the null condition, then for every index $I$
there exists a constant $C(I)>0$ such that
\be
\big|Z^I\Au^{000}\big| \leq C(I) \, \frac{t^2-r^2}{t^2} = C(I) \, (s/t)^2
 \quad \text{ in the cone } \Kcal.
\ee
\end{proposition}

The proof will be based on a homogeneity lemma stated now, 
which concerns homogeneous functions defined ``near'' the light cone in $\big \{ r=t-1 \big\}$.

\begin{lemma}[Homogeneity lemma]
\label{pre lem homor}
Let $f$ be a smooth function defined in the closed set $\{t\geq 1\}\cap \{t/2\leq r\leq t\}$. Assume that $f$ is homogeneous of degree $\eta$ in the sense that
$$
f(pt,px) = p^{\eta}f(1,x/t), \qquad t/2\leq |x|\leq t, \quad pt \geq 1, \quad t\geq 1.
$$
The following estimate holds:
\be
\big|Z^If(t,x)\big|\leq C(n,|I|)t^{\eta} \qquad \text{ in the region } \{r\geq t/2\}\cap \Kcal.
\ee
\end{lemma}

The proof is similar to the one of Lemma \ref{pre lem homo} and is omitted.

\begin{proof}[Proof of Proposition~\ref{pre lem null quadratic}]
We use the notation
$-\omega_a = \omega^a = x^a/|x|$ and $\omega_0 = \omega^0 = 1$, 
so that 
$(\omega_0)^2 - \sum_a (\omega_a)^2 = 0$.
We consider the component $\Tu^{00}$ and write
$$
\aligned
\Tu^{00}= T^{\alphar\betar}\Psi_{\alphar}^0\Psi_{\betar}^0
&= T^{\alphar\betar}\Psi_{\alphar}^0\Psi_{\betar}^0 - T^{\alphar\betar}\omega_{\alphar}\omega_{\betar}
\\
&= T^{\alphar\betar}\big(\Psi_{\alphar}^0\Psi_{\betar}^0 - \omega_{\alphar}\omega_{\betar}\big).
\endaligned
$$

First, we consider the region ``away'' from the light cone. When $r\leq t/2$, we have
$$
Z^I\big(T^{\alphar\betar}\Psi_{\alphar}^0\Psi_{\betar}^0\big)
= \sum_{I_1+I_2=I}T^{\alpha\beta}Z^{I_1}\Psi_{\alpha}^0Z^{I_2}\Psi_{\beta}^0.
$$
Applying Lemma \ref{pre lem homo}, we obtain
$$
\big|Z^I\big(T^{\alphar\betar}\Psi_{\alphar}^0\Psi_{\betar}^0\big)\big|
\leq \sum_{I_1+I_2=I}K\big|Z^{I_1}\Psi_{\alpha}^0\big| \, \big|Z^{I_2}\Psi_{\beta}^0\big|\leq C(n,|I|)K(s/t)^2(t/s)^2
$$
for some constant $K>0$. 
Recall that when $0\leq r\leq t/2$, we have $(t/s)^2\leq 4/3$ so in the region $\{r\leq t/2\} \cap \Kcal$,
$$
\big|Z^I\big(T^{\alphar\betar}\Psi_{\alphar}^0\Psi_{\betar}^0\big)\big|\leq C(n,|I|)K(s/t)^2.
$$

Second, in the region $\{r\geq t/2\}\cap \Kcal$, we have 
$$
\Tu^{00} = T^{\alpha\beta}\big(\Psi^0_{\alpha}\Psi^0_{\beta} - \omega_{\alpha}\omega_{\beta}\big)
$$
and, thus,
$$
Z^I\Tu^{00} = T^{\alpha\beta}Z^I(\Psi^0_{\alpha}\Psi^0_{\beta} - \omega_{\alpha}\omega_{\beta}\big).
$$
When $\alpha=\beta=0$, we have $(\Psi^0_{\alpha}\Psi^0_{\beta} - \omega_{\alpha}\omega_{\beta}\big)=0.$
When $\alpha=a>0, \beta=0$, we have
$$
\aligned
Z^I(\Psi^0_{\alpha}\Psi^0_{\beta} - \omega_{\alpha}\omega_{\beta}\big)
&= - Z^I\big(\omega_a\big(1-(r/t)\big)\big)
\\
&= - \sum_{I_1+I_2=I}Z^{I_1}\omega_a\,Z^{I_2}\bigg(\frac{t-r}{t}\bigg).
\endaligned
$$
When $\alpha = a>0,\,\beta=0>0$, we obtain
$$
\aligned
Z^I\big(\Psi^0_{\alpha}\Psi^0_{\beta} - \omega_{\alpha}\omega_{\beta}\big)
&= Z^I\bigg(\frac{x^ax^b}{t^2} - \frac{x^ax^b}{r^2}\bigg)
\\
&=\sum_{I_1+I_2+I_3+I_4=I}
Z^{I_1}\omega_a\,Z^{I_2}\omega_b \,Z^{I_3}\bigg(1+\frac{r}{t}\bigg)\,Z^{I_4}\bigg(\frac{r-t}{t}\bigg).
\endaligned
$$

We focus on the estimates of $Z^I\omega_a$, $Z^I\big(1+\frac{r}{t}\big)$ and $Z^I\big(\frac{t-r}{t}\big)$.
By Lemma \ref{pre lem homor}, $Z^I\omega_a$ and $Z^I\big(1+\frac{r}{t}\big)$ are bounded by $C(n,|I|)$. For the estimate of $Z^I\big(\frac{t-r}{t}\big)$, we write 
$$
\frac{t-r}{t} = \frac{s^2}{t^2}\frac{t}{t+r}
$$
and then
$$
Z^I\big((t-r)/t\big) = \sum_{I_1+I_2+I_3=I}Z^{I_1}\big(t/(t+r)\big)Z^{I_2}(s/t)Z^{I_3}(s/t).
$$
We observe that $t/(t+r)$ is smooth in $\{t\geq1\}\cap\{t/2\leq r\leq t\}$
 and is homogeneous of degree $0$. We have
$$
\big|Z^{I_1}\big(t/(t+r)\big)\big|\leq C(|I|)
$$
and the term $Z^{I_2}(s/t)$ is bounded by $C(|I|)(s/t)$, thanks to 
Lemma \ref{pre lem lem commutator s/t}. We conclude with
$$
\big|Z^I\big((t/r)/t\big)\big|\leq C(|I|)(s/t)^2. 
$$ 
\end{proof}

Next, we have the following result concerning null quadratic forms.

\begin{proposition}[Estimate of null forms] 
\label{pre lem null 1}
For all null quadratic form $T^{\alpha\beta}\del_\alpha u\del_{\beta}v$
with constant coefficients $T^{\alpha\beta}$ and for any index $I$, one has 
$$
\aligned
\big|Z^I\big(T^{\alpha\beta}\del_\alpha u\del_{\beta}v\big)\big|
& \leq CK(s/t)^2\sum_{|I_1|+|I_2|\leq|I|}\big|Z^{I_1}\del_tuZ^{I_2}\del_tv\big|
\\
&\quad+ CK\sum_{a,\beta,\atop |I_1|+|I_2|\leq|I|}\Big(\big|Z^{I_1}\delu_au\,Z^{I_2}\delu_{\beta}v\big| + \big|Z^{I_1}\delu_{\beta}u\,Z^{I_2}\delu_av\big| \Big),
\endaligned
$$
with $K = \max_{\alpha,\beta}\big|T^{\alpha\beta}\big|$.
\end{proposition}

The importance of this estimate lies on the factor $(s/t)^2$ in front of the component $Z^{I_1}\del_tuZ^{I_2}\del_tv$. 
As we will see later, the derivative $\delu_a$ enjoys better $L^\infty$ and $L^2$ estimates in our framework. The derivatives of direction $\del_t$ do not always have enough decay, and the factor $(s/t)^2$ precisely allows us to overcome
this potential lack of $L^2$ and $L^\infty$ decay.

\begin{proof}
The proof is based on a change of frame. In the semi-hyperboloidal frame we have
$$
\aligned
T^{\alpha\beta}\del_\alpha u\del_{\beta}v
&= \Tu^{\alpha\beta}\delu_\alpha u\delu_{\beta}v
\\
&=\Tu^{00}\del_tu\del_tv + \Tu^{0b}\del_tu\delu_bv + \Tu^{a0}\delu_au\del_tv + \Tu^{ab}\delu_au\delu_bv.
\endaligned
$$
We have
$$
\aligned
&Z^I\big(T^{\alpha\beta}\del_\alpha u\del_{\beta}v\big)
\\
&= Z^I\big(\Tu^{\alpha\beta}\delu_\alpha u\delu_{\beta}v\big)
\\
&= Z^I\big(\Tu^{00}\del_tu\del_tv\big) + Z^I\big(\Tu^{0b}\del_tu\delu_b v\big) + Z^I\big(\Tu^{a0}\delu_au\del_t v\big)
 + Z^I\big(\Tu^{ab}\delu_au\delu_bv\big)
\\
&=: R_1 + R_2 + R_3 + R_4.
\endaligned
$$
Recalling the null condition satisfied by the null form under consideration and by Proposition~\ref{pre lem null quadratic}, we get
$$
\aligned
|R_1|
\leq & \sum_{I_1+I_2+I_3=I}\big|Z^{I_3}\big(\Tu^{00}\big)\, Z^{I_1}\big(\del_tu\big) \,  Z^{I_2}\big(\del_tv\big) \big|
\\
\leq & CK(s/t)^2 \sum_{|I_1|+|I_2|\leq |I|}\big|Z^{I_1}\del_tu\big|\,  \big|Z^{I_2}\del_tv\big|.
\endaligned
$$
The term $R_2$ are estimated directly. Recalling that by Lemma \ref{pre lem frame}, $|Z^I\Tu^{\alpha\beta}|\leq C(I)K$, we find  
$$
\aligned
|R_2|
&\leq \sum_{b \atop I_1+I_2+I_3=I}\big|Z^{I_3}\big(\Tu^{0b}\big)\, Z^{I_1}\big(\del_tu\big)\,  Z^{I_2}\big(\delu_bv\big) \big|
\\
&\leq CK \sum_{b \atop |I_1|+|I_2|\leq |I|}\big|Z^{I_1}\big(\del_tu\big)\,  Z^{I_2}\big(\delu_bv\big) \big|.
\endaligned
$$
The terms $R_3$ and $R_4$ are estimated similarly and the proof is complete. 
\end{proof}


\section{Estimating second-order derivatives}
\label{sec:42}

We can also deal with second-order derivatives
$B^{\alpha\beta}\del_\alpha \del_{\beta} u$,
with constant $B^{\alpha\beta}$. Recall that 
a second-order operator is said to satisfy the null condition
\index{null}
if
$$
B^{\alpha\beta}\xi_\alpha \xi_{\beta} = 0 \quad \text{ when }  \xi_0\xi_0 - \sum_{i=1}^3\xi_i\xi_i = 0.
$$

\begin{proposition}\label{pre lem null 2order0}
Let $B^{\alpha\beta}\del_\alpha \del_{\beta}$ be
a second-order operator satisfying the null condition with constants $B^{\alpha\beta}$ bounded by $K$. One has 
$$
\aligned
& \big|Z^I\big(B^{\alpha\beta}\del_\alpha \del_{\beta}u\big)\big|
\\
& \leq C(n,|I|)K(s/t)^2\sum_{|I_1|\leq |I|}Z^I\big(\del_t\del_t u\big) + C(n,|I|)K\sum_{a,\alpha \atop |I_1|\leq |I|}Z^I\big(\delu_{a}\del_\alpha u\big)
\\
&\quad + C(n,|I|)\frac{K}{t}\sum_{\alpha,|I_1|\leq |I|}|Z^{I_1}\big(\del_\alpha u\big)|.
\endaligned
$$
\end{proposition}

\begin{proof} We have
$$
\aligned
B^{\alpha\beta}\del_\alpha \del_{\beta}u
=& \Bu^{\alpha\beta}\delu_\alpha \delu_{\beta}u
+ B^{\alpha\beta}\del_\alpha \big(\Psi_{\beta}^{\beta'}\big)\delu_{\beta'}u
\\
=& \Bu^{00}\del_t\del_tu + \Bu^{a0}\delu_a\del_tu + \Bu^{0b}\del_t\delu_b u + \Bu^{ab}\delu_a\delu_b u
\\
&
+  B^{\alpha\beta}\del_\alpha \big(\Psi_{\beta}^{\beta'}\big)\delu_{\beta'}u.
\endaligned
$$
Recall that, for a null quadratic form, $|Z^{J}\Bu^{00}|\leq C(n,|J|)K(s/t)^2$ and
$$
\aligned
\big|Z^{I}\Bu^{00}\del_t\del_tu \big|
&\leq \sum_{I_1+I_2=I}\big|Z^{I_2}\Bu^{00}\big|\, \big|Z^{I_1}\del_t\del_tu\big|
\\
&\leq C(n,|I|)K(s/t)^2\sum_{|I_1|\leq|I|}\big|Z^{I_1}\del_t\del_tu\big|.
\endaligned
$$
Also, from Lemma \ref{pre lem frame}, we have
 $\big|Z^I\Bu\big|\leq C(n,|I|)|B|$ and 
$$
\aligned
\big|Z^I\big(\Bu^{a0}\delu_a\del_tu\big)\big|
& \leq \sum_{a\atop I_1+I_2=I}\big|Z^{I_2}\Bu^{a0}\big|\,\big|Z^{I_1}\delu_a\del_tu\big|
\\
& \leq C(n,|I|)K\sum_{a\atop |I_1|\leq |I|}\big|Z^{I_1}\delu_a\del_tu\big|.
\endaligned
$$
For the term $Z^I\big(\Bu^{0b}\del_t\delu_b u\big)$ by applying \eqref{pre commutator base} we have
$$
\aligned
\big|Z^I\big(\Bu^{0b}\del_t\delu_bu\big)\big|
& \leq \sum_{I_1+I_2=I}\big|Z^{I_2}\Bu^{b0}\big|\,\big|Z^{I_1}\del_t\delu_bu\big|
\\
&\leq C(n,|I|)K\sum_{b\atop|I_1|\leq |I|}\big|Z^{I_1}\del_t\delu_bu\big|
\endaligned
$$
and
$$
\aligned
\big|Z^{I_1}\del_t\delu_bu\big|\leq& \big|Z^{I_1}\delu_b\del_t u\big| + \big|Z^{I_1}[\delu_b,\del_t]u\big|
\\
\leq &\big|Z^{I_1}\delu_b\del_t u\big| +\big|Z^{I_1}\big(t^{-1}\Gammau_{0b}^{\gamma}\del_{\gamma}u\big)\big|
\\
\leq &\big|Z^{I_1}\delu_b\del_t u\big| +\sum_{I_3+I_4=I_2}\big|Z^{I_4}\big(t^{-1}\Gammau_{0b}^{\gamma}\big)Z^{I_3}\del_{\gamma}u\big)\big|
\\
\leq &\big|Z^{I_1}\delu_b\del_t u\big| + Ct^{-1}\sum_{\gamma\atop |I_3|\leq|I_1|}\big|Z^{I_3}\del_{\gamma}u\big|.
\endaligned
$$
Here, we have used the estimate $\big|Z^J\big(t^{-1}\Gammau_{0b}^\gamma\big)\big|\leq C(n,|J|)t^{-1}$ because $t^{-1}\Gammau_{0b}^\gamma$ is homogeneous of degree $-1$. So we have proven that
$$
\big|Z^{I_1}\Bu^{0b}\del_t\delu_bu\big|
\leq
 C(n,|I|)K\sum_{b\atop|I_1|\leq |I|}\big|Z^{I_1}\delu_b\del_t u\big|
+ C(n,|I|)Kt^{-1}\sum_{\gamma\atop |I_1|\leq |I|}\big|Z^{I_1}\del_{\gamma}u\big|.
$$

The term $Z^I\big(\Bu^{ab}\delu_a\delu_bu\big)$ is estimated as follows:
$$
\aligned
\big|Z^I\big(\Bu^{ab}\delu_a\delu_bu\big)\big|
& \leq C(n,|I|)\sum_{a,b\atop I_1+I_2=I}\big|Z^{I_2}\Bu^{ab}\big|\,\big|Z^{I_1}\delu_a\delu_bu\big|
\\
& \leq C(n,|I|)K\sum_{\beta,a\atop I_1+I_2=I}\big|Z^{I_1}\delu_a\delu_{\beta}u\big|.
\endaligned
$$
The term $Z^I\big(B^{\alpha\beta}\del_\alpha \big(\Psi_{\beta}^{\beta'}\big)\delu_{\beta'}u\big)$ is estimated by applying the decay rate supplied by $Z^{J}\big|\del_\alpha \big(\Psi_{\beta}^{\beta'}\big)\big|\leq C(n,|J|)t^{-1}$. This completes the proof.
\end{proof}

As a direct corollary, the following estimates hold.

\begin{proposition}\label{pre lem null 2order1}  
Consider a bilinear form  
$B^{\alpha\beta}u\del_{\alpha}\del_{\beta}v$
acting on the function $u$ and the Hessian of $v$
and suppose that the quadratic form $B^{\alpha\beta}\del_{\alpha}\del_{\beta}$ satisfies the null condition.
The following estimates hold:
\begin{subequations}\label{pre null B}
\bel{pre null B1}
\aligned
\big|Z^{I}\big(B^{\alpha\beta}u\del_{\alpha}\del_{\beta}v\big)\big|
\leq
\, 
& C(n,|I|)K(s/t)^2\sum_{|I_1|+|I_2|\leq|I|}\big|Z^{I_1}u\big|\,\big|Z^{I_2}\del_t\del_t v\big|
\\
   &+ C(n,|I|)K\sum_{\alpha,b\atop |I_1|+|I_2|\leq|I|}\big|Z^{I_1}u\big|\,\big|Z^{I_2}\delu_{\alpha}\delu_b v\big|
\\
&    + C(n,|I|)K\sum_{a,\beta\atop |I_1|+|I_2|\leq|I|}\big|Z^{I_1}u\big|\,\big|Z^{I_2}\delu_a\delu_{\beta} v\big|
\\
   &+ C(n,|I|)Kt^{-1}\sum_{\alpha\atop |I_1|+|I_2|\leq|I|}\big|Z^{I_1}u\big|\,\big|Z^{I_2}\delu_{\alpha}v\big|,
\endaligned
\ee
\bel{pre null B2}
\aligned
\big|[Z^I,B^{\alpha\beta}u\del_{\alpha}\del_{\beta}]v\big|
\leq &C(n,|I|)K(s/t)^2\sum_{|I_1|+|I_2|\leq |I|\atop |I_2|<|I|}\big|Z^{I_1}u\big|\,\big|Z^{I_2}\del_t\del_tv\big|
\\
&+C(n,|I|)K\sum_{a,\beta}\sum_{|I_1|+|I_2|\leq |I|\atop |I_2|<|I|}\big|Z^{I_1}u\big|\,\big|Z^{I_2}\delu_a\delu_{\beta}v\big|
\\
&
 +C(n,|I|)K\sum_{\alpha,b}\sum_{|I_1|+|I_2|\leq |I|\atop |I_2|<|I|}\big|Z^{I_1}u\big|\,\big|Z^{I_2}\delu_{\alpha}\delu_bv\big|
\\
&+C(n,|I|)Kt^{-1}\sum_{|I_1|+|I_2|\leq |I|\atop |I_2|<|I|,\beta}\big|Z^{I_1}u\big|\,\big|Z^{I_2}\del_{\beta'}v\big|.
\endaligned
\ee
\end{subequations}
\end{proposition}

\begin{proof} We first establish \eqref{pre null B1}:
$$
\aligned
Z^I\big(B^{\alpha\beta}u\del_{\alpha}\del_{\beta}v\big)
=&Z^I\big(\Bu^{\alpha\beta}u\delu_{\alpha}\delu_{\beta}v + B^{\alpha\beta}u\del_{\alpha}\Psi_{\beta}^{\beta'}\delu_{\beta'}v\big)
\\
=&\sum_{I_1+I_2+I_3=I}Z^{I_3}\Bu^{\alpha\beta}\,Z^{I_1}u\,Z^{I_2}\delu_{\alpha}\delu_{\beta}v
\\
&
 +\sum_{I_1+I_2+I_3=I}B^{\alpha\beta}Z^{I_1}u\,Z^{I_3}\del_{\alpha}\Psi_{\beta}^{\beta'}\,Z^{I_2}\delu_{\beta'}v
=:R_1 + R_2,
\endaligned
$$
in which $R_1$ can be estimated as follows:
$$
\aligned
R_1 = &\sum_{I_1+I_2+I_3=I}Z^{I_3}\Bu^{\alpha\beta}\,Z^{I_1}u\,Z^{I_2}\delu_{\alpha}\delu_{\beta}v
\\
=     &\sum_{I_1+I_2+I_3=I}Z^{I_3}\Bu^{00}\,Z^{I_1}u\,Z^{I_2}\del_t\del_t v
\\
     &+\sum_{I_1+I_2+I_3=I}Z^{I_3}\Bu^{a0}\,Z^{I_1}u\,Z^{I_2}\delu_a\del_t v
\\
&
      +\sum_{I_1+I_2+I_3=I}Z^{I_3}\Bu^{0b}\,Z^{I_1}u\,Z^{I_2}\del_t\delu_b v
\\
     &+\sum_{I_1+I_2+I_3=I}Z^{I_3}\Bu^{ab}\,Z^{I_1}u\,Z^{I_2}\delu_a\delu_b v.
\endaligned
$$
Since $\big|Z^I\Bu^{\alpha\beta}\big|\leq C(n,|I|)K$ and $\big|Z^I\Bu^{00}\big|\leq C(n,|I|)(s/t)^2$, we have 
$$
\aligned
|R_1|
\leq &C(n,|I|)K(s/t)^2\sum_{|I_1|+|I_2|\leq|I|}\big|Z^{I_1}u\big|\,\big|Z^{I_2}\del_t\del_t v\big|
\\
         &+C(n,|I|)K\sum_{\alpha,b\atop |I_1|+|I_2|\leq|I|}\big|Z^{I_1}u\big|\,\big|Z^{I_2}\delu_{\alpha}\delu_b v\big|
 \\
&         +C(n,|I|)K\sum_{a,\beta\atop |I_1|+|I_2|\leq|I|}\big|Z^{I_1}u\big|\,\big|Z^{I_2}\delu_a\delu_{\beta} v\big|.
\endaligned
$$
Observe now that $\big|Z^I\del_{\alpha}\Psi_{\beta}^{\beta'}\big|\leq C(n,|I|)t^{-1}$, then $R_2$ is bounded by
$$
\aligned
|R_2|
& \leq
 \sum_{I_1+I_2+I_3=I}\big|B^{\alpha\beta}\big|\,\big|Z^{I_1}u\big|\,\big|Z^{I_3}\del_{\alpha}\Psi_{\beta}^{\beta'}\big|\,
 \big|Z^{I_2}\delu_{\beta'}v\big|
\\
&\leq C(n,|I|)Kt^{-1}\sum_{\alpha\atop |I_1|+|I_2|\leq|I|}\big|Z^{I_1}u\big|\,\big|Z^{I_2}\delu_{\alpha}v\big|.
\endaligned
$$
In view of 
$Z^{I_2}\delu_{\alpha}v = Z^{I_2}\big(\Psi_{\alpha}^{\alpha'}\del_{\alpha'}v\big) = Z^{I_4}\Psi_{\alpha}^{\alpha'}Z^{I_5}\del_{\alpha'}v$
and that $Z^{I_4}\Psi_{\alpha}^{\alpha'}$ is bounded, 
the estimate \eqref{pre null B1} is proven.

Next, in order to derive \eqref{pre null B2}, we observe that
$$
[Z^I,B^{\alpha\beta}u\del_{\alpha}\del_{\beta}]v
= [Z^I,\Bu^{\alpha\beta}u\delu_{\alpha}\delu_{\beta}]v + [Z^I,B^{\alpha\beta}u\del_{\alpha}\Psi_{\beta}^{\beta'}\delu_{\beta'}]v
=: R_3+R_4.
$$
The term
$R_3$ is decomposed as follows:
$$
\aligned
&[Z^I,\Bu^{\alpha\beta}u\delu_{\alpha}\delu_{\beta}]v
\\
&=\sum_{I_1+I_2+I_3\atop |I_2|<|I|}Z^{I_3}\Bu^{\alpha\beta}\,Z^{I_1}u\,Z^{I_2}\delu_{\alpha}\delu_{\beta}v
 +\Bu^{\alpha\beta}u[Z^I,\delu_{\alpha}\delu_{\beta}]v
\\
&=:R_5 + R_6.
\endaligned
$$
The term $R_5$ is estimated as follows:
$$
\aligned
R_5 =& \sum_{I_1+I_2+I_3\atop |I_2|<|I|}Z^{I_3}\Bu^{\alpha\beta}\,Z^{I_1}u\,Z^{I_2}\delu_{\alpha}\delu_{\beta}v
\\
=&\sum_{I_1+I_2+I_3\atop |I_2|<|I|}Z^{I_3}\Bu^{00}\,Z^{I_1}u\,Z^{I_2}\del_t\del_tv
\\
&+\sum_{I_1+I_2+I_3\atop |I_2|<|I|}Z^{I_3}\Bu^{a0}\,Z^{I_1}u\,Z^{I_2}\delu_a\del_tv
 +\sum_{I_1+I_2+I_3\atop |I_2|<|I|}Z^{I_3}\Bu^{0b}\,Z^{I_1}u\,Z^{I_2}\del_t\delu_bv
\\
&+\sum_{I_1+I_2+I_3\atop |I_2|<|I|}Z^{I_3}\Bu^{ab}\,Z^{I_1}u\,Z^{I_2}\delu_a\delu_bv, 
\endaligned
$$
 thus
$$
\aligned
|R_5|
\leq&\sum_{I_1+I_2+I_3=I\atop |I_2|<|I|}\big|Z^{I_3}\Bu^{00}\big|\,\big|Z^{I_1}u\big|\,\big|Z^{I_2}\del_t\del_tv\big|
\\
&+\sum_{I_1+I_2+I_3=I\atop |I_2|<|I|}\big|Z^{I_3}\Bu^{a0}\big|\,\big|Z^{I_1}u\big|\,\big|Z^{I_2}\delu_a\del_tv\big|
\\
&
 +\sum_{I_1+I_2+I_3=I\atop |I_2|<|I|}\big|Z^{I_3}\Bu^{0b}\big|\,\big|Z^{I_1}u\big|\,\big|Z^{I_2}\del_t\delu_bv\big|
\\
&+\sum_{I_1+I_2+I_3=I\atop |I_2|<|I|}\big|Z^{I_3}\Bu^{ab}\big|\,\big|Z^{I_1}u\big|\,\big|Z^{I_2}\delu_a\delu_bv\big|,
\endaligned
$$
hence
$$
\aligned
|R_5|
\leq& C(n,|I|)K(s/t)^2\sum_{|I_1|+|I_2|\leq|I| \atop |I_2|<|I|}\big|Z^{I_1}u\big|\,\big|Z^{I_2}\del_t\del_tv\big|
\\
&+ C(n,|I|)K\sum_{|I_1|+|I_2|\leq|I| \atop |I_2|<|I|}\big|Z^{I_1}u\big|\,\big|Z^{I_2}\delu_a\delu_{\beta}v\big|
\\
&
 + C(n,|I|)K\sum_{|I_1|+|I_2|\leq|I| \atop |I_2|<|I|}\big|Z^{I_1}u\big|\,\big|Z^{I_2}\delu_{\alpha}\delu_bv\big|.
\endaligned
$$

The term $R_6$ is decomposed as follows:
$$
\aligned
\Bu^{\alpha\beta}u[Z^I,\delu_{\alpha}\delu_{\beta}]v
= & \Bu^{00}u[Z^I,\del_t\del_t]v
+ \Bu^{a0}u[Z^I,\delu_a\del_t]v
\\
&+ \Bu^{0b}u[Z^I,\del_t\delu_b]v + \Bu^{ab}u[Z^I,\delu_a\delu_b]v.
\endaligned
$$
We apply \eqref{pre lem commutator second-order} and \eqref{pre lem commutator second-order bar}
and that $|\Bu^{00}|\leq C(s/t)^2$:
$$
\aligned
&\big|\Bu^{\alpha\beta}u[Z^I,\delu_{\alpha}\delu_{\beta}]v\big|
\\
&= \big|\Bu^{00}u[Z^I,\del_t\del_t]v\big|
+ \big|\Bu^{a0}u[Z^I,\delu_a\del_t]v\big|
\\
& \quad + \big|\Bu^{0b}u[Z^I,\del_t\delu_b]v\big| + \big|\Bu^{ab}u[Z^I,\delu_a\delu_b]v\big|
\\
&\leq C(n,|I|)K(s/t)^2|u|\sum_{\gamma,\gamma'\atop |I'|<|I|}\big|Z^{I_1}\del_{\gamma}\del_{\gamma'}v\big|
\\
    &\quad 
+C(n,|I|)K |u|\sum_{a,\beta\atop |I_1|<|I|}\big|\delu_a\delu_{\beta}Z^{I_1}v\big|
+C(n,|I|)Kt^{-1} |u|\sum_{\gamma\atop |I_1|<|I|}\big|\del_{\gamma}Z^{I_1}v\big|.
\endaligned
$$

The term $R_4$ is estimated as follows:
$$
\aligned
& [Z^I,B^{\alpha\beta}u\del_{\alpha}\Psi_{\beta}^{\beta'}\delu_{\beta'}]v
\\
& =\sum_{I_1+I_2+I_3=I\atop |I_2|<|I|}B^{\alpha\beta}Z^{I_1}u\,Z^{I_3}\del_{\alpha}\Psi_{\beta}^{\beta'}\,Z^{I_2}\delu_{\beta'}v
    +B^{\alpha\beta}u\del_{\alpha}\Psi_{\beta}^{\beta'}[Z^I,\delu_{\beta'}]v.
\endaligned
$$
Thanks to the additional decreasing rate in
 $\big|Z^I\del_{\alpha}\Psi_{\beta}^{\beta'}\big|\leq C(n,|I|)t^{-1}$, the first term is bounded by
$$
C(n,|I|)Kt^{-1}\sum_{|I_1|+|I_2|\leq |I|\atop |I_2|<|I|,\beta'}\big|Z^{I_1}u\big|\,\big|Z^{I_2}\del_{\beta'}v\big|.
$$
The second term is estimated by \eqref{pre lem commutator bar} and the additional decreasing rate supplied by $\big|Z^I\del_{\alpha}\Psi_{\beta}^{\beta'}\big|\leq C(n,|I|)t^{-1}$. It can also be bounded by
$$
C(n,|I|)Kt^{-1}\sum_{|I_1|+|I_2|\leq |I|\atop |I_2|<|I|,\beta}\big|Z^{I_1}u\big|\,\big|Z^{I_2}\del_{\beta'}v\big|.
$$
This completes the proof of \eqref{pre null B2}.
\end{proof}


\section{Products of first-order and second-order derivatives}
\label{sec:43}

The third type of null form we treat is a null quadratic form acting on the gradient and the Hessian, as now stated.

\begin{proposition}\label{pre lem null 2order2}
Consider a quadratic form acting on the gradient of $u$ and the Hessian of $v$, that is, 
$A^{\alpha\beta\gamma}\del_{\gamma}u \del_\alpha \del_{\beta}v$, 
for functions $u,v$ defined in the cone $\Kcal$, and suppose that 
$A^{\alpha\beta\gamma}$ satisfies the null condition.
The following estimate holds for any index $I$:
$$
\aligned
& \big|Z^I\big(A^{\alpha\beta\gamma}\del_{\gamma}u\del_\alpha \del_{\beta}v\big)\big|
\\
&\leq  C(n,|I|)K (s/t)^2\sum_{|I_1|+|I_2|\leq |I|}\big|Z^{I_1}\del_tu\, Z^{I_2}\del_t\del_tv\big|
 \\
 &\quad + C(n,|I|)K\big(\Omega_1(I,u,v) + \Omega_2(I,u,v)\big),
\endaligned
$$
where 
$$
\aligned
\Omega_1(I,u,v) =& \sum_{a,\beta,\gamma\atop|I_1|+|I_2|\leq|I|}
                  \big|Z^{I_1}\delu_a u\big|\,\big|Z^{I_2}\delu_{\beta}\delu_{\gamma}v\big|
                 +\sum_{\alpha,b,\gamma\atop|I_1|+|I_2|\leq|I|}
                  \big|Z^{I_1}\delu_{\alpha}u\big|\,\big|Z^{I_2}\delu_b\delu_{\gamma}v\big|
\\
                 &+\sum_{\alpha\beta,c\atop|I_1|+|I_2|\leq|I|}
                  \big|Z^{I_1}\delu_{\alpha}u\big|\,\big|Z^{I_2}\delu_{\beta}\delu_c v\big|
\endaligned
$$
and
$$
\Omega_2(I,u,v) \leq t^{-1}\sum_{\alpha,\beta,\atop |I_1|+|I_2|\leq |I|}|Z^{I_1}\del_\alpha u\,  Z^{I_2}\del_{\beta}v|.
$$
\end{proposition}

\begin{proof}
We observe the following change of frame formula: 
$$
\aligned
A^{\alpha\beta\gamma}\del_{\gamma}u\del_\alpha \del_{\beta}v
&= \Au^{\alpha\beta\gamma}\delu_{\gamma}u\delu_\alpha \delu_{\beta}v
+ A^{\alpha\beta\gamma}\del_{\gamma}u\del_\alpha \big(\Psi_{\beta}^{\beta'}\big)\delu_{\beta'}u
\\
&=: R_1(u,v) + R_2(u,v).
\endaligned
$$
The term $Z^IR_2$ can be estimated as follows. Recall that $|Z^I \del_\alpha \big(\Psi_{\beta}^{\beta'}\big)|\leq \frac{C}{t}$, then
$$
Z^{I}\big(A^{\alpha\beta\gamma}\del_{\gamma}u\del_\alpha \big(\Psi_{\beta}^{\beta'}\big)\delu_{\beta'}v\big)
= \sum_{I_1+I_2+I_3=I}A^{\alpha\beta\gamma}
Z^{I_1}\del_{\gamma}u\, Z^{I_2}\delu_{\beta'}v\,  Z^{I_3}\del_\alpha \big(\Psi_{\beta}^{\beta'}\big),
$$
which can be estimated as
$$
\big|Z^{I}\big(A^{\alpha\beta\gamma}\del_{\gamma}u\del_\alpha \big(\Psi_{\beta}^{\beta'}\big)\delu_{\beta'}v\big) \big|
\leq CKt^{-1}\sum_{\alpha,\beta\atop |I_1|+|I_2|\leq|I|}|Z^{I_1}\del_\alpha uZ^{I_2}\del_{\beta}v| = \Omega_2.
$$

The term $R_1$ can be estimated as follows:
$$
\aligned
R_1=& \Au^{\alpha\beta\gamma}\delu_{\gamma}u\delu_\alpha \delu_{\beta}v
\\
=& \Au^{000}\del_tu\del_t\del_tv
\\
&+\Au^{00c}\delu_cu\del_t\del_tv
 +\Au^{0b0}\del_tu\del_t\delu_bv
 +\Au^{a00}\del_tu\delu_a\del_tv
\\
&+\Au^{ab0}\del_tu\delu_a\delu_bv
 +\Au^{0bc}\delu_cu\del_t\delu_bv
 +\Au^{a0c}\delu_cu\delu_a\del_tv
 +\Au^{abc}\delu_au\delu_b\del_cv
\\
=:&R_3(u,v)+R_4(u,v), 
\endaligned
$$
where
$$
R_3(u,v) := \Au^{000}\del_tu\del_t\del_tv
$$
and $R_4(u,v)$ denotes the remaining terms. We have $Z^IR_1 = Z^IR_3 + Z^IR_4$. To estimate $|Z^IR_3|$, we observe that $A^{\alpha\beta\gamma}$ is a null cubic form so by Proposition~\ref{pre lem null quadratic}, $|Z^IA^{000}|\leq C(s/t)^2$. We find 
$$
Z^I\big(\Au^{000}\del_tu\del_t\del_tv\big) = \sum_{I_1+I_2+I_3=I}Z^{I_3}A^{000}\,  Z^{I_1}\del_tu\,  Z^{I_2}\del_t\del_tv, 
$$
so that
$$
\big|Z^I\big(\Au^{000}\del_tu\del_t\del_tv\big)\big|
\leq C(n,|I|)K(s/t)^2\sum_{|I_1|+|I_2|\leq|I|}|Z^{I_1}\del_tu\,  Z^{I_2}\del_t\del_tv|. 
$$
To see the estimates on $Z^IR_4$ terms, we just remark that by Lemma \ref{pre lem frame}, $|Z^I\Au^{\alpha\beta\gamma}|\leq C(n,I)K$. 
We can control $Z^IR_3$ by $\Omega_1$.

Finally, by combining the estimate on $R_2, R_3$ and $R_4$, the desired result is established.
\end{proof}


The fourth type of terms we need to control is the commutator between $Z^I$ 
and a null quadratic form.

\begin{proposition}\label{pre lem null 2order3}
Consider a null quadratic form acting on the gradient of function $u$ and 
on the Hessian matrix of $v$, that is, 
$A^{\alpha\beta\gamma}\del_{\gamma}u\del_{\alpha}\del_{\beta}v$, 
for functions $u, v$ defined in the cone $\Kcal$. The following estimate holds:
$$
\aligned
\big|[Z^I,A^{\alpha\beta\gamma}\del_{\gamma}u\del_\alpha \del_{\beta}]v\big|
\leq& CK(s/t)^2\sum_{|I_2|+|I_3|\leq |I|,\atop |I_3|<|I|}\big|Z^{I_2}\del_t uZ^{I_3}\del_t\del_tv\big|
\\
&
       +     CK\sum_{c,\alpha,\beta}\sum_{|I_2|+|I_3|\leq |I|,\atop |I_3|<|I|}
                    \big|Z^{I_2}\delu_c uZ^{I_3}\delu_\alpha \delu_{\beta}v\big|
\\
       &+    CK\sum_{c,\alpha,\beta}\sum_{ |I_2|+|I_3|\leq |I|,\atop |I_3|<|I|}
                    \big|Z^{I_2}\delu_\alpha  uZ^{I_3}\delu_c\delu_{\beta}v\big|
  \\
&      +    CK\sum_{c,\alpha,\beta}\sum_{|I_2|+|I_3|\leq |I|,\atop |I_3|<|I|}
                    \big|Z^{I_2}\delu_\alpha  uZ^{I_3}\delu_{\beta}\delu_cv\big|
\\
       &+t^{-1}CK
       \sum_{\alpha,\gamma}\sum_{|I_2|\leq |I|-1\atop |I_1|+|I_2|\leq |I|}\big|\del_{\gamma}Z^{I_1}u\del_\alpha Z^{I_2}v\big|.
\endaligned
$$
\end{proposition}

\begin{proof}
We use the change of frame formula
$$
A^{\alpha\beta\gamma}\del_{\gamma}u\del_\alpha \del_{\beta}v
= \Au^{\alpha\beta\gamma}\delu_{\gamma}u\, \delu_\alpha \delu_{\beta}v
  +A^{\alpha\beta\gamma}\del_{\gamma}u\, \del_\alpha \big(\Psi_{\beta}^{\beta'}\big)\delu_{\beta'}v
$$
and find
$$
\aligned
&[Z^I,A^{\alpha\beta\gamma}\del_{\gamma}u\del_\alpha \del_{\beta}]v
\\
&= [Z^I,\Au^{\alpha\beta\gamma}\delu_{\gamma}u\delu_\alpha \delu_{\beta}]v + [Z^I,A^{\alpha\beta\gamma}\del_{\gamma}u\del_\alpha \big(\Psi_{\beta}^{\beta'}\big)\delu_{\beta'}]v
\\
&=:  R_1(I,u,v) + R_2(I,u,v).
\endaligned
$$

We first decompose $R_1(I,u,v)$ as
$$
\aligned
&[Z^I,\Au^{\alpha\beta\gamma}\delu_{\gamma}u\delu_\alpha \delu_{\beta}]v
\\
&=\sum_{I_2+I_3+I_4=I,\atop |I_3|<|I|}
Z^{I_4}\Au^{\alpha\beta\gamma} Z^{I_2}\big(\delu_{\gamma}u\big)Z^{I_3}\big(\delu_\alpha \delu_{\beta}v\big)
+ \Au^{\alpha\beta\gamma}\delu_{\gamma}u[Z^I,\delu_\alpha \delu_{\beta}]v
\\
&=: R_3(I,u,v) + R_4(I,u,v),
\endaligned
$$
while
$R_3(I,u,v)$ is decomposed as
$$
\aligned
&R_3(I,u,v)
\\
&= \sum_{I_2+I_3+I_4=I,\atop |I_3|<|I|}Z^{I_4}\Au^{000}Z^{I_2}\del_tuZ^{I_3}\del_t\del_tv
\\
           &
\quad 
+ \sum_{I_2+I_3+I_4=I,\atop |I_3|<|I|}Z^{I_4}\Au^{00c}Z^{I_2}\delu_cuZ^{I_3}\del_t\del_tv
            + \sum_{I_2+I_3+I_4=I,\atop |I_3|<|I|}Z^{I_4}\Au^{0b0}Z^{I_2}\del_tuZ^{I_3}\del_t\delu_cv
\\
           &\quad + \sum_{I_2+I_3+I_4=I,\atop |I_3|<|I|}Z^{I_4}\Au^{a00}Z^{I_2}\del_tuZ^{I_3}\delu_a\del_tv
            + \sum_{I_2+I_3+I_4=I,\atop |I_3|<|I|}Z^{I_4}\Au^{0bc}Z^{I_2}\delu_cuZ^{I_3}\del_t\delu_bv
\\
           &\quad + \sum_{I_2+I_3+I_4=I,\atop |I_3|<|I|}Z^{I_4}\Au^{a0c}Z^{I_2}\delu_cuZ^{I_3}\delu_a\del_tv
            + \sum_{I_2+I_3+I_4=I,\atop |I_3|<|I|}Z^{I_4}\Au^{ab0}Z^{I_2}\del_tuZ^{I_3}\delu_a\delu_bv
\\
           &\quad + \sum_{I_2+I_3+I_4=I,\atop |I_3|<|I|}Z^{I_4}\Au^{abc}Z^{I_2}\delu_cuZ^{I_3}\delu_a\delu_bv.
\endaligned
$$
We observe that $\mathcal{A}$ is a null cubic form and, by Proposition~\ref{pre lem null quadratic}, and therefore 
$$
\big|Z^{I_4}\Au^{000}\big|\leq C(|I|)K(s/t)^2.
$$

By Lemma \ref{pre lem frame}, we have
$$
\big|Z^{I_4}\Au^{\alpha\beta\gamma}\big|\leq K.
$$
The term $R_3(I,u,v)$ is estimated as follows:
$$
\aligned
&\big|R_3(I,u,v)\big|
    \\
&\leq CK(s/t)^2\sum_{|I_2|+|I_3|\leq |I|,\atop |I_3|<|I|}\big|Z^{I_2}\del_t uZ^{I_3}\del_t\del_tv\big|
\\
           &\quad+ CK\sum_{|I_2|+|I_3|\leq |I|,\atop |I_3|<|I|}\big|Z^{I_2}\delu_c uZ^{I_3}\del_t\del_tv\big|
 + CK\sum_{|I_2|+|I_3|\leq |I|,\atop |I_3|<|I|}\big|Z^{I_2}\del_t uZ^{I_3}\del_t\delu_cv\big|
\\
           &\quad+ CK\sum_{|I_2|+|I_3|\leq |I|,\atop |I_3|<|I|}\big|Z^{I_2}\del_t uZ^{I_3}\delu_a\del_tv\big|
  + CK\sum_{|I_2|+|I_3|\leq |I|,\atop |I_3|<|I|}\big|Z^{I_2}\delu_cu Z^{I_3}\del_t\delu_bv\big|
\\
           &\quad+ CK\sum_{|I_2|+|I_3|\leq |I|,\atop |I_3|<|I|}\big|Z^{I_2}\delu_cuZ^{I_3}\delu_a\del_tv\big|
     + CK\sum_{|I_2|+|I_3|\leq |I|,\atop |I_3|<|I|}\big|Z^{I_2}\del_tuZ^{I_3}\delu_a\delu_bv\big|
\\
           &\quad+ CK\sum_{|I_2|+|I_3|\leq |I|,\atop |I_3|<|I|}\big|Z^{I_2}\delu_cuZ^{I_3}\delu_a\delu_bv\big|,
\endaligned
$$
so that
$$
\aligned
\big|R_3(I,u,v)\big|
 \leq& CK(s/t)^2\sum_{|I_2|+|I_3|\leq |I|,\atop |I_3|<|I|}\big|Z^{I_2}\del_t uZ^{I_3}\del_t\del_tv\big|
    \\
&   +     CK\sum_{c,\alpha,\beta}\sum_{|I_2|+|I_3|\leq |I|,\atop |I_3|<|I|}
                    \big|Z^{I_2}\delu_c uZ^{I_3}\delu_\alpha \delu_{\beta}v\big|
\\
       &+    CK\sum_{c,\alpha,\beta}\sum_{ |I_2|+|I_3|\leq |I|,\atop |I_3|<|I|}
                    \big|Z^{I_2}\delu_\alpha  uZ^{I_3}\delu_c\delu_{\beta}v\big|
   \\
&     +    CK\sum_{c,\alpha,\beta}\sum_{|I_2|+|I_3|\leq |I|,\atop |I_3|<|I|}
                    \big|Z^{I_2}\delu_\alpha  uZ^{I_3}\delu_{\beta}\delu_cv\big|.
\endaligned
$$

The term $R_4(I,u,v)$ is also decomposed as:
$$
\aligned
R_4(I,u,v)=& \Au^{000}\del_tu[Z^I,\del_t\del_t]v+ \Au^{00c}\delu_cu[Z^I,\del_t\del_t]v
\\
           &
            + \Au^{0b0}\del_tu[Z^I,\del_t\delu_b]v
            + \Au^{a00}\del_tu[Z^I,\delu_a\del_t]v
\\
           &+ \Au^{0bc}\delu_cu[Z^I,\del_t\delu_b]v
            + \Au^{a0c}\delu_cu[Z^I,\delu_a\del_t]v
\\
           & + \Au^{ab0}\del_tu[Z^I,\delu_a\delu_b]v + \Au^{abc}\delu_cu[Z^I,\delu_a\delu_b]v.
\endaligned
$$
Thanks to the null condition, we have
$$
\aligned
\big|R_4(I,u,v)\big|\leq& CK(s/t)^2\del_tu\big|[Z^I,\del_t\del_t]v\big| + CK\big|\delu_cu[Z^I,\del_t\del_t]v\big|
\\
           &
            + CK\big|\del_tu[Z^I,\del_t\delu_b]v\big|
            + CK\big|\del_tu[Z^I,\delu_a\del_t]v\big|
\\
           &+ CK\big|\delu_cu[Z^I,\del_t\delu_b]v\big|
            + CK\big|\delu_cu[Z^I,\delu_a\del_t]v\big|
\\
           &
            + CK\big|\del_tu[Z^I,\delu_a\delu_b]v\big|+ CK\big|\delu_cu[Z^I,\delu_a\delu_b]v\big|.
\endaligned
$$
Observe now that thanks to the commutator estimates
 \eqref{pre lem commutator second-order} and  \eqref{pre lem commutator second-order bar}, we have
$$
\big|\delu_\alpha u[Z^I,\del_t\del_t]v\big|\leq C\sum_{\alpha,\beta\atop|I'|\leq |I|-1}\big|\delu_\alpha u \del_\alpha \del_{\beta}Z^{I'}v\big|,
$$
$$
\big|\delu_\alpha u[Z^I,\delu_a\delu_{\beta}]v\big|\leq
C\sum_{a,\beta\atop |I'|\leq |I|-1}\big|\delu_\alpha u \delu_a\delu_{\beta}Z^{I'}v\big|
+ Ct^{-1}\sum_{\gamma\atop |I'|<|I|}\big|\delu_\alpha u\del_{\gamma}Z^{I'}v\big|
$$
and
$$
\big|\delu_\alpha u[Z^I,\delu_{\beta}\delu_c]v\big|\leq
C\sum_{\beta,c\atop |I'|\leq |I|-1}\big|\delu_\alpha u \delu_c \delu_{\beta} Z^{I'}v\big|
+ Ct^{-1}\sum_{\gamma\atop |I'|<|I|}\big|\delu_\alpha u \del_{\gamma}Z^{I'}v\big|.
$$
We also note that $\sum_\alpha \big|\delu_\alpha u\big|\leq C\sum_\alpha \big|\del_\alpha u\big|$,
so that $|R_4(I,u,v)|$ is bounded by
$$
\aligned
&\big|R_4(I,u,v)\big|
\\
       \leq& CK(s/t)^2\sum_{\alpha,\beta\atop |I'|\leq |I|-1}\big|\del_tu \del_\alpha \del_{\beta}Z^{I'}v\big|
+ CK\sum_{\alpha,\beta,c\atop|I'|\leq |I|-1}\big|\delu_cu  \del_\alpha \del_{\beta}Z^{I'}v\big|
\\
           &
            + CK\sum_{\alpha,\beta,c\atop|I'|\leq |I|-1}\big|\del_\alpha u \delu_c\delu_{\beta}v\big|
            + CKt^{-1}\sum_{\alpha,\gamma\atop |I'|\leq |I|-1}\big|\del_{\gamma}u\del_\alpha Z^{I'}v\big|.
\endaligned
$$
In conclusion, $R_1(I,u,v)$ is bounded by
$$
\aligned
\big|R_1(I,u,v)\big|
   \leq
&CK(s/t)^2\sum_{|I_2|+|I_3|\leq |I|,\atop |I_3|<|I|}\big|Z^{I_2}\del_t uZ^{I_3}\del_t\del_tv\big|
  \\
&     +     CK\sum_{c,\alpha,\beta}\sum_{|I_2|+|I_3|\leq |I|,\atop |I_3|<|I|}
                    \big|Z^{I_2}\delu_c uZ^{I_3}\delu_\alpha \delu_{\beta}v\big|
\\
       &+    CK\sum_{c,\alpha,\beta}\sum_{ |I_2|+|I_3|\leq |I|,\atop |I_3|<|I|}
                    \big|Z^{I_2}\delu_\alpha  uZ^{I_3}\delu_c\delu_{\beta}v\big|
        \\
&+    CK\sum_{c,\alpha,\beta}\sum_{|I_2|+|I_3|\leq |I|,\atop |I_3|<|I|}
                    \big|Z^{I_2}\delu_\alpha  uZ^{I_3}\delu_{\beta}\delu_cv\big|
\\
       & +t^{-1}CK\sum_{\alpha,\gamma\atop |I'|\leq |I|-1}\big|\del_{\gamma}u\del_\alpha Z^{I'}v\big|.
\endaligned
$$

Next, we turn our attention to the estimate of $R_2(I,u,v)$:
$$
\aligned
R_2(I,u,v) =
& A^{\alpha\beta\gamma}\del_{\gamma}u\,  \del_\alpha \Psi_{\beta}^{\beta'} [Z^I,\delu_{\beta'}]v
\\
&
+\sum_{I_1+I_2+I_3=I,\atop |I_2|<|I|}
 A^{\alpha\beta\gamma}Z^{I_1}\del_{\gamma}u\,  Z^{I_2}\delu_{\beta'}v\,  Z^{I_3}\del_\alpha \Psi_{\beta}^{\beta'}
\\
=: & R_5(I,u,v) + R_6(I,u,v).
\endaligned
$$
To estimate the term $R_5(I,u,v)$, we recall \eqref{pre lem commutator pr2} and the fact that $\big|\del_\alpha \Psi_{\beta}^{\beta'}\big|\leq Ct^{-1}$, so that
$$
\big| A^{\alpha\beta\gamma}\del_{\gamma}u\,  \del_\alpha \Psi_{\beta}^{\beta'} [Z^I,\delu_{\beta'}]v\big|\leq
CKt^{-1}\sum_{\beta,\gamma\atop |I'|\leq |I|-1}\big|\del_{\gamma}u \del_{\beta}Z^{I'}v\big|.
$$
In the same way, for $R_6(I,u,v)$ we have
$$
R_6(I,u,v)\leq
CKt^{-1}\sum_{\beta,\gamma}\sum_{|I_1|+|I_2|\leq |I|\atop |I_2|<|I|}\big|\del_{\beta}Z^{I_1}u\del_{\gamma}Z^{I_2}v\big|.
$$
This establishes the desired result.
\end{proof}

\chapter[Sobolev and Hardy inequalities on hyperboloids]{Sobolev and Hardy inequalities on hyperboloids \label{cha:5}}

\section{A Sobolev inequality on hyperboloids}
\label{sec:51}

To turn $L^2$ energy estimates into $L^\infty$ estimates,
we will rely on the following Sobolev inequality.

\begin{proposition}[Sobolev-type estimate on hyperboloids]
\label{pre lem sobolev} 
Let $u$ be a sufficiently regular function defined in the  cone $\Kcal = \{|x|<t-1\}$, then for all $s>0$ (and with $t = \sqrt{s^2+|x|^2}$)
\begin{equation}
\label{pre ineq sobolev}
\sup_{\Hcal_s} t^{3/2} |u(t,x)| \leq C \sum_L \sum_{|I|\leq 2} \| L^I u \|_{L^2(\Hcal_s)},
\end{equation}
where $C>0$ is a universal constant and the summation in $L$ is over all vector 
fields $L_a = x^a\del_t + t\del_a$, $a=1,2,3$.
\end{proposition}

In comparison to Lemma 7.6.1 in \citet{Hormander97}, observe that the right-hand side of
\eqref{pre ineq sobolev} does not contain the rotation fields $\Omega_{ab} := x^a\del_b - x^b\del_a$.

\begin{proof}
Recall the relation
$t = \sqrt{s^2 + |x|^2}$ on $\Hcal_s$ and consider a function $u$ defined in $\Kcal$, its restriction to the hyperboloid $\Hcal_s$ is, by definition, 
$$
w_s(x) := u(\sqrt{s^2+|x|^2},x).
$$
Fix $s_0$ and a point $(t_0,x_0)$ on the hyperboloid $\Hcal_{s_0}$, with 
$t_0 = \sqrt{s_0^2 + |x_0|^2}$. Observe that
\bel{eq:gh67}
\del_aw_{s_0}(x) = \delu_au\big(\sqrt{s_0^2+|x|^2},x\big) = \delu_au(t,x),
\ee
with $t = \sqrt{s_0^2 + |x|^2}$ and, therefore, 
\begin{equation}
t\del_aw_{s_0}(x) = t\delu_au\big(\sqrt{s_0^2+|x|^2},t\big) = L_au(t,x).
\ee

Introduce the function
$g_{s_0,t_0}(y) := w_{s_0}(x_0 + t_0\,y)$
and note that
$$
g_{s_0,t_0}(0) = w_{s_0}(x_0) = u\big(\sqrt{s_0^2+|x_0|^2},x_0\big)=u(t_0,x_0).
$$
Applying the standard Sobolev inequality to $g_{s_0,t_0}$, we obtain
$$
\big|g_{s_0,t_0}(0)\big|^2\leq C\sum_{|I|\leq 2}\int_{B(0,1/3)}|\del^Ig_{s_0,t_0}(y)|^2
\, dy,
$$
where $B(0, 1/3) \subset \RR^3$ denotes the ball centered at the origin and with radius $1/3$.

Taking into account the identity (with $x = x_0 + t_0y$)  
$$
\aligned
\del_ag_{s_0,t_0}(y)
& = t_0\del_aw_{s_0}(x_0 + t_0y) = t_0\del_aw_{s_0}(x)
\\
& = t_0\delu_au\big(t,x)\big),
\endaligned
$$
in view of \eqref{eq:gh67}, we see that for all $I$ 
$$
\del^Ig_{s_0,t_0}(y) = (t_0\delu)^I u(t,x)
$$
and, therefore,
$$
\aligned
\big|g_{s_0,t_0}(0)\big|^2 \leq& C\sum_{|I|\leq 2}\int_{B(0,1/3)}\big|(t_0\delu)^I u\big(t,x)\big)\big|^2dy
\\
= & C t_0^{-3}\sum_{|I|\leq 2}\int_{B((t_0,x_0),t_0/3)\cap \Hcal_{s_0}}\big|(t_0\delu)^I u\big(t,x)\big)\big|^2dx.
\endaligned
$$

We can check that
$$
\aligned
(t_0\delu_a(t_0\delu_b w_{s_0})) 
& = t_0^2\delu_a\delu_bw_{s_0} 
\\
& = (t_0/t)^2(t\delu_a)(t\delu_b) w_{s_0} - (t_0/t)^2 (x^a/t)L_b w_{s_0}.
\endaligned
$$
We also remark that $x^a/t = x^a_0/t + yt_0/t = (x^a_0/t_0 + y)(t_0/t)$, so 
that, in the region $y\in B(0,1/3)$ of interest, the factor $|x^a/t|$ is bounded by $C(t_0/t)$.  We conclude that for any $|I| \leq 2$,
$$
|(t_0 \delu)^I u| \leq \sum_{|J| \leq |I|} | L^I u| (t_0/t)^I. 
$$

On the other hand, when $|x_0|\leq t_0/2$ then $t_0\leq \frac{2}{\sqrt{3}}s_0$ and thus
$$
t_0 \leq C s_0 \leq C \sqrt{|x|^2 + s_0^2} = Ct, 
$$
$C$ being a universal constant thoughout. 
When $|x_0|\geq t_0/2$ then in the region $B((t_0,x_0),t_0/3)\cap \Hcal_{s_0}$ we have
also 
$$
t_0 \leq C|x| \leq C\sqrt{|x|^2 + s_0^2} =Ct.
$$
Consequently, it follows that
\be
|(t_0 \delu)^I u| \leq C \, \sum_{|J| \leq |I|} | L^I u| 
\ee
and 
$$
\aligned
\big|g_{s_0,t_0}(y_0)\big|^2
\leq& Ct_0^{-3}\sum_{|I|\leq 2}\int_{B(x_0,t_0/3)\cap\Hcal_{s_0}}\big|(t\delu)^I u\big(t,x)\big)\big|^2 \, dx
\\
\leq& Ct_0^{-3}\sum_{|I|\leq 2}\int_{\Hcal_{s_0}}\big|L^I u(t,x)\big|^2 \, dx,
\endaligned
$$
which completes the proof of Proposition~\ref{pre lem sobolev}. 
\end{proof}


\section{Application of the Sobolev inequality on hyperboloids}

Using now 
 the hyperboloidal energy defined in
 \eqref{pre pre expression of energy} 
and combining Proposition~\ref{pre lem sobolev} 
with the technical estimate \eqref{pre lem commutator T/t'},
we can deduce various 
sup-norm estimates, presented now.
For clarity in the presention, we make use of the notation $E_{m, \sigma}$ (defined in Chapter 2) in order to emphasize the dependency of the hyperboloidal energy upon the coeffficient $\sigma$ in the Klein-Gordon equation.

\begin{lemma}[$L^\infty$ estimates on derivatives up to first-order]\label{pre lem decay}
If $u$ is a  sufficiently regular  function supported in $\Kcal$,
then the following estimates hold: 
\bel{pre lem decay eq}
\aligned
&\text{(a)} \quad  \sup_{\Hcal_s}\big|t^{1/2}s\del_\alpha u\big| \leq C\sum_{|I|\leq 2}E_{m}(s,Z^I u)^{1/2},
\\
&\text{(b)} \quad \sup_{\Hcal_s}\big|t^{3/2}\delu_a  u\big| \leq C\sum_{|I|\leq 2}E_{m}(s,Z^I u)^{1/2},
\\
&\text{(c)} \quad  \sup_{\Hcal_s}\big|\sigma \, t^{3/2} u\big| \leq C\sum_{|I|\leq 2}E_{m,\sigma}(s,Z^I u)^{1/2},
\endaligned
\ee
where $Z$ stands for any admissible vector field, that is, any of $\del_{\alpha}, L_a$, 
and $C$ is a universal constant.
\end{lemma}

\begin{remark} Let us illustrate our result with the homogeneous linear wave equation
\be
\Box w = 0,\quad w|_{\Hcal_{B+1}} = w_0,\quad \del_t w|_{\Hcal_{B+1}} = w_1,
\ee
where the solution $w_i$ is defined in $\Hcal_{B+1}\cap \Kcal$. 
Thanks to the energy estimate in Proposition \ref{pre lem energy}, the  energy $E_m(s,Z^Iw)$ is controlled by the initial energy.
By the commutator estimates \eqref{pre lem commutator trivial} and 
\eqref{pre lem commutator} and by 
the Sobolev inequality,  
 we find
\be
\big|\delu_a w\big| \leq Ct^{-3/2},\quad \big|\del_\alpha  w\big| \leq Ct^{-3/2+1}(t^2-r^2)^{-1/2}, 
\ee
which are classical estimates. 
We emphasize that our argument of 
proof is ``robust'' in the sense 
that it 
uses neither the explicit expression of the solution nor the scaling vector field $S=r\del_r + t\del_t$.
\end{remark}

Now, we turn our attention to the energy and decay estimates for the ``good" second-order derivatives,
that is, derivatives such as $\delu_a \del_\alpha  u$.

\begin{lemma}[Bounds on second-order derivatives]\label{pre lem decay high order}
For every  sufficiently regular  func\-tion $u$ supported in the cone $\Kcal$, the following estimates hold:
\bel{pre decay 2order}
\sup_{\Hcal_s}\big|t^{3/2}s\delu_a\delu_\alpha  u\big| + \sup_{\Hcal_s}\big|t^{3/2}s\delu_\alpha  \delu_a u\big| \leq C\sum_{|I|\leq 3} E_m(s,Z^Iu)^{1/2},
\ee
\bel{pre L2 2order}
\int_{\Hcal_s}\big|s\delu_a\delu_\alpha  u\big|^2 dx + \int_{\Hcal_s}\big|s\delu_\alpha \delu_a u\big|^2 dx \leq C\sum_{|I|\leq 1}E_m(s,Z^I u).
\ee
\end{lemma}

\begin{proof} Recalling that 
$\delu_a  = t^{-1}L_a$, 
we obtain
$|\delu_a\del_\alpha u|\leq t^{-1}|L_a\del u|$. 
By Lemma \ref{pre lem decay} and the commutator estimate \eqref{pre lem commutator pr1},  we obtain \eqref{pre decay 2order}. 
The second estimate is immediate in view of the expression \eqref{pre pre expression of energy}.
\end{proof}

\begin{remark} Energy estimates and $L^\infty$ estimates for
the second-order time derivative $\delu_0\delu_0 u$ will be derived
 later from the wave equation itself, thanks to the decomposition in Proposition \ref{decompW}. 
\end{remark}

At the end of this section, we state the $L^\infty$ estimates of the solution of wave equation (i.e. $c_i = 0$).

\begin{lemma}\label{pre lem u it-self}
If $u$ is a  sufficiently regular  function supported in the cone $\Kcal$, then for any multi-index $J$, if $\sum_{|I|\leq |J|+2} E_m(s,Z^I u)^{1/2}\leq C's^\delta$ for 
some $\delta \geq 0$, then one has
\bel{pre lem u it-self decay}
\big|Z^J u\big| \leq CC't^{(-2+\delta)/2}(t-r)^{(1+\delta)/2} \leq C''C't^{-3/2}s^{1+\delta},
\ee
where $C, C',C''$ are universal constants. 
\end{lemma}

\begin{proof}
Using that $\sum_{|I|\leq |J|+2} E_m(s,Z^I u)^{1/2}$ is bounded by $C's^{\delta}$
and recalling by Lemma \ref{pre lem decay}, we find in the cone $\Kcal$ 
$$
|\del_r Z^I u| \leq Ct^{-(2-\delta)/2}(t-r)^{-(1-\delta)/2}.
$$
Then, \eqref{pre lem u it-self decay} follows by integration along radial directions.
\end{proof}


\section{Hardy inequality for the hyperboloidal foliation}
\label{sec:52}

In this section we establish an analogue of the classical Hardy inequality but generalized to hyperboloidal foliations. This inequality will be used in order to control $L^2$ norms of wave components such as $\|Z^I u\|_{L^2(\Hcal_s)}$.

\begin{proposition}[The hyperboloidal Hardy inequality]
\label{pre Hardy hyper}
For all sufficiently regular functions $u$ supported in the cone $\Kcal$, one has 
\bel{pre Hardy ineq}
\aligned
\|s^{-1}u\|_{L^2(\Hcal_s)} \leq&
 C\|s_0^{-1}u\|_{L^2(\Hcal_{s_0})} + C\sum_{a}\|\delu_a u\|_{L^2(\Hcal_s)}
\\
&+ C\sum_a\int_{s_0}^{s}\tau^{-1}
\Big(\|(\tau/t)\del_a u\|_{L^2(\Hcal_\tau)}+\|\delu_a u\|_{L^2(\Hcal_\tau)}\Big) \, d\tau.
\endaligned
\ee
\end{proposition}

Before proving this result, we begin with the following 
modified version of the classical Hardy inequality.

\begin{lemma}
\label{pre Hardy lem1}
For all sufficiently regular functions $u$ supported in the cone $\Kcal$, one has 
$$
\|r^{-1} u\|_{L^2(\Hcal_s)}\leq C\sum_a \|\delu_a u\|_{L^2(\Hcal_s)}.
$$
\end{lemma}

\begin{proof} As in the proof of Proposition~\ref{pre lem sobolev}, we consider the function
$w_s(x):= u\big(\sqrt{s^2+|x|^2},x\big)$, which satisfies 
$$
\del_aw_s(x) = \delu_a u\big(\sqrt{s^2+|x|^2},x\big). 
$$
We then apply the classical Hardy inequality to $w_s$ and obtain 
$$
\aligned
\int_{\mathbb{R}^3}|r^{-1} w_s(x)|^2dx &\leq C\int_{\mathbb{R}^3}|\nabla w_s(x)|^2dx = C\sum_a\int_{\mathbb{R}^3}\big|\delu_au(\sqrt{s^2+r^2},x)\big|^2dx
\\
&\leq C\sum_a\int_{\Hcal_s}\big|\delu_a u(t,x)\big|^2dx.
\endaligned
$$ 
\end{proof}

\begin{proof}[Proof of Proposition \ref{pre Hardy hyper}] 
We introduce a smooth cut-off function $\chi$ satisfying
$$
\chi (r) =
\left\{
\aligned
&1, \quad 2/3\leq r,
\\
&0, \quad 0\leq r\leq 1/3
\endaligned
\right.
$$
and consider the decomposition
$$
\|s^{-1}u\|_{L^2(\Hcal_s)} \leq \|\chi(r/t)s^{-1} u\|_{L^2(\Hcal_s)} + \|(1-\chi(r/t))s^{-1} u\|_{L^2(\Hcal_s)}
$$
which distinguish between the region ``near'' and ``away'' from the light cone. 

The estimate of $\|(1-\chi(r/t))s^{-1} u\|_{L^2(\Hcal_s)}$ is based on the following observation:
$$
\big(1-\chi(r/t)\big)s^{-1}\leq Ct^{-1}  \qquad \text{ in the cone } \Kcal,
$$
so that, by Lemma \ref{pre Hardy lem1}, 
\begin{equation}\label{eq pr2 Hardy}
\aligned
\|(1-\chi(r/t))us^{-1}\|_{L^2(\Hcal_s)}
& \leq \|t^{-1}u\|_{L^2(\Hcal_s)}
\\
& \leq \|r^{-1}u\|_{L^2(\Hcal_s)}\leq C\sum_a\|\delu_a u\|_{L^2(\Hcal_s)}.
\endaligned
\end{equation}

The estimate near the light cone is more delicate and to deal with 
the term $\|\chi(r/t)s^{-1} u\|_{L^2(\Hcal_s)}$, we proceed as follows:
 in the region $\Kcal_{[s_0,s]}$, we can find a positive constant $C$
$$
\chi(r/t)\leq C\frac{\chi(r/t)r}{(1+r^2)^{1/2}}, 
$$
and thus 
$$
\|\chi(r/t)s^{-1}u\|_{L^2(\Hcal_s)}\leq C\|r(1+r^2)^{-1/2}\chi(r/t)s^{-1}u\|_{L^2(\Hcal_s)}.
$$
So, we can focus on controlling this latter term. 

To this end, we consider the vector field 
$$
W = \bigg(0, -x^a  {t (u\chi(r/t))^2\over (1+r^2)s^2}\bigg)
$$
defined in $\Kcal$
and we compute its divergence
$$
\aligned
&\text{div} \hskip.05cm  W
\\
& = s^{-1} \del_a u \frac{r\chi(r/t)u}{(1+r^2)^{1/2}s}\frac{-2x^a t\chi(r/t)}{r(1+r^2)^{1/2}}
- s^{-1}\frac{u}{r} \frac{r\chi(r/t) u}{s(1+r^2)^{1/2}}\frac{2\chi'(r/t)r}{(1+r^2)^{1/2}}
\\
& \quad 
-\bigg(\frac{r^2t+3t}{(1+r^2)^2s^2} + \frac{2r^2t}{(1+r^2)s^4}\bigg)\big(u\chi(r/t)\big)^2.
\endaligned
$$

Next, we integrate the above inequality in the region $\Kcal_{[s_0,s_1]}\subset \Kcal\cap \{s_0\leq \sqrt{t^2-r^2}\leq s_1\}$ with respect to the Lebesgue measure in $\mathbb{R}^4$:
$$
\aligned
&\int_{\Kcal_{[s_0,s_1]}}\!\!\!\!\text{div} \hskip.05cm  W dxdt
\\
&= -2\int_{\Kcal_{[s_0,s_1]}}\!\!\!\!s^{-1} \bigg( \del_a u \frac{r\chi(r/t)u}{(1+r^2)^{1/2}s}\frac{x^a t\chi(r/t)}{r(1+r^2)^{1/2}}\bigg)dxdt
\\
& \quad -2\int_{\Kcal_{[s_0,s_1]}}\!\!\!\! s^{-1}\frac{u}{r}  \frac{r\chi(r/t) u}{s(1+r^2)^{1/2}}\frac{\chi'(r/t)r}{(1+r^2)^{1/2}}dxdt
\\
&\quad -\int_{\Kcal_{[s_0,s_1]}}\!\!\!\!\bigg(\frac{r^2t+3t}{(1+r^2)^2s^2} + \frac{2r^2t}{(1+r^2)s^4}\bigg)\big(u\chi(r/t)\big)^2dxdt
\endaligned
$$
thus 
$$
\aligned
&\int_{\Kcal_{[s_0,s_1]}}\!\!\!\!\text{div} \hskip.05cm  W dxdt
\\
&= -2\int_{s_0}^{s_1}\int_{\Hcal_s}(s/t)s^{-1} \bigg( \del_a u \frac{r\chi(r/t)u}{(1+r^2)^{1/2}s}\frac{x^a t\chi(r/t)}{r(1+r^2)^{1/2}}\bigg)dxds
\\
&\quad -2\int_{s_0}^{s_1}\int_{\Hcal_s}(s/t) s^{-1}\frac{u}{r}  \frac{r\chi(r/t) u}{s(1+r^2)^{1/2}}\frac{\chi'(r/t)r}{(1+r^2)^{1/2}}dxds
\\
&\quad -\int_{s_0}^{s_1}\int_{\Hcal_s}(s/t)\bigg(\frac{r^2t+3t}{(1+r^2)^2s^2} + \frac{2r^2t}{(1+r^2)s^4}\bigg)\big(u\chi(r/t)\big)^2dxds
\\
&=:\int_{s_0}^{s_1}\big( T_1 + T_2 +T_3\big) \, ds,
\endaligned
$$ 
where
$$
\aligned
T_1 =& -2 s^{-1}\int_{\Hcal_s}(s/t) \bigg( \del_a u \frac{r\chi(r/t)u}{(1+r^2)^{1/2}s}\frac{x^a t\chi(r/t)}{r(1+r^2)^{1/2}}\bigg) \, dx
\\
\leq & 2s^{-1}\bigg{\|}\frac{r u\chi(r/t)}{s(1+r^2)^{1/2}}\bigg{\|}_{L^2(\Hcal_s)}
\\
& \qquad \qquad \cdot 
\sum_a \|(s/t)\del_a u\|_{L^2(\Hcal_s)}\big{\|}\chi(r/t)x^atr^{-1}(1+r^2)^{-1/2}\big{\|}_{L^\infty(\Hcal_s)}
\\
\leq &Cs^{-1}\bigg{\|}\frac{r u\chi(r/t)}{s(1+r^2)^{1/2}}\bigg{\|}_{L^2(\Hcal_s)}\sum_a\|(s/t)\del_a u\|_{L^2(\Hcal_s)},
\endaligned
$$
$$
\aligned
T_2 =& -2 s^{-1}\int_{\Hcal_s}(s/t)\frac{u}{r}  \frac{r\chi(r/t) u}{s(1+r^2)^{1/2}}\frac{\chi'(r/t)r}{(1+r^2)^{1/2}}dx
\\
\leq &Cs^{-1}
\bigg{\|}\frac{r u\chi(r/t)}{s(1+r^2)^{1/2}}\bigg{\|}_{L^2(\Hcal_s)}
\|ur^{-1}\|_{L^2(\Hcal_s)}\big{\|}r\chi'(r/t)(1+r^2)^{-1/2}\big{\|}_{L^\infty(\Hcal_s)}
\\
\leq &Cs^{-1}\bigg{\|}\frac{r u\chi(r/t)}{s(1+r^2)^{1/2}}\bigg{\|}_{L^2(\Hcal_s)}\sum_{a}\|\delu_a u\|_{L^2(\Hcal_s)},
\endaligned
$$
where Lemma \ref{pre Hardy lem1} is used. We also observe that 
$T_3\leq 0$.

We write our identity in the form
\be
\frac{d}{ds}\bigg(\int_{\Kcal_{[s_0,s]}}\!\!\!\!\text{div} \hskip.05cm W \, dxdt\bigg) = T_1 + T_2 + T_3
\ee
and obtain
\bel{pre Hardy hyper eq0}
\aligned
& \frac{d}{ds}\bigg(\int_{\Kcal_{[s_0,s_1]}}\!\!\!\!\text{div} \hskip.05cm W \, dxdt\bigg)
\\
& \leq Cs^{-1}\bigg{\|}\frac{r u\chi(r/t)}{s(1+r^2)^{1/2}}\bigg{\|}_{L^2(\Hcal_s)}
\sum_{a}\big(\|(s/t)\del_a\|_{L^2(\Hcal_s)}+\|\delu_a u\|_{L^2(\Hcal_s)}\big).
\endaligned
\ee
On the other hand, we apply Stokes' formula in the region $\Kcal_{[s_0,s_1]}$,
and find 
$$
\aligned
& \int_{\Kcal_{[s_0,s_1]}}\!\!\!\!\text{div} \hskip.05cm W \, dxdt
\\
=&\int_{\Hcal_{s}}W \cdot n \, d\sigma + \int_{\Hcal_{s_0}}W\cdot n \, d\sigma
\\
=&
\int_{\Hcal_{s}}\frac{r^2}{1+r^2}\big|u\chi(r/t)s^{-1}\big|^2 dx
-\int_{\Hcal_{s_0}}\frac{r^2}{1+r^2}\big|u\chi(r/t)s^{-1}\big|^2dx.
\endaligned
$$
By differentiating this identity with respect to $s$, it follows that 
\bel{pre Hardy hyper eq1}
\aligned
&\frac{d}{ds}\bigg(\int_{\Kcal_{[s_0,s_1]}}\!\!\!\!\text{div} \hskip.05cm W \, dxdt\bigg) =
\frac{d}{ds}\bigg(\int_{\Hcal_{s}}\frac{r^2}{1+r^2}\big|u\chi(r/t)s^{-1}\big|^2 dx\bigg)
\\
&= 2 \, \bigg{\|}\frac{r u\chi(r/t)}{s(1+r^2)^{1/2}}\bigg{\|}_{L^2(\Hcal_{s})}\frac{d}{ds}\bigg{\|}\frac{r u\chi(r/t)}{s(1+r^2)^{1/2}}\bigg{\|}_{L^2(\Hcal_{s})}.
\endaligned
\ee
Finally, 
combining \eqref{pre Hardy hyper eq0} and \eqref{pre Hardy hyper eq1} yields us
$$
\frac{d}{ds}\bigg{\|}\frac{r u\chi(r/t)}{s(1+r^2)^{1/2}}\bigg{\|}_{L^2(\Hcal_{s})}
\leq Cs^{-1}\sum_a\big(\|\delu_a u\|_{L^2(\Hcal_{s})} + \|\delu_a u\|_{L^2(\Hcal_{s})}\big)
$$
and, by integration over the interval $[s_0,s]$,
\bel{pre Hardy hyper eq2}
\aligned
&
\big{\|}r(1+r^2)^{-1/2}\chi(r/t)s^{-1}u\big{\|}_{L^2(\Hcal_s)}
\\
&\leq \big{\|}r(1+r^2)^{-1/2}\chi(r/t)s_0^{-1}u\big{\|}_{L^2(\Hcal_{s_0})}
\\
&\quad + C\sum_a\int_{s_0}^s\tau^{-1}\big(\|\delu_a u\|_{L^2(\Hcal_\tau)} + \|\delu_a u\|_{L^2(\Hcal_\tau)}\big).
\endaligned
\ee
In view of Lemma \ref{pre Hardy lem1}, we conclude that 
\begin{equation}\label{eq pr1 ch5 prop 1}
\aligned
&\|\chi(r/t)s^{-1}u\|_{L^2(\Hcal_s)}
\\
&\leq C \, \|r(1+r^2)^{-1/2}\chi(r/t)s^{-1}u\|_{L^2(\Hcal_s)}
\\
&\leq C \, \big{\|}s_0^{-1}u\big{\|}_{L^2(\Hcal_{s_0})} +  C\sum_a\int_{s_0}^s\tau^{-1}\big(\|\delu_a u\|_{L^2(\Hcal_\tau)} + \|\delu_a u\|_{L^2(\Hcal_\tau)}\big).
\endaligned
\end{equation} 
The desired conclusion is reached by combining \eqref{eq pr2 Hardy} with \eqref{eq pr1 ch5 prop 1}.
\end{proof}

\chapter[Revisiting scalar wave equations]{Revisiting scalar wave equations\label{cha:10}}

\section{Background and statement of the main result}

In this chapter, we revisit the classical global existence theory for scalar nonlinear wave equations (with initial data imposed on a hyperboloid $\Hcal_{s_0}$, with $s_0\geq1$): 
\bel{wave eq main}
\aligned
&\Box u = P^{\alpha\beta}\del_{\alpha}u\del_{\beta}u
\\
&u|_{\Hcal_{s_0}} = u_0,\qquad \del_t u|_{\Hcal_{s_0}} = u_1,
\endaligned
\ee
with smooth initial data $u_0, u_1$ compactly supported in the open ball $\mathcal{B}(0,s_0)$.
We denote by $\mathcal{B}(0,s_0)$ the intersection of the spacelike hypersurface
$\Hcal_{s_0}$ and the cone $\Kcal = \big\{ (t,x) \, / \, |x|< t-1 \big\}$.
We impose the classical null condition on the bilinear form $P^{\alpha\beta}$, that is,
\be
\aligned
& P^{\alpha\beta}\xi_{\alpha}\xi_{\beta} = 0 \quad
\text{for all $\xi\in \RR^4$ satisfying $-\xi_0^2 + \sum_{a}\xi_a^2 = 0$.}
\endaligned
\ee
We are going to revisit this classical problem with the hyperboloidal foliation method and 
we establish that the energy of solutions is uniformly bounded, 
while, according to the classical technique of proof, the energy is only known to be at most polynomially increasing.  
The hyperboloidal energy $E_m=E_{m,0}$ was introduced in \eqref{exp14}, while the admissible vector fields $Z \in \Zscr$ were defined in \eqref{vectorsZ}.

\begin{theorem}[Existence theory for scalar wave equations]
\label{wave thm main} 
There exist $\eps_0, C_1>0$ such that for all initial data satisfying 
\be
\label{wave eq 0 thm main}
E_m(s_0,Z^I u)^{1/2}\leq \eps\leq \eps_0 \qquad \text{ for all } |I|\leq 3, \quad Z \in \Zscr, 
\ee
the local-in-time solution $u$ to the Cauchy problem \eqref{wave eq main} extends to arbitrarily large times and, 
furthermore,  
\be\label{wave eq 1 thm main}
E_m(s,Z^I u)^{1/2}\leq C_1\eps  \qquad \text{ for all } |I|\leq 3, \quad Z \in \Zscr, 
\ee
and
\be\label{wave eq 2 thm main}
\big|\del_\alpha u(t,x)\big|\leq C_1\varepsilon t^{-1}(t- |x|)^{-1/2}.
\ee
\end{theorem}


\section{Structure of the proof}
\label{subsec wave 2}

We proceed with a bootstrap strategy and, for some large constant $C_1>1$, we assume that, in a time interval $[s_0, s_1]$  the local solution satisfies the bound 
\be
\label{wave eq 1-2}
 \sum_{|I|\leq 3 \atop Z \in \Zscr} E_m(s,Z^Iu)^{1/2} 
\leq C_1\eps\quad \text{for } s\in [s_0,s_1].
\ee
We take 
$$
s_1:=\sup\Big\{ s\geq s_0 \, / \,  \sum_{|I|\leq 3 \atop Z \in \Zscr}E_m(\tau,Z^Iu)^{1/2} \leq C_1\eps 
\, \text{ for all } \tau \in [s_0, s] \Big\}
$$
to be the largest such time and we suppose that it would be finite.
Since $C_1> 1$, by a continuity argument, we know that $s_1>s_0$. 

Our objective is to establish, for a suitable choice of $\eps_0, C_1>0$, that for all 
$\eps\leq \eps_0$, 
\be\label{wave eq 2-2}
 \sum_{|I|\leq 3 \atop Z \in \Zscr}
E_m(s,Z^Iu)^{1/2} 
\leq \frac{1}{2}C_1\eps \quad \text{for } s\in [s_0,s_1].
\ee
This leads us to
$$
 \sum_{|I|\leq 3 \atop Z \in \Zscr}
E_m(s_1,Z^Iu)^{1/2} \leq \frac{1}{2}C_1\eps, 
$$
and, by continuity, 
\be
s_1 < \sup \Big\{ s\geq s_0 \, / \, 
 \sum_{|I|\leq 3 \atop Z \in \Zscr} E_m(\tau,Z^Iu) \leq C_1\eps
\, \text{ for all } \tau \in [s_0, s] \Big\}, 
\ee
which would be a contradiction. We can then conclude that $s_1=+\infty$. Namely, in view of the local-in-time existence theory (cf.~Theorem~\ref{local semi-hyper}), the solution $u$ extends to all times.

In other words, our task reduces to proving the following result,
and the rest of this chapter is devoted to its proof.

\begin{proposition}\label{wave prop 1}
Let $u$ be a solution to \eqref{wave eq main} defined in $[s_0,s_1]$ and with initial data satisfying
\be\label{wave eq 1 prop 1}
 \sum_{|I|\leq 3 \atop Z \in \Zscr}
E_m(s_0,Z^Iu)^{1/2}\leq \varepsilon.  
\ee
There exist constants $C_1,\varepsilon_0>0$ such that if $u$ satisfies the estimate \eqref{wave eq 1-2} with $\varepsilon\leq \varepsilon_0$,  then the estimate \eqref{wave eq 2-2} holds.
\end{proposition}

\section{Energy estimate}

The following lemma is essentially Proposition \ref{pre lem energy} in the special 
case of \eqref{wave eq main}, but for convenience we provide here a direct proof.

\begin{lemma}\label{wave lem 1}
If $u$ is a solution to \eqref{wave eq main} defined in $[s_0,s_1]$, then the energy estimate 
\be\label{wave eq 1 lem 1}
\aligned
&  \sum_{|I|\leq 3 \atop Z \in \Zscr} E_m(s,Z^I u)^{1/2}
\\
& \leq  \sum_{|I|\leq 3 \atop Z \in \Zscr} E_m(s_0,Z^I u)^{1/2} 
+  \int_{s_0}^s \sum_{|I|\leq 3 \atop Z \in \Zscr} \big\|Z^I\big(P^{\alpha\beta}\del_{\alpha}u\del_\beta u\big)\big\|_{L^2(\Hcal_\tau)}d\tau. 
\endaligned
\ee
\end{lemma}

\begin{proof} 
We apply to the equation \eqref{wave eq main} a product $Z^I$ with $|I|\leq 3$
and, by recalling the commutation relation $[Z^I,\Box] =0$, we find
$$
\Box \big( Z^I u \big)
= Z^I\big(P^{\alpha\beta}\del_{\alpha}u\del_{\beta}u\big).
$$
By multiplying this equation by $\del_t Z^I u$ and performing the standard energy calculation, 
we deduce that the function $\widetilde{u} := Z^I u$ satisfies 
$$
\frac{1}{2}\del_t\bigg((\del_t \widetilde{u})^2 + \sum_a (\del_a \widetilde{u})^2\bigg) - \del_a\big(\del_a \widetilde{u} \del_t \widetilde{u}\big) = Z^I\big(P^{\alpha\beta}\del_{\alpha}\widetilde{u}\del_{\beta}\widetilde{u}\big)\del_t \widetilde{u}.
$$
Integrating this equation in the region $\Kcal_{[s_0,s]}$, we have
\be
\label{wave eq pr1 lem 1}
\aligned
& \int_{\Kcal_{[s_0,s]}}
\Bigg(
\frac{1}{2}\del_t\bigg((\del_t \widetilde{u})^2 + \sum_a (\del_a \widetilde{u})^2\bigg) - \del_a\big(\del_a \widetilde{u} \del_t \widetilde{u} \Bigg)
\, dtdx
\\
& = \int_{\Kcal_{[s_0,s]}}Z^I\big(P^{\alpha\beta}\del_{\alpha}\widetilde{u}\del_{\beta}u\big)\del_t \widetilde{u}\,dtdx.
\endaligned
\ee

By Huygens' principle, the solution is supported in the cone $\Kcal$
and, in a neighborhood of the cone $\{|x| = t-1\}\cap \{s_0\leq \tau \leq s \}$, 
we have $\widetilde{u}=0$. By Stokes' formula, the left-hand side of \eqref{wave eq pr1 lem 1} reduces
 to
$$
\aligned
& \frac{1}{2}\int_{\Hcal_s}
\Big(
\big( \del_t \widetilde{u} \big)^2 + \sum_a \big( \del_a \widetilde{u} \big)^2, 2\del_t \widetilde{u}\del_a \widetilde{u}\Big) \cdot n \,
d\sigma
\\
& - \frac{1}{2}\int_{\Hcal_{s_0}}\big(|\del_t \widetilde{u}|^2 + \sum_a|\del_a \widetilde{u}|^2, 2\del_t \widetilde{u}\del_a \widetilde{u}\big).n d\sigma,
\endaligned
$$
where $n$ is the (future oriented) unit normal vector to the hyperboloids
 and $d\sigma$ is the induced Lebesgue measure on the hyperboloids, with 
$$
n = \big(t^2 + |x|^2\big)^{-1/2}(t,-x^a),\quad d\sigma = \frac{\big(t^2 + |x|^2\big)^{1/2}}{t} dx. 
$$
Hence, the left-hand side of \eqref{wave eq pr1 lem 1} reads
$$
\aligned
& \frac{1}{2}\int_{\Hcal_s}\bigg(|\del_t \widetilde{u}|^2 + \sum_a|\del_a \widetilde{u} |^2 + 2\frac{x^a}{t}\del_a \widetilde{u} \del_t \widetilde{u} \bigg)dx
\\
& -
\frac{1}{2}\int_{\Hcal_{s_0}}\bigg(|\del_t \widetilde{u}|^2 + \sum_a|\del_a \widetilde{u} |^2 + 2\frac{x^a}{t}\del_a \widetilde{u} \del_t \widetilde{u} \bigg)dx,
\endaligned
$$
which is $\frac{1}{2}E_m(s,Z^I u) - \frac{1}{2}E_m(s_0,Z^I u)$.

On the other hand, in the region $\Kcal_{[s_0,s]}$, we use the change of variable $\tau = (t^2 - |x|^2)^{1/2}$ and the identity $dtdx = (\tau/t)d\tau dx$, so that 
the right-hand side of \eqref{wave eq pr1 lem 1} becomes
$$
\int_{s_0}^s 
\int_{\Hcal_\tau}(\tau/t)\del_t u Z^I\big(P^{\alpha\beta}\del_{\alpha}u\del_{\beta}u\big)dxd\tau.
$$
We thus conclude that \eqref{wave eq pr1 lem 1} is equivalent to 
\be\label{wave eq pr2 lem 1}
\frac{1}{2}E_m(s,Z^I u) - \frac{1}{2}E_m(s_0,Z^I u)
 = \int_{s_0}^sds\int_{\Hcal_\tau}(\tau/t)\del_t u Z^I\big(P^{\alpha\beta}\del_{\alpha}u\del_{\beta}u\big)dx.
\ee

Next, we differentiate \eqref{wave eq pr2 lem 1} with respect to the variable $s$ and obtain 
$$
\aligned
E(\tau,Z^I u)^{1/2}\frac{d}{d\tau}E(\tau,Z^I u)^{1/2} =& \int_{\Hcal_\tau}(\tau/t)\del_t u Z^I\big(P^{\alpha\beta}\del_{\alpha}u\del_{\beta}u\big)dx
\\
\leq& \|(\tau/t)\del_t u\|_{L^2(\Hcal_\tau)} \big\|Z^I\big(P^{\alpha\beta}\del_{\alpha}u\del_{\beta}u\big)\big\|_{\Hcal_\tau}.
\endaligned
$$
Recalling the expression of the hyperboloidal energy \eqref{pre pre expression of energy}, we have 
$E(s,u)^{1/2} \geq \|(\tau/t)\del_t u\|_{L^2(\Hcal_s)}$ and therefore 
$$
\frac{d}{ds}  \sum_{|I|\leq 3 \atop Z \in \Zscr}
 E(s,Z^I u)^{1/2} 
\leq 
\sum_{|I|\leq 3 \atop Z \in \Zscr} 
 \big\|Z^I\big(P^{\alpha\beta}\del_{\alpha}u\del_{\beta}u\big)\big\|_{\Hcal_s}. 
$$
The conclusion follows by integrating over $[s_0,s]$.
\end{proof}


\section{Basic $L^2$ and $L^\infty$ estimates}

In this section, from the bootstrap assumption \eqref{wave eq 1-2}, we deduce
$L^2$ and $L^{\infty}$ estimates.
First of all, the following lemma is immediate in view of the expression of the 
hyperboloidal energy 
and \eqref{wave eq 1-2}.

\begin{lemma}[Basic $L^2$ estimates]
\label{wave lem 2}
By relying on \eqref{wave eq 1-2}, 
the following estimate hold for all $s\in [s_0,s_1]$:
\begin{subequations}\label{wave eq 1 lem 2}
\be\label{wave eq 1 lem 2 a}
 \sum_{|I|\leq 3 \atop Z \in \Zscr}
\|\delu_a Z^I u\|_{L^2(\Hcal_s)}\leq CC_1\eps,
\ee
\be\label{wave eq 1 lem 2 b}
 \sum_{|I|\leq 3 \atop Z \in \Zscr} \|(s/t)\delu_0 Z^I u\|_{L^2(\Hcal_s)}\leq CC_1\eps,
\ee
where $C>0$ is a (universal) constant.
\end{subequations}
\end{lemma}


Now, we combine Lemma \ref{wave lem 2} with the commutator estimates in Lemmas \ref{pre lem commutator1} and \ref{pre lem commutator}. 

\begin{lemma}\label{wave lem 2.5}
By relying on \eqref{wave eq 1-2}, the following estimate hold for all $s\in [s_0,s_1]$:
\begin{subequations}\label{wave eq 1 lem 2.5}
\be\label{wave eq 1 lem 2.5 a}
\|Z^{I_1}\delu_a Z^{I_2}u\|_{L^2(\Hcal_s)}\leq CC_1\eps \quad \text{for all } |I_1|+|I_2|\leq 3,
\ee
\be\label{wave eq 1 lem 2.5 b}
\|Z^{I_1}\big((s/t)\delu_0 Z^{I_2}u\big)\|_{L^2(\Hcal_s)}\leq CC_1\eps 
\quad \text{for all } |I_1|+|I_2|\leq 3, 
\ee
where $C>0$ is a universal constant.
\end{subequations}
\end{lemma}

Furthermore, the following decay estimate is immediate in view of \eqref{wave eq 1 lem 2.5} and
the Sobolev estimate on hyperboloids
 \eqref{pre ineq sobolev}.

\begin{lemma}[Basic $L^{\infty}$ estimates]\label{wave lem 3}
By relying on \eqref{wave eq 1-2}, the following estimate hold for all 
$s\in [s_0,s_1]$:
\begin{subequations}\label{wave eq 1 lem 3}
\be\label{wave eq 1 lem 3 a}
\|t^{3/2}\delu_a Z^J u\|_{L^\infty(\Hcal_s)}\leq CC_1\eps  \quad \text{ for all } |J|\leq 1, 
\ee
\be\label{wave eq 1 lem 3 b}
\|t^{1/2}s\delu_0 Z^I u\|_{L^\infty(\Hcal_s)}\leq CC_1\eps  \quad \text{ for all } |J|\leq 1, 
\ee
where $C>0$ is a universal constant.
\end{subequations}
\end{lemma}


\section{Estimate on the interaction term}

We are now in a position to control the interaction term $P^{\alpha\beta}\del_{\alpha}u\del_{\beta}u$ with the help of the $L^2$ bound \eqref{wave eq 1 lem 2.5} and the $L^\infty$ bound \eqref{wave eq 1 lem 3}.

\begin{lemma}\label{wave lem 4}
By reyling on the inequalities \eqref{wave eq 1 lem 2.5} and \eqref{wave eq 1 lem 3} 
and by assuming that the billinear form associated with $P^{\alpha\beta}$ is a null form,
the following estimate holds for all $s \in [s_0, s_1]$: 
\bel{wave eq 1 lem 4}
\big\|Z^I\big(P^{\alpha\beta}\del_{\alpha}u\del_{\beta}u\big)\big\|_{L^2(\Hcal_s)}\leq C K \, (C_1\varepsilon)^2s^{-3/2}, 
\ee
where $C>0$ is a universal constant 
and 
$K = \max_{\alpha,\beta}|P^{\alpha\beta}|$.
\end{lemma}

\begin{proof} We  recall Proposition \ref{pre lem null 1} which tells us how to estimate a null form, that is
$$
\aligned
\big|Z^I\big(P^{\alpha\beta}\del_\alpha u\del_{\beta}u\big)\big|
&\leq CK(s/t)^2\sum_{|I_1|+|I_2|\leq|I|}\big|Z^{I_1}\del_tuZ^{I_2}\del_tu\big|
\\
&\quad+ CK\sum_{a,\beta,\atop |I_1|+|I_2|\leq|I|}\Big(\big|Z^{I_1}\delu_au\,Z^{I_2}\delu_{\beta}u\big| + \big|Z^{I_1}\delu_{\beta}u\,Z^{I_2}\delu_au\big| \Big)
\\
&=: T_1 + T_2,
\endaligned
$$
where $K = \max_{\alpha,\beta}|P^{\alpha\beta}|$.

To estimate the $L^2$ norm of each term in $T_1$, we remark that $|I_1| + |I_2| \leq |I|\leq 3$ implies that $|I_1|\leq 1$ or $|I_2|\leq 1$. Without loss of generality, we assume $|I_2|\leq 1$ and then write 
$$
\aligned
& \big\|(s/t)^2Z^{I_1}\del_tuZ^{I_2}\del_tu\big\|_{L^2(\Hcal_s)}
\\
& \leq \big\|(s/t)\del_tu\big\|_{L^2(\Hcal_s)}\big\|(s/t)t^{-1/2}s^{-1}\big(t^{1/2}sZ^{I_2}\del_tu\big)\big\|_{L^{\infty}(\Hcal_s)}
\\
& \leq CC_1\eps \, s^{-3/2}CC_1\eps\leq C(C_1\eps)^2s^{-3/2}.
\endaligned
$$

The terms in $T_2$ are estimated along the same idea by writing, when $|I_1|\leq 1$,
$$
\aligned
\big\|Z^{I_1}\delu_au\,Z^{I_2}\delu_{\beta}u\big\|_{L^2(\Hcal_s)}=
& \big\|t^{3/2}Z^{I_1}\delu_au\,t^{-3/2}(t/s)(s/t)Z^{I_2}\delu_{\beta}u\big\|_{L^2(\Hcal_s)}
\\
\leq & \big\|t^{3/2}Z^{I_1}\delu_au\,s^{-3/2}(s/t)Z^{I_2}\delu_{\beta}u\big\|_{L^2(\Hcal_s)}
\\
=& s^{-3/2}\big\|t^{3/2}Z^{I_1}\delu_au\big\|_{L^\infty(\Hcal_s)}\big\|(s/t)Z^{I_2}\delu_{\beta}u\big\|_{L^2(\Hcal_s)}
\\
\leq & C(C_1\eps)^2s^{-3/2}
\endaligned
$$
and, on the other hand when $|I_2|\leq 1$,
$$
\aligned
\big\|Z^{I_1}\delu_au\,Z^{I_2}\delu_{\beta}u\big\|_{L^2(\Hcal_s)}
=& \big\|Z^{I_1}\delu_au \,t^{-1/2}s^{-1}\,t^{1/2}sZ^{I_2}\delu_{\beta}u\big\|_{L^2(\Hcal_s)}
\\
\leq & \big\|Z^{I_1}\delu_au\,s^{-3/2}t^{1/2}sZ^{I_2}\delu_{\beta}u\big\|_{L^2(\Hcal_s)}
\\
=& s^{-3/2}\big\|Z^{I_1}\delu_au\big\|_{L^2(\Hcal_s)}\big\|t^{1/2}sZ^{I_2}\delu_{\beta}u\big\|_{L^\infty(\Hcal_s)}
\\
\leq & C(C_1\eps)^2s^{-3/2}.
\endaligned
$$ 
\end{proof}


\section{Conclusion}

\begin{proof}[Proof of Proposition \ref{wave prop 1}]
Now, we combine the null form estimate \eqref{wave eq 1 lem 4} and the energy estimate 
\eqref{wave eq 1 lem 1}, and obtain 
$$
\aligned
\sum_{|I|\leq 3 \atop Z \in \Zscr} 
E_m(s,Z^I u)^{1/2}
\leq& \sum_{|I|\leq 3 \atop Z \in \Zscr}  E_m(s_0,Z^I u)^{1/2} + C(C_1\eps)^2\int_{s_0}^s\tau^{-3/2} \, d\tau
\\
\leq& \sum_{|I|\leq 3 \atop Z \in \Zscr}  E_m(s_0,Z^I u)^{1/2} + C(C_1\eps)^2\int_{1}^{+\infty}\tau^{-3/2} \, d\tau 
\\
\leq& \eps + C(C_1\eps)^2.
\endaligned
$$
In order to conclude, we take $C_1>2$ and $\eps\leq \eps_0 = \frac{C_1-2}{2CC_1^2}$, 
and we find that 
$$
 \sum_{|I|\leq 3 \atop Z \in \Zscr} 
E_m(s,Z^I u)^{1/2} \leq \frac{1}{2}C_1\eps. 
$$ 
\end{proof}

The time-asymptotics of the solutions is also clear: \eqref{wave eq 1 thm main} has already been proved by
Proposition \ref{wave prop 1}. To see \eqref{wave eq 2 thm main}, we remark that by Lemma \ref{wave lem 3},
$$
|\del_t u|\leq C_1\eps t^{-1/2}s^{-1} 
\leq C_1\eps t^{-1}(t-r)^{-1/2},\quad |\delu_a u| \leq C_1\eps t^{-3/2}.
$$
By recalling the relation
$$
\del_a u = \delu_a u - \frac{x^a}{t}\del_t u,
$$
then \eqref{wave eq 2 thm main} follows. In view of the discussion at the beginning of Section \ref{subsec wave 2}, Theorem \ref{wave thm main} is now established.

\chapter[Fundamental $L^\infty$ and $L^2$ estimates]{Fundamental $L^\infty$ and $L^2$ estimates  \label{cha:6}}

\section{Objective of this chapter}
\label{sec:61}

We begin the discussion of the bootstrap arguments and we suppose 
that \eqref{proof energy assumption} holds on some time interval.
Our aim is to derive additional estimates  from these assumptions.
In the present chapter, we are able to deduce several $L^2$ and $L^\infty$ estimates on the solution and its derivatives. These rather immediate estimates will serve as a basis for the following chapters. The estimates in this chapter are classified into two groups: $L^2$ estimates and $L^\infty$ estimates:

\bei

\item The $L^2$-type estimates themselves are classified into two generations:

\bei

\item The estimates of the first generation are immediate consequences of \eqref{proof energy assumption}.

\item The estimates of the second generation are deduced from those of the first generation, by recalling
the commutator estimates (cf.~Lemmas \ref{pre lem commutator1} and \ref{pre lem commutator}). These are the $L^2$ bounds  that will be more often used in the following discussion.

\eei

\item The $L^\infty$-type estimates are also classified into two generations:

\bei

\item The estimates of the first generation follow immediately from \eqref{proof energy assumption} and the
Sobolev inequalities \eqref{pre lem sobolev}.

\item The estimates of the second generation are deduced from the ones of the first generation and the commutator
 estimates in Lemmas \ref{pre lem commutator1} and \ref{pre lem commutator}.

\eei

\eei

In the following, the letter $C$ will be used to represent a constant who depends on the structure of the system \eqref{main eq main}, such as the number of equations $n_0$, the number of wave components $j_0$, and the Klein-Gordon constants $c_i$.


\section{$L^2$ estimates of the first generation}
\label{sec:62}

From the expression of the energy \eqref{pre pre expression of energy} and the energy assumption \eqref{proof energy assumption}, the following estimates hold of all
$|\Is|\leq 5$ and $|\Id|\leq 4$:
\begin{subequations}\label{proof L2 1ge1}
\bel{proof L2 1ge1a}
\sum_{\alpha, i}\big{\|}(s/t)\del_\alpha Z^{\Is}w_i\big{\|}_{L^2(\Hcal_s)} +
\sum_{\alpha,i}\big{\|}(s/t)\delu_{\alpha}Z^{\Is}w_i\big{\|}_{L^2(\Hcal_s)}\leq CC_1\eps s^{\delta},
\ee
\bel{proof L2 1ge1a'}
\sum_{\alpha, i}\big{\|}(s/t)\del_\alpha Z^{\Id}w_i\big{\|}_{L^2(\Hcal_s)} +
\sum_{\alpha,i}\big{\|}(s/t)\delu_{\alpha}Z^{\Id}w_i\big{\|}_{L^2(\Hcal_s)}\leq CC_1\eps s^{\delta/2},
\ee
\bel{proof L2 1ge1b}
\sum_{a, i}\big{\|}\delu_aZ^{\Is}w_i\big{\|}_{L^2(\Hcal_s)}\leq CC_1\eps s^{\delta},
\ee
\bel{proof L2 1ge1b'}
\sum_{a, i}\big{\|}\delu_aZ^{\Id}w_i\big{\|}_{L^2(\Hcal_s)}\leq CC_1\eps s^{\delta/2},
\ee
\bel{proof L2 1ge1c}
\sum_{\ic}\big{\|} Z^{\Is}v_{\ic}\big{\|}_{L^2(\Hcal_s)}\leq CC_1\eps s^{\delta},
\ee
\bel{proof L2 1ge1c'}
\sum_{\ic}\big{\|} Z^{\Id}v_{\ic}\big{\|}_{L^2(\Hcal_s)}\leq CC_1\eps s^{\delta/2}.
\ee
\end{subequations}
By taking 
$$
\aligned
& Z^{\Is} = \del_\alpha Z^{\Id}, 
\quad 
Z^{\Id} = \del_\alpha Z^I, 
\\
& Z^{\Is} = t\delu_aZ^{\Id} = L_aZ^{\Id},
\quad 
Z^{\Id} = t\delu_aZ^I = L_aZ^I
\endaligned
$$
 in  \eqref{proof L2 1ge1c} and \eqref{proof L2 1ge1c'},
we have especially the following estimates on Klein-Gordon components for $|\Id|\leq 4$ and $|I|\leq 3$:
\begin{subequations}\label{proof L2 1ge2}
\bel{proof L2 1ge2a}
\sum_{\ic,\alpha}\big{\|} \del_\alpha Z^{\Id}v_{\ic}\big{\|} _{L^2(\Hcal_s)}
+ \sum_{\ic, a}\big{\|}\delu_{\alpha} Z^{\Id}v_{\ic}\big{\|} _{L^2(\Hcal_s)}
\leq CC_1\eps s^{\delta},
\ee
\bel{proof L2 1ge2a'}
\sum_{\ic,\alpha}\big{\|} \del_\alpha Z^Iv_{\ic}\big{\|} _{L^2(\Hcal_s)}
+ \sum_{\ic, a}\big{\|}\delu_{\alpha} Z^Iv_{\ic}\big{\|} _{L^2(\Hcal_s)}
\leq CC_1\eps s^{\delta/2},
\ee
\bel{proof L2 1ge2b}
\sum_{\ic, a}\big{\|} t \delu_a Z^{\Id}v_{\ic}\big{\|} _{L^2(\Hcal_s)} \leq CC_1\eps s^{\delta},
\ee
\bel{proof L2 1ge2b'}
\sum_{\ic, a}\big{\|} t \delu_a Z^Iv_{\ic}\big{\|} _{L^2(\Hcal_s)} \leq CC_1\eps s^{\delta/2}.
\ee
\end{subequations}
The bound on $\delu_{\alpha}Z^{I^\dag}v_{\ic}$ in \eqref{proof L2 1ge2a} is 
derived from the estimate on $\del_tZ^{I^\dag}v_{\ic}$ and the estimates on $\delu_aZ^{I^\dag}v_{\ic}$. 

For any $|I|\leq 3$, from \eqref{pre pre expression of energy} and \eqref{proof energy assumption e} we control the wave components:
\begin{subequations}\label{proof L2 1ge3}
\bel{proof L2 1ge3a}
\sum_{\alpha,\ih}\big{\|}(s/t)\del_\alpha Z^Iu_{\ih}\big{\|}_{L^2(\Hcal_s)}
+\sum_{\alpha,\ih}\big{\|}(s/t)\delu_\alpha Z^Iu_{\ih}\big{\|}_{L^2(\Hcal_s)}
\leq CC_1\eps,
\ee
\bel{proof L2 1ge3b}
\sum_{a,\ih}\big{\|}\delu_aZ^Iu_{\ih}\big{\|}_{L^2(\Hcal_s)} \leq CC_1\eps.
\ee
\end{subequations}


\section{$L^2$ estimates of the second generation}
\label{sec:63}

By using Lemma \ref{pre lem commutator}, we can commute the vector fields under consideration and, relying on the estimates established in the previous section, we obtain the following result.

The first group of estimates is obtained by \eqref{pre lem commutator trivial}, \eqref{pre lem commutator bar},  and \eqref{proof L2 1ge1}:
\begin{subequations}\label{proof L2 2ge1}
\bel{proof L2 2ge1a}
\sum_{\ih,\alpha}\big{\|}(s/t)Z^{\Is} \del_\alpha  w_i\big{\|}_{L^2(\Hcal_s)}
+\sum_{\ih,\alpha}\big{\|}(s/t)Z^{\Is} \delu_\alpha  w_i\big{\|}_{L^2(\Hcal_s)}\leq C C_1\eps s^{\delta},
\ee
\bel{proof L2 2ge1a'}
\sum_{\ih,\alpha}\big{\|}(s/t)Z^{\Id} \del_\alpha  w_i\big{\|}_{L^2(\Hcal_s)}
+\sum_{\ih,\alpha}\big{\|}(s/t)Z^{\Id} \delu_\alpha  w_i\big{\|}_{L^2(\Hcal_s)}\leq C C_1\eps s^{\delta/2},
\ee
\bel{proof L2 2ge1b}
\sum_{\ih,a}\big{\|}Z^{\Is} \delu_a w_{\ih}\big{\|}_{L^2(\Hcal_s)} \leq C C_1\eps s^{\delta},
\ee
\bel{proof L2 2ge1b'}
\sum_{\ih,a}\big{\|}Z^{\Id} \delu_a w_{\ih}\big{\|}_{L^2(\Hcal_s)} \leq C C_1\eps s^{\delta/2},
\ee
\bel{proof L2 2ge1c}
\sum_{\jc}\big{\|}Z^{\Is}v_{\jc}\big{\|}_{L^2(\Hcal_s)} \leq C C_1\eps s^{\delta},
\ee
\bel{proof L2 2ge1c'}
\sum_{\jc}\big{\|}Z^{\Id}v_{\jc}\big{\|}_{L^2(\Hcal_s)} \leq C C_1\eps s^{\delta/2}.
\ee
\end{subequations}

The second group of estimates follows from \eqref{pre lem commutator trivial},
\eqref{pre lem commutator bar}, and \eqref{proof L2 1ge2}:
\begin{subequations}\label{proof L2 2ge2}
\bel{proof L2 2ge2a}
\sum_{\ic,\alpha }\big{\|} Z^{\Id}\del_\alpha v_{\ic}\big{\|}_{L^2(\Hcal_s)}
+ \sum_{\ic,\alpha}\big{\|} Z^{\Id}\delu_\alpha v_{\ic}\big{\|}_{L^2(\Hcal_s)}
\leq CC_1\eps s^{\delta},
\ee
\bel{proof L2 2ge2a'}
\sum_{\ic,\alpha }\big{\|} Z^I\del_\alpha v_{\ic}\big{\|}_{L^2(\Hcal_s)}
+ \sum_{\ic,\alpha}\big{\|} Z^I\delu_\alpha v_{\ic}\big{\|}_{L^2(\Hcal_s)}
\leq CC_1\eps s^{\delta/2},
\ee
\bel{proof L2 2ge2b}
\sum_{\ic, a}\big{\|} t Z^{\Id} \delu_a v_{\ic}\big{\|}_{L^2(\Hcal_s)}\leq CC_1\eps s^{\delta},
\ee
\bel{proof L2 2ge2b'}
\sum_{\ic, a}\big{\|} t Z^I \delu_a v_{\ic}\big{\|}_{L^2(\Hcal_s)}\leq CC_1\eps s^{\delta/2}.
\ee
\end{subequations}

The third group of estimates follows from
\eqref{pre lem commutator trivial}, \eqref{pre lem commutator bar}, and \eqref{proof L2 1ge3}:
\begin{subequations}\label{proof L2 2ge3}
\bel{proof L2 2ge3a}
\sum_{\ih,\alpha}\big{\|}(s/t)Z^I \del_\alpha  u_{\ih}\big{\|}_{L^2(\Hcal_s)} +
\sum_{\ih,\alpha}\big{\|}(s/t)Z^I \delu_\alpha  u_{\ih}\big{\|}_{L^2(\Hcal_s)}
\leq C C_1\eps,
\ee
\bel{proof L2 2ge3b}
\sum_{\ih,a}\big{\|}Z^I \delu_a  u_{\ih}\big{\|}_{L^2(\Hcal_s)}
 \leq C C_1\eps.
\ee
\end{subequations}

The fourth group of estimates concerns  the second order derivatives. The first is deduced from \eqref{pre L2 2order},
 \eqref{proof energy assumption b}, and \eqref{proof energy assumption a},
 while the second is deduced from 
 \eqref{pre L2 2order},
 \eqref{proof energy assumption c}, and \eqref{proof energy assumption d},
and the last one from
\eqref{pre L2 2order} and \eqref{proof energy assumption e}:
\begin{subequations}\label{proof L2 2ge4}
\bel{proof L2 2ge4a}
\sum_{a,\beta,i}\big{\|}s \delu_a\delu_{\beta} Z^{\Id} w_i\big{\|}_{L^2(\Hcal_s)}
+\sum_{a,\beta,i} \big{\|}s \delu_{\beta}\delu_a Z^{\Id} w_i\big{\|}_{L^2(\Hcal_s)} \leq C C_1\eps s^{\delta},
\ee
\bel{proof L2 2ge4a'}
\sum_{a,\beta,i}\big{\|}s \delu_a\delu_{\beta} Z^I w_i\big{\|}_{L^2(\Hcal_s)}
+\sum_{a,\beta,i} \big{\|}s \delu_{\beta}\delu_a Z^I w_i\big{\|}_{L^2(\Hcal_s)} \leq C C_1\eps s^{\delta/2},
\ee
\bel{proof L2 2ge4b}
\sum_{a,\beta,\ih}\big{\|}s \delu_a\delu_{\beta} Z^{\If} u_{\ih}\big{\|}_{L^2(\Hcal_s)}
+\sum_{a,\beta,\ih} \big{\|}s \delu_{\beta}\delu_a Z^{\If}  u_{\ih}\big{\|}_{L^2(\Hcal_s)} \leq C C_1\eps ,
\ee
\end{subequations}
where the order of $\If$ is less or equal to $2$. We can also use the commutator estimates \eqref{pre lem commutator second-order bar} and obtain 
\begin{subequations}
\bel{proof L2 2ge4'a}
\sum_{a,\beta, i}\big{\|}sZ^{\Id} \delu_a\delu_{\beta} w_i \big{\|}_{L^2(\Hcal_s)}
+\sum_{a,\beta,i} \big{\|}sZ^{\Id} \delu_{\beta}\delu_a w_i\big{\|}_{L^2(\Hcal_s)} \leq C C_1\eps s^{\delta},
\ee
\bel{proof L2 2ge4'a'}
\sum_{a,\beta, i}\big{\|}sZ^I \delu_a\delu_{\beta} w_i \big{\|}_{L^2(\Hcal_s)}
+\sum_{a,\beta,i} \big{\|}sZ^I \delu_{\beta}\delu_a w_i\big{\|}_{L^2(\Hcal_s)} \leq C C_1\eps s^{\delta/2},
\ee
\bel{proof L2 2ge4'b}
\sum_{a,\beta,\ih}\big{\|}sZ^{\If} \delu_a\delu_{\beta} u_{\ih}\big{\|}_{\Hcal_s}
+\sum_{a,\beta,\ih} \big{\|}sZ^{\If} \delu_{\beta}\delu_a u_{\ih}\big{\|}_{\Hcal_s} \leq C C_1\eps. 
\ee
\end{subequations}

Finally, we have the $L^2$ estimates for the wave components themselves\footnote{In \eqref{proof L2 2ge5a},
 if we want to be more precise, $C = C'(C_0/C_1 + 1+\delta^{-1})$ where $C'$ is a constant depending only on the structure of the system. We see that it also depends on $C_0$ and $\delta$. But as we can assume that $1/12\leq \delta< 1/6$ and $C_1\geq C_0$, we denote it again by $C$ and regard it as a constant determined by the structure of the system.} 
for all $|\Is| \leq 5$ and all $|\Id|\leq 4$: 
\begin{subequations}
\bel{proof L2 2ge5a}
\big{\|}s^{-1}Z^{\Is}u_{\ih}\big{\|}_{L^2(\Hcal_s)}\leq CC_1\eps s^{\delta},
\ee
\bel{proof L2 2ge5a'}
\big{\|}s^{-1}Z^{\Id}u_{\ih}\big{\|}_{L^2(\Hcal_s)}\leq CC_1\eps s^{\delta/2},
\ee
\bel{proof L2 2ge5b}
\big{\|}t^{-1}Z^I u_{\ih}\big{\|}_{L^2(\Hcal_s)}\leq CC_1\eps.
\ee
\end{subequations}
The first two inequalities are direct consequences of \eqref{proof L2 1ge1a}, \eqref{proof L2 1ge1a'}, \eqref{proof L2 1ge1b}, and 
\eqref{proof L2 1ge1b'} combined with Proposition \ref{pre Hardy hyper}. The last inequality is a result of Lemma \ref{pre Hardy lem1} combined with \eqref{proof L2 1ge3b}.


\section{$L^\infty$ estimates of the first generation}
\label{sec:64}

For convenience in the presentation, we introduce 
a new convention (which is parallel to the index convention already made in \eqref{index-I}):
\begin{equation}
\aligned
\Js & \quad \text{index of order $\leq 3$,} 
\\
J^{\dag} &\quad \text{index of order $\leq 2$,} 
\\
J &\quad \text{index of order $\leq 1$.} 
\endaligned
\end{equation}
Now we give the decay estimates based on the energy assumption \eqref{proof energy assumption} and the Sobolev-type inequality given in Proposition \ref{pre lem sobolev}.

The first group of inequalities is a direct consequence of \eqref{pre lem decay eq},
 \eqref{proof energy assumption a}, and \eqref{proof energy assumption b}:
\begin{subequations}\label{proof decay 1ge1}
\bel{proof decay 1ge1a}
\sup_{\Hcal_s}\Big(s t^{1/2}\big|\del_\alpha Z^{\Js} w_j\big|\Big) + \sup_{\Hcal_s}\Big( s t^{1/2}\big|\delu_\alpha Z^{\Js} w_j\big|\Big)\leq C C_1\eps s^{\delta},
\ee
\bel{proof decay 1ge1a'}
\sup_{\Hcal_s}\Big(s t^{1/2}\big|\del_\alpha Z^{\Jd} w_j\big|\Big) + \sup_{\Hcal_s}\Big( s t^{1/2}\big|\delu_\alpha Z^{\Jd} w_j\big|\Big)\leq C C_1\eps s^{\delta/2},
\ee
\bel{proof decay 1ge1b}
\sup_{\Hcal_s}\Big(t^{3/2}\big|\delu_a Z^{\Js} w_j\big|\Big) \leq C C_1\eps s^{\delta},
\ee
\bel{proof decay 1ge1b'}
\sup_{\Hcal_s}\Big(t^{3/2}\big|\delu_a Z^{\Jd} w_j\big|\Big) \leq C C_1\eps s^{\delta/2},
\ee
\bel{proof decay 1ge1c}
\sup_{\Hcal_s}\Big(t^{3/2}|Z^{\Js}v_{\jc}|\Big)\leq C C_1\eps s^{\delta},
\ee
\bel{proof decay 1ge1c'}
\sup_{\Hcal_s}\Big(t^{3/2}|Z^{\Jd}v_{\jc}|\Big)\leq C C_1\eps s^{\delta/2}.
\ee
\end{subequations}

The second group of estimates is a special case of the second estimate, by taking
 $$
\aligned
& Z^{\Js} = \del_\alpha Z^{\Jd}, \quad Z^{\Js} = L_aZ^{\Jd} = t\delu_aZ^{\Jd}, 
\\
& Z^{\Jd} = \del_\alpha Z^J, \quad 
Z^{\Jd} = L_aZ^J = t\delu_aZ^J
\endaligned
$$
 in \eqref{proof decay 1ge1c} and \eqref{proof decay 1ge1c'}:
\begin{subequations}\label{proof decay 1ge2}
\bel{proof decay 1ge2a}
\aligned
\sup_{\Hcal_s}\Big(t^{3/2}|\del_\alpha Z^{\Jd}v_{\jc}|\Big) + \sup_{\Hcal_s}\Big(t^{3/2}\big|\delu_{\alpha}Z^{\Jd}v_{\jc}\big|\Big)\leq CC_1\eps s^{\delta},
\endaligned
\ee
\bel{proof decay 1ge2a'}
\aligned
\sup_{\Hcal_s}\Big(t^{3/2}|\del_\alpha Z^Jv_{\jc}|\Big) + \sup_{\Hcal_s}\Big(t^{3/2}\big|\delu_{\alpha}Z^Jv_{\jc}\big|\Big)\leq CC_1\eps s^{\delta/2},
\endaligned
\ee
\bel{proof decay 1ge2b}
\aligned
&\sup_{\Hcal_s}\Big(t^{5/2}\big|\delu_aZ^{\Jd}v_{\jc}\big|\Big)\leq CC_1\eps s^{\delta},
\endaligned
\ee
\bel{proof decay 1ge2b'}
\aligned
&\sup_{\Hcal_s}\Big(t^{5/2}\big|\delu_aZ^Jv_{\jc}\big|\Big)\leq CC_1\eps s^{\delta/2}.
\endaligned
\ee
\end{subequations}

The estimates in the third group follow from
\eqref{pre lem decay eq} (a) and (b) combined with \eqref{proof energy assumption e}:
\begin{subequations}\label{proof decay 1ge3}
\bel{proof decay 1ge3a}
  \sup_{\Hcal_s}\Big(s t^{1/2}\big| \del_\alpha Z^J u_{\kh}\big| \Big)
+ \sup_{\Hcal_s}\Big(s t^{1/2}\big|\delu_{\alpha} Z^J u_{\kh}\big|\Big) \leq C C_1\eps,
\ee
\bel{proof decay 1ge3b}
\sup_{\Hcal_s}\Big(t^{3/2}\big|\delu_a Z^J u_{\kh}\big|\Big)
\leq C C_1\eps.
\ee
\end{subequations}

The fourth group concerns the ``second order derivative" of the solution. They are deduced from \eqref{pre decay 2order} and \eqref{proof energy assumption}:
\begin{subequations}\label{proof decay 1ge4}
\bel{proof decay 1ge4a}
\sup_{\Hcal_s}\big(t^{3/2}s\big|\delu_\alpha \delu_aZ^{\Jd} w_j\big| \big)
+\sup_{\Hcal_s}\big(t^{3/2}s|\delu_a\delu_{\alpha}Z^{\Jd}w_j|\big)
\leq CC_1\eps s^{\delta},
\ee
\bel{proof decay 1ge4a'}
\sup_{\Hcal_s}\big(t^{3/2}s\big|\delu_\alpha \delu_aZ^J w_j\big| \big)
+\sup_{\Hcal_s}\big(t^{3/2}s|\delu_a\delu_{\alpha}Z^Jw_j|\big)
\leq CC_1\eps s^{\delta/2},
\ee
\bel{proof decay 1ge4b}
\sup_{\Hcal_s}\big(t^{3/2}s\big|\delu_\alpha \delu_a u_{\jh}\big| \big)
+\sup_{\Hcal_s}\big(t^{3/2}s\big|\delu_a\delu_\alpha u_{\jh}\big| \big)
 \leq CC_1\eps.
\ee
\end{subequations}


\section{$L^\infty$ estimates of the second generation}
\label{sec:65} 

The estimates in this section follow from the first generation $L^\infty$ estimates established in the previous section,  combined with 
the commutator estimates in Lemmas \ref{pre lem commutator1} 
and  \ref{pre lem commutator}.

The first group of estimates are results of  \eqref{pre lem commutator trivial}, \eqref{pre lem commutator bar} combined with \eqref{proof decay 1ge1}:
\begin{subequations}\label{proof decay 2ge1}
\bel{proof decay 2ge1a}
\sup_{\Hcal_s}\big|s t^{1/2}Z^{\Js}\del_\alpha w_j\big| + \sup_{\Hcal_s}\big|st^{1/2} Z^{\Js} \delu_{\alpha}  w_j\big| \leq C C_1\eps s^{\delta},
\ee
\bel{proof decay 2ge1a'}
\sup_{\Hcal_s}\big|s t^{1/2}Z^{\Jd}\del_\alpha w_j\big| + \sup_{\Hcal_s}\big|st^{1/2} Z^{\Jd} \delu_{\alpha}  w_j\big| \leq C C_1\eps s^{\delta/2},
\ee
\bel{proof decay 2ge1b}
\sup_{\Hcal_s}\big|t^{3/2} Z^{\Js} \delu_a  w_j\big| \leq C C_1\eps s^{\delta},
\ee
\bel{proof decay 2ge1b'}
\sup_{\Hcal_s}\big|t^{3/2} Z^{\Jd} \delu_a  w_j\big| \leq C C_1\eps s^{\delta/2},
\ee
\bel{proof decay 2ge1c}
\sup_{\Hcal_s}\big|t^{3/2}Z^{\Js}v_{\kc}\big| \leq C C_1\eps s^{\delta},
\ee
\bel{proof decay 2ge1c'}
\sup_{\Hcal_s}\big|t^{3/2}Z^{\Jd}v_{\kc}\big| \leq C C_1\eps s^{\delta/2}.
\ee
\end{subequations}

The second group consists of the following estimates. They are results of \eqref{pre lem commutator trivial}, \eqref{pre lem commutator bar} combined with \eqref{proof decay 1ge2}: 
\begin{subequations}\label{proof decay 2ge2}
\bel{proof decay 2ge2a}
\sup_{\Hcal_s}\big|t^{3/2}Z^{\Jd} \del_av_{\jc}\big| + \sup_{\Hcal_s}\big|t^{3/2}Z^{\Jd} \delu_{\alpha} v_{\jc}\big|\leq C C_1\eps s^{\delta},
\ee
\bel{proof decay 2ge2a'}
\sup_{\Hcal_s}\big|t^{3/2}Z^J \del_av_{\jc}\big| + \sup_{\Hcal_s}\big|t^{3/2}Z^J \delu_{\alpha} v_{\jc}\big|\leq C C_1\eps s^{\delta/2},
\ee
\bel{proof decay 2ge2b}
\sup_{\Hcal_s}\big|t^{5/2}Z^{\Jd} \delu_av_{\jc}\big|\leq C C_1\eps s^{\delta},
\ee
\bel{proof decay 2ge2b-DEUX}
\sup_{\Hcal_s}\big|t^{5/2}Z^J \delu_av_{\jc}\big|\leq C C_1\eps s^{\delta/2}.
\ee
\end{subequations}

The third group consists of the following estimates which are direct 
consequences of \eqref{pre lem commutator trivial} and \eqref{pre lem commutator bar} combined with \eqref{proof decay 1ge3}: 
\begin{subequations}\label{proof decay 2ge3}
\bel{proof decay 2ge3a}
\sup_{\Hcal_s}\big|s t^{1/2}Z^J \del_\alpha  u_{\jh}\big| + \sup_{\Hcal_s}\big|s t^{1/2}Z^J\delu_{\alpha}  u_{\jh}\big| \leq C C_1\eps,
\ee
\bel{proof decay 2ge3b}
\sup_{\Hcal_s}\big|t^{3/2}Z^J\delu_a  u_{\jh}\big| \leq C C_1\eps. 
\ee
\end{subequations}

The estimates of the fourth group are deduced from \eqref{pre lem commutator second-order bar}, \eqref{proof decay 1ge4a}, \eqref{proof decay 1ge4a'},
 and \eqref{proof decay 1ge4b}: 
\begin{subequations}
\bel{proof decay 2ge4a}
\sup_{\Hcal_s}\big(t^{3/2}s Z^{\Jd}\delu_\alpha \delu_a w_j \big) + \sup_{\Hcal_s}\big(t^{3/2}s Z^{\Jd}\delu_a\delu_\alpha  w_j \big)
\leq CC_1\eps s^{\delta},
\ee
\bel{proof decay 2ge4a'}
\sup_{\Hcal_s}\big(t^{3/2}s Z^J\delu_\alpha \delu_a w_j \big) + \sup_{\Hcal_s}\big(t^{3/2}s Z^J\delu_a\delu_\alpha  w_j \big)
\leq CC_1\eps s^{\delta/2},
\ee
\bel{proof decay 2ge4b}
\sup_{\Hcal_s}\big|st^{3/2}\delu_a\delu_{\beta} u_{\jh}\big| + \sup_{\Hcal_s}\big|st^{3/2}\delu_\alpha \delu_b u_{\jh}\big|
\leq C C_1\eps s.
\ee
\end{subequations}

We finally can write down the decay estimates for the wave components, which follow from Lemma \ref{pre lem u it-self} combined with \eqref{proof energy assumption}:
\begin{subequations}\label{proof decay 2ge5}
\bel{proof decay 2ge5a}
\sup_{\Hcal_s}\Big(t^{3/2}s^{-1}|Z^{\Js} u_{\ih}|\Big)\leq CC_1\eps s^{\delta},
\ee
\bel{proof decay 2ge5a'}
\sup_{\Hcal_s}\Big(t^{3/2}s^{-1}|Z^{\Jd} u_{\ih}|\Big)\leq CC_1\eps s^{\delta/2},
\ee
\bel{proof decay 2ge5b}
\sup_{\Hcal_s}\big(t^{3/2}s^{-1}|Z^J u_{\ih}|\big) \leq CC_1\eps.
\ee
\end{subequations}


\chapter[Second-order derivatives of the wave components]{Second-order derivatives of the wave components \label{cha:7}}

\section{Preliminaries}
\label{sec:71}

The estimates in this chapter concern ``second-order" terms $\del_\alpha \del_{\beta}Z^I u_{\ih}$ and $Z^{I}\del_\alpha \del_{\beta} u_{\ih}$, which will also be used for the 
control of 
$$
u_{\ih} \, \del_{\alpha}\del_{\beta}u_{\jh}.
$$
This chapter is more technical than the derivation of our earlier estimates 
and our strategy now is as follows. 
We are going to analyze the structure of certain second-order derivatives of the wave components in the semi-hyperboloidal frame and concentrate first on
 the component of the Hessian 
 $\delu_0\delu_0 Z^I u_{\ih}$ or, equivalently, $\del_t\del_tZ^I u_{\ih}$. 
Other components have a main part which can be expressed in terms of this component. 
Next, we will analyze $\del_t\del_tZ^I u_{\ih}$ and give a general sketch of the proof, while postponing the technical aspects to the last two sections. 

More precisely, let us first reduce the derivation of an estimating of the Hessian of $Z^Iu_{\ih}$
to an estimate of its component $\delu_0\delu_0Z^Iu_{\ih}$.

\begin{lemma}\label{proof 2order lem1}
If $u$ is a smooth function compactly supported in the half-cone $\Kcal = \{|x|\leq t-1\}$, then for all index $I$ one has 
$$
|\del_\alpha \del_{\beta}Z^I u|
\leq |\del_t\del_t Z^I u|
 + C\sum_{a,\beta}\big|\delu_a\del_{\beta}Z^I u\big| +\frac{C}{t}\sum_{\gamma}|\del_{\gamma}Z^I u|,
$$
where $C$ is a universal constant.
\end{lemma}

\begin{proof}
We recall the following identity based on a change of frame (between 
the natural frame and the semi-hyperboloidal frame):
$$
\aligned
\del_\alpha \del_{\beta}Z^Iu
=& \Psi_\alpha ^{\alpha'}\Psi_{\beta}^{\beta'}\delu_{\alpha'}\delu_{\beta'}Z^Iu + \del_\alpha \big(\Psi_{\beta}^{\beta'}\big)\delu_{\beta}Z^Iu
\\
=&\Psi_\alpha ^0\Psi_{\beta}^0\delu_0\delu_0Z^Iu+ \sum_{b}\Psi_\alpha ^0\Psi_{\beta}^b\delu_0\delu_bZ^Iu
\\
&
 + \sum_{a}\Psi_\alpha ^a\Psi_{\beta}^0\delu_a\delu_0Z^Iu
 + \sum_{a,b}\Psi_\alpha ^a\Psi_{\beta}^b\delu_a\delu_bZ^I u
 + \del_\alpha \big(\Psi_{\beta}^{\beta'}\big)\delu_{\beta}Z^Iu.
\endaligned
$$
By observing that
$\big|\del_\alpha \Psi_{\beta}^{\beta'}\big| \leq \frac{C}{t}$
and
$$
\big|[\delu_b,\del_t](Z^Iu)\big|\leq Ct^{-1}\sum_\alpha \big|\del_\alpha Z^Iu\big|,
$$
we can write 
$$
\sum_b|\del_t\delu_b Z^I u| + \sum_{a,b}|\delu_a\delu_b Z^I u|
\leq
\sum_{a,\beta}|\delu_a\del_{\beta}Z^I u| + Ct^{-1}\sum_{\gamma}|\del_{\gamma}Z^Iu|.
$$ 
\end{proof}

Now we assume that the energy is controlled.

\begin{lemma}
Under the energy assumption \eqref{proof energy assumption}, the following $L^2$ estimates hold for all $|\Id|\leq 4, |I|\leq 3$ and $|\If|\leq 2$:
\begin{subequations}\label{proof 2order general L2 1}
\bel{proof 2order general L2 1a}
\big{\|}s^3t^{-2}\del_\alpha \del_{\beta}Z^{\Id} u_{\ih}\big{\|}_{L^2(\Hcal_s)}
\leq
\big{\|}s^3t^{-2}\del_t\del_tZ^{\Id}u_{\ih}\big{\|}_{L^2(\Hcal_s)} + C C_1\eps s^{\delta},
\ee
\bel{proof 2order general L2 1a'}
\big{\|}s^3t^{-2}\del_\alpha \del_{\beta}Z^I u_{\ih}\big{\|}_{L^2(\Hcal_s)}
\leq
\big{\|}s^3t^{-2}\del_t\del_tZ^Iu_{\ih}\big{\|}_{L^2(\Hcal_s)} + C C_1\eps s^{\delta/2},
\ee
\bel{proof 2order general L2 1b}
\big{\|}s^3t^{-2}\del_\alpha \del_{\beta}Z^{\If} u_{\ih}\big{\|}_{L^2(\Hcal_s)}
 \leq
\big{\|}s^3t^{-2}\del_t\del_tZ^{\If}u_{\ih}\big{\|}_{L^2(\Hcal_s)} + C C_1\eps.
\ee
\end{subequations}
Furthermore, the following $L^\infty$ estimates hold for all $|\Jd|\leq 2$ and $|J| \leq 1$: 
\begin{subequations}\label{proof 2order general decay1}
\bel{proof 2order general decay1a}
\sup_{\Hcal_s}|s^3t^{-1/2}\del_\alpha \del_{\beta}Z^{\Jd}u|
\leq \sup_{\Hcal_s}|s^3t^{-1/2}\del_t\del_t Z^{\Jd} u| + CC_1\eps s^{\delta},
\ee
\bel{proof 2order general decay1a'}
\sup_{\Hcal_s}|s^3t^{-1/2}\del_\alpha \del_{\beta}Z^Ju|
\leq \sup_{\Hcal_s}|s^3t^{-1/2}\del_t\del_t Z^J u| + CC_1\eps s^{\delta/2},
\ee
\bel{proof 2order general decay1b}
\sup_{\Hcal_s}|s^3t^{-1/2}\del_\alpha \del_{\beta} u| \leq \sup_{\Hcal_s}|s^3t^{-1/2}\del_t\del_t u| + CC_1\eps.
\ee
\end{subequations}
\end{lemma}
 
\begin{proof}
The estimate \eqref{proof 2order general L2 1a} is a combination of Lemma \ref{proof 2order lem1} and the estimate \eqref{proof L2 2ge4a},
while  \eqref{proof 2order general L2 1a'} is a combination of Lemma \ref{proof 2order lem1} and the estimate \eqref{proof L2 2ge4a'}. On the other hand, 
 \eqref{proof 2order general L2 1b} is a combination of Lemma \ref{proof 2order lem1} and the estimate \eqref{proof L2 2ge4b}. For the $L^\infty$ estimates,
we combine Lemma \ref{proof 2order lem1} with \eqref{proof decay 1ge4}.
\end{proof}


\section{Analysis of the algebraic structure}\label{section 2order2}
\label{sec:72}

The aim of this section is to establish a general control on $\del_t\del_tZ^I u_{\ih}$,  by making use of the identity \eqref{pre expression of wave under oneframe}. Our strategy is to express $\del_t\del_tZ^I u_{\ih}$ by other terms arising in \eqref{pre expression of wave under oneframe}.
In this first lemma, we relate $\del_t\del_tZ^I u_{\ih}$ to remaining
terms which have faster decay.

\begin{lemma}\label{proof 2order lem2}
Let $U_{tt}(I) :=\big(\del_t\del_tZ^Iu_1,\del_t\del_tZ^Iu_2,\dots,\del_t\del_tZ^Iu_{j_0}\big)^{T}$ which is a $j_0$-vector and $\Ib_{j_0}$ be the $j_0$- dimensional identity matrix. 
The following identity holds:
\bel{proof 2order eq1}
\aligned
& \big((s/t)^2\Ib_{j_0} + \Gb(w,\del w)\big)U_{tt}(I)
\\
& = \big([Z^I,G_{\ih}^{\jh00}\del_t\del_t]u_{\jh} - {Q_G}_{\ih}(I,w,\del w,\del\del w)  + Z^IF_{\ih} + R(Z^I u_{\ih})\big)_{1\leq\ih\leq j_0},
\endaligned
\ee
where
$\big(\Gb(w,\del w)\big)_{1\leq \ih,\jh\leq j_0} = \big(\Gu_{\ih}^{\jh00}(w,\del w)\big)_{1\leq\ih,\jh\leq j_0}$
is a $j_0\times j_0$ order matrix, 
$$
\aligned
& {Q_G}_{\ih}(I,w,\del w,\del\del w)
\\
& =  Z^I\big( \Gu_{\ih}^{\jh a0}(w,\del w) \delu_a\del_t u_{\jh}
+ \Gu_{\ih}^{\jh 0b}(w,\del w) \del_t\delu_b u_{\jh}
+ \Gu_{\ih}^{\jh ab}(w,\del w) \delu_a\delu_b u_{\jh}\big)
\\
&\quad + Z^I\big(G_{\ih}^{\jh\alpha\beta}(w,\del w)\delu_{\beta'}u_{\jh}\del_\alpha \Psi_{\beta}^{\beta'}\big)
+ Z^I\big(G_{\ih}^{\jc\alpha\beta}(w,\del w)\del_\alpha \del_{\beta}v_{\jc}\big)
\endaligned
$$
and
$$
\aligned
R(Z^I u_{\ih}):=
& - \mbar^{0\ar}\delu_0 \delu_{\ar} Z^Iu_{\ih} - \mbar^{\ar 0}\delu_{\ar} \delu_0 Z^Iu_{\ih} -\mbar^{\ar\br}\delu_{\ar}\del_{\br}Z^Iu_{\ih}
\\
&+ \mbar^{\alphar\betar}\big(\delu_{\alphar}\Phi_{\betar}^{\betar'}\big)\delu_{\betar'}Z^Iu_{\ih}.
\endaligned
$$
\end{lemma}

The terms $Z^IF_{\ih}$ and ${Q_G}_{\ih}$ are bilinear and can be expected to enjoy better estimates. The term $R(Z^I u_{\ih})$ contains only  ``good" second-order derivatives and can also be expected to enjoy better estimates.

\begin{proof}
Recall the identity \eqref{pre expression of wave under oneframe} with the function $Z^Iu_{\ih}$
\bel{proof 2order 00}
\aligned
(s/t)^2\delu_0\delu_0 Z^Iu_{\ih}=&
\Box Z^Iu_{\ih} - \mbar^{0\ar}\delu_0 \delu_{\ar} Z^Iu_{\ih} - \mbar^{\ar 0}\delu_{\ar} \delu_0 Z^Iu_{\ih}
\\
& -\mbar^{\ar\br}\delu_{\ar}\del_{\br}Z^Iu_{\ih}
+ \mbar^{\alphar\betar}\big(\delu_{\alphar}\Phi_{\betar}^{\betar'}\big)\delu_{\betar'}Z^Iu_{\ih}
\\
=:&  \Box Z^Iu_{\ih} + R(Z^I,u_{\ih}).
\endaligned
\ee
By equation \eqref{main eq main}, the first term in the right-hand side can be written as
$$
\Box Z^Iu_{\ih} = Z^I\Box u_{\ih}
= -Z^I\big(G_{\ih}^{j\alpha\beta}(w,\del w)\del_\alpha \del_{\beta}w_j\big) + Z^I(F_{\ih}(w,\del w)), 
$$
where
$$
G_{\ih}^{j\alpha\beta}(w,\del w)\del_\alpha \del_{\beta}w_j
= G_{\ih}^{\jh\alpha\beta}(w,\del w)\del_\alpha \del_{\beta}u_{\jh}
  +G_{\ih}^{\jc\alpha\beta}(w,\del w)\del_\alpha \del_{\beta}v_{\jc}
$$
and
$$
 G_{\ih}^{\jh\alpha\beta}(w,\del w)\del_\alpha \del_{\beta}u_{\jh}
 = \Gu_{\ih}^{\jh\alpha\beta}(w,\del w)\delu_\alpha \delu_{\beta}u_{\jh}
 + G_{\ih}^{\jh\alpha\beta}(w,\del w)\delu_{\beta'}u_{\jh}\del_\alpha \Psi_{\beta}^{\beta'}.
$$
We have
$$
\aligned
&Z^I\big(G_{\ih}^{\jh\alpha\beta}(w,\del w)\del_\alpha \del_{\beta}u_{\jh}\big)
\\
&=Z^I\big(\Gu_{\ih}^{\jh\alpha\beta}(w,\del w)\delu_\alpha \delu_{\beta}u_{\jh}\big)
 + Z^I\big(G_{\ih}^{\jh\alpha\beta}(w,\del w)\delu_{\beta'}u_{\jh}\del_\alpha \Psi_{\beta}^{\beta'}\big)
\\
&=Z^I\big(\Gu_{\ih}^{\jh00}\del_t\del_tu_{\jh} + \Gu_{\ih}^{\jh a0}\delu_a\del_tu_{\jh}
+ \Gu_{\ih}^{\jh 0b}\del_t\delu_bu_{\jh} + \Gu_{\ih}^{\jh ab}\delu_a\delu_bu_{\jh}\big)
\\
&\quad + Z^I\big(G_{\ih}^{\jh\alpha\beta}(w,\del w)\delu_{\beta'}u_{\jh}\del_\alpha \Psi_{\beta}^{\beta'}\big)
\\
&=\Gu_{\ih}^{\jh00}(w,\del w)\del_t\del_tZ^Iu_{\jh}
+ [Z^I,\Gu_{\ih}^{\jh00}(w,\del w)\del_t\del_t]u_{\jh}
\\
&\quad + Z^I\big(\Gu_{\ih}^{\jh0b}(w,\del w)\del_t\delu_bu_{\jh} + \Gu_{\ih}^{\jh a0}(w,\del w)\delu_a\del_tu_{\jh}
+ \Gu_{\ih}^{\jh ab}(w,\del w)\delu_a\delu_bu_{\jh}\big)
\\
&\quad 
+ Z^I\big(G_{\ih}^{\jh\alpha\beta}(w,\del w)\delu_{\beta'}u_{\jh}\del_\alpha \Psi_{\beta}^{\beta'}\big).
\endaligned
$$
So we conclude that 
$$
\aligned
\Box Z^Iu_{\ih}
=& - \Gu_{\ih}^{\jh00}(w,\del w)\del_t\del_tZ^Iu_{\jh}- [Z^I,\Gu_{\ih}^{\jh00}(w,\del w)\del_t\del_t]u_{\jh}
\\
&- Z^I\big(\Gu_{\ih}^{\jh a0}(w,\del w) \delu_a\del_t u_{\jh} + \Gu_{\ih}^{\jh 0b}(w,\del w) \del_t\delu_b u_{\jh}
 + \Gu_{\ih}^{\jh ab}(w,\del w) \delu_a\delu_b u_{\jh}\big)
\\
&- Z^I\big(G_{\ih}^{\jh\alpha\beta}(w,\del w)\delu_{\beta'}u_{\jh}\del_\alpha \Psi_{\beta}^{\beta'}\big)
- Z^I\big(G_{\ih}^{\jc\alpha\beta}(w,\del w)\del_\alpha \del_{\beta}v_{\jc}\big)
\\
&+Z^IF_{\ih}(w,\del w).
\endaligned
$$

Substituting this result into the equation \eqref{proof 2order 00}, we obtain
$$
\aligned
&(s/t)^2\del_t\del_tZ^Iu_{\ih} +  \Gu_{\ih}^{\jh00}(w,\del w)\del_t\del_tZ^Iu_{\jh}
\\
&=- [Z^I,\Gu_{\ih}^{\jh00}(w,\del w)\del_t\del_t]u_{\jh}
\\
&\quad - Z^I\big( \Gu_{\ih}^{\jh a0}(w,\del w) \delu_a\del_t u_{\jh} - \Gu_{\ih}^{\jh 0b}(w,\del w) \del_t\delu_b u_{\jh}
 - \Gu_{\ih}^{\jh ab}(w,\del w) \delu_a\delu_b u_{\jh}\big)
\\
&\quad - Z^I\big(G_{\ih}^{\jh\alpha\beta}(w,\del w)\delu_{\beta'}u_{\jh}\del_\alpha \Psi_{\beta}^{\beta'}\big)
- Z^I\big(G_{\ih}^{\jc\alpha\beta}(w,\del w)\del_\alpha \del_{\beta}v_{\jc}\big)
\\
&\quad + Z^IF_{\ih}(w,\del w) + R(Z^I,u_{\ih}). 
\endaligned
$$ 
\end{proof}

Now we derive estimates from the algebraic relation \eqref{proof 2order eq1}. A first step is to get the estimate of the inverse of the linear operator $(\Ib_{j_0} + (t/s)^2\Gb)$.
We can expect that when $|(t/s)^2\Gb|$ is small, $(\Ib_{j_0} + (t/s)^2\Gb)$ is invertible and
we 
can estimate $U_{tt}(I)$ from \eqref{proof 2order eq1}.

\begin{lemma}\label{proof 2order lem3}
There exists a positive constant $\eps_0''$, such that if the following sup-norm estimates
\bel{proof 2order lem3 1}
\aligned
& |\del u_{\ih}|\leq CC_1\eps t^{-1/2}s^{-1}, \quad |\del v_{\ic}|\leq CC_1\eps t^{-3/2}s^{\delta},
\quad
\\
& |v_{\ic}|\leq CC_1\eps t^{-3/2}s^{\delta},\quad |u_{\ic}|\leq CC_1\eps t^{-3/2}s
\endaligned
\ee
hold for $C_1\eps<\eps_0''$, then the following estimate holds:
\bel{proof 2order eq2}
\aligned
\big|\del_t\del_tZ^Iu_{\ih}\big|
\leq
& C(t/s)^2\max_{1\leq \kh\leq j_0}\big\{\big(\big|[Z^I,G_{\kh}^{\jh00}\del_t\del_t]u_{\jh}\big|
+ \big|{Q_G}_{\kh}(I,w,\del w,\del\del w)\big|
\\
&  + \big|Z^IF_{\kh}\big| + \big|R(Z^I u_{\kh})\big|\big)\big\}.
\endaligned
\ee
\end{lemma}

\begin{remark}
\eqref{proof 2order lem3 1} can be guaranteed by the energy assumption \eqref{proof energy assumption} via the $L^\infty$ estimates \eqref{proof decay 2ge1}, \eqref{proof decay 2ge2}, \eqref{proof decay 2ge3},
 and \eqref{proof decay 2ge5} with $C$ determined by the structure of the system.
\end{remark}

\begin{proof}
By the structure of $\Gb(w, \del w)$
$$
\Gu_{\ih}^{\jh00} = \Au_{\ih}^{\jh00\gamma\kh}\delu_{\gamma}u_{\kh} + \Au_{\ih}^{\jh00\gamma \kc}\delu_{\gamma} v_{\kc}
+ \Bu_{\ih}^{\jh00\kh}u_{\kh} + \Bu_{\ih}^{\jh00\kc}v_{\kc}.
$$
Taking into account the assumption of \eqref{proof 2order lem3 1},
$$
|\Gu_{\ih}^{\jh00}|\leq CC_1\eps(s/t)^2.
$$
When $CC_1\eps\leq C\eps_0''$ sufficiently small, the linear operator $\Ib_{j_0} + (t/s)^2\Gb(w,\del w)$ is invertible (viewed as a linear mapping from $(\mathbb{R}^{j_0},\| \, \|_{\infty})$ to itself) and $ \|(\Ib_{j_0} + (t/s)^2\Gb(w,\del w))^{-1} \|_{\infty,\infty}$ is bounded by a fixed constant. 

By Lemma \ref{proof 2order lem2}, we have
$$
\aligned
&(s/t)^2\del_t\del_tZ^Iu_{\ih}
\\
&=  (\Ib_{j_0} + (t/s)^2\Gb(w,\del w))^{-1}
\\
&\quad
\Big( [Z^I,G_{\ih}^{\jh00}\del_t\del_t]u_{\jh} - {Q_G}_{\ih}(I,w,\del w,\del\del w)
 + Z^IF_{\ih} + R(Z^I u_{\ih})\Big)_{1\leq \ih\leq j_0}
\endaligned
$$
and so
$$
\aligned
\big|(s/t)^2\del_t\del_tZ^Iu_{\ih}\big|
\leq
& C\max_{1\leq \kh\leq j_0}\big(\big|[Z^I,G_{\kh}^{\jh00}\del_t\del_t]u_{\jh}\big|
+ \big|{Q_G}_{\kh}(I,w,\del w,\del\del w)\big|  \big)
\\
& + C \big( \big|Z^IF_{\kh}\big| + \big|R(Z^I u_{\i\kh})\big|\big).
\endaligned
$$
\end{proof}

We need to estimate the terms appearing in the right-hand
side of \eqref{proof 2order eq2}. We observe that the term $[Z^I,G_{\ih}^{\jh00}\del_t\del_t]u_{\jh}$ contains the factor $\del_t\del_tZ^J u_{\jh}$ which is also a second-order derivative but with $|J|<|I|$. This structure leads us to the following induction estimates.
Recall again our convention: $Z^I =0 $ when $|I|<0$.

\begin{lemma}\label{proof 2order lem4}
Followed by the notation of Lemma \ref{proof 2order lem2} and \ref{proof 2order lem3}, the following induction estimate holds
\bel{proof 2order eq3}
\aligned
&|\del_t\del_tZ^Iu_{\ih}|
\\
&\leq C(t/s)^2\sum_{|I_1|+|I_2|\leq |I|\atop |I_1|<|I|,\gamma,\jh,i}
\big(|Z^{I_2}\del_{\gamma} w_i|+|Z^{I_2}w_i|\big)|\del_t\del_tZ^{I_1}u_{\jh}|
\\
& \quad+ C(t/s)^2 \max_{1\leq \kh\leq j_0}\big\{{Q_T}_{\kh}(I,w,\del w,\del\del w) + C(t/s)^2|{Q_G}_{\kh}(I,w,\del w\del\del w)|
\\
&\quad + C(t/s)^2|Z^IF_{\kh}|+ C(t/s)^2|R(Z^Iu_{\kh})|\big\},
\endaligned
\ee
where
$$
{Q_T}_{\ih} = \sum_{|I_1|+|I_2|\leq |I|\atop|I_1|<|I|,\gamma,j,\jh}\sum_{a,\beta,\gamma'}
\big(|Z^{I_2}\del_{\gamma}w_j| + |Z^{I_2}w_j|\big)\big(|\delu_a\delu_{\beta} Z^{I_1}u_{\jh}| + t^{-1}|\del_{\gamma'}Z^{I_1}u_{\jh}|\big)
$$
\end{lemma}

\begin{proof}
This is purely an estimate for the term $[Z^I,\Gu_{\ih}^{\jh00}\del_t\del_t]u_{\jh}$, derived
 as follows:
\bel{proof 2order lem 4 eq1}
[Z^I,\Gu_{\ih}^{\jh00}\del_t\del_t]u_{\jh} =
\sum_{I_1+I_2=I\atop |I_1|< |I|}Z^{I_2}\big(\Gu_{\ih}^{\jh00}\big) \,  Z^{I_1}\big(\del_t\del_tu_{\jh}\big) + \Gu_{\ih}^{\jh00}[Z^I,\del_t\del_t]u_{\jh}.
\ee
We have 
$$
\aligned
\big|Z^{I_2}\big(\Gu_{\ih}^{\jh00}\big)\big|
\leq&
\big|Z^{I_2}\big(\Au_{\ih}^{\jh00\gamma k}\delu_{\gamma}w_k\big)\big| + \big|Z^{I_2}\big(\Bu_{\ih}^{\jh00 k}w_k\big)\big|
\\
\leq& C\sum_{j,\gamma \atop  |I_3|\leq |I_2|}\big(\big|Z^{I_3}\del_{\gamma}w_j\big| + \big|Z^{I_3}w_j\big|\big). 
\endaligned
$$
 
The first term of the right-hand side of \eqref{proof 2order lem 4 eq1} is bounded by
$$
\aligned
&C\sum_{|I_1|+|I_2|\leq |I|\atop|I_1|<|I|,\gamma,j,\jh}\big(|Z^{I_2}\del_{\gamma}w_j| + |Z^{I_2}w_j|\big)|\del_t\del_t Z^{I_1}u_{\jh}|
\\
&+C\sum_{|I_1|+|I_2|\leq |I|\atop|I_1|<|I|,\gamma,j,\jh}\sum_{a,\beta,\gamma'}
\big(|Z^{I_2}\del_{\gamma}w_j| + |Z^{I_2}w_j|\big)\big(|\delu_a\delu_{\beta} Z^{I_1}u_{\jh}| + t^{-1}|\del_{\gamma'}Z^{I_1}u_{\jh}|\big)
\endaligned
$$
Similarly, the second term of the right-hand side of \eqref{proof 2order lem 4 eq1} is estimated as follows, by
\eqref{pre lem commutator second-order} and Lemma \ref{proof 2order lem1}: 
$$
\aligned
\big|[Z^I,\del_t\del_t]u_{\jh}\big| 
&\leq C\sum_{\alpha,\beta\atop |J|<|I|}|\del_{\alpha}\del_{\beta} Z^Ju_{\jh}|
\\
&\leq C\sum_{|J|<|I|}|\del_t\del_t Z^Ju_{\jh}| + C\sum_{a,\beta\atop |J|<|I|}|\delu_a\delu_{\beta} Z^Ju_{\jh}|
\\
&\quad
+ Ct^{-1}\sum_{\gamma\atop |J|<|I|}|\del_{\gamma}Z^Ju_{\jh}|.
\endaligned
$$
Combined with the following estimate, we have 
$$
|\Gu_{\ih}^{\jh00}|\leq C\sum_{\jh,\gamma}\big(|\del_{\gamma}w_j| + |w_j|\big). 
$$ 
\end{proof}


\section{Structure of the quadratic terms} 
\label{sec:73}

The aim of this section is to analyze the structure of the quadratic terms ${Q_T}_{\ih}$, ${Q_G}_{\ih}$ and $Z^IF_{\ih}$. We emphasize that some components of these terms satisfy the null condition and we can
use it to make the estimates a bit simpler, but we prefer to {\sl avoid the use of the null structure} here
in order to show the independence of these analysis on this structure.

First, we note that the terms ${Q_T}_{\ih}$ are 
linear combinations of the following terms with constant coefficients with $|I_1| + |I_2| \leq |I|$
and $|I_1| < |I|$
\bel{proof 2order eq QT}
\aligned
&|Z^{I_2}\del_{\gamma}u_{\kh}|\,|\delu_a\delu_{\beta}Z^{I_1}u_{\jh}|,\quad |Z^{I_2}\del_{\gamma}v_{\kc}|\,|\delu_a\delu_{\beta}Z^{I_1}u_{\jh}|,\quad
\\&
|Z^{I_2}u_{\kh}|\,|\delu_a\delu_{\beta}Z^{I_1}u_{\jh}|,\quad
|Z^{I_2}v_{\kc}|\,|\delu_a\delu_{\beta}Z^{I_1}u_{\jh}|,\quad
\\
&t^{-1}|Z^{I_2}\del_{\gamma}u_{\kh}|\,|\del_{\gamma'}Z^{I_1}u_{\jh}|,\quad
 t^{-1}|Z^{I_2}\del_{\gamma}v_{\kc}|\,|\del_{\gamma'}Z^{I_1}u_{\jh}|,\quad
\\
&
 t^{-1}|Z^{I_2}u_{\kh}|\,|\del_{\gamma'}Z^{I_1}u_{\jh}|,\quad
 t^{-1}|Z^{I_2}v_{\kc}|\,|\del_{\gamma'}Z^{I_1}u_{\jh}|. 
\endaligned
\ee
We consider ${Q_G}_{\ih}$ and write 
$$
\aligned
&{Q_G}_{\ih}(I,w,\del w,\del\del w)
\\
&=  Z^I\big( \Gu_{\ih}^{\jh a0}(w,\del w) \delu_a\del_t u_{\jh}
+ \Gu_{\ih}^{\jh 0b}(w,\del w) \del_t\delu_b u_{\jh}
+ \Gu_{\ih}^{\jh ab}(w,\del w) \delu_a\delu_b u_{\jh}\big)
\\
&\quad + Z^I\big(G_{\ih}^{\jh\alpha\beta}(w,\del w)\delu_{\beta'}u_{\jh}\del_\alpha \Psi_{\beta}^{\beta'}\big)
+ Z^I\big(G_{\ih}^{\jc\alpha\beta}(w,\del w)\del_\alpha \del_{\beta}v_{\jc}\big)
\\
&= Z^I\big(\Au_{\ih}^{\jh a0\gamma \kh}\delu_{\gamma}u_{\kh}\delu_a\del_tu_{\jh}
         + \Bu_{\ih}^{\jh a0 \kh}u_{\kh}\delu_a\del_tu_{\jh}
         + \Au_{\ih}^{\jh a0\gamma \kc}\delu_{\gamma}v_{\kc}\delu_a\del_tu_{\jh}
\\
        &\qquad + \Bu_{\ih}^{\jh a0 \kc}v_{\kc}\delu_a\del_tu_{\jh}\big)
\\
 &\quad +Z^I\big(\Au_{\ih}^{\jh 0b\gamma \kh}\delu_{\gamma}u_{\kh}\del_t\delu_bu_{\jh}
         + \Bu_{\ih}^{\jh 0b \kh}u_{\kh}\del_t\delu_bu_{\jh}
         + \Au_{\ih}^{\jh 0b\gamma \kc}\delu_{\gamma}v_{\kc}\del_t\delu_bu_{\jh}
\\
        &\quad + \Bu_{\ih}^{\jh 0b \kc}v_{\kc}\del_t\delu_bu_{\jh}\big)
\\
 &\quad +Z^I\big(\Au_{\ih}^{\jh ab\gamma \kh}\delu_{\gamma}u_{\kh}\delu_a\delu_bu_{\jh}
         + \Bu_{\ih}^{\jh ab \kh}u_{\kh}\delu_a\delu_bu_{\jh}
         + \Au_{\ih}^{\jh ab\gamma \kc}\delu_{\gamma}v_{\kc}\delu_a\delu_bu_{\jh}
\\
        &\quad + \Bu_{\ih}^{\jh ab \kc}v_{\kc}\delu_a\delu_bu_{\jh}\big)
\\
 &\quad +Z^I\big(A_{\ih}^{\jh\alpha\beta\gamma\kh}\del_{\gamma}u_{\kh}\delu_{\beta'}u_{\jh}\del_\alpha \Psi_{\beta}^{\beta'}
          +B_{\ih}^{\jh\alpha\beta\kh}u_{\kh}\delu_{\beta'}u_{\jh}\del_\alpha \Psi_{\beta}^{\beta'}
\\
& \quad 
          +A_{\ih}^{\jh\alpha\beta\gamma\kc}\del_{\gamma}v_{\kc}\delu_{\beta'}u_{\jh}\del_\alpha \Psi_{\beta}^{\beta'}
          +B_{\ih}^{\jh\alpha\beta\kc}v_{\kc}\delu_{\beta'}u_{\jh}\del_\alpha \Psi_{\beta}^{\beta'} \big)
\\
 &\quad +Z^I\big(A_{\ih}^{\jc \alpha\beta\gamma\kh}\del_{\gamma}u_{\kh}\del_\alpha \del_{\beta}v_{\jc}
          +A_{\ih}^{\jc \alpha\beta\gamma\kc}\del_{\gamma}v_{\kc}\del_\alpha \del_{\beta}v_{\jc}
          +B_{\ih}^{\jc \alpha\beta\kc}v_{\kc}\del_\alpha \del_{\beta}v_{\kc}\big).
\endaligned
$$
Let $K := \max_{\alpha,\beta,\gamma, i,j,k}\{A_i^{j\alpha\beta\gamma k},B_i^{j\alpha\beta}\}$ and remark that by \eqref{pre lem frame},
$$
|Z^I \Au_{\ih}^{\jh \alpha\beta\gamma k}| + |Z^I\Bu_{\ih}^{\jh\alpha\beta k}|\leq C(I)K.
$$

We observe that the terms 
${Q_G}_{\ih}(I,w,\del w,\del\del w)$ are linear combinations with bounded coefficients of the following terms:
\bel{proof 2order eq4 QG}
\aligned
&Z^I\big(u_{\ih}\delu_a\delu_{\beta}u_{\jh}\big),\quad Z^I\big(v_{\ic}\delu_a\delu_{\beta}u_{\jh}\big),
\quad Z^I\big(u_{\ih}\del_t\delu_au_{\jh}\big),\quad Z^I\big(v_{\ic}\del_t\delu_au_{\jh}\big),
\\
&Z^I\big(\delu_{\gamma}u_{\ih}\delu_a\delu_{\beta}u_{\jh}\big),\quad Z^I\big(\delu_{\gamma}v_{\ic}\delu_a\delu_{\beta}u_{\jh}\big),
\quad Z^I\big(\delu_{\gamma}u_{\ih}\del_t\delu_au_{\jh}\big),\quad Z^I\big(\delu_{\gamma}v_{\ic}\del_t\delu_au_{\jh}\big),
\\
&Z^I\big(\del_\alpha \Psi_{\beta}^{\beta'}u_{\ih}\delu_{\beta'}u_{\jh}\big),\quad
Z^I\big(\del_\alpha \Psi_{\beta}^{\beta'}v_{\ic}\delu_{\beta'}u_{\jh}\big),\quad
\\
&
Z^I\big(\del_\alpha \Psi_{\beta}^{\beta'}\del_{\gamma}u_{\ih}\delu_{\beta'}u_{\jh}\big),\quad
Z^I\big(\del_\alpha \Psi_{\beta}^{\beta'}\del_{\gamma}v_{\ic}\delu_{\beta'}u_{\jh}\big),
\\
&Z^I\big(v_{\ic}\del_\alpha \del_{\beta}v_{\jc}\big),
\quad Z^I\big(\del_{\gamma}u_{\ih}\del_\alpha \del_{\beta}v_{\jc}\big),\quad Z^I\big(\del_{\gamma}v_{\ic}\del_\alpha \del_{\beta}v_{\jc}\big). 
\endaligned
\ee

For the term $Z^IF_{\ih}$, recalling its definition
$$
Z^IF_{\ih} = Z^I\big(P_{\ih}^{\alpha\beta jk}\del_\alpha w_j\del_{\beta}w_k + Q_{\ih}^{\alpha j\kc}v_{\kc}\del_\alpha w_j + R_{\ih}^{\jc\kc}v_{\jc}v_{\kc}\big), 
$$
we classify its components into the following groups:
$$
\aligned
\del_\alpha u_{\jh}\del_{\beta}u_{\kh},\quad \del_\alpha v_{\jc}\del_{\beta}w_k, \quad
v_{\kc}\del_\alpha w_j,\quad
v_{\jc}v_{\kc}.
\endaligned
$$
We regard $Z^IF_{\ih}$ as a linear combination with constant coefficients bounded by $K$ of the following terms:
\bel{proof 2order eq5 ZF}
\aligned
&Z^I\big(\del_\alpha u_{\jh}\del_{\beta}u_{\kh}\big),\quad Z^I\big(\del_\alpha v_{\jc}\del_{\beta}w_k\big),
\\
&Z^I\big(v_{\kc}\del_\alpha u_{\jh}\big),\quad
Z^I\big(v_{\jc}v_{\kc}\big).
\endaligned
\ee
The estimates on $\del_t\del_tZ^I u_{\ih}$ turn out to be the estimates on these terms listed in \eqref{proof 2order eq4 QG} and \eqref{proof 2order eq5 ZF} and $R(Z^I u_{\kh})$.


\section{$L^\infty$ estimates}
\label{sec:74}

The purpose of the section is to establish the $L^\infty$ estimates on $\del_t\del_tZ^Ju_{\ih}$ (and then $\del_\alpha \del_{\beta}Z^Ju_{\ih}$) under the energy assumption \eqref{proof energy assumption}. 
We need to combine estimates on the terms ${Q_T}_{\ih}$, ${Q_G}_{\ih}$, $Z^{\Jd}F_{\ih}$ and $R(Z^{\Jd}u_{\ih})$ with Lemma \ref{proof 2order lem4}. So, we first consider these terms in the following two lemmas.

\begin{lemma}\label{proof 2order lem5}
Under the energy assumption \eqref{proof energy assumption}, the following estimates hold for any $|\Jd|\leq 2$
\begin{subequations}\label{proof 2order eq6}
\bel{proof 2order eq6a}
|Z^{\Jd}F_{\ih}|\leq C(C_1\eps)^2 t^{-3/2}s^{-1+\delta},
\ee
\bel{proof 2order eq6b}
|{Q_G}_{\ih}(\Jd,w,\del w,\del\del w)|\leq C(C_1\eps)^2 t^{-3/2}s^{-1+\delta},
\ee
\bel{proof 2order eq6bc}
{Q_T}_{\ih}(\Jd,w,\del w,\del\del w)\leq C(C_1\eps)^2 t^{-3/2}s^{-1+\delta},
\ee
\bel{proof 2order eq6c}
|R(Z^{\Jd}u_{\ih})|\leq CC_1\eps t^{-3/2}s^{-1+\delta}.
\ee
\end{subequations}
\end{lemma}

\begin{lemma}\label{proof 2order lem5'}
Under the energy assumption \eqref{proof energy assumption}, the following estimates hold for any $|J|\leq 1$
\begin{subequations}\label{proof 2order eq6'}
\bel{proof 2order eq6'a}
|Z^JF_{\ih}|\leq C(C_1\eps)^2 t^{-3/2}s^{-1+\delta/2},
\ee
\bel{proof 2order eq6'b}
|{Q_G}_{\ih}(J,w,\del w,\del\del w)|\leq C(C_1\eps)^2 t^{-3/2}s^{-1+\delta/2},
\ee
\bel{proof 2order eq6'bc}
{Q_T}_{\ih}(J,w,\del w,\del\del w)\leq C(C_1\eps)^2 t^{-3/2}s^{-1+\delta/2},
\ee
\bel{proof 2order eq6'c}
|R(Z^Ju_{\ih})|\leq CC_1\eps t^{-3/2}s^{-1+\delta/2}.
\ee
\end{subequations}
\end{lemma}

\begin{lemma}\label{proof 2order lem6}
Under the energy assumption \eqref{proof energy assumption}, the following estimates hold:\
\begin{subequations}\label{proof 2order eq7}
\bel{proof 2order eq7a}
|F_{\ih}|\leq C(C_1\eps)^2 t^{-3/2}s^{-1},
\ee
\bel{proof 2order eq7b}
|{Q_G}_{\ih}(0,w,\del w,\del\del w)|\leq C(C_1\eps)^2 t^{-3/2}s^{-1},
\ee
\bel{proof 2order eq7bc}
{Q_T}_{\ih}(0,w,\del w,\del\del w)\leq C(C_1\eps)^2 t^{-3/2}s,
\ee
\bel{proof 2order eq7c}
|R(u_{\ih})|\leq CC_1\eps t^{-3/2}s^{-1}.
\ee
\end{subequations}
\end{lemma}

The proofs of these three lemmas are essentially the same. We compute the relevant terms and express them as linear combinations of some bilinear terms, then suitably estimate each of them. The difference of the decay rate between these three lemmas is due to the difference of regularity.

\begin{proof}[Proof of Lemma \ref{proof 2order lem5}]
The control of ${Q_G}_{\ih}(\Jd,w,\del w,\del\del w)$ is obtained as follows.
Recall that ${Q_G}_{\ih}(\Jd,w,\del w,\del\del w)$ is a linear combination of the terms listed in \eqref{proof 2order eq4 QG} with the index $I$ replaced by $\Jd$. We estimate 
$Z^{\Jd}\big(u_{\ih}\delu_a\del_tu_{\jh}\big)$ in details:
$$
\aligned
\big|Z^{\Jd}\big(u_{\ih}\delu_a\del_tu_{\jh}\big)\big| 
&\leq
\sum_{J_2+J_3=\Jd}\big|Z^{J_2}u_{\ih}\big| \big|Z^{J_3}\delu_a\del_tu_{\jh}\big|
\\
&= \big|u_{\ih}\big| \big|Z^{\Jd}\delu_a\del_tu_{\jh}\big|
 + \big|Z^{J_2}u_{\ih}\big| \big|Z^{J_3}\delu_a\del_tu_{\jh}\big|
 + \big|Z^{\Jd}u_{\ih}\big| \big|\delu_a\del_tu_{\jh}\big|.
\endaligned
$$
When $J_2=0$ and $J_3=\Jd$, we apply \eqref{proof decay 2ge5b} and \eqref{proof decay 2ge4a};
     when $J_3=0$ and $J_2=\Jd$, we apply \eqref{proof decay 2ge5a} and \eqref{proof decay 2ge4b};
and  when $|J_2|=1$ and $|J_3|=1$, we apply \eqref{proof decay 2ge5b} and \eqref{proof decay 2ge4a};
$$
\big|Z^{\Jd}\big(w_i\delu_a\del_tu_{\jh}\big)\big|\leq C(C_1\eps)^2t^{-3+\delta}.
$$

For the other terms, we will specify the $L^\infty$ estimates to be used but omit the details.

Here the in the different columns represent the different partitions of $\Jd = J_2+J_3$. $(a,\leq b)$ means $|J_2|=a,\,|J_3|\leq b$. In the last column, we specify the decay rate obtained directly by applying the given inequalities (modulo the constant $C(C_1\eps)^2$). The first coefficient is the one we used for
 the first factor, while  the second coefficient is used for the second factor: 
$$
\begin{array}{ccccc}
\text{Products}&(2,\leq 0)&(1,\leq 1)&(0,\leq 2)&\text{Decay rate}
\\
u_{\ih}\delu_a\delu_{\beta}u_{\jh}
&\eqref{proof decay 2ge5a},\eqref{proof decay 2ge4b}
&\eqref{proof decay 2ge5b},\eqref{proof decay 2ge4a}
&\eqref{proof decay 2ge5b},\eqref{proof decay 2ge4a}
&t^{-3}s^{\delta}
\\
v_{\ic}\delu_a\delu_{\beta}u_{\jh}
&\eqref{proof decay 2ge1c},\eqref{proof decay 2ge4b}
&\eqref{proof decay 2ge1c},\eqref{proof decay 2ge4a}
&\eqref{proof decay 2ge1c},\eqref{proof decay 2ge4a}
&t^{-3}s^{-1+2\delta}
\\
u_{\ih}\del_t\delu_au_{\jh}
&\eqref{proof decay 2ge5a},\eqref{proof decay 2ge4b}
&\eqref{proof decay 2ge5b},\eqref{proof decay 2ge4a}
&\eqref{proof decay 2ge5b},\eqref{proof decay 2ge4a}
&t^{-3}s^{\delta}
\\
v_{\ic}\del_t\delu_au_{\jh}
&\eqref{proof decay 2ge1c},\eqref{proof decay 2ge4b}
&\eqref{proof decay 2ge1c},\eqref{proof decay 2ge4a}
&\eqref{proof decay 2ge1c},\eqref{proof decay 2ge4a}
&t^{-3}s^{-1+2\delta}
\\
\delu_{\gamma}v_{\ic}\delu_a\delu_{\beta}u_{\jh}
&\eqref{proof decay 2ge2a},\eqref{proof decay 2ge4b}
&\eqref{proof decay 2ge2a},\eqref{proof decay 2ge4a}
&\eqref{proof decay 2ge2a},\eqref{proof decay 2ge4a}
&t^{-3}s^{-1+2\delta}
\\
\delu_{\gamma}u_{\ih}\delu_a\delu_{\beta}u_{\jh}
&\eqref{proof decay 2ge1a},\eqref{proof decay 2ge4b}
&\eqref{proof decay 2ge3a},\eqref{proof decay 2ge4a}
&\eqref{proof decay 2ge3a},\eqref{proof decay 2ge4a}
&t^{-2}s^{-2+\delta}
\end{array}
$$

$$
\begin{array}{ccccc}
\text{Products} &(2,\leq 0) &(1,\leq 1) &(0,\leq 2) &\text{Decay rate}
\\
\delu_{\gamma}u_{\ih}\del_t\delu_au_{\jh}
&\eqref{proof decay 2ge1a},\eqref{proof decay 2ge4b}
&\eqref{proof decay 2ge3a},\eqref{proof decay 2ge4a}
&\eqref{proof decay 2ge3a},\eqref{proof decay 2ge4a}
&t^{-2}s^{-2+\delta}
\\
\delu_{\gamma}v_{\ic}\del_t\delu_au_{\jh}
&\eqref{proof decay 2ge2a},\eqref{proof decay 2ge4b}
&\eqref{proof decay 2ge2a},\eqref{proof decay 2ge4a}
&\eqref{proof decay 2ge2a},\eqref{proof decay 2ge4a}
&t^{-3}s^{-1+2\delta}
\\
v_{\ic}\del_\alpha \del_{\beta}v_{\jc}
&\eqref{proof decay 2ge1c},\eqref{proof decay 2ge2a}
&\eqref{proof decay 2ge1c},\eqref{proof decay 2ge2a}
&\eqref{proof decay 2ge1c},\eqref{proof decay 2ge2a}
&t^{-3}s^{2\delta}
\\
\del_{\gamma}u_{\ih}\del_\alpha \del_{\beta}v_{\jc}
&\eqref{proof decay 2ge1a},\eqref{proof decay 2ge2a}
&\eqref{proof decay 2ge3a},\eqref{proof decay 2ge2a}
&\eqref{proof decay 2ge3a},\eqref{proof decay 2ge2a}
&t^{-2}s^{-1+2\delta}
\\
\del_{\gamma}v_{\ic}\del_\alpha \del_{\beta}v_{\jc}
&\eqref{proof decay 2ge2a},\eqref{proof decay 2ge2a}
&\eqref{proof decay 2ge2a},\eqref{proof decay 2ge2a}
&\eqref{proof decay 2ge2a},\eqref{proof decay 2ge2a}
&t^{-3}s^{2\delta}
\end{array}
$$
Taking into account the fact that $s\leq Ct\leq Cs^2$ and $\delta<1/6$, we
conclude that these terms are bounded by $C(C_1\eps)^2t^{-3/2}s^{-1+\delta}$.

Again, we have the following four terms from ${Q_G}_{\ih}$, which are estimated separately:
$$
\aligned
& Z^{\Jd}\big(\del_\alpha \Psi_{\beta}^{\beta'}u_{\ih}\delu_{\beta'}u_{\jh}\big),\quad Z^{\Jd}\big(\del_\alpha \Psi_{\beta}^{\beta'}v_{\ic}\delu_{\beta'}u_{\jh}\big),\quad
\\
&
Z^{\Jd}\big(\del_\alpha \Psi_{\beta}^{\beta'}\del_{\gamma}u_{\ih}\delu_{\beta'}u_{\jh}\big),\quad
Z^{\Jd}\big(\del_\alpha \Psi_{\beta}^{\beta'}\del_{\gamma}v_{\ic}\delu_{\beta'}u_{\jh}\big).
\endaligned
$$
By observing that $\big|Z^I\del_\alpha \Psi_{\beta}^{\beta'}\big|\leq C(I)t^{-1}$,
these terms can be  estimated by $C(C_1\eps)^2t^{-3}s^{2\delta}$. We give the proof for $Z^{J_1}\big(\del_\alpha \Psi_{\beta}^{\beta'}u_{\ih}\delu_{\beta'}u_{\jh}\big)$:
$$
\aligned
&\big|Z^{\Jd}\big(\del_\alpha \Psi_{\beta}^{\beta'}u_{\ih}\delu_{\beta'}u_{\jh}\big)\big|
\leq \sum_{J_2+J_3+J_4=J_1}\big|Z^{J_4}\del_\alpha \Psi_{\beta}^{\beta'}\big| \,  \big|Z^{J_2}u_{\ih}\big| \,  \big|Z^{J_3}\delu_{\beta'}u_{\jh}\big|
\\
&\leq Ct^{-1}\sum_{|J_2|+|J_3|\leq |\Jd|} \big|Z^{J_2}u_{\ih}\big| \,  \big|Z^{J_3}\delu_{\beta'}u_{\jh}\big|
\leq Ct^{-1}C_1\eps t^{-3/2}s^{1+\delta}  C_1\eps t^{-1/2}s^{-1+\delta}
\\
&=C(C_1\eps)^2t^{-3}s^{2\delta}
\leq C(C_1\eps)^2t^{-3/2}s^{-1+\delta}, 
\endaligned
$$
where we recall that $\delta<1/6$, and \eqref{proof decay 2ge5a} and \eqref{proof decay 2ge1a} are been used. The other terms are estimated similarly, and we omit the details and only
 list out the inequalities we use for each term and each partition of indices; cf. Table 1.


\protect\begin{landscape}\thispagestyle{empty}
$$
\begin{array}{ccccc}
\text{Terms} &(2,\leq 0) &(1,\leq 1) &(0,\leq 2) &\text{ Decay rate}
\\
t^{-1}Z^{J_2}u_{\ih}Z^{J_3}\delu_{\beta'}u_{\jh}
&\eqref{proof decay 2ge5a},\eqref{proof decay 2ge3a}
&\eqref{proof decay 2ge5b},\eqref{proof decay 2ge3a}
&\eqref{proof decay 2ge5b},\eqref{proof decay 2ge1a}
&t^{-3}s^{\delta}
\\
t^{-1}Z^{J_2}v_{\ic}Z^{J_3}\delu_{\beta'}u_{\jh}
&\eqref{proof decay 2ge1c},\eqref{proof decay 2ge3a}
&\eqref{proof decay 2ge1c},\eqref{proof decay 2ge3a}
&\eqref{proof decay 2ge1c},\eqref{proof decay 2ge1a}
&t^{-3}s^{-1+2\delta}
\\
t^{-1}Z^{J_2}\del_{\gamma}u_{\ih}Z^{J_3}\delu_{\beta'}u_{\jh}
&\eqref{proof decay 2ge1a},\eqref{proof decay 2ge3a}
&\eqref{proof decay 2ge3a},\eqref{proof decay 2ge3a}
&\eqref{proof decay 2ge3a},\eqref{proof decay 2ge1a}
&t^{-2}s^{-2+\delta}
\\
t^{-1}Z^{J_2}\del_{\gamma}v_{\ic}Z^{J_3}\delu_{\beta'}u_{\jh}
&\eqref{proof decay 2ge1a},\eqref{proof decay 2ge3a}
&\eqref{proof decay 2ge1a},\eqref{proof decay 2ge3a}
&\eqref{proof decay 2ge1a},\eqref{proof decay 2ge1a}
&t^{-2}s^{-2+\delta}
\end{array}
$$
\centerline{Table 1} 

\

\ 

$$
\begin{array}{ccccc}
\text{Products} &(2,\leq0) &(1,\leq1) &(0,\leq2) &\text{Decay rate}
\\
\delu_a\delu_{\beta}Z^{J_1}u_{\ih}\,Z^{J_2}\del_{\gamma}u_{\kh}
&\eqref{proof decay 1ge4a},\eqref{proof decay 2ge3a}
&\eqref{proof decay 1ge4a},\eqref{proof decay 2ge3a}
&\eqref{proof decay 1ge4b},\eqref{proof decay 2ge1a}
&t^{-2}s^{-2+\delta}
\\
\delu_a\delu_{\beta}Z^{J_1}u_{\ih}\,Z^{J_2}\del_{\gamma}v_{\kc}
&\eqref{proof decay 1ge4a},\eqref{proof decay 2ge2a}
&\eqref{proof decay 1ge4a},\eqref{proof decay 2ge2a}
&\eqref{proof decay 1ge4b},\eqref{proof decay 2ge2a}
&t^{-3}s^{-1+2\delta}
\\
\delu_a\delu_{\beta}Z^{J_1}u_{\ih}\,Z^{J_2}u_{\kh}
&\eqref{proof decay 1ge4a},\eqref{proof decay 2ge5b}
&\eqref{proof decay 1ge4a},\eqref{proof decay 2ge5b}
&\eqref{proof decay 1ge4b},\eqref{proof decay 2ge5a}
&t^{-3}s^{\delta}
\\
\delu_a\delu_{\beta}Z^{J_1}u_{\ih}\,Z^{J_2}v_{\kc}
&\eqref{proof decay 1ge4a},\eqref{proof decay 2ge1c}
&\eqref{proof decay 1ge4a},\eqref{proof decay 2ge1c}
&\eqref{proof decay 1ge4b},\eqref{proof decay 2ge1c}
&t^{-3}s^{-1+2\delta}
\end{array}
$$
\centerline{Table 2} 

\

\

$$
\begin{array}{ccccc}
\text{Products} &(2,\leq0) &(1,\leq1) &(0,\leq2) &\text{Decay rate}
\\
t^{-1}\del_{\gamma'}Z^{J_1}u_{\jh}\,Z^{J_2}\del_{\gamma}u_{\kh}
&\eqref{proof decay 1ge1a},\eqref{proof decay 2ge3a}
&\eqref{proof decay 1ge3a},\eqref{proof decay 2ge3a}
&\eqref{proof decay 1ge3a},\eqref{proof decay 2ge1a}
&t^{-2}s^{-2+\delta}
\\
t^{-1}\del_{\gamma'}Z^{J_1}u_{\jh}\,Z^{J_2}\del_{\gamma}v_{\kc}
&\eqref{proof decay 1ge1a},\eqref{proof decay 2ge2a}
&\eqref{proof decay 1ge3a},\eqref{proof decay 2ge2a}
&\eqref{proof decay 1ge3a},\eqref{proof decay 2ge2a}
&t^{-3}s^{-1+2\delta}
\\
t^{-1}\del_{\gamma'}Z^{J_1}u_{\jh}\,Z^{J_2}u_{\kh}
&\eqref{proof decay 1ge1a},\eqref{proof decay 2ge5b}
&\eqref{proof decay 1ge3a},\eqref{proof decay 2ge5b}
&\eqref{proof decay 1ge3a},\eqref{proof decay 2ge5a}
&t^{-3}s^{\delta}
\\
t^{-1}\del_{\gamma'}Z^{J_1}u_{\jh}\,Z^{J_2}v_{\kc}
&\eqref{proof decay 1ge1a},\eqref{proof decay 2ge1c}
&\eqref{proof decay 1ge3a},\eqref{proof decay 2ge1c}
&\eqref{proof decay 1ge3a},\eqref{proof decay 2ge1c}
&t^{-3}s^{-1+2\delta}
\end{array}
$$
\centerline{Table 3} 

\end{landscape}

The term ${Q_T}_{\ih}$ a linear combination of the terms presented in \eqref{proof 2order eq QT},
 with $I_i$ replaced by $J_i$ and $i=1,2$. 
Each term is estimated as for $Q_G$. We give the list of inequalities 
we use for every partition $|J_1|+|J_2|\leq |\Jd|$: 
 in Table 2 and Table 3 the symbol $(a,\leq b)$
means $|J_1| = a, |J_2|\leq b$. 
We conclude with \eqref{proof 2order eq6bc}.

The estimates on $Z^{\Jd}F_{\ih}$ are essentially the same. Recall that $Z^{\Jd}F_{\ih}$ is a linear combination of the terms listed in \eqref{proof 2order eq5 ZF} with $I$ replaced by $\Jd$. We write in details the estimate of the term $Z^{\Jd}\big(\del_\alpha u_{\jh}\del_{\beta}u_{\kh}\big)$, as follows:
$$
\aligned
\big|Z^{\Jd}\big(\del_\alpha u_{\jh}\del_{\beta}u_{\kh}\big)\big|
\leq& |Z^{\Jd}\del_\alpha u_{\jh}| \, |\del_{\beta}u_{\kh}|
+|\del_\alpha u_{\jh}| \, |Z^{J_1}\del_{\beta}u_{\kh}|
\\
& \quad
+ \sum_{|J_1|,|J_2| \leq 1  \atop J_2+J_3=\Jd}|Z^{J_2}\del_\alpha u_{\jh}| \, |Z^{J_3}\del_{\beta}u_{\kh}|
\\
\leq& CC_1\eps t^{-1/2}s^{-1+\delta} \,  CC_1\eps t^{-1/2}s^{-1}
\\
& \quad + CC_1\eps t^{-1/2}s^{-1} \,  CC_1\eps t^{-1/2}s^{-1+\delta}
\\
    &
    + CC_1\eps t^{-1/2}s^{-1} \,  CC_1\eps t^{-1/2}s^{-1}
\\
\leq& C(C_1\eps)^2t^{-1}s^{-2+\delta}
\\
\leq& C(C_1\eps)^2t^{-3/2}s^{-1+\delta}.
\endaligned
$$
For the three partitions of $\Jd = J_1+J_2$, the $L^\infty$ estimates 
we use are: when $|J_1|=0$ and $|J_2|\leq |\Jd|$, we apply \eqref{proof decay 2ge3a} and \eqref{proof decay 2ge1a}; when $|J_1|=|J_2|=1$, we
apply \eqref{proof decay 2ge3a} and \eqref{proof decay 2ge3a}; when $|J_1|=|\Jd|$ and $|J_2|\leq 0$,
we apply \eqref{proof decay 2ge1a} and \eqref{proof decay 2ge3a}.

For the remaining terms, for each partition of $J_1$, we just list out 
 the $L^\infty$ estimates we use (with decay rate modulo a factor $C(C_1\eps)^2$):
$$
\begin{array}{ccccc}
\text{Products} &(2,\leq0) &(1,\leq1) &(0,\leq2) &\text{Decay rate}
\\
\del_\alpha u_{\jh}\del_{\beta}u_{\kh}
&\eqref{proof decay 2ge1a},\eqref{proof decay 2ge3a}
&\eqref{proof decay 2ge3a},\eqref{proof decay 2ge3a}
&\eqref{proof decay 2ge3a},\eqref{proof decay 2ge1a}
& t^{-1}s^{-2+\delta}
\\
\del_\alpha v_{\jc}\del_{\beta}w_k
&\eqref{proof decay 2ge2a},\eqref{proof decay 2ge1a}
&\eqref{proof decay 2ge2a},\eqref{proof decay 2ge1a}
&\eqref{proof decay 2ge2a},\eqref{proof decay 2ge1a}
&t^{-2}s^{-1+2\delta}
\\
v_{\kc}\del_\alpha w_j
&\eqref{proof decay 2ge1c},\eqref{proof decay 2ge1a}
&\eqref{proof decay 2ge1c},\eqref{proof decay 2ge1a}
&\eqref{proof decay 2ge1c},\eqref{proof decay 2ge1a}
&t^{-2}s^{-1+2\delta}
\\
v_{\jc}v_{\kc}
&\eqref{proof decay 2ge1c},\eqref{proof decay 2ge1c}
&\eqref{proof decay 2ge1c},\eqref{proof decay 2ge1c}
&\eqref{proof decay 2ge1c},\eqref{proof decay 2ge1c}
&t^{-3+2\delta}
\end{array}
$$
Combined with the condition $\delta<1/6$ and the fact that $s\leq Ct \leq Cs^2$, the estimate \eqref{proof 2order eq6a} is proved.

On the other hand, the estimate of $R(Z^{J_1}u_{\ih})$ is a direct application of \eqref{proof decay 1ge4a}.
\end{proof}

\begin{proof}[Proof of Lemma \ref{proof 2order lem5'}]
The proof is essentially the same as that of \eqref{proof 2order lem5'}, and 
 the main difference lies in
 the inequalities we use for each term and partition of the index. We omit the details 
and present the inequalities we 
use, as follows.

For the estimates of ${Q_G}_{\ih}(J,w,\del w,\del\del w)$, the inequalities we use are listed in:
$$
\begin{array}{cccc}
\text{Products} &(1,\leq 0) &(0,\leq 1) &\text{ Decay rate}
\\
u_{\ih}\delu_a\delu_{\beta}u_{\jh}
&\eqref{proof decay 2ge5b},\eqref{proof decay 2ge4b}
&\eqref{proof decay 2ge5b},\eqref{proof decay 2ge4a'}
&t^{-3}s^{\delta/2}
\\
v_{\ic}\delu_a\delu_{\beta}u_{\jh}
&\eqref{proof decay 2ge1c'},\eqref{proof decay 2ge4b}
&\eqref{proof decay 2ge1c'},\eqref{proof decay 2ge4a'}
&t^{-3}s^{-1+\delta}
\\
u_{\ih}\del_t\delu_au_{\jh}
&\eqref{proof decay 2ge5b},\eqref{proof decay 2ge4b}
&\eqref{proof decay 2ge5b},\eqref{proof decay 2ge4a'}
&t^{-3}s^{\delta/2}
\\
v_{\ic}\del_t\delu_au_{\jh}
&\eqref{proof decay 2ge1c'},\eqref{proof decay 2ge4b}
&\eqref{proof decay 2ge1c'},\eqref{proof decay 2ge4a'}
&t^{-3}s^{-1+\delta}
\\
\delu_{\gamma}v_{\ic}\delu_a\delu_{\beta}u_{\jh}
&\eqref{proof decay 2ge2a'},\eqref{proof decay 2ge4b}
&\eqref{proof decay 2ge2a'},\eqref{proof decay 2ge4a'}
&t^{-3}s^{-1+\delta}
\\
\delu_{\gamma}u_{\ih}\delu_a\delu_{\beta}u_{\jh}
&\eqref{proof decay 2ge3a},\eqref{proof decay 2ge4b}
&\eqref{proof decay 2ge3a},\eqref{proof decay 2ge4a'}
&t^{-2}s^{-2+\delta/2}
\end{array}
$$
$$
\begin{array}{cccc}
\text{Products} &(1,\leq 0) &(0,\leq 1) &\text{Decay rate}
\\
\delu_{\gamma}u_{\ih}\del_t\delu_au_{\jh}
&\eqref{proof decay 2ge3a},\eqref{proof decay 2ge4b}
&\eqref{proof decay 2ge3a},\eqref{proof decay 2ge4a'}
&t^{-2}s^{-2+\delta/2}
\\
\delu_{\gamma}v_{\ic}\del_t\delu_au_{\jh}
&\eqref{proof decay 2ge2a'},\eqref{proof decay 2ge4b}
&\eqref{proof decay 2ge2a'},\eqref{proof decay 2ge4a'}
&t^{-3}s^{-1+\delta}
\\
v_{\ic}\del_\alpha \del_{\beta}v_{\jc}
&\eqref{proof decay 2ge1c'},\eqref{proof decay 2ge2a'}
&\eqref{proof decay 2ge1c'},\eqref{proof decay 2ge2a'}
&t^{-3}s^{\delta}
\\
\del_{\gamma}u_{\ih}\del_\alpha \del_{\beta}v_{\jc}
&\eqref{proof decay 2ge3a},\eqref{proof decay 2ge2a'}
&\eqref{proof decay 2ge3a},\eqref{proof decay 2ge2a'}
&t^{-2}s^{-1+\delta/2}
\\
\del_{\gamma}v_{\ic}\del_\alpha \del_{\beta}v_{\jc}
&\eqref{proof decay 2ge2a'},\eqref{proof decay 2ge2a'}
&\eqref{proof decay 2ge2a'},\eqref{proof decay 2ge2a'}
&t^{-3}s^{\delta}
\end{array}
$$

$$
\begin{array}{ccccc}
\text{Terms} &(1,\leq 0) &(0,\leq 1)  &\text{ Decay rate}
\\
t^{-1}Z^{J_2}u_{\ih}Z^{J_3}\delu_{\beta'}u_{\jh}
&\eqref{proof decay 2ge5b},\eqref{proof decay 2ge3a}
&\eqref{proof decay 2ge5b},\eqref{proof decay 2ge3a}
&t^{-3}
\\
t^{-1}Z^{J_2}v_{\ic}Z^{J_3}\delu_{\beta'}u_{\jh}
&\eqref{proof decay 2ge1c'},\eqref{proof decay 2ge3a}
&\eqref{proof decay 2ge1c'},\eqref{proof decay 2ge3a}
&t^{-3}s^{-1+\delta/2}
\\
t^{-1}Z^{J_2}\del_{\gamma}u_{\ih}Z^{J_3}\delu_{\beta'}u_{\jh}
&\eqref{proof decay 2ge3a},\eqref{proof decay 2ge3a}
&\eqref{proof decay 2ge3a},\eqref{proof decay 2ge3a}
&t^{-2}s^{-2}
\\
t^{-1}Z^{J_2}\del_{\gamma}v_{\ic}Z^{J_3}\delu_{\beta'}u_{\jh}
&\eqref{proof decay 2ge1a'},\eqref{proof decay 2ge3a}
&\eqref{proof decay 2ge1a'},\eqref{proof decay 2ge3a}
&t^{-2}s^{-2+\delta/2}
\end{array}
$$

\

\noindent For the term ${Q_T}_{\ih}(J,w,\del w,\del\del w)$, we have the list: 
$$
\begin{array}{cccc}
\text{Products} &(1,\leq0) &(0,\leq1)  &\text{Decay rate}
\\
\delu_a\delu_{\beta}Z^{J_1}u_{\ih}\,Z^{J_2}\del_{\gamma}u_{\kh}
&\eqref{proof decay 1ge4a'},\eqref{proof decay 2ge3a}
&\eqref{proof decay 1ge4b},\eqref{proof decay 2ge3a}
&t^{-2}s^{-2+\delta/2}
\\
\delu_a\delu_{\beta}Z^{J_1}u_{\ih}\,Z^{J_2}\del_{\gamma}v_{\kc}
&\eqref{proof decay 1ge4a'},\eqref{proof decay 2ge2a'}
&\eqref{proof decay 1ge4b},\eqref{proof decay 2ge2a'}
&t^{-3}s^{-1+\delta}
\\
\delu_a\delu_{\beta}Z^{J_1}u_{\ih}\,Z^{J_2}u_{\kh}
&\eqref{proof decay 1ge4a'},\eqref{proof decay 2ge5b}
&\eqref{proof decay 1ge4b},\eqref{proof decay 2ge5b}
&t^{-3}s^{\delta}
\\
\delu_a\delu_{\beta}Z^{J_1}u_{\ih}\,Z^{J_2}v_{\kc}
&\eqref{proof decay 1ge4a'},\eqref{proof decay 2ge1c'}
&\eqref{proof decay 1ge4b},\eqref{proof decay 2ge1c'}
&t^{-3}s^{-1+\delta}
\end{array}
$$

$$
\begin{array}{cccc}
\text{Terms} &(1,\leq0) &(0,\leq1) &\text{Decay rate}
\\
t^{-1}\del_{\gamma'}Z^{J_1}u_{\jh}\,Z^{J_2}\del_{\gamma}u_{\kh}
&\eqref{proof decay 1ge3a},\eqref{proof decay 2ge3a}
&\eqref{proof decay 1ge3a},\eqref{proof decay 2ge3a}
&t^{-2}s^{-2}
\\
t^{-1}\del_{\gamma'}Z^{J_1}u_{\jh}\,Z^{J_2}\del_{\gamma}v_{\kc}
&\eqref{proof decay 1ge3a},\eqref{proof decay 2ge2a'}
&\eqref{proof decay 1ge3a},\eqref{proof decay 2ge2a'}
&t^{-3}s^{-1+\delta/2}
\\
t^{-1}\del_{\gamma'}Z^{J_1}u_{\jh}\,Z^{J_2}u_{\kh}
&\eqref{proof decay 1ge3a},\eqref{proof decay 2ge5b}
&\eqref{proof decay 1ge3a},\eqref{proof decay 2ge5b}
&t^{-3}
\\
t^{-1}\del_{\gamma'}Z^{J_1}u_{\jh}\,Z^{J_2}v_{\kc}
&\eqref{proof decay 1ge3a},\eqref{proof decay 2ge1c'}
&\eqref{proof decay 1ge3a},\eqref{proof decay 2ge1c'}
&t^{-2}s^{-1+\delta/2}
\end{array}
$$

\

\noindent Finally, for the term $F_{\ih}$ we have 
$$
\begin{array}{cccc}
\text{Products} &(1,\leq0) &(0,\leq1) &\text{Decay rate}
\\
\del_\alpha u_{\jh}\del_{\beta}u_{\kh}
&\eqref{proof decay 2ge3a},\eqref{proof decay 2ge3a}
&\eqref{proof decay 2ge3a},\eqref{proof decay 2ge3a}
& t^{-1}s^{-2}
\\
\del_\alpha v_{\jc}\del_{\beta}w_k
&\eqref{proof decay 2ge2a'},\eqref{proof decay 2ge1a'}
&\eqref{proof decay 2ge2a'},\eqref{proof decay 2ge1a'}
&t^{-2}s^{-1+\delta}
\\
v_{\kc}\del_\alpha w_j
&\eqref{proof decay 2ge1c'},\eqref{proof decay 2ge1a'}
&\eqref{proof decay 2ge1c'},\eqref{proof decay 2ge1a'}
&t^{-2}s^{-1+\delta}
\\
v_{\jc}v_{\kc}
&\eqref{proof decay 2ge1c'},\eqref{proof decay 2ge1c'}
&\eqref{proof decay 2ge1c'},\eqref{proof decay 2ge1c'}
&t^{-3+\delta}
\end{array}
$$
On the other hand, the estimate of $R(Z^J u)$ is a direct result of \eqref{proof decay 1ge4a'}.
\end{proof}

\begin{proof}[Proof of Lemma \ref{proof 2order lem6}]
The proof is essentially the same as the proof of Lemma \ref{proof 2order lem5} but much easier. We analyze
 first ${Q_G}_{\ih}(w,\del w,\del\del w)$. This is a direct application of the $L^\infty$ estimates established earlier
on the terms listed in \eqref{proof 2order eq4 QG} with $I=0$. As done in the proof of Lemma \ref{proof 2order lem5}, this gives the following list (with decay rate modulo $C(C_1\eps)^2$):
$$
\begin{array}{ccc}
\text{Terms} &(0,0) &\text{Decay rate}
\\
u_{\ih}\delu_a\delu_{\beta}u_{\jh}
&\eqref{proof decay 2ge5b},\eqref{proof decay 2ge4b}
&t^{-3}
\\
v_{\ic}\delu_a\delu_{\beta}u_{\jh}
&\eqref{proof decay 2ge1c},\eqref{proof decay 2ge4b}
&t^{-3}s^{-1+\delta}
\\
u_{\ih}\del_t\delu_au_{\jh}
&\eqref{proof decay 2ge5b},\eqref{proof decay 2ge4b}
&t^{-3}
\\
v_{\ic}\del_t\delu_au_{\jh}
&\eqref{proof decay 2ge1c},\eqref{proof decay 2ge4b}
&t^{-3}s^{-1+\delta}
\\
\delu_{\gamma}u_{\ih}\delu_a\delu_{\beta}u_{\jh}
&\eqref{proof decay 2ge3a},\eqref{proof decay 2ge4b}
&t^{-2}s^{-2}
\\
\delu_{\gamma}v_{\ic}\delu_a\delu_{\beta}u_{\jh}
&\eqref{proof decay 2ge2a},\eqref{proof decay 2ge4b}
&t^{-3}s^{-1+\delta}
\end{array}
$$

$$
\begin{array}{ccc}
\text{Terms} &(0,0) &\text{Decay rate}
\\
\delu_{\gamma}u_{\ih}\del_t\delu_au_{\jh}
&\eqref{proof decay 2ge3a},\eqref{proof decay 2ge4b}
&t^{-2}s^{-2}
\\
\delu_{\gamma}v_{\ic}\del_t\delu_au_{\jh}
&\eqref{proof decay 2ge2a},\eqref{proof decay 2ge4b}
&t^{-3}s^{-1+\delta}
\\
v_{\ic}\del_\alpha \del_{\beta}v_{\jc}
&\eqref{proof decay 2ge1c},\eqref{proof decay 2ge2a}
&t^{-3}s^{2\delta}
\\
\del_{\gamma}u_{\ih}\del_\alpha \del_{\beta}v_{\jc}
&\eqref{proof decay 2ge3a},\eqref{proof decay 2ge2a}
&t^{-2}s^{-1+\delta}
\\
\del_{\gamma}v_{\ic}\del_\alpha \del_{\beta}v_{\jc}
&\eqref{proof decay 2ge2a},\eqref{proof decay 2ge2a}
&t^{-3}s^{2\delta}
\end{array}
$$
The following four terms
$$
\del_\alpha \Psi_{\beta}^{\beta'} u_{\ih}\delu_{\beta'}u_{\jh},
\quad
\del_\alpha \Psi_{\beta}^{\beta'}v_{\ic}\delu_{\beta'}u_{\jh},
\quad
\del_\alpha \Psi_{\beta}^{\beta'}\del_{\gamma}u_{\ih}\delu_{\beta'}u_{\jh},
\quad
\del_\alpha \Psi_{\beta}^{\beta'}\del_{\gamma}v_{\ic}\delu_{\beta'}u_{\jh}
$$
are estimated by taking into account the additional decay supplied by the factor $|\Psi_{\beta}^{\beta'}|\leq Ct^{-1}$. The inequalities we use for each term are listed as follows:
$$
\begin{array}{ccc}
\text{Terms} &(0,0) &\text{Decay rate}
\\
t^{-1}u_{\ih}\delu_{\beta'}u_{\jh}
&\eqref{proof decay 2ge5b},\eqref{proof decay 2ge3a}
&t^{-3}
\\
t^{-1}v_{\ic}\delu_{\beta'}u_{\jh}
&\eqref{proof decay 2ge1c},\eqref{proof decay 2ge3a}
&t^{-3}s^{-1+\delta}
\\
t^{-1}\del_{\gamma}u_{\ih}\delu_{\beta'}u_{\jh}
&\eqref{proof decay 2ge3a},\eqref{proof decay 2ge3a}
&t^{-2}s^{-2}
\\
t^{-1}\del_{\gamma}v_{\ic}\delu_{\beta'}u_{\jh}
&\eqref{proof decay 2ge2a},\eqref{proof decay 2ge3a}
&t^{-3}s^{-1+\delta}
\end{array}
$$

The estimates of ${Q_T}_i$ are similar. We establish the following list and omit the details:
$$
\begin{array}{ccc}
\text{Products} &(0,0) &\text{Decay rate}
\\
\del_{\gamma}u_{\kh}\,\delu_a\delu_{\beta}u_{\ih}
&\eqref{proof decay 2ge3a},\eqref{proof decay 2ge4b}
&t^{-2}s^{-2}
\\
\del_{\gamma}v_{\kc}\,\delu_a\delu_{\beta}u_{\ih}
&\eqref{proof decay 2ge2a},\eqref{proof decay 2ge4b}
&t^{-3}s^{-1+\delta}
\\
u_{\kh}\,\delu_a\delu_{\beta}u_{\ih}
&\eqref{proof decay 2ge5b},\eqref{proof decay 2ge4b}
&t^{-3}
\\
v_{\kc}\,\delu_a\delu_{\beta}u_{\ih}
&\eqref{proof decay 2ge1c},\eqref{proof decay 2ge4b}
&t^{-3}s^{-1+\delta}
\end{array}
$$
$$
\begin{array}{ccc}
\text{Products} &(0,0) &\text{Decay rate}
\\
t^{-1}\del_{\gamma}u_{\kh}\,\del_{\gamma'}u_{\jh}
&\eqref{proof decay 2ge3a},\eqref{proof decay 2ge3a}
&t^{-2}s^{-2}
\\
t^{-1}\del_{\gamma}v_{\kc}\,\del_{\gamma'}u_{\jh}
&\eqref{proof decay 2ge2a},\eqref{proof decay 2ge3a}
&t^{-3}s^{-1+\delta}
\\
t^{-1}u_{\kh}\,\del_{\gamma'}u_{\jh}
&\eqref{proof decay 2ge5b},\eqref{proof decay 2ge3a}
&t^{-3}
\\
t^{-1}v_{\kc}\,\del_{\gamma'}u_{\jh}
&\eqref{proof decay 2ge1c},\eqref{proof decay 2ge3a}
&t^{-3}s^{-1+\delta}
\end{array}
$$
We conclude with \eqref{proof 2order eq7bc}.

The estimate of $F_{\ih}$ is as follows: recall the structure of $F_{\ih}$ described by \eqref{proof 2order eq5 ZF},
we need to estimate these terms with $I=0$. 
As in the proof of Lemma \ref{proof 2order lem5}, the following list is established (with decay rate modulo $C(C_1\eps)^2$):
$$
\begin{array}{ccc}
\text{Products} &(0,0) &\text{Decay rate}
\\
\del_\alpha u_{\ih}\del_{\beta}u_{\jh}
&\eqref{proof decay 2ge3a},\eqref{proof decay 2ge3a}
&t^{-1}s^{-2}
\\
\del_\alpha v_{\ic}\del_{\beta}w_j
&\eqref{proof decay 2ge2a},\eqref{proof decay 2ge1a}
&t^{-2}s^{-1+2\delta}
\\
v_{\ic}\del_\alpha w_j
&\eqref{proof decay 2ge1c},\eqref{proof decay 2ge1a}
&t^{-2}s^{-1+2\delta}
\\
v_{\ic}v_{\jc}
&\eqref{proof decay 2ge1c},\eqref{proof decay 2ge1c}
&t^{-3}s^{2\delta}
\end{array}
$$
By taking into account the condition $\delta<1/6$ and the fact that $C^{-1} s \leq t \leq C s^2$,
we conclude with \eqref{proof 2order eq7a}.

The estimate of $R(u_{\ih})$ is a direct application of \eqref{proof decay 2ge4b}.
\end{proof}

We can now prove the main result of this section.

\begin{proposition}\label{proof 2order prop decay}
Let $w_i$ the solution of \eqref{main eq main} and assume that \eqref{proof energy assumption} holds with $C_1\eps\leq \min\{1,\eps_0''\}$. 
For the wave components, the following decay estimates hold for any $\Jd\leq 2$ and $|J| \leq 1$: 
\begin{subequations}\label{proof 2order eq8}
\bel{proof 2order eq8a}
\sup_{\Hcal_s}\big(s^3t^{-1/2}|\del_t\del_t Z^{\Jd}u_{\ih}|\big)
\leq CC_1\eps s^{\delta},
\ee
\bel{proof 2order eq8a'}
\sup_{\Hcal_s}\big(s^3t^{-1/2}|\del_t\del_t Z^Ju_{\ih}|\big)
\leq CC_1\eps s^{\delta/2},
\ee
\bel{proof 2order eq8b}
\sup_{\Hcal_s}\big(s^3t^{-1/2}|\del_t\del_t u_{\ih}|\big)\leq CC_1\eps,
\ee
\end{subequations}
and more generally
\begin{subequations}\label{proof 2order general decay}
\bel{proof 2order general decay a}
\sup_{\Hcal_s}\big(s^3t^{-1/2}|\del_\alpha \del_{\beta}Z^{\Jd}u|\big)
\leq CC_1\eps s^{\delta},
\ee
\bel{proof 2order general decay a'}
\sup_{\Hcal_s}\big(s^3t^{-1/2}|\del_\alpha \del_{\beta}Z^Ju|\big)
\leq CC_1\eps s^{\delta/2},
\ee
\bel{proof 2order general decay b}
\sup_{\Hcal_s}\big(s^3t^{-1/2}|\del_\alpha \del_{\beta} u|\big) \leq  CC_1\eps.
\ee
\end{subequations}
\end{proposition}
\begin{proof}
\eqref{proof 2order eq8b} will be proved first.
Substitute \eqref{proof 2order eq7} into \eqref{proof 2order eq3} with $|I|=0$. Recall the convention $|I|<0 \Rightarrow Z^I = 0$. 
The desired result is proven.

Observe that for \eqref{proof 2order eq8a'}, the case $|J|\leq 0$ is proved by \eqref{proof 2order eq8b}. For the case $|J|=1$, we recall \eqref{proof 2order eq3}:
$$
\aligned
|\del_t\del_tZ^Ju_{\ih}|
\leq& C(t/s)^2\sum_{|J_2|=1\atop \gamma,\jh,i}
\big(|Z^{J_2}\del_{\gamma} w_i|+|Z^{J_2}w_i|\big)|\del_t\del_tu_{\jh}|
\\
& + C(t/s)^2 {Q_T}_{\ih}(J_2,w,\del w,\del\del w) + C(t/s)^2|{Q_G}_{\ih}(J_2,w,\del w\del\del w)|
\\
&+ C(t/s)^2|Z^{J_2}F_{\ih}|+ C(t/s)^2|R(Z^{J_2}u_{\ih})|
\endaligned
$$
substitutes \eqref{proof 2order eq8b} (an estimate on $\del_t\del_t u_{\jh}$) and \eqref{proof 2order eq6'} into \eqref{proof 2order eq3},
 together with the following estimate:
$$
\sum_{\gamma\,\jh,\ih \atop |J|\leq 1}\big(|Z^J\del_{\gamma}w_i| + |Z^Jw_i|\big)\leq CC_1\eps \big(t^{-1/2}s^{-1} + t^{-3/2}s + t^{-3/2}s^{\delta/2} \big)\leq CC_1\eps t^{-3/2}s.
$$

For \eqref{proof 2order eq8a}, remark that the case $|\Jd|\leq 1$ is guaranteed by \eqref{proof 2order eq8a'} and \eqref{proof 2order eq8b}.

The case $|\Jd|=2$ is done by substituting \eqref{proof 2order eq8a'} (with $|J|= 1$), \eqref{proof 2order eq8b} and \eqref{proof 2order eq6} into  \eqref{proof 2order eq3} with $|I|=2$ together with the following estimate which is a direct result of \eqref{proof decay 2ge1a}, \eqref{proof decay 2ge5b} and \eqref{proof decay 2ge1c}:
$$
\aligned
\sum_{\gamma\,\jh,\ih}\big(|Z^{\Jd}\del_{\gamma}w_i| + |Z^{\Jd}w_i|\big)
& \leq CC_1\eps \big(t^{-1/2}s^{-1+\delta/2} + t^{-3/2}s + t^{-3/2}s^{\delta/2} \big)
\\
& \leq CC_1\eps t^{-3/2}s^{1+\delta/2}.
\endaligned
$$
We recall \eqref{proof 2order eq3} with $I=\Jd$:
$$
\aligned
&|\del_t\del_tZ^{\Jd}u_{\ih}|
\\
\leq& C(t/s)^2\sum_{|J_1|+|J_2|\leq 2\atop |J_1|\leq 1 ,\gamma,\jh,i}
\big(|Z^{J_2}\del_{\gamma} w_i|+|Z^{J_2}w_i|\big)\,|\del_t\del_tZ^{J_1}u_{\jh}|
\\
& + C(t/s)^2 {Q_T}_{\ih}(J_2,w,\del w,\del\del w) + C(t/s)^2|{Q_G}_{\ih}(J_2,w,\del w\del\del w)|
\\
& + C(t/s)^2|Z^{J_2}F_{\ih}|+ C(t/s)^2|R(Z^{J_2}u_{\ih})|.
\endaligned
$$
Substituting  \eqref{proof 2order eq8a'} (with $|J|= 1$), \eqref{proof 2order eq8b} and \eqref{proof 2order eq6}, we find that 
$$
\aligned
&|\del_t\del_tZ^{\Jd}u_{\ih}|
\\
&\leq C(t/s)^2\sum_{|J_1|+|J_2|\leq 2\atop |J_1|\leq 1 ,\gamma,\jh,i}
\big(|Z^{J_2}\del_{\gamma} w_i|+|Z^{J_2}w_i|\big)\,|\del_t\del_tZ^{J_1}u_{\jh}|
 + CC_1\eps t^{1/2}s^{-3+\delta}
\\
&\leq C(t/s)^2\sum_{|J_1|+|J_2|\leq 2\atop |J_1| = 1 ,\gamma,\jh,i}
\big(|Z^{J_2}\del_{\gamma} w_i|+|Z^{J_2}w_i|\big)\,|\del_t\del_tZ^{J_1}u_{\jh}|
\\
    &\quad 
+ C(t/s)^2\sum_{|J_1|+|J_2|\leq 2\atop |J_1| = 0 ,\gamma,\jh,i}
\big(|Z^{J_2}\del_{\gamma} w_i|+|Z^{J_2}w_i|\big)\,|\del_t\del_tZ^{J_1}u_{\jh}| + CC_1\eps 2t^{1/2}s^{-3+\delta}
\\
&\leq CC_1\eps (t/s)^2 t^{-3/2}s \,  t^{1/2}s^{-3+\delta/2} + CC_1\eps (t/s)^2t^{-3/2}s^{1+\delta/2} \,  t^{1/2}s^{-3}
\\
& \quad+ CC_1\eps t^{1/2}s^{-3+\delta}
\\
&\leq 
CC_1\eps ts^{-4+\delta/2} + CC_1\eps t^{1/2}s^{-3+\delta} \leq  CC_1\eps t^{1/2}s^{-3+\delta}.
\endaligned
$$
The bound \eqref{proof 2order general decay} are direct result of \eqref{proof 2order eq8} combined with \eqref{proof 2order general decay1}.
\end{proof}

We can give the complete $L^\infty$ estimates of the second-order derivatives. 

\begin{proposition}
By relying on \eqref{proof energy assumption} with $C_1\eps\leq \min\{1,\eps_0''\}$,
 the following estimates hold for all $|\Jd|\leq 2$ and $|J|\leq 1$:
\begin{subequations}\label{proof decay 2order}
\bel{proof decay 2order a}
\aligned
&
\sup_{\Hcal_s}|s^3t^{-1/2}\del_\alpha \del_{\beta}Z^{\Jd}u| + \sup_{\Hcal_s}|s^3t^{-1/2}Z^{\Jd}\del_\alpha \del_{\beta}u|
\leq  CC_1\eps s^{\delta},
\endaligned
\ee
\bel{proof decay 2order a'}
\aligned
& \sup_{\Hcal_s}|s^3t^{-1/2}\del_\alpha \del_{\beta}Z^Ju| + \sup_{\Hcal_s}|s^3t^{-1/2}Z^J\del_\alpha \del_{\beta}u|
\leq CC_1\eps s^{\delta/2},
\endaligned
\ee
\bel{proof decay 2order b}
\aligned
& \sup_{\Hcal_s}|s^3t^{-1/2}\del_\alpha \del_{\beta} u| \leq CC_1\eps.
\endaligned
\ee
\end{subequations}
\end{proposition}

\begin{proof}
These estimates are a consequence of \eqref{pre lem commutator second-order} and \eqref{proof 2order general decay}.
\end{proof}


\section{$L^2$ estimates}
\label{sec:75}

The aim of this section is to get the $L^2$ estimates on $\del_t\del_tZ^Iu_{\jh}$. As in the last section, the strategy is to make use of \eqref{proof 2order eq3}. First, we estimate the terms ${Q_T}_{\ih}$, ${Q_G}_{\ih}$, $Z^I F_{\ih}$ and $R(Z^Iu_{\ih})$.

\begin{lemma}\label{proof 2order lem7}
Under the energy assumption \eqref{proof energy assumption}, the following estimates hold for all $|\Id|\leq 4$:
\begin{subequations}\label{proof 2order eq9}
\bel{proof 2order eq9a}
\big{\|}s{Q_G}_{\ih}(\Id,w,\del w,\del\del w)\big{\|}_{L^2(\Hcal_s)} \leq C(C_1\eps)^2s^{\delta},
\ee
\bel{proof 2order eq9ab}
\big{\|}s{Q_T}_{\ih}(\Id,w,\del w,\del\del w)\big{\|}_{L^2(\Hcal_s)} \leq C(C_1\eps)^2s^{\delta},
\ee
\bel{proof 2order eq9b}
\big{\|}sZ^{\Id}F_{\ih}\big{\|}_{L^2(\Hcal_s)} \leq C(C_1\eps)^2s^{\delta},
\ee
\bel{proof 2order eq9c}
\big{\|}sR(Z^{\Id}u_{\ih})\big{\|}_{L^2(\Hcal_s)}\leq CC_1\eps s^{\delta}.
\ee
\end{subequations}
\end{lemma}

\begin{lemma}\label{proof 2order lem7'}
Under the energy assumption \eqref{proof energy assumption}, the following estimates hold for all $|I|\leq 3$:
\begin{subequations}\label{proof 2order eq9-DEUX}
\bel{proof 2order eq9'a}
\big{\|}s{Q_G}_{\ih}(I,w,\del w,\del\del w)\big{\|}_{L^2(\Hcal_s)} \leq C(C_1\eps)^2s^{\delta/2},
\ee
\bel{proof 2order eq9'ab}
\big{\|}s{Q_T}_{\ih}(I,w,\del w,\del\del w)\big{\|}_{L^2(\Hcal_s)} \leq C(C_1\eps)^2s^{\delta/2},
\ee
\bel{proof 2order eq9'b}
\big{\|}sZ^{I}F_{\ih}\big{\|}_{L^2(\Hcal_s)} \leq C(C_1\eps)^2s^{\delta/2},
\ee
\bel{proof 2order eq9'c}
\big{\|}sR(Z^{I}u_{\ih})\big{\|}_{L^2(\Hcal_s)}\leq CC_1\eps s^{\delta/2}.
\ee
\end{subequations}
\end{lemma}

\begin{lemma}\label{proof 2order lem8}
Under the energy assumption \eqref{proof energy assumption}, the following estimates hold for all $|\If|\leq 2$:
\begin{subequations}\label{proof 2order eq10}
\bel{proof 2order eq10a}
\big{\|}s{Q_G}_{\ih}(\If,w,\del w,\del\del w)\big{\|}_{L^2(\Hcal_s)} \leq C(C_1\eps)^2,
\ee
\bel{proof 2order eq10ab}
\big{\|}s{Q_G}_{\ih}(\If,w,\del w,\del\del w)\big{\|}_{L^2(\Hcal_s)} \leq C(C_1\eps)^2,
\ee
\bel{proof 2order eq10b}
\big{\|}sZ^{\If}F_{\ih}\big{\|}_{L^2(\Hcal_s)} \leq C(C_1\eps)^2,
\ee
\bel{proof 2order eq10c}
\big{\|}sR(Z^{\If}u_{\ih})\big{\|}_{L^2(\Hcal_s)}\leq CC_1\eps.
\ee
\end{subequations}
\end{lemma}

\begin{proof}[Proof of Lemma \ref{proof 2order lem7}]
We consider first ${Q_G}_{\ih}$ and 
recall the structure of ${Q_G}_{\ih}$ expressed by \eqref{proof 2order eq4 QG}. 
We do an $L^2$ estimate on each term of \eqref{proof 2order eq4 QG} with $I = \Id, |\Id|\leq 4$. We take $u_{\ih}\delu_a\del_tu_{\jh}$  as an example and we write down of the argument:
$$
\aligned
&\big{\|}sZ^{\Id}\big(u_{\ih}\delu_a\del_tu_{\jh}\big)\big{\|}_{L^2(\Hcal_s)}
\\
&\leq \sum_{I_1+I_2 = \Id} \!\!\!\!
\big{\|}s\big( Z^{I_1}u_{\ih} \,  Z^{I_2}\delu_a\del_tu_{\jh}\big)\big{\|}_{L^2(\Hcal_s)}
\\
&\leq  \big{\|}s\big(u_{\ih} \,  Z^{\Id}\delu_a\del_tu_{\jh}\big)\big{\|}_{L^2(\Hcal_s)}
   +\!\!\!\!\sum_{|I_1|=1 \atop I_1+I_2 = \Id}\!\!\!\!
   \big{\|}s\big( Z^{I_1}u_{\ih} \,  Z^{I_2}\delu_a\del_tu_{\jh}\big)\big{\|}_{L^2(\Hcal_s)}
\\
 &\quad +\!\!\!\!\sum_{|I_1|=2 \atop I_1+I_2 = \Id}\!\!\!\!
   \big{\|}s\big( Z^{I_1}u_{\ih} \,  Z^{I_2}\delu_a\del_tu_{\jh}\big)\big{\|}_{L^2(\Hcal_s)}
   +\sum_{|I_1|=3 \atop I_1+I_2 = \Id}
   \big{\|}s\big( Z^{I_1}u_{\ih} \,  Z^{I_2}\delu_a\del_tu_{\jh}\big)\big{\|}_{L^2(\Hcal_s)}
\\
 &\quad 
+\big{\|}s\big(Z^{\Id}u_{\ih} \, \delu_a\del_tu_{\jh}\big)\big{\|}_{L^2(\Hcal_s)}
\\
&=: T_0 + T_1 + T_2 + T_3 + T_4.
\endaligned
$$
Observe now that in term $T_k$, $|I_2|\leq 4-k$.

The term $T_0$ is estimated by \eqref{proof decay 2ge5b} and \eqref{proof L2 2ge4'a}:
$$
\aligned
T_0=& \big{\|}s\big(u_{\ih} \,  Z^{I_1}\delu_a\del_tu_{\jh}\big)\big{\|}_{L^2(\Hcal_s)}
\leq CC_1\eps\big{\|}s\big(t^{-3/2}\, sZ^{I_1}\delu_a\del_tu_{\jh}\big)\big{\|}_{L^2(\Hcal_s)}
\\
\leq &  CC_1\eps s^{-1/2}\big{\|}sZ^{I_1}\delu_a\del_tu_{\jh}\big{\|}_{L^2(\Hcal_s)}
\\
\leq& C(C_1\eps)^2s^{-1/2+\delta}.
\endaligned
$$

The term $T_1$ is estimated by \eqref{proof decay 2ge5b} and \eqref{proof L2 2ge4'a}:
$$
\aligned
T_1 =& \sum_{|I_1|=1 \atop I_1+I_2 = \Id}\!\!\!\!
   \big{\|}s\big( Z^{I_1}u_{\ih} \,  Z^{I_2}\delu_a\del_tu_{\jh}\big)\big{\|}_{L^2(\Hcal_s)}
\\
\leq& CC_1\eps\!\!\!\!\sum_{|I_1|=1 \atop  I_1+I_2 = \Id}\!\!\!\!
   \big{\|}s\big(t^{-3/2} \,  sZ^{I_2}\delu_a\del_tu_{\jh}\big)\big{\|}_{L^2(\Hcal_s)}
\\
\leq& CC_1\eps s^{-1/2}\!\!\!\!\sum_{|I_1|=1 \atop  I_1+I_2 = \Id}\!\!\!\!
\big{\|}sZ^{I_3}\delu_a\del_tu_{\jh}|\big|_{L^2(\Hcal_s)}
\leq C(C_1\eps)^2s^{-1/2+\delta}.
\endaligned
$$

The term $T_2$ is estimated by \eqref{proof decay 2ge5a} and \eqref{proof L2 2ge4'b} (remark that $|I_2|\leq 4-2=2$):
$$
\aligned
T_2 =& \sum_{|I_1|=2 \atop  I_1+I_2 = \Id}\!\!\!\!
   \big{\|}s\big( Z^{I_1}u_{\ih} \,  Z^{I_2}\delu_a\del_tu_{\jh}\big)\big{\|}_{L^2(\Hcal_s)}
\\
\leq& CC_1\eps\!\!\!\!\sum_{|I_1|=2 \atop  I_1+I_2 = \Id}\!\!\!\!
   \big{\|}s\big( t^{-3/2}s^{1+\delta} \,  Z^{I_2}\delu_a\del_tu_{\jh}\big)\big{\|}_{L^2(\Hcal_s)}
\\
=&CC_1\eps s^{-1/2+\delta}\!\!\!\!\sum_{|I_1|=2 \atop  I_1+I_2 = \Id}\!\!\!\!
   \big{\|} sZ^{I_2}\delu_a\del_tu_{\jh}\big{\|}_{L^2(\Hcal_s)}
\leq C(C_1\eps)^2s^{-1/2+\delta}.
\endaligned
$$

The term $T_3$ is estimated by \eqref{proof decay 2ge5a} and \eqref{proof L2 2ge4'b} (remark that $|I_2|\leq 4-3=1$):
$$
\aligned
T_3 =& \sum_{|I_1|=3 \atop  I_1+I_2 = \Id}\!\!\!\!
   \big{\|}s\big( Z^{I_1}u_{\ih} \,  Z^{I_2}\delu_a\del_tu_{\jh}\big)\big{\|}_{L^2(\Hcal_s)}
\\
\leq& CC_1\eps\!\!\!\!\sum_{|I_1|=3 \atop  I_1+I_2 = \Id}\!\!\!\!
   \big{\|}s\big( t^{-3/2}s^{\delta} \,  sZ^{I_2}\delu_a\del_tu_{\jh}\big)\big{\|}_{L^2(\Hcal_s)}
\\
=&CC_1\eps s^{-1/2+\delta}\!\!\!\!\sum_{|I_1|=3 \atop  I_1+I_2 = \Id}\!\!\!\!
   \big{\|} sZ^{I_2}\delu_a\del_tu_{\jh}\big{\|}_{L^2(\Hcal_s)}
\leq C(C_1\eps)^2 s^{-1/2+\delta}.
\endaligned
$$
The term 
$T_4$ is estimated by \eqref{proof L2 2ge5a} and \eqref{proof decay 2ge4b}:
$$
\aligned
T_4 &=
   \big{\|}s\big(Z^{\Id}u_{\ih} \, \delu_a\del_tu_{\jh}\big)\big{\|}_{L^2(\Hcal_s)}
   \leq CC_1\eps
   \big{\|}s Z^{\Id}u_{\ih} \,  t^{-3/2}s^{-1}\big{\|}_{L^2(\Hcal_s)}
\\
& = CC_1\eps s^{-1/2}\big{\|}t^{-1}Z^{\Id}u_{\ih}\big{\|}_{L^2(\Hcal_s)} \leq CC_1\eps s^{-1/2} \big{\|}s^{-1}Z^{\Id}u_{\ih}\big{\|}_{L^2(\Hcal_s)}
\\
&
\leq C(C_1\eps)^2s^{-1/2+\delta}.
\endaligned
$$
So we  conclude with
$$
\big{\|}sZ^{\Id}\big(u_{\ih}\delu_a\del_tu_{\jh}\big)\big{\|}_{L^2(\Hcal_s)}\leq C(C_1\eps)^2s^{-1/2+\delta}.
$$

For the remaining terms, we will not write the details, but
 for each partition of $\Id = I_1 + I_2$ we give 
the $L^2$ and $L^\infty$ estimates in Chapter~\ref{cha:6}.
As in the estimate of $u_{\ih}\delu_a\del_tu_{\jh}$,
we denote by $(k,\leq 4-k)$ the terms with which $|I_1|=k,\, |I_2|\leq 4-k$. This leads us to Table 4 and Table 5. 

There are four terms to be estimated separately:
$$
\aligned
& sZ^{\Id}\big(\del_\alpha \Psi_{\beta}^{\beta'}u_{\ih}\delu_{\beta'}u_{\jh}\big),\quad
sZ^{\Id}\big(\del_\alpha \Psi_{\beta}^{\beta'}v_{\ic}\delu_{\beta'}u_{\jh}\big),\quad
\\
&
sZ^{\Id}\big(\del_\alpha \Psi_{\beta}^{\beta'}\del_{\gamma}u_{\ih}\delu_{\beta'}u_{\jh}\big),\quad
sZ^{\Id}\big(\del_\alpha \Psi_{\beta}^{\beta'}\del_{\gamma}v_{\ic}\delu_{\beta'}u_{\jh}\big).
\endaligned
$$
We will use the additional decay supplied by $|Z^I\del_{\alpha}\Psi_{\beta}^{\beta'}|\leq C(I)t^{-1}$.
Let us take
$$
Z^{\Id}\big(\del_\alpha \Psi_{\beta}^{\beta'}u_{\ih}\delu_{\beta'}u_{\jh}\big)
$$
as an example and write its estimate in details:
$$
\aligned
&\big{\|}sZ^{\Id}\big(\del_\alpha \Psi_{\beta}^{\beta'}u_{\ih}\delu_{\beta'}u_{\jh}\big)\big{\|}_{L^2(\Hcal_s)}
\leq  \sum_{I_1+I_2+I_3=\Id}\!\!\!\!
 \big{\|}sZ^{I_3}\del_{\alpha}\Psi_{\beta}^{\beta'} \,  Z^{I_1}u_{\ih} \,  Z^{I_2}\delu_{\beta'}u_{\jh}\big{\|}_{L^2(\Hcal_s)}
\\
&\leq C\sum_{|I_1|+|I_2|\leq |\Id|}\big{\|}st^{-1}Z^{I_1}u_{\ih} \,  Z^{I_2}\delu_{\beta'}u_{\jh}\big{\|}_{L^2(\Hcal_s)}
\\
&\leq  C\sum_{|I_2|\leq |\Id|}\big{\|}st^{-1} u_{\ih} \,  Z^{I_2}\delu_{\beta'}u_{\jh}\big{\|}_{L^2(\Hcal_s)}
   +C\sum_{|I_1|=1 \atop  |I_2|\leq |\Id|-1}
   \big{\|}st^{-1}Z^{I_1}u_{\ih} \,  Z^{I_2}\delu_{\beta'}u_{\jh}\big{\|}_{L^2(\Hcal_s)}
\\
  &\quad 
+C\sum_{|I_1|=2 \atop  |I_2|\leq |\Id|-2}
  \big{\|}st^{-1}Z^{I_1}u_{\ih} \,  Z^{I_2}\delu_{\beta'}u_{\jh}\big{\|}_{L^2(\Hcal_s)}
\\
  &\quad 
+C\sum_{|I_1|=3 \atop  |I_2|\leq |\Id|-3}
   \big{\|}st^{-1}Z^{I_1}u_{\ih} \,  Z^{I_2}\delu_{\beta'}u_{\jh}\big{\|}_{L^2(\Hcal_s)}
   +C\sum_{|I_1|=4}\big{\|}st^{-1}Z^{I_1}u_{\ih} \,  \delu_{\beta'}u_{\jh}\big{\|}_{L^2(\Hcal_s)}
\\
&=: T_0 + T_1 + T_2 + T_3 + T_4.
\endaligned
$$


\protect\begin{landscape}\thispagestyle{empty}
$$
\begin{array}{cccccc}
\text{Terms}&(4,\leq 0)&(3,\leq 1)&(2,\leq 2)&(1,\leq 3) &(0,\leq 4)
\\
u_{\ih}\delu_a\delu_{\beta} u_{\jh}
&\eqref{proof L2 2ge5a'},\eqref{proof decay 2ge4b}
&\eqref{proof decay 2ge5a},\eqref{proof L2 2ge4'b}
&\eqref{proof decay 2ge5a},\eqref{proof L2 2ge4'b}
&\eqref{proof decay 2ge5b},\eqref{proof L2 2ge4'a}
&\eqref{proof decay 2ge5b},\eqref{proof L2 2ge4'a}
\\
v_{\ic}\delu_a\delu_{\beta} u_{\jh}
&\eqref{proof L2 2ge1c},\eqref{proof decay 2ge4b}
&\eqref{proof decay 2ge1c},\eqref{proof L2 2ge4'b}
&\eqref{proof decay 2ge1c},\eqref{proof L2 2ge4'b}
&\eqref{proof decay 2ge1c},\eqref{proof L2 2ge4'a}
&\eqref{proof decay 2ge1c},\eqref{proof L2 2ge4'a}
\\
u_{\ih}\del_t\delu_au_{\jh}
&\eqref{proof L2 2ge5a'},\eqref{proof decay 2ge4b}
&\eqref{proof decay 2ge5a},\eqref{proof L2 2ge4'b}
&\eqref{proof decay 2ge5a},\eqref{proof L2 2ge4'b}
&\eqref{proof decay 2ge5b},\eqref{proof L2 2ge4'a}
&\eqref{proof decay 2ge5b},\eqref{proof L2 2ge4'a}
\\
v_{\ic}\del_t\delu_au_{\jh}
&\eqref{proof L2 2ge1c},\eqref{proof decay 2ge4b}
&\eqref{proof decay 2ge1c},\eqref{proof L2 2ge4'b}
&\eqref{proof decay 2ge1c},\eqref{proof L2 2ge4'b}
&\eqref{proof decay 2ge1c},\eqref{proof L2 2ge4'a}
&\eqref{proof decay 2ge1c},\eqref{proof L2 2ge4'a}
\\
\delu_{\gamma}u_{\ih}\delu_a\delu_{\beta} u_{\jh}
&\eqref{proof L2 2ge1a},\eqref{proof decay 2ge4b}
&\eqref{proof decay 2ge1a},\eqref{proof L2 2ge4'b}
&\eqref{proof decay 2ge1a},\eqref{proof L2 2ge4'b}
&\eqref{proof decay 2ge3a},\eqref{proof L2 2ge4'a}
&\eqref{proof decay 2ge3a},\eqref{proof L2 2ge4'a}
\\
\delu_{\gamma}v_{\ic}\delu_a\delu_{\beta} u_{\jh}
&\eqref{proof L2 2ge2a},\eqref{proof decay 2ge4b}
&\eqref{proof decay 2ge2a},\eqref{proof L2 2ge4'b}
&\eqref{proof decay 2ge2a},\eqref{proof L2 2ge4'b}
&\eqref{proof decay 2ge2a},\eqref{proof L2 2ge4'a}
&\eqref{proof decay 2ge2a},\eqref{proof L2 2ge4'a}
\end{array}
$$
\centerline{Table 4} 

\

\

$$
\begin{array}{cccccc}
\text{Terms} &(4,\leq 0) &(3,\leq 1) &(2,\leq 2) &(1,\leq 3) &(0,\leq 4)
\\
\delu_{\gamma}u_{\ih}\del_t\delu_au_{\jh}
&\eqref{proof L2 2ge1a},\eqref{proof decay 2ge4b}
&\eqref{proof decay 2ge1a},\eqref{proof L2 2ge4'b}
&\eqref{proof decay 2ge1a},\eqref{proof L2 2ge4'b}
&\eqref{proof decay 2ge3a},\eqref{proof L2 2ge4'a}
&\eqref{proof decay 2ge3a},\eqref{proof L2 2ge4'a}
\\
\delu_{\gamma}v_{\ic}\del_t\delu_au_{\jh}
&\eqref{proof L2 2ge2a},\eqref{proof decay 2ge4b}
&\eqref{proof decay 2ge2a},\eqref{proof L2 2ge4'b}
&\eqref{proof decay 2ge2a},\eqref{proof L2 2ge4'b}
&\eqref{proof decay 2ge2a},\eqref{proof L2 2ge4'a}
&\eqref{proof decay 2ge2a},\eqref{proof L2 2ge4'a}
\\
v_{\ic}\del_\alpha \del_{\beta}v_{\jc}
&\eqref{proof L2 2ge1c},\eqref{proof decay 2ge2a}
&\eqref{proof decay 2ge1c},\eqref{proof L2 2ge2a}
&\eqref{proof decay 2ge1c},\eqref{proof L2 2ge2a}
&\eqref{proof decay 2ge1c},\eqref{proof L2 2ge2a}
&\eqref{proof decay 2ge1c},\eqref{proof L2 2ge2a}
\\
\del_{\gamma}u_{\ih}\del_\alpha \del_{\beta}v_{\jc}
&\eqref{proof L2 2ge1a},\eqref{proof decay 2ge2a}
&\eqref{proof decay 2ge1a},\eqref{proof L2 2ge2a}
&\eqref{proof decay 2ge1a},\eqref{proof L2 2ge2a}
&\eqref{proof decay 2ge3a},\eqref{proof L2 2ge2a}
&\eqref{proof decay 2ge3a},\eqref{proof L2 2ge2a}
\\
\del_{\gamma}v_{\ic}\del_\alpha \del_{\beta}v_{\jc}
&\eqref{proof L2 2ge2a},\eqref{proof decay 2ge2a}
&\eqref{proof decay 2ge2a},\eqref{proof L2 2ge2a}
&\eqref{proof decay 2ge2a},\eqref{proof L2 2ge2a}
&\eqref{proof decay 2ge2a},\eqref{proof L2 2ge2a}
&\eqref{proof decay 2ge2a},\eqref{proof L2 2ge2a}
\end{array}
$$
\centerline{Table 5} 

\end{landscape}


For the term $T_0$ by \eqref{proof decay 2ge5b} and \eqref{proof L2 2ge1a}, we have 
$$
\aligned
T_0 =& C\sum_{|I_2|\leq |\Id|}\big{\|}st^{-1}u_{\ih} \,  Z^{I_2}\delu_{\beta'}u_{\jh}\big{\|}_{L^2(\Hcal_s)}
\\
\leq&  CC_1\eps\sum_{|I_2|\leq |\Id|}\big{\|}st^{-1}t^{-3/2}s \,  Z^{I_2}\delu_{\beta'}u_{\jh}\big{\|}_{L^2(\Hcal_s)}
\\
= & CC_1\eps\sum_{|I_2|\leq |\Id|}\big{\|}t^{-5/2}s^2 \,  Z^{I_2}\delu_{\beta'}u_{\jh}\big{\|}_{L^2(\Hcal_s)}
\\
\leq& CC_1\eps s^{-1/2}\sum_{|I_2|\leq |\Id|}\big{\|}Z^{I_2}\delu_{\beta'}u_{\jh}\big{\|}_{L^2(\Hcal_s)}
\\
\leq & C(C_1\eps)^2s^{-1/2+\delta}.
\endaligned
$$
The term
$T_1$ is estimated by \eqref{proof decay 2ge5b} and \eqref{proof L2 2ge3a}: 
$$
\aligned
T_1 =& C\!\!\!\!\sum_{|I_1|=1 \atop  |I_2|\leq |\Id|-1}\!\!\!\!
\big{\|}st^{-1}Z^{I_1}u_{\ih} \,  Z^{I_2}\delu_{\beta'}u_{\jh}\big{\|}_{L^2(\Hcal_s)}
\\
\leq&  CC_1\eps\!\!\!\! \sum_{|I_1|=1 \atop  |I_2|\leq |I^\dag|-1}\!\!\!\!
\big{\|}st^{-1} t^{-3/2}s \,  Z^{I_2}\delu_{\beta'}u_{\jh}\big{\|}_{L^2(\Hcal_s)}
\\
\leq & CC_1\eps s^{-1/2}\!\!\!\! \sum_{|I_1|=1 \atop  |I_2|\leq |\Id|-1}\!\!\!\!
\big{\|} Z^{I_2}\delu_{\beta'}u_{\jh}\big{\|}_{L^2(\Hcal_s)}
\\
\leq & C(C_1\eps)^2s^{-1/2+\delta}. 
\endaligned
$$
The term $T_2$ is estimated by \eqref{proof decay 2ge5a} and \eqref{proof L2 2ge3a}: 
$$
\aligned
T_2 =& C\!\!\!\!\sum_{|I_1|=2 \atop  |I_2|\leq |\Id|-2}\!\!\!\!
\big{\|}st^{-1}Z^{I_1}u_{\ih} \,  Z^{I_2}\delu_{\beta'}u_{\jh}\big{\|}_{L^2(\Hcal_s)}
\\
\leq&  CC_1\eps\!\!\!\! \sum_{|I_1|=2 \atop  |I_2|\leq |\Id|-2}\!\!\!\!
\big{\|}st^{-1} t^{-3/2}s^{1+\delta} \,  Z^{I_2}\delu_{\beta'}u_{\jh}\big{\|}_{L^2(\Hcal_s)}
\\
\leq & CC_1\eps s^{-1/2+\delta}\!\!\!\! \sum_{|I_1|=2 \atop  |I_2|\leq |\Id|-2}\!\!\!\!
\big{\|} Z^{I_2}\delu_{\beta'}u_{\jh}\big{\|}_{L^2(\Hcal_s)}
\\
\leq&  C(C_1\eps)^2s^{-1/2+\delta}. 
\endaligned
$$
The term $T_3$ is estimated by \eqref{proof decay 2ge5a} and \eqref{proof L2 2ge3a}:
$$
\aligned
T_3 =& C\!\!\!\!\sum_{|I_1|=3 \atop  |I_1|\leq |\Id|-3}\!\!\!\!
\big{\|}st^{-1}Z^{I_1}u_{\ih} \,  Z^{I_2}\delu_{\beta'}u_{\jh}\big{\|}_{L^2(\Hcal_s)}
\\
\leq&  CC_1\eps\!\!\!\! \sum_{|I_1|=3 \atop  |I_2|\leq |\Id|-3}\!\!\!\!
\big{\|}st^{-1} t^{-3/2}s^{1+\delta} \,  Z^{I_2}\delu_{\beta'}u_{\jh}\big{\|}_{L^2(\Hcal_s)}
\\
\leq & CC_1\eps s^{-1/2+\delta} \!\!\!\!\sum_{|I_1|=3 \atop  |I_2|\leq |\Id|-3}\!\!\!\!
\big{\|} Z^{I_2}\delu_{\beta'}u_{\jh}\big{\|}_{L^2(\Hcal_s)}
\\
\leq&  C(C_1\eps)^2s^{-1/2+\delta}
\endaligned
$$
The term
$T_4$ is estimated by \eqref{proof L2 2ge5a} and \eqref{proof decay 2ge3a}:
$$
\aligned
T_4 =& C\sum_{|I_2|=4}\big{\|}st^{-1}Z^{I_2}u_{\ih} \,  \delu_{\beta'}u_{\jh}\big{\|}_{L^2(\Hcal_s)}
\\
\leq& CC_1\eps \sum_{|I_2|=4}\big{\|}st^{-1}Z^{I_2}u_{\ih} \,  t^{-1/2}s^{-1}\big{\|}_{L^2(\Hcal_s)}
\\
=&CC_1\eps s^{-1/2}\sum_{|I_2|=4}\big{\|}t^{-1}Z^{I_2}u_{\ih}\big{\|}_{L^2(\Hcal_s)}
\\
\leq& C(C_1\eps)^2s^{-1/2+\delta}.
\endaligned
$$
For the remaining three terms, we list out the $L^2$ and $L^\infty$ estimates to be used for each term and each partition of $\Id$,  recall here the notation $(a,\leq b)$ means $|I_1|=a, |I_2|\leq b$: see Table 6. 
Finally we conclude with \eqref{proof 2order eq9a}.

We estimate the term ${Q_T}_{\ih}$ which is similar to that of ${Q_G}_{\ih}$. Recall the structure of ${Q_T}_{\ih}$ presented in \eqref{proof 2order eq QT}. We take the term $Z^{I_2}\del_{\gamma}u_{\kh}\,\delu_a\delu_{\beta}Z^{I_1}u_{\jh}$ as an example:
$$
\aligned
&\sum_{|I_1|+|I_2|
\leq |\Id|}\big{\|}sZ^{I_2}\del_{\gamma}u_{\kh}\,\delu_a\delu_{\beta}Z^{I_1}u_{\jh}\big{\|}_{L^2(\Hcal_s)}
\\
= &\sum_{|I_2|\leq |\Id|}\big{\|}sZ^{I_2}\del_{\gamma}u_{\kh}\,\delu_a\delu_{\beta}u_{\jh}\big{\|}_{L^2(\Hcal_s)}
+\sum_{|I_1|+|I_2|\leq |\Id|\atop |I_1|=1}
\big{\|}sZ^{I_2}\del_{\gamma}u_{\kh}\,\delu_a\delu_{\beta}Z^{I_1}u_{\jh}\big{\|}_{L^2(\Hcal_s)}
\\
&+\sum_{|I_1|+|I_2|\leq |\Id|\atop |I_1|=2}
\big{\|}sZ^{I_2}\del_{\gamma}u_{\kh}\,\delu_a\delu_{\beta}Z^{I_1}u_{\jh}\big{\|}_{L^2(\Hcal_s)}
\\&
+\sum_{|I_1|+|I_2|\leq |\Id|\atop |I_1|=3}
\big{\|}sZ^{I_2}\del_{\gamma}u_{\kh}\,\delu_a\delu_{\beta}Z^{I_1}u_{\jh}\big{\|}_{L^2(\Hcal_s)}
\\
&+\sum_{|I_1|\leq |\Id|}
\big{\|}s\del_{\gamma}u_{\kh}\,\delu_a\delu_{\beta}Z^{I_1}u_{\jh}\big{\|}_{L^2(\Hcal_s)}
\\
=:& T_0 + T_1 + T_2 + T_3 + T_4.
\endaligned
$$


\protect\begin{landscape}\thispagestyle{empty}
$$
\begin{array}{cccccc}
\text{Terms}&(4,\leq 0)&(3,\leq 1)&(2,\leq 2)&(1,\leq 3)&(0,\leq 4)
\\
u_{\ih}\delu_{\beta'}u_{\jh}
&\eqref{proof L2 2ge5a'},\eqref{proof decay 2ge3a}
&\eqref{proof decay 2ge5a},\eqref{proof L2 2ge3a}
&\eqref{proof decay 2ge5a},\eqref{proof L2 2ge3a}
&\eqref{proof decay 2ge5b},\eqref{proof L2 2ge3a}
&\eqref{proof decay 2ge5b},\eqref{proof L2 2ge1a}
\\
v_{\ic}\delu_{\beta'}u_{\jh}
&\eqref{proof L2 2ge1c},\eqref{proof decay 2ge3a}
&\eqref{proof decay 2ge1c},\eqref{proof L2 2ge3a}
&\eqref{proof decay 2ge1c},\eqref{proof L2 2ge3a}
&\eqref{proof decay 2ge1c},\eqref{proof L2 2ge3a}
&\eqref{proof decay 2ge1c},\eqref{proof L2 2ge1a}
\\
\del_{\gamma}u_{\ih}\delu_{\beta'}u_{\jh}
&\eqref{proof L2 2ge1a},\eqref{proof decay 2ge3a}
&\eqref{proof decay 2ge1a},\eqref{proof L2 2ge3a}
&\eqref{proof decay 2ge1a},\eqref{proof L2 2ge3a}
&\eqref{proof decay 2ge3a},\eqref{proof L2 2ge3a}
&\eqref{proof decay 2ge3a},\eqref{proof L2 2ge1a}
\\
\del_{\gamma}v_{\ic}\delu_{\beta'}u_{\jh}
&\eqref{proof L2 2ge2a},\eqref{proof decay 2ge3a}
&\eqref{proof decay 2ge2a},\eqref{proof L2 2ge3a}
&\eqref{proof decay 2ge2a},\eqref{proof L2 2ge3a}
&\eqref{proof decay 2ge2a},\eqref{proof L2 2ge3a}
&\eqref{proof decay 2ge2a},\eqref{proof L2 2ge1a}
\end{array}
$$
\centerline{Table 6} 

\

\

$$
\begin{array}{cccccc}
\text{Products} &(4,\leq 0)&(3,\leq 1)&(2,\leq2)&(1,\leq 3)&(0,\leq 4)
\\
\del_\alpha u_{\jh}\del_{\beta}u_{\kh}
&\eqref{proof L2 2ge1a},\eqref{proof decay 2ge3a}
&\eqref{proof L2 2ge3a},\eqref{proof decay 2ge3a}
&\eqref{proof decay 2ge1a},\eqref{proof L2 2ge3a}
&\eqref{proof decay 2ge3a},\eqref{proof L2 2ge3a}
&\eqref{proof decay 2ge3a},\eqref{proof L2 2ge1a}
\\
\del_\alpha v_{\jh}\del_{\beta}w_k
&\eqref{proof L2 2ge2a},\eqref{proof decay 2ge1a}
&\eqref{proof L2 2ge2a},\eqref{proof decay 2ge1a}
&\eqref{proof decay 2ge2a},\eqref{proof L2 2ge1a}
&\eqref{proof decay 2ge2a},\eqref{proof L2 2ge1a}
&\eqref{proof decay 2ge2a},\eqref{proof L2 2ge1a}
\\
v_{\kh}\del_\alpha w_j
&\eqref{proof L2 2ge1c},\eqref{proof decay 2ge1a}
&\eqref{proof L2 2ge1c},\eqref{proof decay 2ge1a}
&\eqref{proof decay 2ge1c},\eqref{proof L2 2ge1a}
&\eqref{proof decay 2ge1c},\eqref{proof L2 2ge1a}
&\eqref{proof decay 2ge1c},\eqref{proof L2 2ge1a}
\\
v_{\kh}v_{\jh}
&\eqref{proof L2 2ge1c},\eqref{proof decay 2ge1c}
&\eqref{proof L2 2ge1c},\eqref{proof decay 2ge1c}
&\eqref{proof decay 2ge1c},\eqref{proof L2 2ge1c}
&\eqref{proof decay 2ge1c},\eqref{proof L2 2ge1c}
&\eqref{proof decay 2ge1c},\eqref{proof L2 2ge1c}
\end{array}
$$
\centerline{Table 7} 

\end{landscape}


We will estimate each partition of $\Id$: 
$$
\aligned
T_0 &= \sum_{|I_2|\leq |\Id|}\big{\|}sZ^{I_2}\del_{\gamma}u_{\kh}\,\delu_a\delu_{\beta}u_{\jh}\big{\|}_{L^2(\Hcal_s)}
\leq \sum_{|I_2|\leq 4}\big{\|}sZ^{I_2}\del_{\gamma}u_{\kh}\,\delu_a\delu_{\beta}u_{\jh}\big{\|}_{L^2(\Hcal_s)}
\\
&\leq  \sum_{|I_2|\leq 4}\big{\|}(s/t)Z^{I_2}\del_{\gamma}u_{\kh}\big{\|}_{L^2(\Hcal_s)}\,
\big{\|}s (t/s)\delu_a\delu_{\beta}u_{\jh} \big{\|}_{L^\infty(\Hcal_s)}\leq CC_1\eps s^{\delta} \,  CC_1\eps s^{-3/2}
\\
&\leq C(C_1\eps)^2s^{-3/2+\delta},
\endaligned
$$
where \eqref{proof L2 2ge1a} and \eqref{proof decay 1ge4b} are used.
$$
\aligned
T_1 &\leq \sum_{|I_1|=1,|I_2|\leq 3}
\big{\|}sZ^{I_2}\del_{\gamma}u_{\kh}\,\delu_a\delu_{\beta}Z^{I_1}u_{\jh}\big{\|}_{L^2(\Hcal_s)}
\\
&\leq \sum_{|I_1|=1,|I_2|\leq 3}
\big{\|}(s/t)Z^{I_2}\del_{\gamma}u_{\kh}\big{\|}_{L^2(\Hcal_s)}
\big{\|}s(t/s)\delu_a\delu_{\beta}Z^{I_1}u_{\jh}\big{\|}_{L^\infty(\Hcal_s)}
\\
&\leq CC_1\eps  \,  CC_1\eps s^{-3/2+\delta} \leq C(C_1\eps)^2s^{-3/2+\delta},
\endaligned
$$
where \eqref{proof L2 2ge3a} and \eqref{proof decay 1ge4a} are used.

Similarly, the term $T_2$, $T_3$ and $T_4$ are estimated by applying respectively \eqref{proof L2 2ge3a}\eqref{proof decay 1ge4a}, \eqref{proof L2 2ge4a}, 
\eqref{proof decay 2ge3a}, \eqref{proof L2 2ge4a}, and \eqref{proof decay 2ge3a}.

For the remaining terms, we will  not write in details the proof but list out the inequalities to be used on each term and each partition of the index ${|I_1|+|I_2|\leq |\Id|}$, in the following list:
$$
\begin{array}{cccccc}
\text{Products} &(4,\leq0) &(3,\leq1) &(2,\leq2)
\\
s\delu_a\delu_{\beta}Z^{I_1}u_{\ih}Z^{I_2}\del_{\gamma}u_{\kh}
&\eqref{proof L2 2ge4a},\eqref{proof decay 2ge3a}
&\eqref{proof L2 2ge4a},\eqref{proof decay 2ge3a}
&\eqref{proof L2 2ge4b},\eqref{proof decay 2ge1a}
\\
s\delu_a\delu_{\beta}Z^{I_1}u_{\ih}Z^{I_2}\del_{\gamma}v_{\kc}
&\eqref{proof L2 2ge4a},\eqref{proof decay 2ge2a}
&\eqref{proof L2 2ge4a},\eqref{proof decay 2ge2a}
&\eqref{proof L2 2ge4b},\eqref{proof decay 2ge2a}
\\
s\delu_a\delu_{\beta}Z^{I_1}u_{\ih}Z^{I_2}u_{\kh}
&\eqref{proof L2 2ge4a},\eqref{proof decay 2ge5b}
&\eqref{proof L2 2ge4a},\eqref{proof decay 2ge5b}
&\eqref{proof L2 2ge4b},\eqref{proof decay 2ge5a}
\\
s\delu_a\delu_{\beta}Z^{I_1}u_{\ih}Z^{I_2}v_{\kc}
&\eqref{proof L2 2ge4a},\eqref{proof decay 2ge1c}
&\eqref{proof L2 2ge4a},\eqref{proof decay 2ge1c}
&\eqref{proof L2 2ge4b},\eqref{proof decay 2ge1c}
\end{array}
$$

$$
\begin{array}{ccc}
\text{Products} &(1,\leq3) &(0,\leq4)
\\
s\delu_a\delu_{\beta}Z^{I_1}u_{\ih}Z^{I_2}\del_{\gamma}u_{\kh}
&\eqref{proof decay 1ge4a},\eqref{proof L2 2ge3a}
&\eqref{proof decay 1ge4b},\eqref{proof L2 2ge1a}
\\
s\delu_a\delu_{\beta}Z^{I_1}u_{\ih}Z^{I_2}\del_{\gamma}v_{\kc}
&\eqref{proof decay 1ge4a},\eqref{proof L2 2ge2a}
&\eqref{proof decay 1ge4b},\eqref{proof L2 2ge2a}
\\
s\delu_a\delu_{\beta}Z^{I_1}u_{\ih}Z^{I_2}u_{\kh}
&\eqref{proof decay 1ge4a},\eqref{proof L2 2ge5b}
&\eqref{proof decay 1ge4b},\eqref{proof L2 2ge5a}
\\
s\delu_a\delu_{\beta}Z^{I_1}u_{\ih}Z^{I_2}v_{\kc}
&\eqref{proof decay 1ge4a},\eqref{proof L2 2ge1c}
&\eqref{proof decay 1ge4b},\eqref{proof L2 2ge1c}
\end{array}
$$

$$
\begin{array}{cccccc}
\text{Products} &(4,\leq0) &(3,\leq1) &(2,\leq2)
\\
st^{-1}\del_{\gamma'}Z^{I_1}u_{\jh}Z^{I_2}\del_{\gamma}u_{\kh}
&\eqref{proof L2 1ge1a},\eqref{proof decay 2ge3a}
&\eqref{proof L2 1ge3a},\eqref{proof decay 2ge3a}
&\eqref{proof L2 1ge3a},\eqref{proof decay 2ge1a}
\\
st^{-1}\del_{\gamma'}Z^{I_1}u_{\jh}Z^{I_2}\del_{\gamma}v_{\kc}
&\eqref{proof L2 1ge1a},\eqref{proof decay 2ge2a}
&\eqref{proof L2 1ge3a},\eqref{proof decay 2ge2a}
&\eqref{proof L2 1ge3a},\eqref{proof decay 2ge2a}
\\
st^{-1}\del_{\gamma'}Z^{I_1}u_{\jh}Z^{I_2}u_{\kh}
&\eqref{proof L2 1ge1a},\eqref{proof decay 2ge5b}
&\eqref{proof L2 1ge3a},\eqref{proof decay 2ge5b}
&\eqref{proof L2 1ge3a},\eqref{proof decay 2ge5a}
\\
st^{-1}\del_{\gamma'}Z^{I_1}u_{\jh}Z^{I_2}v_{\kc}
&\eqref{proof L2 1ge1a},\eqref{proof decay 2ge1c}
&\eqref{proof L2 1ge3a},\eqref{proof decay 2ge1c}
&\eqref{proof L2 1ge3a},\eqref{proof decay 2ge1c}
\end{array}
$$

$$
\begin{array}{ccc}
\text{Products} &(1,\leq3) &(0,\leq4)
\\
st^{-1}\del_{\gamma'}Z^{I_1}u_{\jh}Z^{I_2}\del_{\gamma}u_{\kh}
&\eqref{proof decay 1ge3a},\eqref{proof L2 2ge3a}
&\eqref{proof decay 1ge3a},\eqref{proof L2 2ge1a}
\\
st^{-1}\del_{\gamma'}Z^{I_1}u_{\jh}Z^{I_2}\del_{\gamma}v_{\kc}
&\eqref{proof decay 1ge3a},\eqref{proof L2 2ge2a}
&\eqref{proof decay 1ge3a},\eqref{proof L2 2ge2a}
\\
st^{-1}\del_{\gamma'}Z^{I_1}u_{\jh}Z^{I_2}u_{\kh}
&\eqref{proof decay 1ge3a},\eqref{proof L2 2ge5b}
&\eqref{proof decay 1ge3a},\eqref{proof L2 2ge5a}
\\
st^{-1}\del_{\gamma'}Z^{I_1}u_{\jh}Z^{I_2}v_{\kc}
&\eqref{proof decay 1ge3a},\eqref{proof L2 2ge1c}
&\eqref{proof decay 1ge3a},\eqref{proof L2 2ge1c}
\end{array}
$$

\

Now, we estimate the term $F_{\ih}$ which is also similar: recall the structure of $Z^{\Id}F_{\ih}$ presented in \eqref{proof 2order eq5 ZF} with $I$ replaced by $\Id$. As before, we consider the term $Z^{\Id}\big(\del_\alpha u_{\jh}\del_{\beta}u_{\kh}\big)$ as an example and we write down the details of the analysis.
For the rest terms, we just give the $L^2$ and $L^\infty$ estimates to be used for each factor:
$$
\aligned
&\big{\|}sZ^{\Id}\big(\del_\alpha u_{\jh}\del_{\beta}u_{\kh}\big)\big{\|}_{L^2(\Hcal_s)}
\\
&\leq \big{\|}s\del_\alpha u_{\jh} \,  Z^{\Id}\del_{\beta}u_{\kh}\big{\|}_{L^2(\Hcal_s)}
  +\!\!\!\!\!\sum_{I_1+I_2=\Id\atop|I_1|=1}\!\!\!\!\!
  \big{\|}sZ^{I_1}\del_\alpha u_{\jh} \,  Z^{I_2}\del_{\beta}u_{\kh}\big{\|}_{L^2(\Hcal_s)}
\\
&\quad
 +\!\!\!\!\!\sum_{I_1+I_2=\Id\atop|I_1|=2}\!\!\!\!\!
  \big{\|}sZ^{I_1}\del_\alpha u_{\jh} \,  Z^{I_2}\del_{\beta}u_{\kh}\big{\|}_{L^2(\Hcal_s)}
\\
  &\quad +\!\!\!\!\sum_{I_1+I_2=\Id\atop|I_1|=3}\!\!\!\!
  \big{\|}sZ^{I_1}\del_\alpha u_{\jh} \,  Z^{I_2}\del_{\beta}u_{\kh}\big{\|}_{L^2(\Hcal_s)}
 +\big{\|}s Z^{\Id}\del_\alpha u_{\jh} \, \del_{\beta}u_{\kh}\big{\|}_{L^2(\Hcal_s)}
\\
&=: T_0 + T_1 +T_2 + T_3 + T_4.
\endaligned
$$
The term $T_0$ is estimated by \eqref{proof decay 2ge3a} and \eqref{proof L2 2ge1a}:
$$
\aligned
T_0 =& \big{\|}s\del_\alpha u_{\jh} \,  Z^{\Id}\del_{\beta}u_{\kh}\big{\|}_{L^2(\Hcal_s)}
\leq  CC_1\eps  \big{\|}s \,  t^{-1/2}s^{-1} \,  Z^{\Id}\del_{\beta}u_{\kh}\big{\|}_{L^2(\Hcal_s)}
\\
=& CC_1\eps  \big{\|}s\, t^{-1/2}s^{-1} \,  (t/s)\big(s/tZ^{\Id}\del_{\beta}u_{\kh}\big)\big{\|}_{L^2(\Hcal_s)}
\leq CC_1\eps  \big{\|}s/tZ^{\Id}\del_{\beta}u_{\kh}\big{\|}_{L^2(\Hcal_s)}
\\
\leq& C(C_1\eps)^2s^{\delta}.
\endaligned
$$
Here the relation $t\leq Cs^2$ (in $\Kcal$) is used. 
The term $T_4$ is estimated in the same way by exchanging
 the role of $u_{\jh}$ and $u_{\kh}$.

The terms $T_1$ and $T_3$ are estimated by \eqref{proof decay 2ge3a} and \eqref{proof L2 2ge3a}. We estimate $T_1$ as follows. The estimate of $T_3$ is done by exchanging $u_{\jh}$ and $u_{\kh}$ in the following argument: 
$$
\aligned
T_1 =& \sum_{I_1+I_2=\Id \atop  |I_1|=1}\!\!\!\!
\big{\|}sZ^{I_1}\del_\alpha u_{\jh} \,  Z^{I_2}\del_{\beta}u_{\kh}\big{\|}_{L^2(\Hcal_s)}
\\
\leq& CC_1\eps\big{\|}s\, t^{-1/2}s^{-1} \, (t/s) (s/t)Z^{I_2}\del_{\beta}u_{\kh}\big{\|}_{L^2(\Hcal_s)}
\\
\leq & CC_1\eps \big{\|}(s/t)Z^{I_2}\del_{\beta}u_{\kh}\big{\|}_{L^2(\Hcal_s)}
\leq  C(C_1\eps)^2.
\endaligned
$$


\protect\begin{landscape}\thispagestyle{empty}
$$
\begin{array}{ccccc}
\text{Products} &(3,\leq 0) &(2,\leq 1) &(1,\leq 2)&(0,\leq 3)
\\
u_{\ih}\delu_a\delu_{\beta} u_{\jh}
&\eqref{proof L2 2ge5b},\eqref{proof decay 2ge4b}
&\eqref{proof L2 2ge5b},\eqref{proof decay 2ge4a'}
&\eqref{proof decay 2ge5b},\eqref{proof L2 2ge4'b}
&\eqref{proof decay 2ge5b},\eqref{proof L2 2ge4'a'}
\\
v_{\ic}\delu_a\delu_{\beta} u_{\jh}
&\eqref{proof L2 2ge1c'},\eqref{proof decay 2ge4b}
&\eqref{proof L2 2ge1c'},\eqref{proof decay 2ge4a'}
&\eqref{proof decay 2ge1c'},\eqref{proof L2 2ge4'b}
&\eqref{proof decay 2ge1c'},\eqref{proof L2 2ge4'a'}
\\
u_{\ih}\del_t\delu_au_{\jh}
&\eqref{proof L2 2ge5b},\eqref{proof decay 2ge4b}
&\eqref{proof L2 2ge5b},\eqref{proof decay 2ge4a'}
&\eqref{proof decay 2ge5b},\eqref{proof L2 2ge4'b}
&\eqref{proof decay 2ge5b},\eqref{proof L2 2ge4'a'}
\\
v_{\ic}\del_t\delu_au_{\jh}
&\eqref{proof L2 2ge1c'},\eqref{proof decay 2ge4b}
&\eqref{proof L2 2ge1c'},\eqref{proof decay 2ge4a'}
&\eqref{proof decay 2ge1c'},\eqref{proof L2 2ge4'b}
&\eqref{proof decay 2ge1c'},\eqref{proof L2 2ge4'a'}
\\
\delu_{\gamma}u_{\ih}\delu_a\delu_{\beta} u_{\jh}
&\eqref{proof L2 2ge3a},\eqref{proof decay 2ge4b}
&\eqref{proof L2 2ge3a},\eqref{proof decay 2ge4a'}
&\eqref{proof decay 2ge3a},\eqref{proof L2 2ge4'b}
&\eqref{proof decay 2ge3a},\eqref{proof L2 2ge4'a'}
\\
\delu_{\gamma}v_{\ic}\delu_a\delu_{\beta} u_{\jh}
&\eqref{proof L2 2ge2a'},\eqref{proof decay 2ge4b}
&\eqref{proof L2 2ge2a'},\eqref{proof decay 2ge4a'}
&\eqref{proof decay 2ge2a'},\eqref{proof L2 2ge4'b}
&\eqref{proof decay 2ge2a'},\eqref{proof L2 2ge4'a'}
\end{array}
$$
\centerline{Table 8}

\

$$
\begin{array}{ccccc}
\text{Products}&(3,\leq 0)&(2,\leq 1)&(1,\leq 2) &(0,\leq 3)
\\
\delu_{\gamma}u_{\ih}\del_t\delu_au_{\jh}
&\eqref{proof L2 2ge3a},\eqref{proof decay 2ge4b}
&\eqref{proof L2 2ge3a},\eqref{proof decay 2ge4a'}
&\eqref{proof decay 2ge3a},\eqref{proof L2 2ge4'b}
&\eqref{proof decay 2ge3a},\eqref{proof L2 2ge4'a'}
\\
\delu_{\gamma}v_{\ic}\del_t\delu_au_{\jh}
&\eqref{proof L2 2ge2a'},\eqref{proof decay 2ge4b}
&\eqref{proof L2 2ge2a'},\eqref{proof decay 2ge4a'}
&\eqref{proof decay 2ge2a'},\eqref{proof L2 2ge4'b}
&\eqref{proof decay 2ge2a'},\eqref{proof L2 2ge4'a'}
\\
v_{\ic}\del_\alpha \del_{\beta}v_{\jc}
&\eqref{proof L2 2ge1c'},\eqref{proof decay 2ge2a'}
&\eqref{proof L2 2ge1c'},\eqref{proof decay 2ge2a'}
&\eqref{proof decay 2ge1c'},\eqref{proof L2 2ge2a'}
&\eqref{proof decay 2ge1c'},\eqref{proof L2 2ge2a}
\\
\del_{\gamma}u_{\ih}\del_\alpha \del_{\beta}v_{\jc}
&\eqref{proof L2 2ge3a},\eqref{proof decay 2ge2a'}
&\eqref{proof L2 2ge3a},\eqref{proof decay 2ge2a'}
&\eqref{proof decay 2ge3a},\eqref{proof L2 2ge2a'}
&\eqref{proof decay 2ge3a},\eqref{proof L2 2ge2a}
\\
\del_{\gamma}v_{\ic}\del_\alpha \del_{\beta}v_{\jc}
&\eqref{proof L2 2ge2a'},\eqref{proof decay 2ge2a'}
&\eqref{proof L2 2ge2a'},\eqref{proof decay 2ge2a'}
&\eqref{proof decay 2ge2a'},\eqref{proof L2 2ge2a'}
&\eqref{proof decay 2ge2a'},\eqref{proof L2 2ge2a}
\end{array}
$$
\centerline{Table 9}

\

\newpage 

$$
\begin{array}{ccccc}
\text{Terms} &(3,\leq 0) &(2,\leq 1) &(1,\leq 2) &(0,\leq 3)
\\
\del_\alpha \Psi_{\beta}^{\beta'}u_{\ih}\delu_{\beta'}u_{\jh}
&\eqref{proof L2 2ge5b},\eqref{proof decay 2ge3a}
&\eqref{proof L2 2ge5b},\eqref{proof decay 2ge3a}
&\eqref{proof decay 2ge5b},\eqref{proof L2 2ge3a}
&\eqref{proof decay 2ge5b},\eqref{proof L2 2ge3a}
\\
\del_\alpha \Psi_{\beta}^{\beta'}v_{\ic}\delu_{\beta'}u_{\jh}
&\eqref{proof L2 2ge1c'},\eqref{proof decay 2ge3a}
&\eqref{proof L2 2ge1c'},\eqref{proof decay 2ge3a}
&\eqref{proof decay 2ge1c'},\eqref{proof L2 2ge3a}
&\eqref{proof decay 2ge1c'},\eqref{proof L2 2ge3a}
\\
\del_\alpha \Psi_{\beta}^{\beta'}\del_{\gamma}u_{\ih}\delu_{\beta'}u_{\jh}
&\eqref{proof L2 2ge3a},\eqref{proof decay 2ge3a}
&\eqref{proof L2 2ge3a},\eqref{proof decay 2ge3a}
&\eqref{proof decay 2ge3a},\eqref{proof L2 2ge3a}
&\eqref{proof decay 2ge3a},\eqref{proof L2 2ge3a}
\\
\del_\alpha \Psi_{\beta}^{\beta'}\del_{\gamma}v_{\ic}\delu_{\beta'}u_{\jh}
&\eqref{proof L2 2ge2a'},\eqref{proof decay 2ge3a}
&\eqref{proof L2 2ge2a'},\eqref{proof decay 2ge3a}
&\eqref{proof decay 2ge2a'},\eqref{proof L2 2ge3a}
&\eqref{proof decay 2ge2a'},\eqref{proof L2 2ge3a}
\end{array}
$$
\centerline{Table 10}

\

$$
\begin{array}{ccccc}
\text{Products} &(3,\leq0)&(2,\leq1)&(1,\leq2)&(0,\leq 3)
\\
s\delu_a\delu_{\beta}Z^{J_1}u_{\ih}\,Z^{J_2}\del_{\gamma}u_{\kh}
&\eqref{proof L2 2ge4a'},\eqref{proof decay 2ge3a}
&\eqref{proof L2 2ge4b},\eqref{proof decay 2ge3a}
&\eqref{proof decay 1ge4a'},\eqref{proof L2 2ge3a}
&\eqref{proof decay 1ge4b},\eqref{proof L2 2ge3a}
\\
s\delu_a\delu_{\beta}Z^{J_1}u_{\ih}\,Z^{J_2}\del_{\gamma}v_{\kc}
&\eqref{proof L2 2ge4a'},\eqref{proof decay 2ge2a'}
&\eqref{proof L2 2ge4b},\eqref{proof decay 2ge2a'}
&\eqref{proof decay 1ge4a'},\eqref{proof L2 2ge2a'}
&\eqref{proof decay 1ge4b},\eqref{proof L2 2ge2a'}
\\
s\delu_a\delu_{\beta}Z^{J_1}u_{\ih}\,Z^{J_2}u_{\kh}
&\eqref{proof L2 2ge4a'},\eqref{proof decay 2ge5b}
&\eqref{proof L2 2ge4b},\eqref{proof decay 2ge5b}
&\eqref{proof decay 1ge4a'},\eqref{proof L2 2ge5b}
&\eqref{proof decay 1ge4b},\eqref{proof L2 2ge5b}
\\
s\delu_a\delu_{\beta}Z^{J_1}u_{\ih}\,Z^{J_2}v_{\kc}
&\eqref{proof L2 2ge4a'},\eqref{proof decay 2ge1c'}
&\eqref{proof L2 2ge4b},\eqref{proof decay 2ge1c'}
&\eqref{proof decay 1ge4a'},\eqref{proof L2 2ge1c'}
&\eqref{proof decay 1ge4b},\eqref{proof L2 2ge1c'}
\end{array}
$$
\centerline{Table 11}

\

$$
\begin{array}{cccccc}
\text{Products} &(3,\leq0)&(2,\leq1)&(1,\leq2) &(0,\leq 3)
\\
st^{-1}\del_{\gamma'}Z^{J_1}u_{\jh}\,Z^{J_2}\del_{\gamma}u_{\kh}
&\eqref{proof L2 1ge3a},\eqref{proof decay 2ge3a}
&\eqref{proof L2 1ge3a},\eqref{proof decay 2ge3a}
&\eqref{proof decay 1ge3a},\eqref{proof L2 2ge3a}
&\eqref{proof decay 1ge3a},\eqref{proof L2 2ge3a}
\\
st^{-1}\del_{\gamma'}Z^{J_1}u_{\jh}\,Z^{J_2}\del_{\gamma}v_{\kc}
&\eqref{proof L2 1ge3a},\eqref{proof decay 2ge2a'}
&\eqref{proof L2 1ge3a},\eqref{proof decay 2ge2a'}
&\eqref{proof decay 1ge3a},\eqref{proof L2 2ge2a'}
&\eqref{proof decay 1ge3a},\eqref{proof L2 2ge2a'}
\\
st^{-1}\del_{\gamma'}Z^{J_1}u_{\jh}\,Z^{J_2}u_{\kh}
&\eqref{proof L2 1ge3a},\eqref{proof decay 2ge5b}
&\eqref{proof L2 1ge3a},\eqref{proof decay 2ge5b}
&\eqref{proof decay 1ge3a},\eqref{proof L2 2ge5b}
&\eqref{proof decay 1ge3a},\eqref{proof L2 2ge5b}
\\
st^{-1}\del_{\gamma'}Z^{J_1}u_{\jh}\,Z^{J_2}v_{\kc}
&\eqref{proof L2 1ge3a},\eqref{proof decay 2ge1c'}
&\eqref{proof L2 1ge3a},\eqref{proof decay 2ge1c'}
&\eqref{proof decay 1ge3a},\eqref{proof L2 2ge1c'}
&\eqref{proof decay 1ge3a},\eqref{proof L2 2ge1c'}
\end{array}
$$
\centerline{Table 12}

\thispagestyle{empty}
\end{landscape}

The term $T_2$ is estimated by \eqref{proof decay 2ge1a} and \eqref{proof L2 2ge3a}:
$$
\aligned
T_2 =& \sum_{I_1+I_2=\Id \atop  |I_1|=2}\!\!\!\!
\big{\|}sZ^{I_1}\del_\alpha u_{\jh} \,  Z^{I_2}\del_{\beta}u_{\kh}\big{\|}_{L^2(\Hcal_s)}
\\
\leq& CC_1\eps\big{\|}
s\, t^{-1/2}s^{-1+\delta} \, (t/s) (s/t)Z^{I_2}\del_{\beta}u_{\kh}\big{\|}_{L^2(\Hcal_s)}
\\
= & CC_1\eps s^{\delta} \big{\|}(s/t)Z^{I_2}\del_{\beta}u_{\kh}\big{\|}_{L^2(\Hcal_s)}
\leq  C(C_1\eps)^2s^{\delta}.
\endaligned
$$
The estimate of other terms are presented in Table 7.  
We conclude \eqref{proof 2order eq9b}.

The estimate of $R(Z^Iu_{\ih})$ is a direct result of \eqref{proof L2 2ge4a}.
\end{proof}


\begin{proof}[Proof of Lemma \ref{proof 2order lem7'}]
The proof is essentially the same as the one of Lemma \ref{proof 2order lem7}. The 
main difference lies in the inequalities we use for each term and partition of the index. We will list out the relevant inequalities and skip the details.

For the proof of \eqref{proof 2order eq9'a}, we list out the inequalities in Table 8.

The following four terms are estimated by apply the additional decay rate supplied by $\del_{\alpha}\Phi_{\beta}^{\gamma}$:
$$
\aligned
& sZ^I\big(\del_\alpha \Psi_{\beta}^{\beta'}u_{\ih}\delu_{\beta'}u_{\jh}\big),\quad
sZ^I\big(\del_\alpha \Psi_{\beta}^{\beta'}v_{\ic}\delu_{\beta'}u_{\jh}\big),\quad
\\
& sZ^I\big(\del_\alpha \Psi_{\beta}^{\beta'}\del_{\gamma}u_{\ih}\delu_{\beta'}u_{\jh}\big),\quad
sZ^I\big(\del_\alpha \Psi_{\beta}^{\beta'}\del_{\gamma}v_{\ic}\delu_{\beta'}u_{\jh}\big).
\endaligned
$$
See Table 10.

For the term ${Q_T}_{\ih}$, we find Table 11 and Table 12. 


For the estimates on $Z^IF_{\ih}$, the inequalities we use are presented in the following list:
$$
\begin{array}{ccccc}
\text{Products}&(3,\leq 0)&(2,\leq 1) &(1,\leq 2) &(0,\leq3)
\\
\del_\alpha u_{\jh}\del_{\beta}u_{\kh}
&\eqref{proof L2 2ge3a},\eqref{proof decay 2ge3a}
&\eqref{proof L2 2ge3a},\eqref{proof decay 2ge3a}
&\eqref{proof decay 2ge3a},\eqref{proof L2 2ge3a}
&\eqref{proof decay 2ge3a},\eqref{proof L2 2ge3a}
\\
\del_\alpha v_{\jc}\del_{\beta}w_k
&\eqref{proof L2 2ge2a'},\eqref{proof decay 2ge1a'}
&\eqref{proof L2 2ge2a'},\eqref{proof decay 2ge1a'}
&\eqref{proof decay 2ge2a'},\eqref{proof L2 2ge1a'}
&\eqref{proof decay 2ge2a'},\eqref{proof L2 2ge1a'}
\\
v_{\kc}\del_\alpha w_j
&\eqref{proof L2 2ge1c'},\eqref{proof decay 2ge1a'}
&\eqref{proof L2 2ge1c'},\eqref{proof decay 2ge1a'}
&\eqref{proof decay 2ge1c'},\eqref{proof L2 2ge1a'}
&\eqref{proof decay 2ge1c'},\eqref{proof L2 2ge1a'}
\\
v_{\kc}v_{\jc}
&\eqref{proof L2 2ge1c'},\eqref{proof decay 2ge1c'}
&\eqref{proof L2 2ge1c'},\eqref{proof decay 2ge1c'}
&\eqref{proof decay 2ge1c'},\eqref{proof L2 2ge1c'}
&\eqref{proof decay 2ge1c'},\eqref{proof L2 2ge1c'}
\end{array}
$$

The estimate on $R(Z^Iu_{\ih})$ is a direct result of \eqref{proof L2 2ge4a'}
\end{proof}


\begin{proof}[Proof of Lemma \ref{proof 2order lem8}]
The proof is essentially the same as
the one of Lemma \ref{proof 2order lem7}. The main difference is the level of
 regularity under consideration. We will not give the details and only list the inequalities we use for each term and partition $I = I_2+I_3$:
$$
\begin{array}{cccc}
\text{Products}&(2,\leq 0) &(1,\leq 1)&(0,\leq 2)
\\
u_{\ih}\delu_a\del_{\beta}u_{\jh}
&\eqref{proof L2 2ge5b},\eqref{proof decay 2ge4b}
&\eqref{proof decay 2ge5b},\eqref{proof L2 2ge4'b}
&\eqref{proof decay 2ge5b},\eqref{proof L2 2ge4'b}
\\
v_{\ic}\delu_a\del_{\beta}u_{\jh}
&\eqref{proof L2 2ge1c},\eqref{proof decay 2ge4b}
&\eqref{proof decay 2ge1c},\eqref{proof L2 2ge4'b}
&\eqref{proof decay 2ge1c},\eqref{proof L2 2ge4'b}
\\
u_{\ih}\del_t\delu_au_{\jh}
&\eqref{proof L2 2ge5b},\eqref{proof decay 2ge4b}
&\eqref{proof decay 2ge5b},\eqref{proof L2 2ge4'b}
&\eqref{proof decay 2ge5b},\eqref{proof L2 2ge4'b}
\\
v_{\ic}\del_t\delu_au_{\jh}
&\eqref{proof L2 2ge1c},\eqref{proof decay 2ge4b}
&\eqref{proof decay 2ge1c},\eqref{proof L2 2ge4'b}
&\eqref{proof decay 2ge1c},\eqref{proof L2 2ge4'b}
\end{array}
$$

$$
\begin{array}{cccc}
\text{Products} &(2,\leq 0) &(1,\leq 1) &(0,\leq 2)
\\
\delu_{\gamma}u_{\ih}\delu_a\del_{\beta}u_{\jh}
&\eqref{proof L2 2ge3a},\eqref{proof decay 2ge4b}
&\eqref{proof decay 2ge3a},\eqref{proof L2 2ge4'b}
&\eqref{proof decay 2ge3a},\eqref{proof L2 2ge4'b}
\\
\delu_{\gamma}v_{\ic}\delu_a\del_{\beta}u_{\jh}
&\eqref{proof L2 2ge2a},\eqref{proof decay 2ge4b}
&\eqref{proof decay 2ge2a},\eqref{proof L2 2ge4'b}
&\eqref{proof decay 2ge2a},\eqref{proof L2 2ge4'b}
\\
\delu_{\gamma}u_{\ih}\del_t\delu_au_{\jh}
&\eqref{proof L2 2ge3a},\eqref{proof decay 2ge4b}
&\eqref{proof decay 2ge3a},\eqref{proof L2 2ge4'b}
&\eqref{proof decay 2ge3a},\eqref{proof L2 2ge4'b}
\\
\delu_{\gamma}v_{\ic}\del_t\delu_au_{\jh}
&\eqref{proof L2 2ge2a},\eqref{proof decay 2ge4b}
&\eqref{proof decay 2ge2a},\eqref{proof L2 2ge4'b}
&\eqref{proof decay 2ge2a},\eqref{proof L2 2ge4'b}
\\
v_{\ic}\del_\alpha \del_{\beta}v_{\jc}
&\eqref{proof L2 2ge1c},\eqref{proof decay 2ge2a}
&\eqref{proof decay 2ge1c},\eqref{proof L2 2ge2a}
&\eqref{proof decay 2ge1c},\eqref{proof L2 2ge2a}
\\
\del_{\gamma}u_{\ih}\del_\alpha \del_{\beta}v_{\jc}
&\eqref{proof L2 2ge3a},\eqref{proof decay 2ge2a}
&\eqref{proof decay 2ge3a},\eqref{proof L2 2ge2a}
&\eqref{proof decay 2ge3a},\eqref{proof L2 2ge2a}
\\
\del_{\gamma}v_{\ic}\del_\alpha \del_{\beta}v_{\jc}
&\eqref{proof L2 2ge2a},\eqref{proof decay 2ge2a}
&\eqref{proof decay 2ge2a},\eqref{proof L2 2ge2a}
&\eqref{proof decay 2ge2a},\eqref{proof L2 2ge2a}
\end{array}
$$

There are four terms to be estimated separately:
$$
\aligned
& Z^{I}\big(\del_\alpha \Psi_{\beta}^{\beta'}u_{\ih}\delu_{\beta'}u_{\jh}\big),\quad Z^{I}\big(\del_\alpha \Psi_{\beta}^{\beta'}v_{\ic}\delu_{\beta'}u_{\jh}\big),\quad
\\
&
Z^{I}\big(\del_\alpha \Psi_{\beta}^{\beta'}\del_{\gamma}u_{\ih}\delu_{\beta'}u_{\jh}\big),\quad
Z^{I}\big(\del_\alpha \Psi_{\beta}^{\beta'}\del_{\gamma}v_{\ic}\delu_{\beta'}u_{\jh}\big).
\endaligned
$$
As before, these terms are to be estimated by the additional decay supplied by $|Z^I\Psi_{\beta}^{\beta'}|\leq C(I)t^{-1}$. We omit the details but list out the inequalities to be used for each term and each partition of $I=I_2+I_3$:
$$
\begin{array}{cccc}
\text{Products}  &(2,0) &(1,1) &(0,2)
\\
\del_\alpha \Psi_{\beta}^{\beta'}u_{\ih}\delu_{\beta'}u_{\jh}
&\eqref{proof L2 2ge5b},\eqref{proof decay 2ge3a}
&\eqref{proof decay 2ge5b},\eqref{proof L2 2ge3a}
&\eqref{proof decay 2ge5b},\eqref{proof L2 2ge3a}
\\
\del_\alpha \Psi_{\beta}^{\beta'}v_{\ic}\delu_{\beta'}u_{\jh}
&\eqref{proof L2 2ge1c},\eqref{proof decay 2ge3a}
&\eqref{proof decay 2ge1c},\eqref{proof L2 2ge3a}
&\eqref{proof decay 2ge1c},\eqref{proof L2 2ge3a}
\\
\del_\alpha \Psi_{\beta}^{\beta'}\del_{\gamma}u_{\ih}\delu_{\beta'}u_{\jh}
&\eqref{proof L2 2ge3a},\eqref{proof decay 2ge3a}
&\eqref{proof decay 2ge3a},\eqref{proof L2 2ge3a}
&\eqref{proof decay 2ge3a},\eqref{proof L2 2ge3a}
\\
\del_\alpha \Psi_{\beta}^{\beta'}\del_{\gamma}v_{\ic}\delu_{\beta'}u_{\jh}
&\eqref{proof L2 2ge2a},\eqref{proof decay 2ge3a}
&\eqref{proof decay 2ge2a},\eqref{proof L2 2ge3a}
&\eqref{proof decay 2ge2a},\eqref{proof L2 2ge3a}
\end{array}
$$
And we conclude with \eqref{proof 2order eq10a}.

We turn our attention to the estimates for ${Q_T}_{\ih}$. As before the details are omitted. The inequalities we use for each term and each partition of the indices are listed:
$$
\begin{array}{cccc}
\text{Products} &(2,\leq 0) &(1,\leq 1) &(0,\leq 2)
\\
s\delu_a\delu_{\beta}Z^{I_1}u_{\ih}Z^{I_2}\del_{\gamma}u_{\kh}
&\eqref{proof L2 2ge4b},\eqref{proof decay 2ge3a}
&\eqref{proof L2 2ge4b},\eqref{proof decay 2ge3a}
&\eqref{proof decay 1ge4b},\eqref{proof L2 2ge3a}
\\
s\delu_a\delu_{\beta}Z^{I_1}u_{\ih}Z^{I_2}\del_{\gamma}v_{\kc}
&\eqref{proof L2 2ge4b},\eqref{proof decay 2ge2a}
&\eqref{proof L2 2ge4b},\eqref{proof decay 2ge2a}
&\eqref{proof decay 1ge4b},\eqref{proof L2 2ge2a}
\\
s\delu_a\delu_{\beta}Z^{I_1}u_{\ih}Z^{I_2}u_{\kh}
&\eqref{proof L2 2ge4b},\eqref{proof decay 2ge5b}
&\eqref{proof L2 2ge4b},\eqref{proof decay 2ge5b}
&\eqref{proof decay 1ge4b},\eqref{proof L2 2ge5b}
\\
s\delu_a\delu_{\beta}Z^{I_1}u_{\ih}Z^{I_2}v_{\kc}
&\eqref{proof L2 2ge4b},\eqref{proof decay 2ge1c}
&\eqref{proof L2 2ge4b},\eqref{proof decay 2ge1c}
&\eqref{proof decay 1ge4b},\eqref{proof L2 2ge1c}
\end{array}
$$

$$
\begin{array}{cccc}
\text{Products} &(2,\leq 0) &(1,\leq 1) &(0,\leq 2)
\\
st^{-1}\del_{\gamma'}Z^{I_1}u_{\jh}Z^{I_2}\del_{\gamma}u_{\kh}
&\eqref{proof L2 1ge3a},\eqref{proof decay 2ge3a}
&\eqref{proof L2 1ge3a},\eqref{proof decay 2ge3a}
&\eqref{proof decay 1ge3a},\eqref{proof L2 2ge3a}
\\
st^{-1}\del_{\gamma'}Z^{I_1}u_{\jh}Z^{I_2}\del_{\gamma}v_{\kc}
&\eqref{proof L2 1ge3a},\eqref{proof decay 2ge2a}
&\eqref{proof L2 1ge3a},\eqref{proof decay 2ge2a}
&\eqref{proof decay 1ge3a},\eqref{proof L2 2ge2a}
\\
st^{-1}\del_{\gamma'}Z^{I_1}u_{\jh}Z^{I_2}u_{\kh}
&\eqref{proof L2 1ge3a},\eqref{proof decay 2ge5b}
&\eqref{proof L2 1ge3a},\eqref{proof decay 2ge5b}
&\eqref{proof decay 1ge3a},\eqref{proof L2 2ge5b}
\\
st^{-1}\del_{\gamma'}Z^{I_1}u_{\jh}Z^{I_2}v_{\kc}
&\eqref{proof L2 1ge3a},\eqref{proof decay 2ge1c}
&\eqref{proof L2 1ge3a},\eqref{proof decay 2ge1c}
&\eqref{proof decay 1ge3a},\eqref{proof L2 2ge1c}
\end{array}
$$

The estimate of $Z^IF_{\ih}$ is essentially the same. As before, we omit the details but list out the inequalities to be used:
$$
\begin{array}{cccc}
\text{Products}&(2,0) &(1,1) &(0,2)
\\
\del_\alpha u_{\jh}\del_{\beta}u_{\kh}
&\eqref{proof L2 2ge3a},\eqref{proof decay 2ge3a}
&\eqref{proof decay 2ge3a},\eqref{proof L2 2ge3a}
&\eqref{proof decay 2ge3a},\eqref{proof L2 2ge3a}
\\
\del_\alpha v_{\jc}\del_{\beta}w_k
&\eqref{proof L2 2ge2a},\eqref{proof decay 2ge1a}
&\eqref{proof decay 2ge2a},\eqref{proof L2 2ge1a}
&\eqref{proof decay 2ge2a},\eqref{proof L2 2ge1a}
\\
v_{\kc}\del_\alpha w_j
&\eqref{proof L2 2ge1c},\eqref{proof decay 2ge1a}
&\eqref{proof decay 2ge1c},\eqref{proof L2 2ge1a}
&\eqref{proof decay 2ge1c},\eqref{proof L2 2ge1a}
\\
v_{\kc}v_{\jc}
&\eqref{proof L2 2ge1c},\eqref{proof decay 2ge1c}
&\eqref{proof decay 2ge1c},\eqref{proof L2 2ge1c}
&\eqref{proof decay 2ge1c},\eqref{proof L2 2ge1c}
\end{array}
$$

The estimate of $R(Z^I u_{\ih})$ is a direct result of \eqref{proof L2 2ge4b}.
\end{proof}

Now we are ready to prove the second main result of this section.

\begin{proposition}\label{proof 2order prop L2}
Let $u_{\ih}$ be wave components of a sufficiently regular, local-in-time solution to \eqref{main eq main} and assume that \eqref{proof energy assumption} holds with $C_1\eps \leq \min\{1,\eps_0''\}$. The following estimates hold for $|\Id|\leq 4$, $|I|\leq 3$ and $|\If|\leq 2$:
\begin{subequations}\label{proof 2order eq 11}
\bel{proof 2order eq 11a}
\|s^{3}t^{-2}\del_t\del_tZ^{\Id}u_{\ih}\|_{L^2(\Hcal_s)}\leq CC_1\eps s^{\delta},
\ee
\bel{proof 2order eq 11a'}
\|s^{3}t^{-2}\del_t\del_tZ^Iu_{\ih}\|_{L^2(\Hcal_s)}\leq CC_1\eps s^{\delta/2},
\ee
\bel{proof 2order eq 11b}
\|s^{3}t^{-2}\del_t\del_tZ^{\If}u_{\ih}\|_{L^2(\Hcal_s)}\leq CC_1\eps
\ee
\end{subequations}
and, furthermore,
\begin{subequations}\label{proof 2order general L2}
\bel{proof 2order general L2a}
\|s^{3}t^{-2}\del_\alpha\del_\beta Z^{\Id}u_{\ih}\|_{L^2(\Hcal_s)}\leq CC_1\eps s^{\delta},
\ee
\bel{proof 2order general L2a'}
\|s^{3}t^{-2}\del_\alpha\del_\beta Z^Iu_{\ih}\|_{L^2(\Hcal_s)}\leq CC_1\eps s^{\delta/2},
\ee
\bel{proof 2order general L2b}
\|s^{3}t^{-2}\del_\alpha\del_\beta Z^{\If}u_{\ih}\|_{L^2(\Hcal_s)}\leq CC_1\eps.
\ee
\end{subequations}
\end{proposition}

\begin{proof}
The proof is a combination of \eqref{proof 2order eq9}, \eqref{proof 2order eq9-DEUX},
 and \eqref{proof 2order eq10} with \eqref{proof 2order eq3}.
We will prove first \eqref{proof 2order eq 11b}.

The proof is done by induction. We first prove \eqref{proof 2order eq 11b} with $|\If|=0$ and 
recall our convention that $|I|<0$ implies $Z^I=0$, 
and \eqref{proof 2order eq3} implies (with $|I|=0$):
$$
\aligned
& \big{\|}s^{3}t^{-2}\del_t\del_t u_{\ih}\big{\|}_{L^2(\Hcal_s)}
 \\
&\leq C\big{\|}s{Q_G}_{\ih}\big{\|}_{L^2(\Hcal_s)}
 + C\big{\|}s{Q_T}_{\ih}\big{\|}_{L^2(\Hcal_s)}
 + C\big{\|}sF_{\ih}\big{\|}_{L^2(\Hcal_s)}
 + C\big{\|}sR(u_{\ih})\big{\|}_{L^2(\Hcal_s)}.
\endaligned
$$
By the group of inequalities \eqref{proof 2order eq10},
$$
\big{\|}s^3t^{-2}\del_t\del_t u_{\ih}\big{\|}_{L^2(\Hcal_s)} \leq CC_1\eps.
$$

Suppose that \eqref{proof 2order eq 11b} holds with $|\If|\leq m$, we will prove \eqref{proof 2order eq 11b} with $|\If|\leq m+1$. By \eqref{proof 2order eq3},
$$
\aligned
&\big{\|}s^3t^{-2}\del_t\del_tZ^{\If}u_{\ih}\big{\|}_{L^2(\Hcal_s)}
\\
&\leq  CK\sum_{|I_2|+|I_3|\leq |\If| \atop  |I_2|<|\If|}
\sum_{\gamma,\jh,i}
\big{\|}s(|Z^{I_3}\del_{\gamma}w_i| + |Z^{I_3}w_i|) \,  \del_t\del_tZ^{I_2}u_{\jh}\big{\|}_{L^2(\Hcal_s)}
\\
&\quad 
+C\|s{Q_G}_{\ih}(\If,w,\del w,\del\del w)_{\ih}\|_{L^2(\Hcal_s)}
+ C\|s{Q_T}_{\ih}(\If,w,\del w,\del\del w)\|_{\ih}\|_{L^2(\Hcal_s)}
\\
&
\quad + C\|sR(Z^{\If}u_{\ih})\|_{L^2(\Hcal_s)}
\endaligned
$$
when $|\If|\leq m+1\leq 2$, and we  apply \eqref{proof 2order eq10}:
$$
\aligned
& \big{\|}s^3t^{-2}\del_t\del_tZ^{\If}u_{\ih}\big{\|}_{L^2(\Hcal_s)}
\\
& \leq CK\sum_{|I_2|+|I_3|\leq |\If| \atop  |I_2|<|\If|}
\sum_{\gamma,\jh,i}
\big{\|}s(|Z^{I_3}\del_{\gamma}w_i| + |Z^{I_3}w_i|) \,  \del_t\del_tZ^{I_2}u_{\jh}\big{\|}_{L^2(\Hcal_s)}
+CC_1\eps
\\
& \leq CK\sum_{|I_2|+|I_3|\leq |\If| \atop  |I_3|=1}
\sum_{\gamma,\jh,i}
\big{\|}s(|Z^{I_3}\del_{\gamma}w_i| + |Z^{I_3}w_i|) \,  \del_t\del_tZ^{I_2}u_{\jh}\big{\|}_{L^2(\Hcal_s)}
\\
 & \quad +CK\sum_{|I_2|+|I_3|\leq |\If| \atop  |I_3|=2}
\sum_{\gamma,\jh,i}
\big{\|}s(|Z^{I_3}\del_{\gamma}w_i| + |Z^{I_3}w_i|) \,  \del_t\del_tZ^{I_2}u_{\jh}\big{\|}_{L^2(\Hcal_s)}
+CC_1\eps.
\endaligned
$$
Observe now that when $I_3=1$, by \eqref{proof decay 2ge3a}, \eqref{proof decay 2ge1b} and \eqref{proof decay 2ge5a} we have 
$$
\aligned
&|Z^{I_3}\del_{\gamma}u_{\ih}|\leq CC_1\eps t^{-1/2}s^{-1},
\quad
|Z^{I_3}\del_{\gamma}v_{\ic}|\leq CC_1\eps t^{-3/2}s^{\delta},
\\
&|Z^{I_3}u_{\ih}|\leq CC_1\eps t^{-3/2}s,
\quad
|Z^{I_3}v_{\jc}|\leq CC_1\eps t^{-3/2}s^{\delta}.
\endaligned
$$
And by the induction assumption \big(\eqref{proof 2order eq 11b} for $|I_2|\leq |\If|-1$\big):
$$
\aligned
&\sum_{|I_2|+|I_3|\leq |\If| \atop  |I_3|=1}
\sum_{\gamma,\jh,i}
\big{\|}s(|Z^{I_3}\del_{\gamma}w_i| + |Z^{I_3}w_i|) \,  \del_t\del_tZ^{I_2}u_{\jh}\big{\|}_{L^2(\Hcal_s)}
\\
&\leq C(C_1\eps)\sum_{|I_2|+|I_3|\leq |\If|\atop \gamma,\jh,i,|I_3|=1}
\big{\|}t^2s^{-2}(t^{-1/2}s^{-1}+t^{-3/2}s+t^{-3/2}s^{\delta}) \,  s^3t^{-2}\del_t\del_tZ^{I_2}u_{\jh}\big{\|}_{L^2(\Hcal_s)}
\\
&\leq C(C_1\eps)\sum_{|I_2|+|I_3|\leq |\If| \atop  |I_3|=1}
\sum_{\gamma,\jh,i}\big{\|}s^3t^{-2}\del_t\del_tZ^{I_2}u_{\jh}\big{\|}_{L^2(\Hcal_s)}
\\
&\leq C(C_1\eps)^2.
\endaligned
$$
When $|I_3|=2$, we observe that $|I_2|\leq 0$. 
By \eqref{proof 2order eq8b}, we have 
$$
\aligned
&\sum_{\gamma,\jh,i}\big{\|}s(|Z^{I_3}\del_{\gamma}w_i| + |Z^{I_3}w_i|) \,  \del_t\del_tZ^{I_2}u_{\jh}\big{\|}_{L^2(\Hcal_s)}
\\
&\leq CC_1\eps\sum_{\gamma,\jh,i}\big{\|}s(|Z^{I_3}\del_{\gamma}w_i| + |Z^{I_3}w_i|) \,  t^{1/2}s^{-3}\big{\|}_{L^2(\Hcal_s)}
\\
&\leq CC_1\eps\sum_{\gamma,\ih}\big{\|}t^{1/2}s^{-2}Z^{I_3}\del_{\gamma}u_{\ih}\big{\|}_{L^2(\Hcal_s)}
     +CC_1\eps\sum_{\gamma,\ic}\big{\|}t^{1/2}s^{-2}Z^{I_3}\del_{\gamma}v_{\ic}\big{\|}_{L^2(\Hcal_s)}
\\
            &\quad 
+CC_1\eps\sum_{\ih}\big{\|}t^{1/2}s^{-2}Z^{I_3}u_{\ih}\big{\|}_{L^2(\Hcal_s)}
             +CC_1\eps\sum_{\ic}\big{\|}t^{1/2}s^{-2}Z^{I_3}v_{\ic}\big{\|}_{L^2(\Hcal_s)}
\endaligned
$$
These four terms can be bounded by $C(C_1\eps)^2$ by applying \eqref{proof L2 2ge3a}, \eqref{proof L2 2ge2a} and \eqref{proof L2 2ge5b}. So for $|\If|\leq m+1\leq 2$, \eqref{proof 2order eq 11b} is proved. 
By induction, \eqref{proof 2order eq 11b} is proved for $|\If|\leq 2$.

We turn to the proof of \eqref{proof 2order eq 11a'} and observe that in \eqref{proof 2order eq 11a'}, the case $|I|\leq 2$ is already proved by \eqref{proof 2order eq 11b}. We need only treat the cases $|I|=3$.

When $|I|=3$, we apply \eqref{proof 2order eq3}:
$$
\aligned
& \big{\|}s^3t^{-2}\del_t\del_tZ^Iu_{\ih}\big{\|}_{L^2(\Hcal_s)}
\\
&\leq CK\sum_{|I_2|+|I_3|\leq 3 \atop  |I_2|<3}
\sum_{\gamma,\jh,i}
\big{\|}s(|Z^{I_3}\del_{\gamma}w_i| + |Z^{I_3}w_i|) \, \del_t\del_tZ^{I_2}u_{\jh}\big{\|}_{L^2(\Hcal_s)}
\\
&\quad+C\|s{Q_G}_{\ih}(I,w,\del w,\del\del w)_{\ih}\|_{L^2(\Hcal_s)} + C\|Z^{I}F_{\ih}\| + C\|sZ^{I}u_{\ih}\|_{L^2(\Hcal_s)}.
\endaligned
$$
We observe that, in view of \eqref{proof 2order eq9-DEUX}, 
$$
\aligned
& \big{\|}s^3t^{-2}\del_t\del_tZ^Iu_{\ih}\big{\|}_{L^2(\Hcal_s)}
\\
&\leq  CK\sum_{|I_2|+|I_3|\leq 3 \atop  |I_2|<3}
\sum_{\gamma,\jh,i}
\big{\|}s(|Z^{I_3}\del_{\gamma}w_i| + |Z^{I_3}w_i|) \, \del_t\del_tZ^{I_2}u_{\jh}\big{\|}_{L^2(\Hcal_s)}
 +CC_1\eps s^{\delta/2}. 
\endaligned
$$

We focus on
$$
\aligned
&\sum_{|I_2|+|I_3|\leq 3 \atop  |I_2|<3}
\sum_{\gamma,\jh,i}
\big{\|}s(|Z^{I_3}\del_{\gamma}w_i| + |Z^{I_3}w_i|) \,  \del_t\del_tZ^{I_2}u_{\jh}\big{\|}_{L^2(\Hcal_s)}
\\
&\leq\sum_{|I_2|+|I_3|\leq 3 \atop  |I_2|=0}
\sum_{\gamma,\jh,\ih}
\big{\|}sZ^{I_3}\del_{\gamma}u_{\ih}\, \del_t\del_tZ^{I_2}u_{\jh}\big{\|}_{L^2(\Hcal_s)}
\\
& \quad
    +\sum_{|I_2|+|I_3|\leq 3 \atop  |I_2|=1}
\sum_{\gamma,\jh,\ic}
\big{\|}sZ^{I_3}\del_{\gamma}v_{\ic}\,\del_t\del_tZ^{I_2}u_{\jh}\big{\|}_{L^2(\Hcal_s)}
\\
    &\quad +\sum_{|I_2|+|I_3|\leq 3 \atop  |I_2|=2}
\sum_{\gamma,\jh,\ih}
\big{\|}sZ^{I_3}u_{\ih} \,  \del_t\del_tZ^{I_2}u_{\jh}\big{\|}_{L^2(\Hcal_s)}.
\endaligned
$$
For each sum, the possible choice of $(|I_2|,|I_3|)$ are
$$
(0,\leq 3),\quad (1,\leq 2),\quad (2,\leq 1).
$$
We list out, for each sum, the relevant inequalities for each possible choice $(|I_2|,|I_3|)$. This makes the following list:
$$
\begin{array}{cccc}
\text{Products} &(0,\leq 3) &(1,\leq 2) &(2,\leq 1)
\\
s\del_t\del_tZ^{I_2}u_{\jh}\, Z^{I_3}\del_{\gamma}u_{\ih}
&\eqref{proof decay 2order b},\eqref{proof L2 2ge3a}
&\eqref{proof 2order eq 11b},\eqref{proof decay 2ge2a'}
&\eqref{proof 2order eq 11b},\eqref{proof decay 2ge3a}
\\
s\del_t\del_tZ^{I_2}u_{\jh}\, Z^{I_3}\del_{\gamma}v_{\ic}
&\eqref{proof decay 2order b},\eqref{proof L2 2ge2a'}
&\eqref{proof 2order eq 11b},\eqref{proof decay 2ge2a}
&\eqref{proof 2order eq 11b},\eqref{proof decay 2ge2a}
\\
s\del_t\del_tZ^{I_2}u_{\jh} \,  Z^{I_3}u_{\ih}
&\eqref{proof decay 2order b},\eqref{proof L2 2ge5a'}
&\eqref{proof 2order eq 11b},\eqref{proof decay 2ge5a'}
&\eqref{proof 2order eq 11b},\eqref{proof decay 2ge5b}
\\
s\del_t\del_tZ^{I_2}u_{\jh} \,  Z^{I_3}v_{\ic}
&\eqref{proof decay 2order b},\eqref{proof L2 2ge2a'}
&\eqref{proof 2order eq 11b},\eqref{proof decay 2ge2a}
&\eqref{proof 2order eq 11b},\eqref{proof decay 2ge2a}
\end{array}
$$
We conclude with the fact that these four terms can be bounded by $C(C_1\eps)^2 s^{\delta/2}$. So \eqref{proof 2order eq 11a'} is proved for $|\Id|=3$.

Now we turns to the proof of \eqref{proof 2order eq 11}. When $|\Id|=4$, we apply \eqref{proof 2order eq 11a'} with $|\Id|=3$ and, more precisely,
\bel{proof 2order eq 12}
\|s^{3}t^{-2}\del_t\del_tZ^{I_1}u_{\ih}\|_{L^2(\Hcal_s)}\leq CC_1\eps s^{\delta/2}\leq CC_1\eps s^{\delta}. 
\ee
As in previous cases, by \eqref{proof 2order eq3} and \eqref{proof 2order eq9},
$$
\aligned
& \big{\|}s^3t^{-2}\del_t\del_tZ^{\Is}u_{\ih}\big{\|}_{L^2(\Hcal_s)}
\\
&\leq  CK\sum_{|I_2|+|I_3|\leq 4 \atop  |I_2|<3}
\sum_{\gamma,\jh,i}
\big{\|}s(|Z^{I_3}\del_{\gamma}w_i| + |Z^{I_3}w_i|) \, \del_t\del_tZ^{I_2}u_{\jh}\big{\|}_{L^2(\Hcal_s)}
 +C(C_1\eps)^2 s^{\delta}. 
\endaligned
$$
The first sum is also decomposed into four parts:
$$
\aligned
&\sum_{|I_2|+|I_3|\leq 4 \atop  |I_2|<4}
\sum_{\gamma,\jh,i}
\big{\|}s(|Z^{I_3}\del_{\gamma}w_i| + |Z^{I_3}w_i|) \,  \del_t\del_tZ^{I_2}u_{\jh}\big{\|}_{L^2(\Hcal_s)}
\\
&\leq\sum_{|I_2|+|I_3|\leq 4 \atop  |I_2|<4}
\sum_{\gamma,\jh,\ih}
\big{\|}sZ^{I_3}\del_{\gamma}u_{\ih}\, \del_t\del_tZ^{I_2}u_{\jh}\big{\|}_{L^2(\Hcal_s)}
 \\
&\quad
   +\sum_{|I_2|+|I_3|\leq 4 \atop  |I_2|<4}
\sum_{\gamma,\jh,\ic}
\big{\|}sZ^{I_3}\del_{\gamma}v_{\ic}\,\del_t\del_tZ^{I_2}u_{\jh}\big{\|}_{L^2(\Hcal_s)}
\\
    &\quad 
+\sum_{|I_2|+|I_3|\leq 4 \atop  |I_2|<4}
\sum_{\gamma,\jh,\ih}
\big{\|}sZ^{I_3}u_{\ih} \,  \del_t\del_tZ^{I_2}u_{\jh}\big{\|}_{L^2(\Hcal_s)}
  \\
& \quad   +\sum_{|I_2|+|I_3|\leq 4 \atop  |I_2|<4}
\sum_{\gamma,\jh,\ic}
\big{\|}sZ^{I_3}v_{\ic} \,  \del_t\del_tZ^{I_2}u_{\jh}\big{\|}_{L^2(\Hcal_s)}.
\endaligned
$$
The possible choices of $(|I_2|,|I_3|)$ are
$$
(3,\leq 1),\quad (2,\leq 2),\quad (1,\leq 3),\quad(0,\leq 4). 
$$
We list out the inequalities to be used for each term and each choice of $(|I_2|,|I_3|)$. Cf. Table 13. 
We conclude with \eqref{proof 2order eq 11a}.
On the other hand, \eqref{proof 2order general L2} are combinations of \eqref{proof 2order eq 11} and 
Lemma~\ref{proof 2order general L2 1}.
\end{proof}


\protect\begin{landscape}\thispagestyle{empty}
$$
\begin{array}{ccccc}
\text{Products} &(3,\leq 1) &(2,\leq 2) &(1,\leq 3) &(0,\leq 4)
\\
s\del_t\del_tZ^{I_2}u_{\jh}\, Z^{I_3}\del_{\gamma}u_{\ih}
&\eqref{proof 2order eq 12},\eqref{proof decay 2ge3a}
&\eqref{proof 2order eq8a},\eqref{proof L2 2ge3a}
&\eqref{proof 2order eq8a},\eqref{proof L2 2ge3a}
&\eqref{proof 2order eq8b},\eqref{proof L2 2ge1a}
\\
s\del_t\del_tZ^{I_2}u_{\jh}\, Z^{I_3}\del_{\gamma}v_{\ic}
&\eqref{proof 2order eq 12},\eqref{proof decay 2ge2a}
&\eqref{proof 2order eq8a},\eqref{proof L2 2ge2a}
&\eqref{proof 2order eq8a},\eqref{proof L2 2ge2a}
&\eqref{proof 2order eq8b},\eqref{proof L2 2ge2a}
\\
s\del_t\del_tZ^{I_2}u_{\jh} Z^{I_3}u_{\ih}
&\eqref{proof 2order eq 12},\eqref{proof decay 2ge5b}
&\eqref{proof 2order eq8a},\eqref{proof L2 2ge5b}
&\eqref{proof 2order eq8a},\eqref{proof L2 2ge5b}
&\eqref{proof 2order eq8b},\eqref{proof L2 2ge5a}
\\
s\del_t\del_tZ^{I_2}u_{\jh} Z^{I_3}v_{\ic}
&\eqref{proof 2order eq 12},\eqref{proof decay 2ge1c}
&\eqref{proof 2order eq8a},\eqref{proof L2 2ge1b}
&\eqref{proof 2order eq8a},\eqref{proof L2 2ge1b}
&\eqref{proof 2order eq8b},\eqref{proof L2 2ge1b}
\end{array}
$$
\centerline{Table 13}

\

\

$$
\begin{array}{ccccc}
\text{Terms} &(3,\leq 0) &(2,\leq 1) &(1,\leq 2) &(0,\leq 3)
\\
(s/t)^2Z^{I_1}\del_tu_{\ih}Z^{I_2}\del_tu_{\jh}
&\eqref{proof L2 2ge3a},\eqref{proof decay 2ge3a}
&\eqref{proof L2 2ge3a},\eqref{proof decay 2ge3a}
&\eqref{proof decay 2ge3a},\eqref{proof L2 2ge3a}
&\eqref{proof decay 2ge3a},\eqref{proof L2 2ge3a}
\\
Z^{I_1}\delu_au_{\ih}Z^{I_2}\delu_{\beta}u_{\jh}
&\eqref{proof L2 2ge3b},\eqref{proof decay 2ge3a}
&\eqref{proof L2 2ge3b},\eqref{proof decay 2ge3a}
&\eqref{proof decay 2ge3b},\eqref{proof L2 2ge3a}
&\eqref{proof decay 2ge3b},\eqref{proof L2 2ge3a}
\\
Z^{I_1}\delu_{\beta}u_{\ih}Z^{I_2}\delu_a u_{\jh}
&\eqref{proof L2 2ge3a},\eqref{proof decay 2ge3b}
&\eqref{proof L2 2ge3a},\eqref{proof decay 2ge3b}
&\eqref{proof decay 2ge3a},\eqref{proof L2 2ge3b}
&\eqref{proof decay 2ge3a},\eqref{proof L2 2ge3b}
\end{array}
$$
\centerline{Table 14}
\end{landscape}


Now, we derive the complete $L^2$ estimates of the second-order derivatives.

\begin{lemma}
Under the assumption \eqref{proof energy assumption}, the following $L^2$ estimates hold for $|\Id|\leq 4, |I|\leq 3|$ and $\If|\leq 2$:
\begin{subequations}\label{proof L2 2order}
\bel{proof L2 2order a}
\aligned
\big{\|}s^3t^{-2}\del_\alpha \del_{\beta}Z^{\Id} u_{\ih}\big{\|}_{L^2(\Hcal_s)}
+ \big{\|}s^3t^{-2}Z^{\Id}\del_\alpha \del_{\beta} u_{\ih}\big{\|}_{L^2(\Hcal_s)}
\leq C C_1\eps s^{\delta},
\endaligned
\ee
\bel{proof L2 2order a'}
\aligned
& \big{\|}s^3t^{-2}\del_\alpha \del_{\beta}Z^I u_{\ih}\big{\|}_{L^2(\Hcal_s)}
+ \big{\|}s^3t^{-2}Z^I \del_\alpha \del_{\beta} u_{\ih}\big{\|}_{L^2(\Hcal_s)}
\leq C C_1\eps s^{\delta/2},
\endaligned
\ee
\bel{proof L2 2order b}
\aligned
& \big{\|}s^3t^{-2}\del_\alpha \del_{\beta}Z^{\If} u_{\ih}\big{\|}_{L^2(\Hcal_s)}
+\big{\|}s^3t^{-2}Z^{\If}\del_\alpha \del_{\beta} u_{\ih}\big{\|}_{L^2(\Hcal_s)}
\leq C C_1\eps.
\endaligned
\ee
\end{subequations}
\end{lemma}
 
\begin{proof}
These inequalities are based on \eqref{pre lem commutator second-order} and \eqref{proof 2order general L2}.
\end{proof}

\chapter[Null forms and decay in time]{Null forms and decay in time\label{cha:8}}

\section{Bounds that are independent of second-order estimates}
\label{sec:81}
 
In this chapter, we are going to estimate the terms 
$$
\mathcal{T}(\del u_{\ih},\del u_{\jh}) := T^{\alpha\beta}\del_{\alpha}u_{\ih}\del_{\beta}u_{\jh},\quad
[Z^I,A^{\alpha\beta\gamma}\del_{\gamma}u_{\ih}\del_{\alpha}\del_{\beta}]u_{\jh}.
$$
where $T^{\alpha\beta}$ and $A^{\alpha\beta\gamma}$ are null quadratic forms. 
Concerning the term
$[Z^I,A^{\alpha\beta\gamma}\del_{\gamma}u_{\ih}\del_{\alpha}\del_{\beta}]u_{\jh}$,
we can apply \eqref{proof 2order eq 11b} and \eqref{proof 2order eq8b} for better decay rates, but the following rates will be sufficient for our main result in this monograph.

\begin{lemma}\label{proof null lem1}
By reyling on the energy assumption \eqref{proof energy assumption e}, 
the following estimates hold for all $|I|\leq 3$ and $|J|\leq 1$:
\begin{subequations}\label{proof null eq1}
\bel{proof null eq1a}
\big{\|}Z^I\mathcal{T}(\del u_{\ih},\del u_{\jh})\big{\|}_{L^2(\Hcal_s)}\leq C(C_1\eps)^2s^{-3/2},
\ee
\bel{proof null eq1c}
\big{\|}[Z^I,A^{\alpha\beta\gamma}\del_{\gamma}u_{\ih}\del_{\alpha}\del_{\beta}]u_{\jh}\big{\|}_{L^2(\Hcal_s)}\leq C(C_1\eps)^2s^{-3/2+\delta},
\ee
\end{subequations}
\be
\sup_{\Hcal_s}\big(t^2s \mathcal{T}(\del u_{\ih},\del u_{\jh})\big)\leq C(C_1\eps)^2.
\ee
\end{lemma}

\begin{proof}
The proof of \eqref{proof null eq1a} is a combination of Proposition \ref{pre lem null 1} with \eqref{proof decay 2ge3} and \eqref{proof L2 2ge3}. By recalling the decomposition of $\mathcal{T}$ presented in Proposition \ref{pre lem null 1},
we have the list in Table 14. 
This completes the argument.

 
The estimate on the term $[Z^I,A^{\alpha\beta\gamma}\del_{\gamma}u_{\ih}\del_{\alpha}\del_\beta]u_{\jh}$ is also proved by using the inequalities presented in the following list. Observe that by Proposition \ref{pre lem null 2order3}, some partition of the index do not exist. We have: 
$$
\begin{array}{cccc}
\text{Terms} &(3,\leq 0) &(2,\leq 1) &(1,\leq 2)
\\
(s/t)^2Z^{I_1}\del_tu_{\kh}Z^{I_2}\del_t\del_tu_{\jh}
&\eqref{proof L2 2ge3a},\eqref{proof decay 2ge3a}
&\eqref{proof L2 2ge3a},\eqref{proof decay 2ge1a}
&\eqref{proof decay 2ge3a},\eqref{proof L2 2ge3a}
\\
Z^{I_1}\delu_au_{\kh}Z^{I_2}\delu_{\beta}\delu_{\gamma}u_{\jh}
&\eqref{proof L2 2ge3b},\eqref{proof decay 2ge3a}
&\eqref{proof L2 2ge3b},\eqref{proof decay 2ge1a}
&\eqref{proof decay 2ge3b},\eqref{proof L2 2ge3a}
\\
Z^{I_1}\delu_{\alpha}u_{\kh}Z^{I_2}\delu_b\delu_{\gamma}u_{\jh}
&\eqref{proof L2 2ge3a},\eqref{proof decay 2ge4b}
&\eqref{proof L2 2ge3a},\eqref{proof decay 2ge4a}
&\eqref{proof decay 2ge3a},\eqref{proof L2 2ge4'b}
\\
Z^{I_1}\delu_{\alpha}u_{\kh}Z^{I_2}\delu_{\beta}\delu_c u_{\jh}
&\eqref{proof L2 2ge3a},\eqref{proof decay 2ge4b}
&\eqref{proof L2 2ge3a},\eqref{proof decay 2ge4a}
&\eqref{proof decay 2ge3a},\eqref{proof L2 2ge4'b}
\\
t^{-1}Z^{I_1}\del_{\alpha}u_{\kh}\delu_{\beta}u_{\jh}
&\eqref{proof L2 2ge3a},\eqref{proof decay 2ge3a}
&\eqref{proof L2 2ge3a},\eqref{proof decay 2ge3a}
&\eqref{proof decay 2ge3a},\eqref{proof L2 2ge3a}
\end{array}
$$
Finally, we derive the $L^\infty$ estimates of $\mathcal{T}$ by the following list: 
$$
\begin{array}{cccc}
\text{Terms} &(1,\leq 0) &(0,\leq1) &\text{Decay rate}
\\
(s/t)^2Z^{J_1}\del_tu_{\ih}Z^{J_2}\del_tu_{\jh}
&\eqref{proof decay 2ge3a},\eqref{proof decay 2ge3a}
&\eqref{proof decay 2ge3a},\eqref{proof decay 2ge3a}
&t^{-3}
\\
Z^{J_1}\delu_au_{\ih}Z^{J_2}\delu_{\beta}u_{\jh}
&\eqref{proof decay 2ge3b},\eqref{proof decay 2ge3a}
&\eqref{proof decay 2ge3b},\eqref{proof decay 2ge3a}
&t^{-2}s^{-1}
\\
Z^{J_1}\delu_{\beta}u_{\ih}Z^{J_2}\delu_a u_{\jh}
&\eqref{proof decay 2ge3a},\eqref{proof decay 2ge3b}
&\eqref{proof decay 2ge3a},\eqref{proof decay 2ge3b}
&t^{-2}s^{-1}
\end{array}
$$
\end{proof}


\section{Bounds that depend on second-order estimates}
\label{sec:82}

In this section we estimate the terms $[Z^I,B^{\alpha\beta}u_{\ih}\del_{\alpha}\del_{\beta}]u_{\jh}$
by 
essentially relying on the second-order estimates \eqref{proof decay 2order}.

\begin{lemma}\label{proof null lem2}
By relying on the energy assumption \eqref{proof energy assumption e},
 the following estimates hold for $|I|\leq 3$:
\begin{subequations}
\label{proof null eq2}
\bel{proof null eq2b3}
\big{\|}[Z^I,B^{\alpha\beta}u_{\ih}\del_{\alpha}\del_{\beta}]u_{\jh}\big{\|}_{L^2(\Hcal_s)}\leq C(C_1\eps)^2s^{-3/2+\delta}.
\ee
\end{subequations}
\end{lemma}

\begin{proof}
Recall the structure of
$[Z^I,B^{\alpha\beta}u_{\ih}\del_{\alpha}\del_{\beta}]u_{\jh}$ presented in \eqref{pre null B2}. It is a linear combinations of several terms and
we need to control each term in each possible partition of the index $|I_1|+|I_2|\leq |I|, |I_2|<|I|$.
As before, we list out the inequalities we use for each term and each partition of the indices:
$$
\begin{array}{cccc}
\text{terms} &(3,\leq 0) &(2,\leq 1) &(1,\leq 2)
\\
(s/t)^2Z^{I_1}u_{\kh}\,Z^{I_2}\del_t\del_tu_{\jh}
&\eqref{proof L2 2ge5b},\eqref{proof 2order eq8b}
&\eqref{proof L2 2ge5b},\eqref{proof 2order eq8a}
&\eqref{proof L2 2ge5b},\eqref{proof 2order eq8a}
\\
Z^{I_1}u_{\kh}Z^{I_2}\delu_a\delu_{\beta}u_{\jh}
&\eqref{proof L2 2ge5b},\eqref{proof decay 2ge4b}
&\eqref{proof L2 2ge5b},\eqref{proof decay 2ge4a}
&\eqref{proof L2 2ge5b},\eqref{proof decay 2ge4a}
\\
Z^{I_1}u_{\kh}Z^{I_2}\delu_{\alpha}\delu_bu_{\jh}
&\eqref{proof L2 2ge5b},\eqref{proof decay 2ge4b}
&\eqref{proof L2 2ge5b},\eqref{proof decay 2ge4a}
&\eqref{proof L2 2ge5b},\eqref{proof decay 2ge4a}
\\
t^{-1}Z^{I_1}u_{\kh}Z^{I_1}\del_{\alpha}u_{\jh}
&\eqref{proof L2 2ge5b},\eqref{proof decay 2ge3a}
&\eqref{proof L2 2ge5b},\eqref{proof decay 2ge3a}
&\eqref{proof decay 2ge5b},\eqref{proof L2 2ge3a}
\end{array}
$$
\end{proof}


\section{Decay estimates}
\label{sec:83}

We are now in a position to prove \eqref{main bootstrap 1}, \eqref{main bootstrap 2}, \eqref{main bootstrap 2'},  and \eqref{main bootstrap 4}.
The proof of \eqref{main bootstrap 1} requires only the $L^{\infty}$ estimate established in Chapter~\ref{cha:6}.

\begin{lemma}\label{proof lem curved energy is big}
Let $\{w_i\}$ be the local-in-time solution of \eqref{main eq main} and suppose that the energy assumption \eqref{proof energy assumption a} and \eqref{proof energy assumption e} hold. 
Then there exists a constant ${\Kcoef}>0$ 
(depending upon  
the constants 
$A_i^{j\alpha\beta\gamma k}$ and $B_i^{j\alpha\beta k}$) 
such that, if $C_1\eps$ is sufficiently small,  
$$
\Kcoef^{-2}\sum_iE_{m,c_i}(s,Z^I w_i) \leq \sum_iE_{G,c_i}(s,Z^I w_i) \leq \Kcoef^2 \sum_iE_{m,c_i}(s,Z^I w_i).
$$
\end{lemma}

\begin{proof} Note that
$$
\aligned
\sum_{i,j,\alpha,\beta} \big|G_i^{j\alpha\beta}\big| &\leq CK\sum_i\big(|\del w_i| + |w_i|\big)
\\
&\leq CK\sum_{\ih,\ic,\alpha}\big(|\del_{\alpha} v_{\ic}| + |\del_{\beta} u_{\ih}| + |v_{\ic}| + |u_{\ih}|\big). 
\endaligned
$$
Applying \eqref{proof decay 2ge2a}, \eqref{proof decay 2ge3a}, \eqref{proof decay 2ge1c}, 
 and \eqref{proof decay 2ge5b}, and recalling $0<\delta<1/6$, we have 
$$
\sum_{i,j,\alpha,\beta} \big|G_i^{j\alpha\beta}\big|\leq CK(t^{-3/2}s + t^{-1/2}s^{-1}).
$$

$$
\aligned
&\sum_i\big|E_{G,c_i}(s,w_i) - E_{m,c_i}(s,w_i) \big|
\\
&=  \bigg|2\int_{\Hcal_s} \big(\del_t w_i \del_{\beta}w_j G_i^{j\alpha\beta}\big)_{0\leq \alpha\leq 1}\, (1,-x/t) dx
- \int_{\Hcal_s} \big(\del_\alpha w_i\del_{\beta}w_j G_i^{j\alpha\beta}\big)dx\bigg|
\\
&\leq 2\int_{\Hcal_s}\bigg(\sum_{i,j,\alpha,\beta}\big|G_i^{j\alpha\beta}\big|\bigg) \, \bigg(\sum_{\alpha,k}|\del_\alpha w_k|^2\bigg) dx
\\
&\leq 2CK\int_{\Hcal_s}\sum_i\big(|\del w_i| + |w_i|\big) \, \bigg(\sum_{\alpha,k}|\del_\alpha w_k|^2\bigg)dx
\\
&\leq 2CKC_1\eps \int_{\Hcal_s}\big(t^{-1/2}s^{-1}+t^{-3/2}s\big)(t/s)^2 \, \bigg(\sum_{\alpha,k}|(s/t)\del_\alpha w_k|^2\bigg)dx
\\
&\leq CKC_1\eps \sum_iE_m(s,w_i) \leq CKC_1\eps \sum_iE_{m,c_i}(s,w_i),
\endaligned
$$
where the relation $t\leq Cs^2$ in $\Kcal$ is taken into consideration. Here we take $CKC_1\eps \leq 2/3$ with $C$ a universal constant, then the lemma is proved by fixing $\Kcoef = \sqrt{3}$.
\end{proof}

The proof of \eqref{main bootstrap 2} will be related to the energy estimate \eqref{proof L2 1ge1a}.

\begin{lemma}\label{proof lem Mv 4}
\label{quasilinear lem energy curveterm is small}
By relying on \eqref{proof energy assumption b}, \eqref{proof energy assumption a},
 and \eqref{proof energy assumption e},
 for any $|\Is|\leq 5$ the following estimate holds:
\bel{main bootstrap 2 pr}
\aligned
&\bigg|\int_{\Hcal_s}\frac{s}{t}\bigg(\big(\del_{\alphar}G_i^{\jr\alphar\betar}\big)\del_t Z^{\Is} w_i \del_{\beta}Z^{\Is} w_{\jr}
- \frac{1}{2}\big(\del_t G_i^{\jr\alphar\betar}\big)\del_{\alphar}Z^{\Is} w_i \del_{\betar}Z^{\Is} w_{\jr}\bigg)dx\bigg|
\\
& \leq C(C_1\eps)^2s^{-1+\delta}E_m(s,Z^{\Is} w_i)^{1/2} \leq C(C_1\eps)^2s^{-1+\delta}E_{m,\sigma}(s,Z^{\Is} w_i)^{1/2}.
\endaligned
\ee
\end{lemma}

\begin{proof} By \eqref{proof decay 2ge1a} and \eqref{proof decay 2ge2a}, 
we note that
$$
\aligned
\big|\del_\alpha G_i^{j\alpha\beta}\big|\leq& C\sum_j|\del_\alpha w_j| + C\sum_{j,\beta}|\del_\alpha \del_{\beta}w_j|
\\
\leq&  C\sum_{\jc,\alpha,\beta}\big(|\del_\alpha v_{\jc}| + |\del_{\alpha}\del_{\beta}v_{\jc}|\big) + C\sum_{\jh,\alpha,\beta}\big(|\del_{\alpha}u_{\jh}| + |\del_\alpha \del_{\beta}u_{\jh}|\big)
\\
\leq& C C_1\eps t^{-3/2}s^{\delta} + C C_1\eps t^{-1/2}s^{-1+\delta}.
\endaligned
$$
By substituting this result into the expression, the first term in the left-hand side of \eqref{main bootstrap 2 pr} is bounded as follows:
$$
\aligned
&\big{\|}(s/t)\big(\del_{\alphar}G_i^{\jr\alphar\betar}\big)\del_t Z^{\Is} w_i \del_{\betar}Z^{\Is} w_{\jr}\big{\|}_{L^1(\Hcal_s)}
\\
&=\big{\|}\big((t/s)\del_{\alphar}G_i^{\jr\alphar\betar}\big)(s/t)\del_t Z^{\Is} w_i (s/t)\del_{\betar}Z^{\Is} w_{\jr}\big{\|}_{L^1(\Hcal_s)}
\\
&\leq \sum_{j,\beta}\big{\|}C(C_1\eps)(t^{1/2}s^{-2} + t^{-1/2}s^{-1+\delta})\,(s/t)\del_t Z^{\Is} w_i (s/t)\del_{\beta}Z^{\Is} w_j\big{\|}_{L^1(\Hcal_s)}
\\
&\leq CC_1\eps s^{-1}\sum_{j,\beta}\big{\|}(s/t)\del_t Z^{\Is} w_i (s/t)\del_{\beta}Z^{\Is} w_j\big{\|}_{L^1(\Hcal_s)}
\\
&\leq
CC_1\eps s^{-1}\sum_{j,\beta}\|(s/t)\del_{\beta}Z^{\Is} w_j\big\|_{L^2(\Hcal_s)}\,
\big\|(s/t)\del_t Z^{\Is} w_i \big\|_{L^2(\Hcal_s)}.
\endaligned
$$
We apply \eqref{proof L2 1ge1a} and find
$$
\aligned
& \big{\|}(s/t)\big(\del_{\alphar}G_i^{\jr\alphar\betar}\big)\del_t Z^{\Is} w_i \del_{\betar}Z^{\Is} w_{\jr}\big{\|}_{L^1(\Hcal_s)}
\\
&\leq C(C_1\eps)^2 s^{-1+\delta}\bigg(\int_{\Hcal_s}\big|(s/t)\del_t Z^{\Is} w_i \big|^2dx\bigg)^{1/2}
\\
&\leq C(C_1\eps)^2 s^{-1+\delta}E_m(s,Z^{\Is}w_i)^{1/2}.
\endaligned
$$
The second term is estimated in the same way and we omit the details.
\end{proof}

The proof of \eqref{main bootstrap 2'} is similar to that of \eqref{main bootstrap 2}.

\begin{lemma}
By relying on
 \eqref{proof energy assumption d}, \eqref{proof energy assumption c},
 and \eqref{proof energy assumption e}, the following estimate holds:
\begin{equation}
\aligned
&\bigg|\int_{\Hcal_s}\frac{s}{t}\bigg(\big(\del_{\alphar}G_i^{\jr\alphar\betar}\big)\del_t Z^{\Id} w_i \del_{\beta}Z^{\Id} w_{\jr}
- \frac{1}{2}\big(\del_t G_i^{\jr\alphar\betar}\big)\del_{\alphar}Z^{\Id} w_i \del_{\betar}Z^{\Id} w_{\jr}\bigg)dx\bigg|
\\
\leq& C(C_1\eps)^2s^{-1+\delta/2}E_m(s,Z^{\Id} w_i)^{1/2}\leq C(C_1\eps)^2s^{-1+\delta/2}E_{m,\sigma}(s,Z^{\Id} w_i)^{1/2}.
\endaligned
\ee
\end{lemma}

\begin{proof}
As in the proof of Lemma \ref{proof lem Mv 4},  
 we can apply  \eqref{proof L2 1ge1a'} together with \eqref{proof decay 2ge1a'} and \eqref{proof decay 2ge2a'}: 
$$
\aligned
& \big{\|}(s/t)\big(\del_{\alphar}G_i^{\jr\alphar\betar}\big)\del_t Z^{\Is} w_i \del_{\betar}Z^{\Is} w_{\jr}\big{\|}_{L^1(\Hcal_s)}
\\
&\leq C(C_1\eps)^2 s^{-1+\delta/2}\bigg(\int_{\Hcal_s}\big|(s/t)\del_t Z^{\Is} w_i \big|^2dx\bigg)^{1/2}
\\
&\leq C(C_1\eps)^2 s^{-1+\delta/2}E_m(s,Z^{\Is}w_i)^{1/2}.
\endaligned
$$
The second term is estimated in the same way and we omit the details.
\end{proof}

The proof of \eqref{main bootstrap 4} is quite similar, but the null structure will be taken into consideration for the sharp decay rate $s^{-3/2+2\delta}$. This is one of the only two places where the null structure is taken into account.

\begin{lemma}\label{proof lem Mv 3}
Suppose \eqref{proof energy assumption a} and \eqref{proof energy assumption e} hold,
then for any $|I|\leq 3$ the following estimates are true:
\begin{equation}
\aligned
&\bigg|\int_{\Hcal_s}\frac{s}{t}\bigg(\big(\del_{\alphar}G_{\ih}^{\jh\alphar\betar}\big)\del_t Z^I u_{\ih} \del_{\betar}Z^I u_{\jh}
- \frac{1}{2}\big(\del_t G_{\ih}^{\jh\alphar\betar}\big)\del_{\alphar}Z^I u_{\ih} \del_{\betar}Z^I u_{\jh}\bigg)dx\bigg|
\\
&\leq C(C_1\eps)^2s^{-3/2+\delta}E_m(s,Z^I u_{\ih})^{1/2}.
\endaligned
\ee
\end{lemma}

\begin{proof}
We decompose the term $G_{\ih}^{\jh\alpha\beta}$ as follows:
$$
\aligned
&(s/t)\del_{\alpha}G_{\ih}^{\jh\alpha\beta} \del_{\betar}Z^I u_{\jh}\del_tZ^Iu_{\ih}
\\
&= (s/t)\big(A_{\ih}^{\jh\alpha\beta\gamma\kc}\del_{\alpha}\del_{\gamma}v_{\kc} + B_{\ih}^{\jh\alpha\beta\kc}\del_{\alpha}v_{\kc}\big)\del_{\betar}Z^I u_{\jh}\del_tZ^Iu_{\ih}
\\
&\quad 
+ (s/t)\big(A_{\ih}^{\jh\alpha\beta\gamma\kh}\del_{\alpha}\del_{\gamma}u_{\kh}
+ B_{\ih}^{\jh\alpha\beta\kh}\del_{\alpha}u_{\kh}\big)\del_{\betar}Z^I u_{\jh}\del_tZ^Iu_{\ih}
\\
&:= R_1 + R_2. 
\endaligned
$$
The estimate of $R_1$ is direct, and we simply apply the inequalities \eqref{proof decay 2ge1c}: 
$$
\aligned
&\|R_1\|_{L^1(\Hcal_s)} = \|(s/t)\big(A_{\ih}^{\jh\alpha\beta\gamma\kc}\del_{\alpha}\del_{\gamma}v_{\kc} + B_{\ih}^{\jh\alpha\beta\kc}\del_{\alpha}v_{\kc}\big)\,(\del_{\betar}Z^I u_{\jh}\del_tZ^Iu_{\ih})\|_{L^1(\Hcal_s)}
\\
&\leq C C_1 \eps \, \|(t^{-3/2}s^{\delta}(t/s))\, (s/t)\del_{\betar}Z^I u_{\jh}\, (s/t)\del_tZ^Iu_{\ih}\|_{L^1(\Hcal_s)}
\\
&\leq C C_1 \eps \,  \|(t^{-1/2}s^{-1+\delta}(s/t)\del_{\betar}Z^I u_{\jh}\|_{L^2(\Hcal_s)} \,  \|(s/t)\del_tZ^Iu_{\ih}\|_{L^2(\Hcal_s)}
\\
&\leq C C_1\eps s^{3/2+\delta} 
\| (s/t) \del_{\betar}Z^I u_{\jh}\|_{L^2(\Hcal_s)} \,  E_m(s,Z^Iu_{\ih})^{1/2}.
\endaligned
$$
In view of \eqref{proof L2 1ge3a}, we get 
\bel{proof lem Mv 3 eq1}
\|R_1\|_{L^1(\Hcal_s)} \leq C(C_1\eps)^2 s^{-3/2+\delta} E_m(s,Z^Iu_{\ih})^{1/2}.
\ee

The estimates on $R_2$ is more involved: we have the following decomposition:
$$
\aligned
&(s/t)\big(A_{\ih}^{\jh\alpha\beta\gamma\kh}\del_{\alpha}\del_{\gamma}u_{\kh}
+ B_{\ih}^{\jh\alpha\beta\kh}\del_{\alpha}u_{\kh}\big)\del_{\betar}Z^I u_{\jh}
\\
= &(s/t)A_{\ih}^{\jh\alpha\beta\gamma\kh}\del_{\alpha}\del_{\gamma}u_{\kh}\del_{\beta}Z^I u_{\jh}
+ (s/t)B_{\ih}^{\jh\alpha\beta\kh}\del_{\alpha}u_{\kh}\del_{\betar}Z^I u_{\jh}
\\
=:&(s/t)\mathcal{A}_{\ih}^{\jh\kh}(\del Z^I u_{\jh},\del\del u_{\kh}) + (s/t)\mathcal{T}_{\ih}^{\jh\kh}(\del u_{\kh}, \del Z^I u_{\jh})
\endaligned
$$
Recall that $\mathcal{A}_{\ih}^{\jh\kh}$ and $\mathcal{T}_{\ih}^{\jh\kh}$ are null forms.

Let us consider first the terms of $\mathcal{A}_{\ih}^{\jh\kh}$. By Proposition \ref{pre lem null 2order2} (with $I=0$), we have
$$
\aligned
&\mathcal{A}_{\ih}^{\jh\kh}(\del Z^I u_{\jh},\del\del u_{\kh})
\\
&\leq CK(s/t)^2|\del_t\del_tu_{\kh}|\,|\del_tZ^Iu_{\jh}|
+CK\Omega_1(0,Z^I u_{\jh},u_{\kh}) + CKt^{-1}\Omega_2(0,Z^I u_{\jh},u_{\kh})
\\
&\leq CK(s/t)^2|\del_t\del_tu_{\kh}|\,|\del_tZ^Iu_{\jh}|
\\
&\quad +CK\sum_{a,\beta,\gamma}|\delu_a\delu_{\beta}u_{\kh}|\,|\delu_{\gamma}Z^Iu_{\jh}|
+ CK\sum_{\alpha,b,\gamma}|\delu_{\alpha}\delu_bu_{\kh}|\,|\delu_{\gamma}Z^Iu_{\jh}|
\\
&\quad 
+ CK\sum_{\alpha,\beta,c}|\delu_{\alpha}\delu_{\beta}u_{\kh}|\,|\delu_cZ^Iu_{\jh}|
+CKt^{-1}\sum_{\alpha,\beta}|\del_{\alpha}Z^Iu_{\jh}|\,|\del_{\beta}u_{\kh}|.
\endaligned
$$
Each term is estimated as follows:
$$
\aligned
& \|(s/t)^2 \del_t\del_tu_{\kh}\,\del_tZ^Iu_{\jh}\|_{L^2(\Hcal_s)} =
\|(s/t)\del_t\del_tu_{\kh} \, (s/t)\del_tZ^Iu_{\jh}\|_{L^2(\Hcal_s)}
\\
&\leq CC_1\eps\|t^{-3/2}\,(s/t)\del_tZ^Iu_{\jh}\|_{L^2(\Hcal_s)}
\\
&\leq CC_1\eps s^{-3/2}\|(s/t)\del_tZ^Iu_{\jh}\|_{L^2(\Hcal_s)}\leq C(C_1\eps)^2 s^{-3/2},
\endaligned
$$
where we applied \eqref{proof decay 2ge3a} and \eqref{proof L2 1ge3a}. We have also 
$$
\aligned
& \|\delu_a\delu_{\beta}u_{\kh}\,\delu_{\gamma}Z^Iu_{\jh}\|_{L^2(\Hcal_s)}
= \|(t/s)\delu_a\delu_{\beta}u_{\kh}\,(s/t)\del_tZ^Iu_{\jh}\|_{L^2(\Hcal_s)}
\\
&\leq CC_1\eps \|(t/s)t^{-3/2}s^{-1}\,(s/t)\del_tZ^Iu_{\jh}\|_{L^2(\Hcal_s)}
\\
&\leq CC_1\eps s^{-3/2} \|(s/t)\del_tZ^Iu_{\jh}\|_{L^2(\Hcal_s)}
\leq C(C_1\eps)^2 s^{-3/2}. 
\endaligned
$$ 

The term $\|\delu_{\alpha}\delu_bu_{\kh}\,\delu_{\gamma}Z^Iu_{\jh}\|_{L^2(\Hcal_s)}$ is estimated in the same manner and we omit the details.
$$
\aligned
\|\delu_{\alpha}\delu_{\beta}u_{\kh}\,\delu_cZ^Iu_{\jh}\|_{L^2(\Hcal_s)}
\leq \|\delu_{\alpha}\delu_{\beta}u_{\kh}\|_{L^\infty(\Hcal_s)}\|\delu_cZ^Iu_{\jh}\|_{L^2(\Hcal_s)}
\leq C(C_1\eps)^2 s^{-3/2}
\endaligned
$$
where \eqref{proof decay 2ge3a} and \eqref{proof L2 1ge3b} are used.
$$
\aligned
\|t^{-1}\del_{\alpha}Z^Iu_{\jh}\,\del_{\beta}u_{\kh}\|_{L^2(\Hcal_s)}
&\leq \|t^{-1}(t/s)\del_{\beta}u_{\kh}\|_{L^\infty(\Hcal_s)}\|(s/t)\del_{\alpha}Z^Iu_{\jh}\|_{L^2(\Hcal_s)}
\\
& \leq C(C_1\eps)^2s^{-5/2}. 
\endaligned
$$
 We conclude with
$$
\big{\|}(s/t)\mathcal{A}_{\ih}^{\jh\kh}(\del Z^Iu_{\ih},\del\del u_{\jh})\big{\|}_{L^2(\Hcal_s)}\leq CK(C_1\eps)^2s^{-3/2} 
$$
and 
\bel{proof lem Mv 3 eq2}
\aligned
& \big{\|}(s/t)A_{\ih}^{\jh\kh}(\del Z^Iu_{\ih},\del\del u_{\jh})\del_tZ^Iu_{\ih}\big{\|}_{L^1(\Hcal_s)}
\\
&\leq \big{\|}A_{\ih}^{\jh\kh}(\del Z^Iu_{\ih},\del\del u_{\jh})\|_{L^2(\Hcal)}\|(s/t)\del_tZ^Iu_{\ih}\|_{L^2(\Hcal_s)}
\\
&\leq C(C_1\eps)^2s^{-3/2}E_m(s,Z^Iu_{\ih})^{1/2}.
\endaligned
\ee

Next, we consider the terms of $\mathcal{T}_{\ih}^{\jh\kh}$. By \eqref{pre null B1}, $\mathcal{T}_{\ih}^{\jh\kh}(\del Z^Iu_{\kh}, \del u_{\jh})$ can be bounded by a linear combination of the terms presented in Proposition \ref{pre lem null 1}. All the estimates are based on the inequalities \eqref{proof decay 2ge3a}, \eqref{proof decay 2ge3b}, \eqref{proof L2 1ge1a}, and \eqref{proof L2 1ge1b}.
$$
\aligned
\big{\|}(s/t)^2\del_tZ^Iu_{\jh}\del_tu_{\kh}\big{\|}_{L^2(\Hcal_s)}
&\leq \big{\|}(s/t)\del_tu_{\kh}\big{\|}_{L^\infty(\Hcal_s)} \big{\|}(s/t)\del_tZ^Iu_{\jh}\big{\|}_{L^2(\Hcal_s)}
\\
& \leq C(C_1\eps)^2 s^{-3/2},
\endaligned
$$
$$
\aligned
\big{\|}\delu_au_{\kh}\delu_{\beta}Z^Iu_{\jh}\big{\|}_{L^2(\Hcal_s)}
& \leq \big{\|}(t/s)\delu_au_{\kh}\big{\|}_{L^\infty(\Hcal_s)} \big{\|}(s/t)\delu_{\beta}Z^Iu_{\jh}\big{\|}_{L^2(\Hcal_s)}
\\
&\leq C(C_1\eps)^2s^{-3/2},
\endaligned
$$
$$
\aligned
\big{\|}\delu_{\alpha}u_{\kh}\delu_bZ^Iu_{\jh}\big{\|}_{L^2(\Hcal_s)}
& \leq \big{\|}\delu_{\alpha}u_{\kh}\big{\|}_{L^\infty(\Hcal_s)} \big{\|}\delu_bZ^Iu_{\jh}\big{\|}_{L^2(\Hcal_s)}
\\
& \leq C(C_1\eps)^2s^{-3/2}.
\endaligned
$$
We conclude with
\bel{proof lem Mv 3 eq3}
\aligned
\big{\|}\mathcal{T}(\del u_{\kh},\del Z^I u_{\jh}) \del_tZ^Iu_{\ih} \big{\|}_{L^1(\Hcal_s)}
&\leq \big{\|}\mathcal{T}(\del u_{\kh},\del Z^I u_{\jh})\big{\|}_{L^2(\Hcal_s)}
\big{\|}(s/t)\del_tZ^Iu_{\ih}\big{\|}_{L^2(\Hcal_s)}
\\
&\leq C(C_1\eps)^2s^{-3/2}E_m(s,Z^Iu_{\ih})^{1/2}.
\endaligned
\ee

Combining \eqref{proof lem Mv 3 eq1}, \eqref{proof lem Mv 3 eq2}, and \eqref{proof lem Mv 3 eq3}, we conclude that
$$
\big{\|}\del_{\alpha}\big(G_{\ih}^{\jh\alpha\beta}(w,\del w)\big)\del_tZ^Iu_{\ih}\big{\|}_{L^1(\Hcal_s)}\leq C(C_1\eps)^2s^{-3/2+\delta}E_m(s,Z^Iu_{\ih})^{1/2}.
$$
By a similar argument, the term $\big{\|}\del_t\big(G_{\ih}^{\jh\alpha\beta}(w,\del w)\big)\del_{\alpha}Z^Iu_{\ih}\big{\|}_{L^1(\Hcal_s)}$ is also bounded by
 $C(C_1\eps)^2s^{-3/2+\delta}$, and we omit the details.
\end{proof}

\chapter[$L^2$ estimate on the interaction terms]{$L^2$ estimates of the interaction terms \label{cha:9}}

\section{$L^2$ estimates on higher-order interaction terms}
\label{sec:91}

In this chapter,  we establish three groups of energy estimates mentioned in the proof of Proposition \ref{main prop2}, that is, \eqref{main bootstrap 3}, \eqref{main bootstrap 3'}, and \eqref{main bootstrap 5} which are derived under the assumption \eqref{proof energy assumption}.
We emphasize that three groups of inequalities correspond to different decreasing rates in time. The proof of \eqref{main bootstrap 3} and \eqref{main bootstrap 3'} 
is much easier than that of \eqref{main bootstrap 5},
 since, roughly speaking, it does not require the null structure, while 
the decreasing rate in \eqref{main bootstrap 5} 
is a consequence of the null structure.
Interestingly, this is one of the two 
places in our proof where the null structure is used in a fundamental way.

\begin{lemma}\label{proof source lem higher}
Under the assumption of \eqref{proof energy assumption} (with $C_1\eps\leq \min(1,\eps_0'')$), the following estimates hold for all $|\Is|\leq 5$ and $s_0\leq s\leq s_1$: 
\begin{subequations}\label{proof source higher}
\bel{proof source higher F}
\big{\|}Z^{\Is} F_{\ih}\big{\|}_{L^2(\Hcal_s)} \leq C(C_1\eps)^2s^{-1+\delta},
\ee
\bel{proof source higher G}
\big{\|}[Z^I,G_i^{j\alpha\beta}(w,\del w)\del_{\alpha}\del_{\beta}]w_j\big{\|}\leq C(C_1\eps)^2s^{-1+\delta}.
\ee
\end{subequations}
\end{lemma}

\begin{proof}
We begin with \eqref{proof source higher F}. This concerns only the basic $L^2$ and $L^\infty$ estimates established
earlier.
Recall that $Z^{\Is}F_i$ is decomposed  as follows:
$$
\aligned
Z^{\Is}F_i &= Z^{\Is}\big(P_i^{\alpha\beta j k}\del_\alpha w_j\del_{\beta}w_k\big) + Z^{\Is}\big(Q_i^{\alpha j \kc}v_{\kc}\del_\alpha w_j\big) + Z^{\Is}\big(R_i^{\jc\kc}v_{\jc}v_{\kc}\big)
\endaligned
$$
We see that $Z^{\Is}F_i$ is a linear combination of the following terms:
$$
\aligned
&Z^{\Is}\big(\del_\alpha u_{\ih}\del_{\beta}u_{\jh}\big),\quad
Z^{\Is}\big(\del_\alpha v_{\ic}\del_{\beta}u_{\jh}\big),\quad
Z^{\Is}\big(\del_\alpha v_{\ic}\del_{\beta}v_{\jc}\big),\quad
\\
&Z^{\Is}\big(v_{\ic}\del_\alpha w_j\big), \quad
Z^{\Is}\big(v_{\jc}v_{\kc}\big).
\endaligned
$$
As done before, 
for each term and each partition of ${\Is = \Is_1 + \Is_2}$, we write the relevant 
inequalities in the following two lists (Recall our convention that, on each term, $\Is_1$ acts on the first factor and $\Is_2$ acts on the second, while 
the symbol $(a,\leq b)$ means $|\Is_1|=a, |\Is_2|\leq b$.):
$$
\begin{array}{cccc}
\text{Products} &(5,\leq 0) &(4,\leq 1) &(3,\leq 2)
\\
\del_{\alpha}u_{\ih}\del_{\beta}u_{\jh}
&\eqref{proof L2 2ge1a},\eqref{proof decay 2ge3a}
&\eqref{proof L2 2ge1a},\eqref{proof decay 2ge3a}
&\eqref{proof L2 2ge3a},\eqref{proof decay 2ge1a}
\\
\del_\alpha v_{\ic}\del_{\beta}u_{\jh}
&\eqref{proof L2 2ge1a},\eqref{proof decay 2ge3a}
&\eqref{proof L2 2ge2a},\eqref{proof decay 2ge3a}
&\eqref{proof L2 2ge2a},\eqref{proof decay 2ge1a}
\\
\del_\alpha v_{\ic}\del_{\beta}v_{\jc}
&\eqref{proof L2 2ge1a},\eqref{proof decay 2ge2a}
&\eqref{proof L2 2ge2a},\eqref{proof decay 2ge2a}
&\eqref{proof L2 2ge2a},\eqref{proof decay 2ge2a}
\\
v_{\ic}\del_\alpha w_j
&\eqref{proof L2 2ge1c},\eqref{proof decay 2ge1a}
&\eqref{proof L2 2ge1c},\eqref{proof decay 2ge1a}
&\eqref{proof L2 2ge1c},\eqref{proof decay 2ge1a}
\\
v_{\jc}v_{\kc}
&\eqref{proof L2 2ge1c},\eqref{proof decay 2ge1c}
&\eqref{proof L2 2ge1c},\eqref{proof decay 2ge1c}
&\eqref{proof L2 2ge1c},\eqref{proof decay 2ge1c}
\end{array}
$$

$$
\begin{array}{cccc}
\text{Products} &(2,\leq 3) &(1,\leq 4) &(0,\leq 5)
\\
\del_{\alpha}u_{\ih}\del_{\beta}u_{\jh}
&\eqref{proof decay 2ge1a},\eqref{proof L2 2ge3a}
&\eqref{proof decay 2ge3a},\eqref{proof L2 2ge1a}
&\eqref{proof decay 2ge3a},\eqref{proof L2 2ge1a}
\\
\del_\alpha v_{\ic}\del_{\beta}u_{\jh}
&\eqref{proof decay 2ge2a},\eqref{proof L2 2ge3a}
&\eqref{proof decay 2ge2a},\eqref{proof L2 2ge1a}
&\eqref{proof decay 2ge2a},\eqref{proof L2 2ge1a}
\\
\del_\alpha v_{\ic}\del_{\beta}v_{\jc}
&\eqref{proof decay 2ge2a},\eqref{proof L2 2ge1a}
&\eqref{proof decay 2ge2a},\eqref{proof L2 2ge1a}
&\eqref{proof decay 2ge2a},\eqref{proof L2 2ge1a}
\\
v_{\ic}\del_\alpha w_j
&\eqref{proof L2 2ge1c},\eqref{proof decay 2ge1a}
&\eqref{proof decay 2ge1c},\eqref{proof L2 2ge1a}
&\eqref{proof decay 2ge1c},\eqref{proof L2 2ge1a}
\\
v_{\jc}v_{\kc}
&\eqref{proof L2 2ge1c},\eqref{proof decay 2ge1c}
&\eqref{proof decay 2ge1c},\eqref{proof L2 2ge1c}
&\eqref{proof decay 2ge1c},\eqref{proof L2 2ge1c}
\end{array}
$$
From these inequalities, we conclude with \eqref{proof source higher F}.

We turn our attention to \eqref{proof source higher G} and, by 
recalling the decomposition of $[Z^{\Is},G_i^{j\alpha\beta}(w,\del w)\del_{\alpha}\del_{\beta}]w_j$, we find
\bel{proof source higher eq0}
\aligned
&\,[Z^{\Is},G_i^{j\alpha\beta}(w,\del w)\del_{\alpha}\del_{\beta}]w_j
\\
=& [Z^{\Is},A_i^{j\alpha\beta\gamma k}\del_{\gamma}w_k \del_{\alpha}\del_{\beta}]w_j
 + [Z^{\Is},B_i^{j\alpha\beta \kc}v_{\kc} \del_{\alpha}\del_{\beta}]w_j
\\
& \quad
 + [Z^{\Is},B_i^{j\alpha\beta \kh}u_{\kh}\del_{\alpha}\del_{\beta}]u_{\jh}. 
\endaligned
\ee
The first term of the right-hand side is decomposed as follows:
\bel{proof source higher eq1}
\aligned
& [Z^{\Is},A_i^{j\alpha\beta\gamma k}\del_{\gamma}w_k \del_{\alpha}\del_{\beta}]w_j
\\
& =\sum_{\Is_1+\Is_2=\Is\atop |\Is_2|\leq |\Is|-1}
A_i^{j\alpha\beta\gamma k}Z^{\Is_1}\del_{\gamma}w_kZ^{\Is_2}\del_\alpha \del_{\beta}w_j
+A_i^{j\alpha\beta\gamma k}\del_{\gamma}w_k [Z^{\Is},\del_{\alpha}\del_{\beta}]w_j, 
\endaligned
\ee
in which we have 
$$
\aligned
&\sum_{\Is_1+\Is_2=\Is\atop |\Is_2|\leq |\Is|-1}
A_i^{j\alpha\beta\gamma k}Z^{\Is_1}\del_{\gamma}w_kZ^{\Is_2}\del_\alpha \del_{\beta}w_j
\\
&=\sum_{\Is_1+\Is_2=\Is\atop |\Is_2|\leq |\Is|-1}\!\!\!\!
\bigg(
 A_i^{j\alpha\beta\gamma \kh}Z^{\Is_1}\del_{\gamma}u_{\kh}Z^{\Is_2}\del_\alpha \del_{\beta}w_j
+A_i^{j\alpha\beta\gamma \kc}Z^{\Is_1}\del_{\gamma}v_{\kc}Z^{\Is_2}\del_\alpha \del_{\beta}w_j
\bigg).
\endaligned
$$
This term is a linear combination of the following terms with constant coefficients:
$$
\aligned
&Z^{\Is_1}\del_{\gamma}u_{\kh}Z^{\Is_2}\del_{\alpha}\del_{\beta}u_{\jh},\quad
Z^{\Is_1}\del_{\gamma}u_{\kh}Z^{\Is_2}\del_{\alpha}\del_{\beta}v_{\jc}
\\
&Z^{\Is_1}\del_{\gamma}v_{\kc}Z^{\Is_2}\del_{\alpha}\del_{\beta}u_{\jh},\quad
Z^{\Is_1}\del_{\gamma}v_{\kc}Z^{\Is_2}\del_{\alpha}\del_{\beta}v_{\jc}
\endaligned
$$
with $|\Is_1| + |\Is_2|\leq |\Is|$ and $\Is_2\leq |\Is|-1$.

The second term in the right-hand side of \eqref{proof source higher eq1} is estimated as follows. We see
 that it is a linear combination of the following terms with constant coefficients:
$$
\aligned
&\del_{\gamma}u_{\kh}[Z^{\Is},\del_{\alpha}\del_{\beta}]u_{\jh},\quad
\del_{\gamma}u_{\kh}[Z^{\Is},\del_{\alpha}\del_{\beta}]v_{\jc},
\\
&\del_{\gamma}v_{\kc}[Z^{\Is},\del_{\alpha}\del_{\beta}]u_{\jh},\quad
\del_{\gamma}v_{\kc}[Z^{\Is},\del_{\alpha}\del_{\beta}]v_{\jc}
\endaligned
$$
By the commutator estimates \eqref{pre lem commutator second-order}, these terms are bounded by
$$
\aligned
& \sum_{|\Is_2|\leq |\Is|-1}\big|\del_{\gamma}u_{\kh}\del_{\alpha}\del_{\beta}Z^{\Is_2}u_{\jh}\big|,\quad
\sum_{|\Is_2|\leq |\Is|-1}\big|\del_{\gamma}u_{\kh}\del_{\alpha}\del_{\beta}Z^{\Is_2}v_{\jc}\big|,
\\
& \sum_{|\Is_2|\leq |\Is|-1}\big|\del_{\gamma}v_{\kc}\del_{\alpha}\del_{\beta}Z^{\Is_2}u_{\jh}\big|,\quad
\sum_{|\Is_2|\leq |\Is|-1}\big|\del_{\gamma}v_{\kc}\del_{\alpha}\del_{\beta}Z^{\Is_2}v_{\jc}\big|.
\endaligned
$$
We observe that $[Z^{\Is},A_i^{j\alpha\beta\gamma k}\del_{\gamma}w_k \del_{\alpha}\del_{\beta}]w_j$ is bounded by a linear combination of the following terms with constant coefficients:
$$
\aligned
&Z^{\Is_1}\del_{\gamma}u_{\kh}Z^{\Is_2}\del_{\alpha}\del_{\beta}u_{\jh},\quad
 Z^{\Is_1}\del_{\gamma}u_{\kh}Z^{\Is_2}\del_{\alpha}\del_{\beta}v_{\jc},\quad
 Z^{\Is_1}\del_{\gamma}v_{\kc}Z^{\Is_2}\del_{\alpha}\del_{\beta}u_{\jh},
\\
&Z^{\Is_1}\del_{\gamma}v_{\kc}Z^{\Is_2}\del_{\alpha}\del_{\beta}v_{\jc},\quad
 Z^{\Is_1}\del_{\gamma}u_{\kh}\del_{\alpha}\del_{\beta}Z^{\Is_2}u_{\jh},\quad
 Z^{\Is_1}\del_{\gamma}u_{\kh}\del_{\alpha}\del_{\beta}Z^{\Is_2}v_{\jc},
\\
&Z^{\Is_1}\del_{\gamma}v_{\kc}\del_{\alpha}\del_{\beta}Z^{\Is_2}u_{\jh},\quad
Z^{\Is_1}\del_{\gamma}v_{\kc}\del_{\alpha}\del_{\beta}Z^{\Is_2}v_{\jc}
\endaligned
$$
with $|\Is_1|+|\Is_2|\leq |\Is|$ and $|\Is_2|\leq 4$.

We give the inequalities we use for each term and each partition of the index:
$$
\begin{array}{cccc}
\text{Products} &(5,\leq 0) &(4,\leq 1) &(3,\leq2)
\\
Z^{\Is_1}\del_{\gamma}u_{\kh}Z^{\Is_2}\del_\alpha \del_{\beta}u_{\jh}
&\eqref{proof L2 2ge1a},\eqref{proof decay 2order b}
&\eqref{proof L2 2ge1a'},\eqref{proof decay 2order a'}
&\eqref{proof L2 2ge3a},\eqref{proof decay 2order a}
\\
Z^{\Is_1}\del_{\gamma}u_{\kh}Z^{\Is_2}\del_{\alpha} \del_{\beta}v_{\jc}
&\eqref{proof L2 2ge1a},\eqref{proof decay 2ge2a}
&\eqref{proof L2 2ge1a'},\eqref{proof decay 2ge2a}
&\eqref{proof L2 2ge3a},\eqref{proof decay 2ge2a}
\\
Z^{\Is_1}\del_{\gamma}v_{\kc}Z^{\Is_2}\del_\alpha \del_{\beta}u_{\jh}
&\eqref{proof L2 2ge1a},\eqref{proof decay 2order b}
&\eqref{proof L2 2ge2a},\eqref{proof decay 2order a'}
&\eqref{proof L2 2ge2a},\eqref{proof decay 2order a}
\\
Z^{\Is_1}\del_{\gamma}v_{\kc}Z^{\Is_2}\del_\alpha \del_{\beta}v_{\jc}
&\eqref{proof L2 2ge1a},\eqref{proof decay 2ge1c}
&\eqref{proof L2 2ge2a},\eqref{proof decay 2ge1c}
&\eqref{proof L2 2ge2a},\eqref{proof decay 2ge1a}
\\
Z^{\Is_1}\del_{\gamma}u_{\kh}\del_{\alpha}\del_{\beta}Z^{\Is_2}u_{\jh}
&\eqref{proof L2 2ge1a},\eqref{proof decay 2order b}
&\eqref{proof L2 2ge1a'},\eqref{proof decay 2order a'}
&\eqref{proof L2 2ge3a},\eqref{proof decay 2order a}
\\
Z^{\Is_1}\del_{\gamma}u_{\kh}\del_{\alpha}\del_{\beta}Z^{\Is_2}v_{\jc}
&\eqref{proof L2 2ge1a},\eqref{proof decay 1ge2a}
&\eqref{proof L2 2ge1a'},\eqref{proof decay 1ge2a}
&\eqref{proof L2 2ge3a},\eqref{proof decay 1ge2a}
\\
Z^{\Is_1}\del_{\gamma}v_{\kc}\del_{\alpha}\del_{\beta}Z^{\Is_2}u_{\jh}
&\eqref{proof L2 2ge1a},\eqref{proof decay 2order b}
&\eqref{proof L2 2ge2a},\eqref{proof decay 2order a'}
&\eqref{proof L2 2ge2a},\eqref{proof decay 2order a}
\\
Z^{\Is_1}\del_{\gamma}v_{\kc}\del_{\alpha}\del_{\beta}Z^{\Is_2}v_{\jc}
&\eqref{proof L2 2ge1a},\eqref{proof decay 1ge1c}
&\eqref{proof L2 2ge2a},\eqref{proof decay 1ge1c}
&\eqref{proof L2 2ge2a},\eqref{proof decay 1ge1a}
\end{array}
$$

$$
\begin{array}{ccc}
\text{Products} &(2,\leq3) &(1,\leq 4)
\\
Z^{\Is_1}\del_{\gamma}u_{\kh}Z^{\Is_2}\del_\alpha \del_{\beta}u_{\jh}
&\eqref{proof decay 2ge1a'},\eqref{proof L2 2order a'}
&\eqref{proof decay 2ge3a},\eqref{proof L2 2order a}
\\
Z^{\Is_1}\del_{\gamma}u_{\kh}Z^{\Is_2}\del_{\alpha} \del_{\beta}v_{\jc}
&\eqref{proof decay 2ge1a},\eqref{proof L2 2ge2a}
&\eqref{proof decay 2ge1a},\eqref{proof L2 2ge2a}
\\
Z^{\Is_1}\del_{\gamma}v_{\kc}Z^{\Is_2}\del_\alpha \del_{\beta}u_{\jh}
&\eqref{proof decay 2ge2a},\eqref{proof L2 2order a}
&\eqref{proof decay 2ge2a},\eqref{proof L2 2order a}
\\
Z^{\Is_1}\del_{\gamma}v_{\kc}Z^{\Is_2}\del_\alpha \del_{\beta}v_{\jc}
&\eqref{proof decay 2ge2a},\eqref{proof L2 2ge1c}
&\eqref{proof decay 2ge2a},\eqref{proof L2 2ge1c}
\\
Z^{\Is_1}\del_{\gamma}u_{\kh}\del_{\alpha}\del_{\beta}Z^{\Is_2}u_{\jh}
&\eqref{proof decay 2ge1a'},\eqref{proof L2 2order a'}
&\eqref{proof decay 2ge3a},\eqref{proof L2 2order a}
\\
Z^{\Is_1}\del_{\gamma}u_{\kh}\del_{\alpha}\del_{\beta}Z^{\Is_2}v_{\jc}
&\eqref{proof decay 2ge1a'},\eqref{proof L2 1ge2a}
&\eqref{proof decay 2ge3a},\eqref{proof L2 1ge2a}
\\
Z^{\Is_1}\del_{\gamma}v_{\kc}\del_{\alpha}\del_{\beta}Z^{\Is_2}u_{\jh}
&\eqref{proof decay 2ge2a},\eqref{proof L2 2order a}
&\eqref{proof decay 2ge2a},\eqref{proof L2 2order a}
\\
Z^{\Is_1}\del_{\gamma}v_{\kc}\del_{\alpha}\del_{\beta}Z^{\Is_2}v_{\jc}
&\eqref{proof decay 2ge2a},\eqref{proof L2 1ge1c}
&\eqref{proof decay 2ge2a},\eqref{proof L2 1ge1c}
\end{array}
$$

The remaining terms in the right-hand side of \eqref{proof source higher eq0} are estimated similarly. These two terms can be bounded by the combination of the following terms with constant coefficients:
$$
\aligned
&Z^{\Is_1}v_{\kc}Z^{\Is_2}\del_{\alpha}\del_{\beta}w_j,\quad
 Z^{\Is_1}u_{\kh}Z^{\Is_2}\del_{\alpha}\del_{\beta}u_{\jh},
 \\
&Z^{\Is_1}v_{\kc}\del_{\alpha}\del_{\beta}Z^{\Is_2}w_j,\quad
 Z^{\Is_1}u_{\kh}\del_{\alpha}\del_{\beta}Z^{\Is_2}u_{\jh} 
\endaligned
$$
with $|\Is_1|+|\Is_2|\leq |\Is|$ and $|\Is_2|\leq |\Is|-1$. We will write in details the estimate for the most critical term $Z^{\Is_1}u_{\kh}Z^{\Is_2}\del_{\alpha}\del_{\beta}u_{\jh}$:
$$
\aligned
&\sum_{|\Is_1|+|\Is_2|\leq|\Is|\atop |\Is_2|\leq |\Is|-1} \big|Z^{\Is_1}u_{\kh}Z^{\Is_2}\del_{\alpha}\del_{\beta}u_{\jh}\big|
\\
&\leq \sum_{|\Is_1|\leq 5}\big|Z^{\Is_1}u_{\kh}\del_{\alpha}\del_{\beta}u_{\jh}\big|
 +\sum_{|\Is_1|=4,|\Is_2|\leq 1}\big|Z^{\Is_1}u_{\kh}Z^{\Is_2}\del_{\alpha}\del_{\beta}u_{\jh}\big|
\\
& \quad
 +\sum_{|\Is_1|=3,|\Is_2|\leq 2}\big|Z^{\Is_1}u_{\kh}Z^{\Is_2}\del_{\alpha}\del_{\beta}u_{\jh}\big|
\\
&\quad 
+\sum_{|\Is_1|=2,|\Is_2|\leq 3}\big|Z^{\Is_1}u_{\kh}Z^{\Is_2}\del_{\alpha}\del_{\beta}u_{\jh}\big|
 +\sum_{|\Is_1|=1,|\Is_2|\leq 4}\big|Z^{\Is_1}u_{\kh}Z^{\Is_2}\del_{\alpha}\del_{\beta}u_{\jh}\big|
\\
&=: T_5 + T_4 + T_3 +T_2 + T_1.
\endaligned
$$
The term $T_5$ is estimated by applying \eqref{proof L2 2ge5a} on the first factor and \eqref{proof decay 2order b} on the second factor:
$$
\aligned
\|T_5\|_{L^2(\Hcal_s)}
&= \sum_{|\Is_1|\leq 5}
\big{\|}s^{-1}Z^{\Is_1}u_{\kh}\,s\del_{\alpha}\del_{\beta}u_{\jh}\big{\|}_{L^2(\Hcal_s)}
\\
& \leq \big{\|}s^{-1}Z^{\Is_1}u_{\kh}\big{\|}_{L^2(\Hcal_s)}\big{\|}s\del_{\alpha}\del_{\beta}u_{\jh}\big{\|}_{L^\infty(\Hcal_s)}
\\
&\leq  C(C_1\eps)^2 s^{\delta}\big{\|}t^{1/2}s^{-2}\big{\|}_{L^\infty(\Hcal_s)} \leq C(C_1\eps)^2s^{-1+\delta}, 
\endaligned
$$
where we observe that $s\leq t\leq Cs^2$ in the half-cone $\Kcal$.  

The term $T_4$ is bounded by applying \eqref{proof L2 2ge5a'} and \eqref{proof decay 2order a'}:
$$
\aligned
\|T_4\|_{L^2(\Hcal_s)}
\leq& \sum_{|\Is_1|=4,|\Is_2|\leq 1}
\big{\|}s^{-1}Z^{\Is_1}u_{\kh}\,sZ^{\Is_2}\del_{\alpha}\del_{\beta}u_{\jh}\big{\|}_{L^2(\Hcal_s)}
\\
\leq& C(C_1\eps)^2 s^{\delta/2} \big{\|}s t^{1/2}s^{-3+\delta/2}\big{\|}_{L^\infty(\Hcal_s)}
\\
\leq& C(C_1\eps)s^{-1+\delta}.
\endaligned
$$

The term $T_3$ is bounded by applying \eqref{proof decay 2ge5a} and \eqref{proof L2 2order b}:
$$
\aligned
\|T_3\|_{L^2(\Hcal_s)}
\leq& \sum_{|\Is_1|=3,|\Is_2|\leq 2} \big{\|}s^{-1}t^{3/2}Z^{\Is_1}u_{\kh}\,st^{-3/2}Z^{\Is_2}\del_{\alpha}\del_{\beta}u_{\jh}\big{\|}_{L^2(\Hcal_s)}
\\
\leq& C C_1\eps s^{\delta} \big{\|} s^3t^{-2}Z^{\Is_2}\del_{\alpha}\del_{\beta}u_{\jh}\big{\|}_{L^2(\Hcal_s)}\big{\|}s^{-2}t^{1/2}\big{\|}_{L^{\infty}(\Hcal_s)}
\\
\leq& C(C_1\eps)^2s^{-1+\delta}.
\endaligned
$$

The term $T_2$ is estimated by applying \eqref{proof decay 2ge5a'}, \eqref{proof L2 2order a'}:
$$
\aligned
\|T_2\|_{L^2(\Hcal_s)}
\leq& \sum_{|\Is_1|\leq 2,|\Is_2|\leq 3}
\big{\|}s^{-3}t^2Z^{\Is_1}u_{\kh}\,s^3t^{-2}Z^{\Is_2}\del_{\alpha}\del_{\beta}u_{\jh}\big{\|}_{L^2(\Hcal_s)}
\\
\leq& \sum_{|\Is_1|\leq 2,|\Is_2|\leq 3}
\big{\|}s^{-3}t^2Z^{\Is_1}u_{\kh}\big{\|}_{L^\infty(\Hcal_s)}
\big{\|}s^3t^{-2}Z^{\Is_2}\del_{\alpha}\del_{\beta}u_{\jh}\big{\|}_{L^2(\Hcal_s)}
\\
\leq& CC_1\eps \big{\|}s^{-3}t^2t^{-3/2}s^{1+\delta/2}\big{\|}_{L^\infty(\Hcal_s)}CC_1\eps s^{\delta/2}
\\
=& C(C_1\eps)^2s^{-1+\delta}.
\endaligned
$$

The term $T_1$ is bounded by applying \eqref{proof decay 2ge5b} and \eqref{proof L2 2order a}:
$$
\aligned
\|T_1\|_{L^2(\Hcal_s)} \leq & \sum_{|\Is_1|\leq 1,|\Is_2|\leq 4}
\big{\|}s^{-3}t^2Z^{\Is_1}u_{\kh}\,s^3t^{-2}Z^{\Is_2}\del_{\alpha}\del_{\beta}u_{\jh}\big{\|}_{L^2(\Hcal_s)}
\\
\leq & \sum_{|\Is_1|\leq 1,|\Is_2|\leq 4}
\big{\|}s^{-3}t^2Z^{\Is_1}u_{\kh}\big{\|}_{L^\infty(\Hcal_s)}
\big{\|}s^3t^{-2}Z^{\Is_2}\del_{\alpha}\del_{\beta}u_{\jh}\big{\|}_{L^2(\Hcal_s)}
\\
\leq &CC_1\eps \|s^{-3}t^2t^{-3/2}s\|_{L^\infty(\Hcal_s)} \cdot CC_1\eps s^{\delta}
\leq C(C_1\eps)^2s^{-1+\delta}.
\endaligned
$$

For the remaining terms, we list out the inequalities to be used on each term and each partition of the index as follows:
$$
\begin{array}{cccc}
\text{Products} &(5,\leq 0) &(4,\leq 1) &(3,\leq 2)
\\
Z^{\Is_1}v_{\kc}Z^{\Is_2}\del_{\alpha}\del_{\beta}w_j
&\eqref{proof L2 2ge1c},\eqref{proof decay 2ge1a}
&\eqref{proof L2 2ge1c},\eqref{proof decay 2ge1a}
&\eqref{proof L2 2ge1c},\eqref{proof decay 2ge1a}
\\
Z^{\Is_1}u_{\kh}Z^{\Is_2}\del_{\alpha}\del_{\beta}u_{\jh}
&\eqref{proof L2 2ge5a},\eqref{proof decay 2order b}
&\eqref{proof L2 2ge5a'},\eqref{proof decay 2order a'}
&\eqref{proof decay 2ge5a},\eqref{proof L2 2order b}
\\
Z^{\Is_1}v_{\kc}\del_{\alpha}\del_{\beta}Z^{\Is_2}w_j
&\eqref{proof L2 2ge1c},\eqref{proof decay 1ge1a}
&\eqref{proof L2 2ge1c},\eqref{proof decay 1ge1a}
&\eqref{proof L2 2ge1c},\eqref{proof decay 1ge1a}
\\
Z^{\Is_1}u_{\kh}\del_{\alpha}\del_{\beta}Z^{\Is_2}u_{\jh}
&\eqref{proof L2 2ge5a},\eqref{proof decay 2order b}
&\eqref{proof L2 2ge5a'},\eqref{proof decay 2order a'}
&\eqref{proof decay 2ge5a},\eqref{proof L2 2order b}
\end{array}
$$

$$
\begin{array}{ccc}
\text{Products} &(2,\leq 3) &(1,\leq 4)
\\
Z^{\Is_1}v_{\kc}Z^{\Is_2}\del_{\alpha}\del_{\beta}w_j
&\eqref{proof decay 2ge1c},\eqref{proof L2 2ge1a}
&\eqref{proof decay 2ge1c},\eqref{proof L2 2ge1a}
\\
Z^{\Is_1}u_{\kh}Z^{\Is_2}\del_{\alpha}\del_{\beta}u_{\jh}
&\eqref{proof decay 2ge5a'},\eqref{proof L2 2order a'}
&\eqref{proof decay 2ge5b},\eqref{proof L2 2order a}
\\
Z^{\Is_1}v_{\kc}\del_{\alpha}\del_{\beta}Z^{\Is_2}w_j
&\eqref{proof decay 2ge1c},\eqref{proof L2 1ge1a}
&\eqref{proof decay 2ge1c},\eqref{proof L2 1ge1a}
\\
Z^{\Is_1}u_{\kh}\del_{\alpha}\del_{\beta}Z^{\Is_2}u_{\jh}
&\eqref{proof decay 2ge5a'},\eqref{proof L2 2order a'}
&\eqref{proof decay 2ge5b},\eqref{proof L2 2order a}
\end{array}
$$
\end{proof}

We estimate the source with fourth-order derivatives.

\begin{lemma}\label{proof source lem higher'}
By relying on the assumption \eqref{proof energy assumption}, the following estimates hold for all $|\Id|\leq 4$ and $s_0\leq s\leq s_1$: 
\begin{subequations}\label{proof source higher'}
\bel{proof source higher' F}
\big{\|}Z^{\Id} F_{\ih}\big{\|}_{L^2(\Hcal_s)} \leq C(C_1\eps)^2s^{-1+\delta/2},
\ee
\bel{proof source higher' G}
\big{\|}[Z^{\Id},G_i^{j\alpha\beta}(w,\del w)\del_{\alpha}\del_{\beta}]w_j\big{\|}_{L^2(\Hcal_s)} 
\leq C(C_1\eps)^2s^{-1+\delta/2}.
\ee
\end{subequations}
\end{lemma}

\begin{proof}
The proof is essentially the same as the one of Lemma \ref{proof source lem higher}. We will not write the proof in detail, but we list out the inequalities 
we use for each term and each partition of the index.

First we list out the estimate on $Z^{\Id}F_i$ in Table 15. 
The terms $[Z^{\Id}, B_i^{j\alpha\beta k}w_k\del_{\alpha}\del_{\beta}]w_j$ are estimated   by the inequalities listed in Table 16. 
\end{proof}

\

\section{$L^2$ estimates on third-order terms}
\label{sec:92}

We now derive the $L^2$ estimates for the source terms \eqref{main bootstrap 5}. This is the second place where the null structure is taken into consideration after Lemma \ref{proof lem Mv 3}.

\begin{lemma}\label{proof source lower}
By relying on the assumption of \eqref{proof energy assumption} with
$C_1 \eps < \min (1, \eps_0'')$, the following estimates hold for all $|I|\leq 3$:
\begin{subequations}
\bel{proof source lower F}
\|Z^IF_{\ih}\|_{L^2(\Hcal_s)}\leq C(C_1\eps)^2s^{-3/2+2\delta},
\ee
\bel{proof source lower G1}
\|Z^IG_{\ih}^{\jc\alpha\beta}(w,\del w)\del_{\alpha}\del_{\beta}v_{\jc}\|_{L^2(\Hcal_s)}\leq C(C_1\eps)^2s^{-3/2+2\delta},
\ee
\bel{proof source lower G2}
\big{\|}[Z^I,G_{\ih}^{\jh\alpha\beta}(w,\del w)\del_{\alpha}\del_{\beta}]u_{\jh}\big{\|}_{L^2(\Hcal_s)}\leq C(C_1\eps)^2s^{-3/2+2\delta}.
\ee
\end{subequations}
\end{lemma}


\protect\begin{landscape}\thispagestyle{empty}

$$
\begin{array}{cccccc}
\text{Products} &(4,\leq 0) &(3,\leq 1) &(2,\leq 2) &(1,\leq 3) &(0,\leq 4)
\\
\del_{\alpha}u_{\ih}\del_{\beta}u_{\jh}
&\eqref{proof L2 2ge1a'},\eqref{proof decay 2ge3a}
&\eqref{proof L2 2ge3a},\eqref{proof decay 2ge3a}
&\eqref{proof L2 2ge3a},\eqref{proof decay 2ge1a'}
&\eqref{proof decay 2ge3a},\eqref{proof L2 2ge3a}
&\eqref{proof decay 2ge3a},\eqref{proof L2 2ge1a'}
\\
\del_{\alpha}v_{\ic}\del_{\beta}w_j
&\eqref{proof L2 2ge2a},\eqref{proof decay 2ge1a}
&\eqref{proof L2 2ge2a},\eqref{proof decay 2ge1a}
&\eqref{proof L2 2ge2a},\eqref{proof decay 2ge1a}
&\eqref{proof decay 2ge2a},\eqref{proof L2 2ge1a}
&\eqref{proof decay 2ge2a},\eqref{proof L2 2ge1a}
\\
v_{\ic}\del_{\alpha}w_j
&\eqref{proof L2 2ge1c},\eqref{proof decay 2ge1a}
&\eqref{proof L2 2ge1c},\eqref{proof decay 2ge1a}
&\eqref{proof L2 2ge1c},\eqref{proof decay 2ge1a}
&\eqref{proof decay 2ge1c},\eqref{proof L2 2ge1a}
&\eqref{proof decay 2ge1c},\eqref{proof L2 2ge1a}
\\
v_{\jc}v_{\kc}
&\eqref{proof L2 2ge1c},\eqref{proof decay 2ge1c}
&\eqref{proof L2 2ge1c},\eqref{proof decay 2ge1c}
&\eqref{proof L2 2ge1c},\eqref{proof decay 2ge1c}
&\eqref{proof decay 2ge1c},\eqref{proof L2 2ge1c}
&\eqref{proof decay 2ge1c},\eqref{proof L2 2ge1c}
\end{array}
$$
\centerline{Table 15}

\

\ 

$$
\begin{array}{ccccc}
\text{Products} &(4,\leq 0)&(3,\leq 1)&(2,\leq2)&(1,\leq 3)
\\
Z^{\Id_1}\del_{\gamma}u_{\kh}Z^{\Id_2}\del_{\alpha}\del_{\beta}u_{\jh}
&\eqref{proof L2 2ge1a'},\eqref{proof decay 2order b}
&\eqref{proof L2 2ge3a},\eqref{proof decay 2order a'}
&\eqref{proof decay 2ge1a'},\eqref{proof L2 2order b}
&\eqref{proof decay 2ge3a},\eqref{proof L2 2order a'}
\\
Z^{\Id_1}\del_{\gamma}u_{\kh}\del_{\alpha}\del_{\beta}Z^{\Id_2}u_{\jh}
&\eqref{proof L2 2ge1a'},\eqref{proof decay 2order b}
&\eqref{proof L2 2ge3a},\eqref{proof decay 2order a'}
&\eqref{proof decay 2ge1a'},\eqref{proof L2 2order b}
&\eqref{proof decay 2ge3a},\eqref{proof L2 2order a'}
\\
Z^{\Id_1}\del_{\gamma}u_{\kh}Z^{\Id_2}\del_{\alpha}\del_{\beta}v_{\jc}
&\eqref{proof L2 2ge1a'},\eqref{proof decay 2ge2a}
&\eqref{proof L2 2ge3a},\eqref{proof decay 2ge2a}
&\eqref{proof decay 2ge1a'},\eqref{proof L2 2ge2a'}
&\eqref{proof decay 2ge3a},\eqref{proof L2 2ge2a}
\\
Z^{\Id_1}\del_{\gamma}u_{\kh}\del_{\alpha}\del_{\beta}Z^{\Id_2}v_{\jc}
&\eqref{proof L2 2ge1a'},\eqref{proof decay 1ge2a}
&\eqref{proof L2 2ge3a},\eqref{proof decay 1ge2a}
&\eqref{proof decay 2ge1a'},\eqref{proof L2 1ge2a'}
&\eqref{proof decay 2ge3a},\eqref{proof L2 1ge2a}
\\
Z^{\Id_1}\del_{\gamma}v_{\kc}Z^{\Id_2}\del_{\alpha}\del_{\beta}w_j
&\eqref{proof L2 2ge2a},\eqref{proof decay 2ge1a}
&\eqref{proof L2 2ge2a},\eqref{proof decay 2ge1a}
&\eqref{proof decay 2ge2a},\eqref{proof L2 2ge1a}
&\eqref{proof decay 2ge2a},\eqref{proof L2 2ge1a}
\\
Z^{\Id_1}\del_{\gamma}v_{\kc}\del_{\alpha}\del_{\beta}Z^{\Id_2}w_j
&\eqref{proof L2 2ge2a},\eqref{proof decay 2ge1a}
&\eqref{proof L2 2ge2a},\eqref{proof decay 2ge1a}
&\eqref{proof decay 2ge2a},\eqref{proof L2 1ge1a}
&\eqref{proof decay 2ge2a},\eqref{proof L2 1ge1a}
\\
Z^{\Id_1}v_{\kc}Z^{\Id_2}\del_{\alpha}\del_{\beta}w_j
&\eqref{proof L2 2ge1c},\eqref{proof decay 2ge1a}
&\eqref{proof L2 2ge1c},\eqref{proof decay 2ge1a}
&\eqref{proof decay 2ge1c},\eqref{proof L2 2ge1a}
&\eqref{proof decay 2ge1c},\eqref{proof L2 2ge1a}
\\
Z^{\Id_1}v_{\kc}\del_{\alpha}\del_{\beta}Z^{\Id_2}w_j
&\eqref{proof L2 2ge1c},\eqref{proof decay 2ge1c}
&\eqref{proof L2 2ge1c},\eqref{proof decay 2ge1c}
&\eqref{proof decay 2ge1c},\eqref{proof L2 1ge1c}
&\eqref{proof decay 2ge1c},\eqref{proof L2 1ge1c}
\\
Z^{\Id_1}u_{\kh}Z^{\Id_2}\del_{\alpha}\del_{\beta}u_{\jh}
&\eqref{proof L2 2ge5a'},\eqref{proof decay 2order b}
&\eqref{proof L2 2ge5b},\eqref{proof decay 2order a'}
&\eqref{proof decay 2ge5a'},\eqref{proof L2 2order b}
&\eqref{proof decay 2ge5b},\eqref{proof L2 2order a'}
\\
Z^{\Id_1}u_{\kh}\del_{\alpha}\del_{\beta}Z^{\Id_2}u_{\jh}
&\eqref{proof L2 2ge5a'},\eqref{proof decay 2order b}
&\eqref{proof L2 2ge5b},\eqref{proof decay 2order a'}
&\eqref{proof decay 2ge5a'},\eqref{proof L2 2order b}
&\eqref{proof decay 2ge5b},\eqref{proof L2 2order a'}
\end{array}
$$
\centerline{Table 16}

\

\end{landscape}

\begin{proof} We first prove \eqref{proof source lower F} and recall the structure of $Z^IF_{\ih}$:
$$
\aligned
Z^IF_{\ih} 
&=Z^I \big(P_{\ih}^{\alpha\beta\jh\kh}\del_{\alpha}u_{\jh}\del_{\beta}u_{\kh} \big)
\\
&\quad 
+ P_{\ih}^{\alpha\beta\jc\kh}Z^I\big(\del_{\alpha}v_{\jc}\del_{\beta}u_{\kh}\big)
 + P_{\ih}^{\alpha\beta\jh\kc}Z^I\big(\del_{\alpha}u_{\jh}\del_{\beta}v_{\kc}\big)
 + P_{\ih}^{\alpha\beta\jc\kc}Z^I\big(\del_{\alpha}v_{\jc}\del_{\beta}v_{\kc}\big)
\\
&\quad 
+Q_{\ih}^{\alpha j\kc}Z^I\big(v_{\kc}\del_{\alpha}w_j\big) + R_{\ih}^{\jh\kh}Z^I\big(v_{\jh}v_{\kh}\big).
\endaligned
$$
The first term in the right-hand side is a null term and we can apply directly \eqref{proof null eq1a}:
$$
\big{\|}P_{\ih}^{\alpha\beta\jh\kh}\del_{\alpha}u_{\jh}\del_{\beta}u_{\kh}\big{\|}_{L^2(\Hcal_s)}\leq C(C_1\eps)^2s^{-3/2}.
$$
The remaining terms are linear combinations of the following terms with constant coefficients:
$$
Z^I\big(\del_{\alpha}v_{\jc}\del_{\beta}w_k\big),\quad
Z^I\big(v_{\kc}\del_{\alpha}w_j\big),\quad
Z^I\big(v_{\jh}v_{\kh}\big).
$$
We will not give in details the estimates of each term, but we list out the inequalities
we use for each term and each partition:
$$
\begin{array}{ccccc}
\text{Products} &(3,\leq 0) &(2,\leq 1) &(1,\leq 2) &(0,\leq 3)
\\
\del_{\alpha}v_{\jc}\del_{\beta}w_k
&\eqref{proof L2 2ge2a},\eqref{proof decay 2ge1a}
&\eqref{proof L2 2ge2a},\eqref{proof decay 2ge1a}
&\eqref{proof decay 2ge2a},\eqref{proof L2 2ge1a}
&\eqref{proof decay 2ge2a},\eqref{proof L2 2ge1a}
\\
v_{\kc}\del_{\alpha}w_j
&\eqref{proof L2 2ge1c},\eqref{proof decay 2ge1a}
&\eqref{proof L2 2ge1c},\eqref{proof decay 2ge1a}
&\eqref{proof decay 2ge1c},\eqref{proof L2 2ge1a}
&\eqref{proof decay 2ge1c},\eqref{proof L2 2ge1a}
\\
v_{\jh}v_{\kh}
&\eqref{proof L2 2ge1c},\eqref{proof decay 2ge1c}
&\eqref{proof L2 2ge1c},\eqref{proof decay 2ge1c}
&\eqref{proof decay 2ge1c},\eqref{proof L2 2ge1c}
&\eqref{proof decay 2ge1c},\eqref{proof L2 2ge1c}
\end{array}
$$

To establish the second inequality \eqref{proof source lower G1}, we decompose it as follows:
$$
Z^I\big(G_{\ih}^{\jc\alpha\beta}(w,\del w)\del_{\alpha}\del_{\beta}v_{\jc}\big)
 = A_{\ih}^{\jc\alpha\beta\gamma k}Z^I\big(\del_{\gamma}w_k\del_{\alpha}\del_{\beta}v_{\jc}\big)
 + B_{\ih}^{\jc\alpha\beta \kc}Z^I\big(v_{\kc}\del_{\alpha}\del_{\beta}v_{\jc}\big). 
$$
So, it is to be bounded by a linear combination of the following terms with constant coefficients:
$$
Z^I\big(\del_{\gamma} w_k\del_{\alpha}\del_{\beta}v_{\jc}\big),\quad
Z^I\big(v_{\kc}\del_{\alpha}\del_{\beta} v_{\jc}\big).
$$
As before, we list out the inequalities 
we use for each term and each partition of the index: 
$$
\begin{array}{ccccc}
\text{Products}  &(3,\leq 0) &(2,\leq 1) &(1,\leq 2) &(0,\leq 3)
\\
\del_{\gamma} w_k\del_{\alpha}\del_{\beta}v_{\jc}
&\eqref{proof L2 2ge1a},\eqref{proof decay 2ge2a}
&\eqref{proof L2 2ge1a},\eqref{proof decay 2ge2a}
&\eqref{proof decay 2ge1a},\eqref{proof L2 2ge2a}
&\eqref{proof decay 2ge1a},\eqref{proof L2 2ge2a}
\\
v_{\kc}\del_{\alpha}\del_{\beta}v_{\jc}
&\eqref{proof L2 2ge1c},\eqref{proof decay 2ge2a}
&\eqref{proof L2 2ge1c},\eqref{proof decay 2ge2a}
&\eqref{proof decay 2ge1c},\eqref{proof L2 2ge2a}
&\eqref{proof decay 2ge1c},\eqref{proof L2 2ge2a}
\end{array}
$$

The proof of the inequality \eqref{proof source lower G2} is also related to the null structure. Recall the following decomposition:
$$
\aligned
&[Z^I,G_{\ih}^{\jh\alpha\beta}(w,\del w)\del_{\alpha}\del_{\beta}]u_{\jh}
\\
 =&[Z^I,A_{\ih}^{\jh\alpha\beta\gamma\kh}\del_{\gamma}u_{\kh}\del_{\alpha}\del_{\beta}]u_{\jh} + [Z^I,B_{\ih}^{\jh\alpha\beta\kh}u_{\kh}\del_{\alpha}\del_{\beta}]u_{\jh}
\\
&+[Z^I,A_{\ih}^{\jh\alpha\beta\gamma\kc}\del_{\gamma}v_{\kc}\del_{\alpha}\del_{\beta}]u_{\jh}
+ [Z^I,B_{\ih}^{jh\alpha\beta\kc}v_{\kc}\del_{\alpha}\del_{\beta}]u_{\jh}
\endaligned
$$
Observing the null structure of the first two terms and applying \eqref{proof null eq1c} and \eqref{proof null eq2b3}, we obtain
$$
\aligned
& \big{\|}[Z^I,A_{\ih}^{\jh\alpha\beta\gamma\kh}\del_{\gamma}u_{\kh}\del_{\alpha}\del_{\beta}]u_{\jh}\big{\|}_{L^2(\Hcal_s)} + \big{\|}[Z^I,B_{\ih}^{\jh\alpha\beta\kh}u_{\kh}\del_{\alpha}\del_{\beta}]u_{\jh}\big{\|}_{L^2(\Hcal_s)}
\\
& \leq C(C_1\eps)^2s^{-3/2+\delta}.
\endaligned
$$
The remaining two terms are linear combinations of the following term:
$$
\aligned
&
Z^{I_1}\del_{\gamma}v_{\kc}Z^{I_2}\del_{\alpha}\del_{\beta}u_{\jh},\quad
Z^{I_1}v_{\kc}Z^{I_2}\del_{\alpha}\del_{\beta}u_{\jh},
\\
&
Z^{I_1}\del_{\gamma}v_{\kc}\del_{\alpha}\del_{\beta}Z^{I_2}u_{\jh},\quad
Z^{I_1}v_{\kc}\del_{\alpha}\del_{\beta}Z^{I_2}u_{\jh}
\endaligned$$
with $|I_1|+|I_2|\leq |I|$ and $|I_2|\leq 2$.

We omit here the details and list out the inequalities applied to each term and each partition of the index:
$$
\begin{array}{cccc}
\text{Products} &(3,\leq 0) &(2,\leq 1) &(1,\leq 2)
\\
Z^{I_1}\del_{\gamma}v_{\kc}Z^{I_2}\del_{\alpha}\del_{\beta}u_{\jh}
&\eqref{proof L2 2ge2a},\eqref{proof decay 2ge1a}
&\eqref{proof L2 2ge2a},\eqref{proof decay 2ge1a}
&\eqref{proof decay 2ge2a},\eqref{proof L2 2ge1a}
\\
Z^{I_1}v_{\kc}Z^{I_2}\del_{\alpha}\del_{\beta}u_{\jh}
&\eqref{proof L2 2ge1c},\eqref{proof decay 2ge1a}
&\eqref{proof L2 2ge1c},\eqref{proof decay 2ge1a}
&\eqref{proof decay 2ge1c},\eqref{proof L2 2ge1a}
\\
Z^{I_1}\del_{\gamma}v_{\kc}\del_{\alpha}\del_{\beta}Z^{I_2}u_{\jh}
&\eqref{proof L2 2ge2a},\eqref{proof decay 1ge1a}
&\eqref{proof L2 2ge2a},\eqref{proof decay 1ge1a}
&\eqref{proof decay 2ge2a},\eqref{proof L2 1ge1a}
\\
Z^{I_1}v_{\kc}\del_{\alpha}\del_{\beta}Z^{I_2}u_{\jh}
&\eqref{proof L2 2ge1c},\eqref{proof decay 1ge1a}
&\eqref{proof L2 2ge1c},\eqref{proof decay 1ge1a}
&\eqref{proof decay 2ge1c},\eqref{proof L2 1ge1a}
\end{array}
$$
\end{proof}

\chapter[The local well-posedness theory]{The local well-posedness theory \label{cha:11}}

\section{Construction of the initial data}
\label{sec:101} 

In this chapter, we establish the local-in-time existence theory for general systems
of the form
\bel{appendix eq nonlinear}
\aligned
&\Box w_i + G_i^{j\alpha\beta}(w,\del w)\del_{\alpha}\del_{\beta} w_j + c_i^2 w_i = F_i(w,\del w),
\\
&w_i(B+1,x) = {w_i}_0,\quad \del_t w_i(B+1,x) = {w_i}_1,
\endaligned
\ee
where
\be
\aligned
&G_i^{j\alpha\beta}(w,\del w) = A^{j\alpha\beta\gamma k}_i\del_{\gamma}w_k + B^{j\alpha\beta k}w_k,
\\
&F_i(w,\del w) = P^{\alpha\beta jk}_i \del_\alpha w_j \del_{\beta}w_k + Q^{\alpha jk}_i \del_\alpha w_j w_k + R^{jk}_i w_j w_k, 
\endaligned
\ee
with constants $A^{\alpha\beta\gamma j}_i,B^{j\alpha\beta k},P^{\alpha\beta jk}_i,Q^{\alpha jk}_i,R^{jk}_i$. 
To guarantee the hyperbolicity property, the following symmetry conditions are assumed:
\be
G_i^{j\alpha\beta} = G_j^{i\alpha\beta},\qquad G_i^{j\alpha\beta} = G_i^{j\beta\alpha}.
\ee

Initial data by $({w_i}_0,{w_i}_1) \in \Hf^{l+1}(\RR^3)\times \Hf^l(\RR^3)$ of sufficiently high regularity 
are prescribed on the hyperplane $\{t = B+1\}$ and we assume the smallness condition:
\be
\|{w_i}_0\|_{\Hf^{l+1}(\RR^3)}
+
 \|{w_i}_1\|_{\Hf^l(\RR^3)}\leq \eps'.
\ee
We also assume that ${w_i}_0$ and ${w_i}_1$ are supported in the ball $\{|x| \leq B\}$.
The following theorem was essentially already established in \citet{Sogge08}.

\begin{theorem}\label{appendix Thm A}
For any integer $l\geq 5$, 
there exists a time interval $[B+1,T(\eps')+B+1]$ on which
the Cauchy problem \eqref{appendix eq nonlinear} admits a (unique) solution
$w_i=w_i(t,x)$ and one has 
\be
w_i(t,x) \in C\big([B+1,T(\eps)+B+1], \Hf^{l+1}(\RR^3)\big)
\ee
and there exists a constant $A>0$ such that
\be
\sum_i \sum_{\alpha \atop |I|\leq 5}\|\del_\alpha Z^I w_i\|_{L^\infty([B+1,T(\eps)],L^2(\RR^3))}\leq A\eps'.
\ee
Moreover, the time of existence approaches infinity when the size of the data approaches zero, that is, 
\be
\lim_{\eps'\rightarrow 0^+}T(\eps') = +\infty. 
\ee
If $T^*$ denotes the supremum of all such times $T(\eps')$ (for fixed initial data), then either $T^*=+\infty$ or else 
\be
\sup_{t \in [0, T^*] \atop x \in \RR} \sum_i 
\sum_{|\beta|\leq 5}|\del^\beta w_i(t,x)| = +\infty.  
\ee
In addition, for $C_0>1$ sufficiently large and some $\eps_0'>0$ (depending only on the structure of the system), the uniform bound 
\be
\sum_{i}E_G(B+1, w_i) \leq C_0\eps'
\ee
holds for all initial data satisfying $\eps'\leq \eps'_0$, 
where $E_G(\cdot,w_i) := E_{G,0}(\cdot,w_i)$ is the hyperboloidal energy
in \eqref{pre expression of curved energy}. 
\end{theorem}

\begin{proof} We only sketh the proof.
The local existence theory and the blow-up criteria are 
discussed in \citet{Sogge08} (cf.~by Theorem 4.1 in 
Section 1.4 therein).  The property of $T(\eps')$ is deduced from our proof of 
global existence in this monograph. Therefore, 
 we can focus here on estimating the energy $E_G(B+1,Z^I w_i)$.

We consider the region $\Kcal_{B+1} : = \{(t,x): t\geq B+1, t^2-|x|^2 \leq (B+1)^2\}\cap \Kcal$ where $\Kcal = \{|x|< t-1, t\geq 0\}$. We observe that in $\Kcal_{B+1}$, $t\leq \frac{(B+1)^2+1}{2}$. 
We can fix some $\eps'_0$ sufficiently
small, so that $T(\eps')\geq \frac{(B+1)^2+1}{2}$. 
The local solution is well defined in the region $\Kcal_{B+1}$.

To estimate $E_G(B+1,Z^I w_i)$, we choose $\del_t w_i$ as a multiplier and, with the notation 
\be
\aligned
&E_G^*(B+1,w_i)
\\
& :=  \int_{\RR^3} \Big(
\sum_\alpha |\del_\alpha w_i|^2 + 2 G_i^{j \alpha\beta} \del_t w_i \del_\beta w_j - G_i^{j\alpha\beta} \del_\alpha w_i \del_\beta w_j \Big) (B+1, \cdot) \, dx
\endaligned
\ee
in order to denote the standard energy defined on the flat hypersurface $t=B+1$, 
we obtain the energy estimate 
$$
\aligned
& \sum_i E_G(B+1,Z^I w_i) - \sum_i  E^*_G(B+1,Z^I w_i)
\\
&= \sum_i  \int_{\Kcal_{B+1}} \big(Z^I F_i(w, \del w) \, \del_t Z^I w_i - [Z^I,G_i^{j\alpha\beta}\del_{\alpha\beta}] Z^Iw_j\, \del_t Z^Iw_i\big) dxdt
\\
&\quad 
+ \sum_i  \int_{\Kcal_{B+1}}\bigg(\del_\alpha G_i^{j\alpha\beta}\del_tZ^Iw_i\del_{\beta}Z^Iw_j - \frac{1}{2}\del_tG_i^{j\alpha\beta}\del_\alpha Z^Iw_i\del_{\beta}Z^Iw_j\bigg) dxdt
\\
&\leq CK\sum_{i,j,k\atop \alpha,\beta}
\int_{B+1}^{\frac{(B+1)^2+1}{2}}|\del_\alpha Z^Iw_i| \, |\del_{\beta}Z^Jw_j| \, |\del_{\gamma}Z^{J'}w_k|dxdt
 \leq  C A \eps'.
\endaligned
$$
Here, $C$ is a constant depending only on the structure of the system. We have
$$
E_G(B+1,Z^I w_i)\leq CA\eps' +C\eps'.
$$
Thus, for $C_0>1$ sufficiently large and $\eps'>0$, we can find some $\eps'_0$ sufficiently small such that
$$
CA\eps' +C\eps'\leq C_0 \eps'.
$$
\end{proof}


\section{Local well-posedness  theory in the hyperboloidal foliation}
\label{sec:102}

For convenience, we introduce the following notation in the cone $\Kcal$ 
$$
\aligned
&\xb^0 := s = \sqrt{t^2-|x|^2},
\\
&\xb^a: = x^a.
\endaligned
$$
The natural frame associated with these variables is
$$
\aligned
&\delb_0 = \del_s = \frac{s}{t}\del_t,
\\
&\delb_a = \delu_a = \frac{x^a}{t}\del_t + \del_a.
\endaligned
$$
By an easy calculation, we express the wave operator in this frame:
\be
\Box u = \delb_0\delb_0 u +\frac{2\xb^a}{s}\delb_0\delb_a u - \sum_a\delb_a\delb_a u + \frac{3}{s}\delu_a u.
\ee
The  symmetric second-order quasi-linear operator $G_i^{j\alpha\beta}(w,\del w)\del_{\alpha}\del_{\beta}w_j$  can also be expressed in this frame:
\be
\aligned
& G_i^{j\alpha\beta}(w,\del w)\del_{\alpha}\del_{\beta}w_j
\\
& = \Gbar_i^{j\alpha\beta}(w,\delb w)\del_{\alpha}\del_{\beta}w_j
+G_i^{j\alpha\beta}(w,\del w)\del_{\alpha}\big(\Psi_{\beta}^{\beta'}\big)\delb_{\beta'}w_j
\endaligned
\ee
with
$$
\Gbar_i^{j\alpha\beta} = G_i^{j\alpha'\beta'}\Psi_{\alpha'}^{\alpha}\Psi_{\beta'}^{\beta}, 
$$
which is again symmetric with respect to $i,j$.

Now, we can transform the system \eqref{appendix eq nonlinear}: 
\bel{appendix eq nonlinear-trans}
\aligned
\delb_0\delb_0w_i + \Gbar_i^{j00}\delb_0\delb_0 w_j
& + \frac{\xb^a}{s}\delb_0\delb_a w_i + \Gbar_i^{j0a}(w,\del w)\delb_0\delb_aw_j
\\
&
-\sum_a\delb_a\delb_a w_i + \Gbar_i^{jab}\delb_a\delb_bw_j = \tilde{F}_i(w,\del w)
\endaligned
\ee
with initial data posed on ${s = B+1}$:
$$
w_i(B+1,\cdot) = {w'_i}_0 , \quad \del_s w_i(B+1,\cdot) = {w'_i}_1.
$$
We observe that this system is again symmetric and
that the perturbation terms $\Gbar$ are small,
so this system is again hyperbolic with respect to the variables $(s,\xb^a)$. 
By the
standard theory of second-order symmetric hyperbolic system, the following theorem holds.
(See, for example, \citet{Taylor}, Proposition 3.1 in Chap. 16, Section~3.)

\begin{theorem}
\label{local semi-hyper}
Let $l \geq 5$ be an integer and $B>0$. Then, the following property holds for a sufficiently small  $\eps>0$.
Let $({w_i}_0,{w_i}_1) \in \Hf^{l+1}\times \Hf^l$ 
be data supported in the ball centered at the origin and with radius $\frac{(B+1)^2-1}{2}$. Then, the Cauchy problem 
\be
\aligned
&\Box w_i + G_i^{j\alpha\beta}(w,\del w)\del_{\alpha\beta} w_j + c_i^2 w_i 
= F_i(w,\del w),
\\
&w_i|_{\Hcal_{B+1}} = {w_i}_0,\quad \del_t w_i|_{\Hcal_{B+1}} = {w_i}_1
\endaligned
\ee
under the condition (on the initial slice) 
\be
\sum_i 
\sum_{|\beta| \leq l+1}\|\del^{\beta} w_i\|_{L^2(\Hcal_{B+1})} 
+ \sum_i \sum_{|\beta|\leq l}\|\del^{\beta} \del_t w_i\|_{L^2(\Hcal_{B+1})}\leq \eps
\ee
admits a unique local-in-time solution, which satisfies 
\be
\sum_i 
\sum_{|\beta| \leq l+1}\|\del^{\beta}w_i\|_{L^2(\Hcal_s)}\leq A\eps 
\ee
in the time interval $[B+1,T(\eps)+B+1]$ and for soe constant $A>0$. 
Furthermore, if $T^*$ denotes the supremum of all such times $T(\eps)$ (for $\eps$ fixed), then either $T^*=+\infty$ or else 
\be
\lim_{s\rightarrow {T^*} \atop s < T^*} 
\sum_i \sum_{\alpha\leq l+1}\|\del^{\alpha}w_i\|_{L^2(\Hcal_s)} = +\infty.
\ee
\end{theorem}





\end{document}